\documentclass[reqno,12pt]{amsbook}
\usepackage{amsmath, latexsym, amsfonts, amssymb, amsthm, amscd}
\usepackage[latin1]{inputenc} 

\usepackage{epsfig}
\usepackage{graphics,epsf,psfrag}

\numberwithin{equation}{chapter}
\numberwithin{figure}{chapter}

\newtheorem{theorem}{Theorem}[chapter]
\newtheorem{lemma}[theorem]{Lemma}
\newtheorem{proposition}[theorem]{Proposition}

\newtheorem{rem}[theorem]{Remark}
\newtheorem{definition}[theorem]{Definition}

\newtheorem{assumption}[theorem]{Assumption}

\newcommand{\bbE}{{\ensuremath{\mathbb E}} }

\newcommand{\bbL}{{\ensuremath{\mathbb L}} }

\newcommand{\bbN}{{\ensuremath{\mathbb N}} }

\newcommand{\bbP}{{\ensuremath{\mathbb P}} }

\newcommand{\bbS}{{\ensuremath{\mathbb S}} }

\newcommand{\bbZ}{{\ensuremath{\mathbb Z}} }

\newcommand{\cA}{{\ensuremath{\mathcal A}} }

\newcommand{\cC}{{\ensuremath{\mathcal C}} }
\newcommand{\cD}{{\ensuremath{\mathcal D}} }

\newcommand{\cF}{{\ensuremath{\mathcal F}} }
\newcommand{\cG}{{\ensuremath{\mathcal G}} }
\newcommand{\cH}{{\ensuremath{\mathcal H}} }

\newcommand{\cL}{{\ensuremath{\mathcal L}} }
\newcommand{\cM}{{\ensuremath{\mathcal M}} }
\newcommand{\cN}{{\ensuremath{\mathcal N}} }

\newcommand{\cP}{{\ensuremath{\mathcal P}} }

\newcommand{\cR}{{\ensuremath{\mathcal R}} }

\newcommand{\cU}{{\ensuremath{\mathcal U}} }

\newcommand{\cZ}{{\ensuremath{\mathcal Z}} }

\newcommand{\ga}{\alpha}
\newcommand{\gb}{\beta}
\newcommand{\gga}{\gamma}            

\newcommand{\gd}{\delta}
\newcommand{\gD}{\Delta}
\newcommand{\gep}{\varepsilon}       

\newcommand{\gk}{\kappa}
\newcommand{\gl}{\lambda}

\newcommand{\gr}{\rho}

\newcommand{\gs}{\sigma}
\newcommand{\gS}{\Sigma}
\newcommand{\gp}{\varphi}

\newcommand{\go}{\omega}
\newcommand{\gO}{\Omega}

\renewcommand{\tilde}{\widetilde}          
\DeclareMathSymbol{\leqslant}{\mathalpha}{AMSa}{"36} 
\DeclareMathSymbol{\geqslant}{\mathalpha}{AMSa}{"3E} 
\DeclareMathSymbol{\eset}{\mathalpha}{AMSb}{"3F}     
\newcommand{\dd}{\text{\rm d}}             

\DeclareMathOperator*{\union}{\bigcup}       
\DeclareMathOperator*{\inter}{\bigcap}       
\newcommand{\sumtwo}[2]{\sum_{\substack{#1 \\ #2}}} 
\newcommand{\limtwo}[2]{\lim_{\substack{#1 \\ #2}}}     

\newcommand{\R}{\mathbb{R}}
\newcommand{\Z}{\mathbb{Z}}
\newcommand{\N}{\mathbb{N}}
\newcommand{\Q}{\mathbb{Q}}

\DeclareMathOperator{\sign}{sign}

\def\bP{\ensuremath{\bs{\mathrm{P}}}}
\def\bE{\ensuremath{\bs{\mathrm{E}}}}

\newcommand{\ind}{\bs{1}}

\def\bs{\boldsymbol}
\def\proof{\noindent\textbf{Proof.} \ }

\def\free{{f}}
\def\freea{\tilde {f}}

\newcommand{\tf}{\textsc{f}}

\newcommand{\M}{\textsf{M}}

\def\oT{\overline{T}}
\def\oH{\overline{H}}
\def\hG{\widehat{G}}
\def\tC{\mathtt{C}}
\newcommand{\mysim}{\overset{\star}{\sim}}

\def\free{{f}}
\def\freea{\tilde {f}}

\def\rc{\mathrm{c}}
\def\rf{\mathrm{f}}

\def\tB{\tilde{B}}
\def\tL{\tilde{L}}

\def\tf{\textsc{f}}
\def\tPhi{\tilde \Phi}
\def\tL{\tilde L}

\newcommand{\obbS}{\overline{\ensuremath{\mathbb S}} }
\newcommand{\obbN}{\overline{\ensuremath{\mathbb N}} }

\title[Random walk models for inhomogeneous polymer chains]{Random walk models
and\\ probabilistic techniques\\for inhomogeneous polymer chains\\
{\large \hfill\\
Ph.D. Thesis}}

\author{Francesco Caravenna
\break}

\address{Universit\`a degli Studi di Milano -- Bicocca (Italy)\hfill\break
\indent\textit{and}\hfill\break
\indent{}Universit{\'e} Paris~7 -- Denis Diderot (France)}
\email{f.caravenna\@@sns.it}
\urladdr{http://www.matapp.unimib.it/\raisebox{0.15ex}{\tiny$\sim$}fcaraven/
\hfill\break
\hfill\break
\hfill\break
\indent\textnormal{$\,$Supervisors: \ Prof. Alberto Gandolfi
\ and \ Prof. Giambattista Giacomin\hfill\break}}

\date{\today}

\begin{document}

\frontmatter


\thispagestyle{empty}

\begin{center}

\large \bf

UNIVERSIT\'E PARIS 7 -- DENIS DIDEROT

UNIVERSIT\`A DI MILANO -- BICOCCA

\vskip 0.7cm

\large \sc

Thèse de Doctorat en Cotutelle

Tesi di Dottorato in Cotutela

\vskip 0.8cm

\normalsize \rm

\hskip -0.3cm
\begin{tabular}{cc}
pour l'obtention du Diplôme de \ \ \ \ &\hskip 0.7cm \ \ \ \ per il conseguimento del Diploma di\\
\rule{0pt}{16pt}\sc \large Docteur de l'Université \ \ &\hskip 1cm  \ \ \sc \large Dottore di Ricerca\\
\hskip -0.5cm
\begin{tabular}{l}
\rule{0pt}{15pt}Discipline~: \ Mathématiques \\
\rule{0pt}{12pt}Spécialité~: \ Probabilités
\end{tabular}
& \hskip 1cm
\begin{tabular}{l}
\rule{0pt}{15pt} \ Disciplina~: \ Matematica\\
\rule{0pt}{12pt} \ Specialità~: \ Probabilità\\
\end{tabular}
\end{tabular}

\vskip 0.9 cm

Présentée par~: / Presentata da~:

\large \bf Francesco CARAVENNA

\vskip 0.3 cm

\rule{5cm}{1pt}
\hrule
\vskip 0.5cm
\LARGE \bf
Random walk models and\\
probabilistic techniques\\
for inhomogeneous polymer chains
\vskip 0.5cm
\hrule
\rule{5cm}{1pt}

\vskip 0.8 cm
\normalsize \rm

Thèse dirigée par~: / Tesi diretta da~:

\large
Giambattista GIACOMIN \ \ \ \normalsize et$\,$/$\,$e \large \ \ \ Alberto GANDOLFI

\vskip 0.6 cm
\normalsize \rm

Soutenue le 21 Octobre 2005 devant le jury composé de~:

Difesa il 21 Ottobre 2005 davanti alla commissione composta da~:


\vskip 0.5cm

\begin{tabular}{ll}
M. Jean BERTOIN & Examinateur / Esaminatore\\
M. Francis COMETS & Examinateur / Esaminatore\\
M. Alberto GANDOLFI & Directeur de Thèse / Relatore di Tesi\\
M. Giambattista GIACOMIN \ \ \ \ \ \ \ \ \ \ \ \ \ \ & Directeur de Thèse / Relatore di Tesi\\
M. Yueyun HU & Rapporteur / Referee
\end{tabular}

\end{center}

\newpage

\thispagestyle{empty}
\mbox{}

\newpage

\thispagestyle{empty}

\mbox{}
\vskip 4.95em

\begin{center}
\large \textbf{Acknowledgements}
\end{center}

\normalsize
\vskip 1.1 em

I would like to thank my supervisors Alberto Gandolfi and
Giambattista Gia\-co\-min for their guidance and support.
I owe my deep gratitude to Giambattista
\mbox{Giacomin} for introducing me to this fascinating field, for his
contagious and enthusiastic way of working and for all the things
he has taught me. I thank Erwin Bolthausen and Yueyun Hu for agreeing
to act as referees. I am very grateful to Jean Bertoin,
to Francis Comets and to Yueyun Hu for accepting to be part of
the thesis committee.

\vskip 0.18cm

I want to express my special thanks
to all my Ph.D. colleagues at the University of Milano--Bicocca
and at the University of Paris 7--Denis Diderot for the great time we have
spent together. I am also grateful to the personnel of both departments
for the stimulating ambience. I am very much indebted to Lorenzo Zambotti for
all the discussions we have had, for his constant encouragement and for his sympathy.
Je tiens aussi à remercier Madame Wasse pour son aide généreuse, sans laquelle je
n'aurais jamais surmonté les difficultés bureaucratiques engendrées par la cotutelle.

\vskip 0.18cm

Last but not least, I want to express all my gratitude to my friends,
to my family and to Laura. Grazie per essermi stati vicini in questi anni
meravigliosi.

\mbox{}
\newpage


\tableofcontents

\mainmatter

\chapter*{Introduction}

Modeling of polymer chains, that is long linear molecules
made up of a sequence of simpler units called monomers,
has, for a lot of time, received a lot of attention in
physics, chemistry, biology, \ldots
Mathematics belongs to this list too. For example,
probabilistic models that naturally arise in statistical
mechanics have been widely studied by mathematicians
for the very challenging and novel problems that
they pose.  This is true to the extent that, in probability,
the word {\sl polymer} has become synonymous with
{\sl self--avoiding walk}, a basic and extremely difficult
mathematical entity. The interaction of a polymer with
the environment leads to even more challenging questions:
these are often tackled in the framework of {\sl directed walks}.
Restricting attention to directed trajectories is a way
of enforcing the self--avoiding constraint that leads to
much more tractable models. Still, the interaction with the environment
may quickly lead to extremely difficult questions.

A particularly interesting situation is that of an
inhomogeneous polymer (or {\sl copolymer}) in the proximity
of an interface between two selective solvents. The polymer is inhomogeneous in
that its monomers may differ in some characteristics and, consequently, the
interaction with the solvents and the interface may vary from monomer to
monomer. In interesting cases there can be a phase transition between a state
in which the polymer sticks very close to the interface (localized regime)
and a state in which it wanders away from it (delocalized regime).
The typical mechanism underlying such phase transitions is an
{\sl energy/entropy competition}.

\vskip 3mm

The main task of this Ph.D.~thesis is to introduce and
study {\sl random walk models} of polymer chains
with the purpose of understanding
this competition in a deep and quantitative way.
Since a random walk can be regarded as
an example of an abstract polymer,
the idea of modeling real polymers using random
walks is quite natural and it has proved to be very successful.

The models we are going to consider are modifications of a basic
model introduced in the late eighties by T.~Garel, D.~A.~Huse,
S.~Leibler and H.~Orland~\cite{cf:GHLO} that in turn had translated
into the language of theoretical physics ideas that were developing in
the applied sciences. Despite the fact that the definition of these
models is extremely elementary, their analysis is not simple at all.
For a number of interesting issues there is still no agreement in
the physical literature. From a mathematical viewpoint it has taken
quite a lot of time and effort to rigorously derive their basic
properties, and several interesting questions are still open.

In this Ph.D.~thesis we present new results that answer some of these questions.
The approach taken here is essentially probabilistic, and
it is interesting to note how the analysis performed has required the application of a
wide range of techniques including Large Deviations
and Concentration Inequalities (Ch.~\ref{ch:cgg}), Perron--Frobenius Theory
(Ch.~\ref{ch:cg}), Renewal Theory (Ch.~\ref{ch:cgz}) and Fluctuation Theory for random walks
(Ch.~\ref{ch:continuous} and~\ref{ch:llt}).
A numerical and statistical study has also been performed (Ch.~\ref{ch:cgg}).
Reciprocally, the study of the models stimulates the extension
of these techniques, see, for instance, the Local Limit Theorem for random walks
conditioned to stay positive presented in Chapter~\ref{ch:llt}.

\vskip 3mm

The thesis is organized as follows.
The definition of the models we consider is given in detail in
Chapter~\ref{ch:first}, where we also give some motivation and we
collect the known results from the literature. The following five chapters
contain original results. A detailed outline of the thesis
may be found in Section~\ref{sec:out} of Chapter~\ref{ch:first}.

%
%
\chapter{Inhomogeneous polymer chains}
\label{ch:first}


In this first chapter, we introduce the class of models that are the center
of our analysis, providing some motivations for their study and recalling the
known results in the literature. The exposition takes inspiration from~\cite{cf:G}.


\smallskip
\section{Introduction and motivations}

\smallskip
\subsection{Polymers and random walks}

The notion of polymer has
originated in the field of chemistry to indicate
a natural or synthetic compound consisting of large molecules which are made up
of a linked series of repeated simple molecules called monomers.
However this concept has spread out and nowadays polymers appear in a variety of different
fields with the broader meaning of \textsl{linear structures} which are built up by
\textsl{joining together} a large number of \textsl{simpler structures} (everything
possibly in an abstract sense). A very relevant example of an abstract polymer is
given by a \textsl{random walk}, where the increments are thought of as monomers.

As a matter of fact both the above meanings of polymers, the concrete one
and the abstract one, are of interest to us.
In fact our main purpose is to build probabilistic models based on random walks
that try to mimic the phenomenology of true chemical polymers in some interesting
situations. We stress from
now that the models we are going to consider are {\sl very simple} and nevertheless
they pose very challenging problems.

\smallskip
\subsection{Copolymers at selective interfaces}

\begin{figure}[t]
\bigskip
\bigskip
\psfragscanon
\psfrag{Oil}[c][c]{\Large \bf Oil}
\psfrag{Water}[c][c]{\Large \bf Water}
\psfrag{Interface}[c][c]{\large Interface}
\centerline{\psfig{file=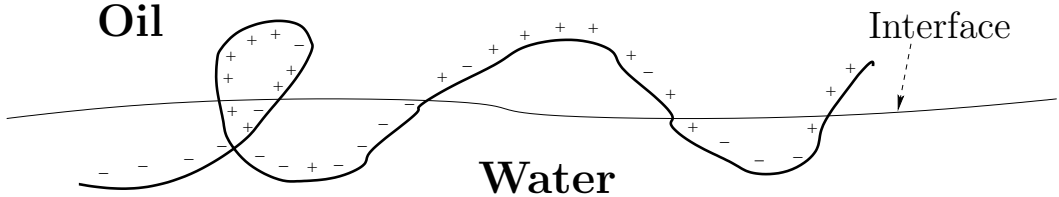,width=14 cm}}
\bigskip
\caption{ \label{fig:interface}
A copolymer in the proximity of an interface between two selective solvents.}
\medskip
\end{figure}

The polymer set--up to which we mainly dedicate our attention is the problem of a
{\sl copolymer in the proximity of an interface between two selective solvents},
say oil and water. ``Copolymer'' is simply a synonym of
inhomogeneous polymer, that is a polymer whose monomers are not identical but
can be of different types. In our case we suppose that the monomers can differ
in only one characteristics: they may be hydrophobic ($+$) or hydrophilic ($-$)
(see Figure~\ref{fig:interface}) and to be definite we assume
that there is a majority of hydrophobic monomers.

At first one could be led to think that the copolymer should prefer to live
in oil, because of the substantial hydrophobicity of the chain. However a second
scenario is also possible: the polymer could decide to stick very close to the
interface, in order to place each monomer (or at least a big part of them)
in the right solvent, that is~$+$ in oil and~$-$ in water.
Observe that this second strategy produces an \textsl{energetic gain},
arising from the fact that a greater fraction of monomers
is placed in the preferred solvent, but it also entails an~\textsl{entropic cost},
because the polymer has access to a much smaller portion of the configuration space
(the trajectories that stay close to the interface are much less than those
who are free to wander in oil).

It should be clear that we are facing a typical
energy/entropy competition: our aim is to build a probabilistic model of
this situation following the paradigms of Statistical Mechanics. This will
allow us a quantitative study of this competition, in order to
decide --in function of the characteristic of the polymer chain and of other
physical quantities, such as the temperature-- which one is the winning strategy, that
is the strategy followed by the polymer.

We stress from now that our interest is in describing the thermodynamic behavior of
the copolymer \textsl{at a fixed time}: nothing will be said in this thesis
about the problem of dynamical evolution.

\smallskip
\subsection{Random walk models}

Let us be a bit more specific about how to build a random walk model
for our copolymer (the precise definitions
will be given in the next section). We first take a random walk with values
in~$\R^d$ (or a sublattice of it) and we fix a large integer~$N$,
the size of the polymer.
The idea is to look at the random walk trajectories up to epoch~$N$
as describing the configurations of the polymer chain when there are no interactions.
Then we modify the law of the walk in the way prescribed by
Statistical Mechanics, that is by giving to each trajectory
an exponential weight (Boltzmann factor) which takes into account the
interaction of the copolymer with the solvents. This new law is
the \textsl{copolymer measure}, which describes the statistical behavior
of the copolymer in thermodynamic equilibrium.

A basic issue is how to choose the random walk. To avoid trivialities
we assume that the space in which the random walk lives
is at least two--dimensional, that is~$d\ge 2$. Moreover, since real polymers
do occupy a physical space, one would rather like to deal with {\sl self--avoiding walks}.
In the lattice case, by this we mean a random walk which is
conditioned not to visit again the sites it has already visited (defining self--avoiding
walks in the continuum case requires some more care,
but we don't want to get in details at this point).

However the point is that self--avoiding walks are a very difficult
object to deal with. One possibility to bypass the problem is to
impose a much simpler excluded volume constraint, by working with
\textsl{directed walks}. By this we mean walks in which one of the
coordinates is forced to be strictly increasing: a typical example
is the case of~$(1+m)$--dimensional directed walks, that
is~$\{(n,S_n)\}_n$ where~$\{S_n\}_n$ is a random walk in~$\R^{m}$.
Although this may appear a too drastic solution, it has been widely
used in the literature and it is the one that we will adopt too:
more precisely, we will work with an~$(1+1)$--dimensional directed
walk.

Of course another
possibility could be to give up any excluded volume constraint and to
work with genuine random walks.
However we stress that, for the model we consider,
working with a~$d$--dimensional random
walk~$\{S_n\}_n$ is equivalent to work with a suitable~$(1+1)$--dimensional
directed walk (this point will be clarified in Chapter~\ref{ch:continuous}).
This consideration gives somehow more value to the directed walk approach.


\smallskip
\section{Copolymers at selective interfaces}
\label{sec:main_model}

We are going to define a random walk model for the copolymer near a selective
interface, that will be the main object of this work.

\smallskip
\subsection{Definition of the model}

Let $S=\{S_n \}_{n=0,1,\ldots}$ be a simple symmetric random walk on~$\Z$, that is
\begin{equation*}
    S_0=0 \qquad \qquad S_n=\sum_{j=1}^n X_j\,,
\end{equation*}
where $\{ X_j\}_j$ is a sequence
of IID random variables with $\bP \left( X_1=  1\right)=\bP \left( X_1=  -1\right)=1/2$.
We take the directed walk point of view, looking at the trajectories of~$\{(n,S_n)\}_n$
as the configurations of our polymer chain.

\begin{figure}[t]
\medskip
\begin{center}
\epsfxsize =13 cm
\psfragscanon
\psfrag{0}[c][l]{$0$}
\psfrag{n}[c][l]{ $n$}
\psfrag{S}[c][l]{$S_n$}
\epsfbox{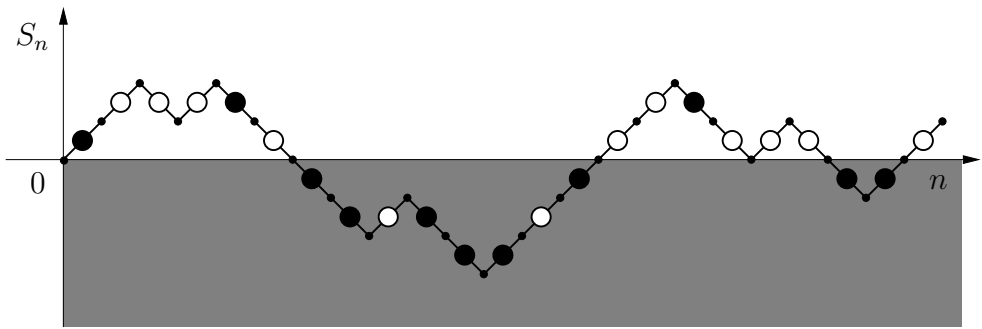}
\end{center}
\caption{\label{fig:cop_bonds}
The process we have introduced
is a model for a non--homogeneous polymer, or {\sl copolymer},
near an interface, the horizontal axis, between two selective solvents, say oil (white) and water (grey).
 In the drawing the monomer {\sl junctions} are the small black rounds and the monomers are the bonds of the random walk. The big round in the middle of each monomer
gives the sign of the charge  (white = positive charge = hydrophobic monomer, black = negative charge = hydrophilic monomer). When $h > 0$ water is  the unfavorable solvent and the question
is whether the polymer is {\sl delocalized} in oil or if it is still more profitable
to place a large number of monomers in the preferred solvent, leading in such a way
to the {\sl localization at the interface} phenomenon.
The conventional  choice of $\sign(0)$ we have made reflects the fact that
the charge is assigned to bonds rather than points.  }
\end{figure}

For $\gl\ge 0$,  $h \ge 0$, $N\in 2\N$  and
$\go =\{ \go_j\}_{j=1,2, \ldots} \in \R ^\N$
we introduce the \textsl{copolymer measure} $\bP _{N, \go}^{\gl,h}$ by giving
the density w.r.t.~$\bP$:
\begin{align}
\frac {\dd \bP_{N, \go}^{\gl, h}} {\dd \bP} (S)
&\, = \,
\frac 1
{{\tilde Z}_{N,\go}^{\gl,h}} \,
\exp\Big( \cH_{N,\go}^{\gl,h} \big( S \big) \Big) \nonumber \\
&\,:=\,
\frac 1
{{\tilde Z}_{N,\go}^{\gl,h}} \,
{\exp\left( \gl \sum_{n=1}^N \left( \go_n +h\right) \sign \left(S_n\right)
\right)} \;,
\label{eq:Boltzmann}
\end{align}
where $\sign \left(S_{2n}\right)$ is set to be equal to $\sign \left(S_{2n-1}\right)$
for any $n$ such that $S_{2n}=0$ (this is a natural choice, as it is explained in
the caption of Fig.~\ref{fig:cop_bonds}). The term ${\tilde Z}_{N,\go}^{\gl, h}$
is simply a normalization constant to make
$\bP_{N,\go}^{\gl,h}$ a probability measure, that is
\begin{equation*}
    {\tilde Z}_{N,\go}^{\gl, h} \;=\; \bE \left[ \exp\left( \gl \sum_{n=1}^N \left( \go_n +h\right) \sign \left(S_n\right)
\right) \right]\,,
\end{equation*}
and it is called the {\sl partition function} of the model.

We refer to the caption of Fig.~\ref{fig:cop_bonds} for a visual
interpretation of the copolymer measure~\eqref{eq:Boltzmann}. The expression
in the exponential, in the r.h.s. of~\eqref{eq:Boltzmann} is called
the \textsl{Hamiltonian} (we stress that there is a minus sign of difference
with respect to the physicists' conventions). Let us discuss
the meaning of the parameters appearing in the definition of~$\bP_{N,\go}^{\gl,h}$:
\begin{itemize}
\smallskip
\item $N$ is of course the size of the copolymer;
\smallskip
\item the parameter $\gl$ tunes the overall strength of the interaction,
and physically it corresponds to the inverse of the
\textsl{temperature};
\smallskip
\item $\sign (\go_n + h)$ tells whether the $n$--th monomer
is hydrophobic ($+$) or hydrophilic ($-$), and $|\go_n + h|$ gives the intensity
of the hydrophobicity (or hydrophily) of the monomer. The reason
for writing~$(\go_n + h)$ is to isolate in the parameter~$h$
the overall asymmetry of the hydrophobicity/hydrophily, as it will be clear
in a moment.
\smallskip
\end{itemize}
It remains to specify how to choose the sequence~$\go$, that will be referred to
as the {\sl charges} or the {\sl environment} of our copolymer. Two possibilities
will be considered in this thesis:
\begin{itemize}
\smallskip
\item {\sl periodic set--up:} $\go$ is a fixed periodic sequence, that is for some
$T\in\N$ we have~$\go_{2T + n} = \go_n$ for all~$n\in\N$: the least such~$T$ will be
denoted by~$T_\go$ and will be called the half--period of the sequence~$\go$
(the choice of an even period is due to the
periodicity of the simple random walk). Up to a redefinition of the parameter~$h$,
we can (and will) assume that the sequence is centered, namely~$\sum_{n=1}^{2T} \go_n = 0$.
Moreover to avoid trivialities
we suppose that $\go_{2n-1} + \go_{2n} \neq 0$ for some~$n$ (remember the
periodicity of the walk);
\smallskip
\item {\sl random set--up}:
$\go$ is a typical realization of an IID sequence of
random variables, whose law is denoted by~$\bbP$. We suppose that
\begin{equation}
\label{eq:hypM}
\M (\ga) := \bbE\left[ \exp\left( \ga \go_1 \right)\right]<\infty \qquad \forall \ga \in\R \, ,
\end{equation}
that $\bbE\left[ \go_1 \right]=0$ (which simply amounts to
redefine~$h$) and we also fix $\bbE[ {\go_1}^2]=1$. We stress that we are dealing with
{\sl quenched randomness}, that is the sequence~$\go$ is chosen at the beginning,
according to~$\bbP$, and then is kept fixed to define the copolymer
measure~$\bP_{N,\go}^{\gl,h}$.
\smallskip
\end{itemize}
The differences of the two set-up will be discussed in detail in the sequel.

\smallskip
\subsection{The free energy approach}
\label{sec:free_approach}

We suppose that the sequence of charges~$\go$, periodic or random,
has been fixed and we turn to the study of the copolymer measure~$\bP_{N,\go}^{\gl,h}$
when the size~$N$ of the copolymer is {\sl very large} (that is we are interested in
asymptotic results as~$N\to\infty$, the so--called \textsl{thermodynamic limit}).
More precisely, we would like to understand, in function of the parameters~$\gl \ge 0$
and~$h\ge 0$, whether the typical trajectories of the copolymer stay close to the interface
({\sl localized} regime) or if they rather prefer to wander away in the
solvents ({\sl delocalized} regime).

\smallskip

In order to have a quantitative criterion to decide between the two situations, it is
convenient to introduce the specific \textsl{free energy} of the system, defined by
\begin{equation}
\label{eq:free_energy_intro}
f_\go(\gl ,h)\, = \, \limtwo{N\to \infty}{N\in 2\N} \; \frac 1N \log {\tilde Z}^{\gl, h}_{N,\go} \,.
\end{equation}
Let us be more precise:
\begin{itemize}
\smallskip
\item when the sequence~$\go$ is periodic, the existence of such a limit follows
by standard superadditive arguments, see e.g.~\cite{cf:G};
\smallskip
\item in the random setting, the existence of the above limit in the
$\bbP\left( \dd \go \right)$--almost sure sense and in $\bbL_1\left( \bbP\right)$
follows by Kingman's Superadditive Ergodic Theorem, see~\cite{cf:G}.
Moreover we stress that in this case the limit does not depend on $\go$, a phenomenon
called {\sl self--averaging}. Therefore in the sequel when treating the random case
the $\go$--dependence of the free energy will be omitted.
\end{itemize}
In both the periodic and the random setting, by convexity arguments
one easily sees that the free energy is a continuous function of
$\gl$ and~$h$.

\medskip

The basic observation is that
\begin{equation}
\label{eq:delocfe}
f_\go(\gl , h) \, \ge \, \gl h.
\end{equation}
In fact if we set $\gO_N^+= \{S:\, S_n>0$ for $ n=1, 2, \ldots , N\}$,
by restricting the integration over~$\gO_N^+$ (for even values of~$N$) we get
\begin{equation}
\label{eq:step_deloc}
\begin{split}
\frac 1N \log {\tilde Z}_{N,\go}^{\gl, h}  &\;\ge\;
\frac 1N \log \bE
\left[
\exp\left(
\gl \sum_{n=1}^N
\left( \go_n +h\right) \sign \left(S_n\right)
\right)
;  \gO_N^+
\right]
\\
&\;=\; \frac {\gl} {N} \sum_{n=1}^N \left( \go_n +h\right) \, + \,
\frac 1N
\log \bP \left( \gO_N^+\right)\, \stackrel{N \to \infty}{\longrightarrow}\, \gl h ,
\end{split}
\end{equation}
where in the random case the limit has to be understood in the
$\bbP(\dd \go)$--almost sure sense, having used the law of large numbers.
We have also applied the well known fact that $\bP\left( \gO_N^+\right)$ behaves like
$N^{-1/2}$ for $N$ large \cite[Ch. III]{cf:Feller}.

The steps in \eqref{eq:step_deloc} show that~$\gl h$ is the contribution
to the free energy coming from paths delocalized in oil. This consideration leads
to the following partition of the phase diagram:
\smallskip
\begin{itemize}
\item the localized region: $\cL = \left\{ (\gl , h): \, f_\go(\gl, h)>\gl h\right\}$;
\smallskip
\item the delocalized region:  $\cD = \left\{ (\gl , h): \, f_\go(\gl, h) = \gl h\right\}$.
\end{itemize}
\smallskip
This definition of (de)localization in terms of the free energy may seem
a bit indirect and it is not a priori obvious whether it corresponds to a really
(de)localized behavior of the typical paths of the polymer measure: we will come back
in \S~\ref{sec:pathwise} to this important issue.

\smallskip

Now the program is to study in detail the phase diagram, both in the periodic and
in the random setting. Notice that a priori it is not even obvious
that~$\cL \neq \emptyset$, while of course $\cD \supseteq \{(\gl,h):\, \gl=0\}$.
We start with a basic result, valid in both settings, which says that
indeed~$\cL \neq \emptyset$ and gives the existence of a {\sl critical line},
which will be a central object of our analysis.

\smallskip
\begin{proposition}
\label{prop:prel}
Both in the periodic and in the random setting, there exists a continuous {\sl increasing}
function $h_c: [0,\infty) \to [0,\infty)$ with~$h_c(0)=0$ such that
\begin{equation*}
    \cD \,=\, \big\{(\gl,h):\: h \ge h_c(\gl) \big\} \qquad \cL \,=\,
    \big\{(\gl,h):\: h < h_c(\gl)\big\} \,.
\end{equation*}
\end{proposition}
\smallskip

\noindent
In particular we have the interesting result that for~$h=0$ and~$\gl > 0$ the copolymer
is localized (a fact that was first proven by Sinai in~\cite{cf:Sinai}).

About the proof of Proposition~\ref{prop:prel}, we point out that just by simple convexity
arguments one can prove the existence of the critical line~$h_c(\cdot)$ and the fact that
for~$\gl > 0$ it can be written as $h_c(\gl) = U(\gl)/\gl$, with~$U(\cdot)$ a convex function
such that~$U(0)=0$, cf.~\cite[\S~1.2]{cf:BG2}.
From this representation some elementary properties of the critical line
follow easily, like for instance the fact that there exists~$\ell \in (0,\infty]$
such that~$h_c(\cdot)$ is continuous and nondecreasing in~$(0,\ell)$
while~$h_c(\gl) = \infty$ for~$\gl > \ell$. It remains to prove that
$\ell = \infty$ and that~$\gl \mapsto h_c(\gl)$ is actually increasing
and continuous also at~$\gl=0$: the easiest way to get these results is to
combine convexity arguments with the bounds on~$h_c(\cdot)$
described in~\S~\ref{sec:per_diagram} and~\S~\ref{sec:ph_random}.

\smallskip

In the following sections we are going to study the
properties of the critical line~$h_c(\cdot)$, and we will see that
a closer look shows important differences between the periodic setting
and the random one. However before proceeding
it is convenient to make some preliminary transformations on our model.

\smallskip
\subsection{A new partition function}
\label{rem:Z}

The content of this section is valid both for the periodic and for the random
setting. From~\eqref{eq:delocfe} it is natural to introduce the {\sl excess free energy}~$\tf_\go$
defined by
\begin{equation*}
    \tf_\go(\gl,h) \;:=\; f_\go(\gl,h) - \gl h\,,
\end{equation*}
so that the condition for localization (resp. delocalization)
becomes $\tf_\go(\gl,h)>0$ (resp. $\tf_\go(\gl,h)=0$). It is clear
that we can obtain~$\tf_\go$ as the free energy of our copolymer,
once we redefine the Hamiltonian $\cH_{N,\go}^{\gl,h} \to
\cH_{N,\go}^{\gl,h} - \gl h N$ (observe that adding to the
Hamiltonian a term that does not depend on~$S$ has no influence on
the copolymer measure). However it is more convenient to redefine
the Hamiltonian in a slightly different way, by subtracting the
term~$\gl \sum_{n=1}^N(\go_n + h)$ instead of just~$\gl h N$. As
this term does not depend on~$S$ too, we can write
\begin{align}
\frac {\dd \bP_{N, \go}^{\gl, h}} {\dd \bP} (S)
&\, = \,
\frac 1
{{Z}_{N,\go}^{\gl,h}} \;
{ \exp\left( \gl \sum_{n=1}^N \left( \go_n +h\right) \big( \sign(S_n) - 1 \big)
\right)} \nonumber \\
&\, = \,
\frac 1
{{Z}_{N,\go}^{\gl,h}} \;
{ \exp\left( -2 \gl \sum_{n=1}^N \left( \go_n +h\right) \Delta_n
\right)}
\label{eq:Boltzmann1}
\end{align}
with $\Delta_n = \left(1-\sign (S_n)\right)/2 = \ind_{\{\sign(S_n) = -1\}}$ and
with a new partition function $Z_{N,\go}^{\gl, h}$ given by
\begin{align}
    \label{eq:Znew}
    Z_{N, \go}^{\gl, h} & \;:=\; \bE\left[ \exp\left( -2 \gl \sum_{n=1}^N
    \left( \go_n +h\right) \Delta_n \right) \right] \\
    & \;=\; \tilde Z_{N, \go}^{\gl, h} \exp \bigg( - \gl \sum_{n=1}^N
    (\go_n +h) \bigg)\,. \nonumber
\end{align}
Hence from~\eqref{eq:free_energy_intro} we get
\begin{equation}
\label{eq:felim2}
\limtwo{N \to \infty}{N\in 2\N}
\frac 1N \log Z_{N, \go}^{\gl, h}= f(\gl, h)-\gl h = \tf_\go (\gl, h) \,,
\end{equation}
where in the random case this limit has to be interpreted in the $\bbP (\dd \go)$--a.s.
or in the~$\bbL_1(\bbP)$ sense.

We will see that the new partition function $Z_{N , \go}$ turns out to be
substantially more useful than~$\tilde Z_{N, \go}$ (this fact had been
already realized in \cite{cf:BdH}).
For this reason, in the following with partition function we will always
mean~$Z_{N, \go}$, and in the same way $\tf_\go (\gl, h)$
will be for us the free energy \textsl{tout court}.

We will use repeatedly also the partition function associated to the
model {\sl pinned} at the right endpoint:
\begin{equation}
\label{eq:pinned}
Z_{N, \go}^{\gl, h}(x) \, :=\,
\bE\left[
 \exp\left(
 -2 \gl \sum_{n=1}^N \left( \go_n +h\right) \Delta_n
\right); \, S_N= x
\right].
\end{equation}
It is worth recalling that one can substitute $Z_{N, \go}^{\gl, h}$ with $Z_{N, \go}^{\gl, h}(x)$,
any fixed even $x$ (with the same parity of~$N$), in
\eqref{eq:felim2} and the limit is unchanged, see e.g.~\cite{cf:BdH} or~\cite{cf:G}.

\smallskip
\subsection{The phase diagram in the periodic case}
\label{sec:per_diagram}

As a matter of fact, the periodic case is essentially simpler than the random case. The reason
is that by expressing the partition function in terms of the random walk excursions, the
problem can be reduced to a finite--dimensional setting, as it has been first point out in~\cite{cf:BG}
(this approach will be exploited in detail in Chapter~\ref{ch:cgz}). The net result is that
the free energy of the model is expressible as the solution of a finite--dimensional
Perron--Frobenius problem, from which sharp estimates for the critical line can be obtained.

To express the results, we introduce the Abelian group $\bbS := \Z/(T_\go \Z)$ (we recall that
$2T_\go$ is the period of the sequence~$\go$), and we define the following $\bbS \times \bbS$ matrix
\begin{equation*}
    \Xi_{\ga,\gb} \;:=\; \sum_{i=2a + 1}^{2b} \go_i\,,
\end{equation*}
which is well defined by choosing representatives~$a \in \ga$ and~$b\in\gb$ with~$a<b$. We also
introduce for~$x\in\N$, $\ga,\gb \in\bbS$ and~$\gl,h\ge 0$
\begin{equation*}
    \Phi_{\ga,\gb}^{\gl, h} (x) \;:=\; \log \bigg( \frac{1 + 2\exp \big( -2 ( \gl \Xi_{\ga,\gb}
    + \gl h x ) \big)}{2} \bigg) \,.
\end{equation*}
Then, denoting by~$K(x) := \bP(\tau_1 = 2x)$ where $\tau_1:=\inf \{n>0:\, S_n=0\}$ is the first
return time to zero of the random walk, we define for~$b\ge 0$ the following $\bbS \times \bbS$
matrix with nonnegative entries:
\begin{equation*}
    A_{\ga,\gb}(b;\gl,h) \;:=\; \sum_{x \in \N} \exp \big( \Phi_{\ga,\gb}^{\gl,h} (2x) - b (2x)
    \big) \, K(2x) \, \ind_{(x \in \gb - \ga)}\,,
\end{equation*}
and we denote by~${\tt Z}(b; \gl, h)$ its Perron--Frobenius eigenvalue, cf.~\cite{cf:Asm}.
Observe that~$\tt Z$ is a decreasing function of~$b$ and~$h$, since~$A_{\ga,\gb}$ are so for
all~$\ga,\gb$. Then the free energy of the model is given by the following theorem
(cf.~\cite[Th.~1.2]{cf:BG}):
\smallskip
\begin{theorem} \label{th:per1}
For $\gl, h \ge 0$ we denote by~$b=\tilde b (\gl, h)$ the unique solution of the implicit
equation ${\tt Z}(b;\gl,h)=1$, if such a solution exists,
and we set~$\tilde b(\gl, h)= 0$ otherwise. Then~$\tilde b(\gl, h)$
is exactly the free energy of the model:
\begin{equation*}
    \tf_\go (\gl,h) \;=\; \tilde b (\gl,h)\,.
\end{equation*}
\end{theorem}
\smallskip
\noindent
It follows in particular that the critical line~$h=h_c(\gl)$ is determined by the implicit equation
${\tt Z}(0, \gl, h_c(\gl)) = 1$, and from this relation one can extract the asymptotic behavior
of~$h_c(\gl)$ both for $\gl \to 0$ and for~$\gl \to \infty$ (cf.~\cite[Th.~1.3]{cf:BG}):
\smallskip
\begin{theorem} \label{th:per2}
There exist two positive constants $m_\go>0$, $M_\go > 0$ such that:
\begin{align*}
    \text{as } \gl \to 0 \qquad h_c(\gl) &\,=\, m_\go \gl^3 \big(1 + o(1) \big) \\
    \text{as } \gl \to \infty \qquad h_c(\gl) &\,=\, \max_{n=1,\ldots , T} \Big( - \frac{\go_{2n-1} + \go_{2n}}{2} \Big) - \frac{\big( M_\go + o(1) \big)}{\gl}\,.
\end{align*}
\end{theorem}
\smallskip

These results give a satisfactory characterization of the phase diagram of the
copolymer in the periodic setting. We point out that the original proof of
Theorem~\ref{th:per1}, cf.~\cite{cf:CG}, is based on Large Deviations techniques.
We do not report it here because in Chapter~\ref{ch:cgz} we present an approach based on
Renewal Theory that allows a much more detailed analysis
for a wide class of periodic inhomogeneous
polymer models, including the periodic copolymer near a selective interface, and
Theorem~\ref{th:per1} will come as a byproduct of our main results (cf. Theorem~\ref{th:as_Z}
in Chapter~\ref{ch:cgz}).

\smallskip
\subsection{The phase diagram in the random case}
\label{sec:ph_random}

\begin{figure}[t]
\medskip
\psfragscanon
\psfrag{hu}[r][c]{\small $\overline{h}(\gl) $}
\psfrag{hc}[c][c]{$h_c(\gl)$}
\psfrag{hd}[c][c]{$\underline{h}( \gl)$}
\psfrag{minset}[l][c]{$h_c^\prime (0)\in[2/3,1]$}
\psfrag{l}[r][c]{$\gl$}
\psfrag{h}[r][c]{$h$}
\psfrag{0}[c][c]{$0$}
\centerline{\psfig{file=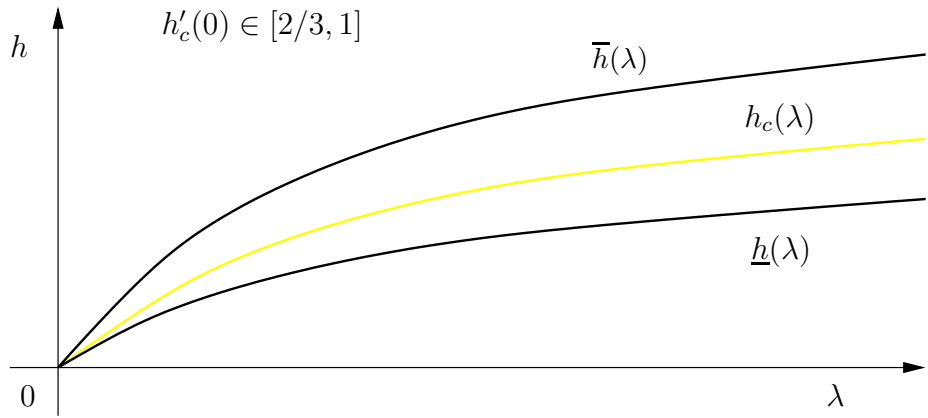,height=5.5 cm}}
\caption{ \label{fig:sumup1}
The phase diagram in the random case.\hfill\break}
\smallskip
\end{figure}

From now on, when speaking of the random case we will omit the $\go$--dependency on the free
energy, that will be simply denoted by~$\tf(\gl,h)$. We sum up in the following theorem what is
known about the critical line of the random model (see Fig.~\ref{fig:sumup1} for a graphical
representation).
\smallskip
\begin{theorem}
\label{th:sumup} For every~$\gl > 0$ the following bounds hold true:
\begin{equation}
\label{eq:sumupq}
\underline h (\gl) \, :=\,
\frac 1{4\gl /3} \log \emph\M \left( -4\gl /3\right) \;\le\;
h_c (\gl) \;\le\; \frac 1{2\gl } \log \emph\M \left( -2\gl \right)
\, =: \, \overline{h} (\gl).
\end{equation}
In particular the slope at the origin of~$h_c(\cdot)$ belongs to $[{2}/{3},1]$,
in the sense  that the inferior limit of $h_c(\gl)/\gl$ as $\gl \searrow 0$
is not smaller than $2/3$ and the superior limit is not larger than~$1$.
\end{theorem}
\smallskip

We recall that $\M(\cdot)$ it the moment generating function of~$\go_1$, see~\eqref{eq:hypM}, and
we observe that the last statement in the theorem follows easily from~\eqref{eq:sumupq}
applying the asymptotic expansion
$\M(\ga) = 1 + \ga^2/2 + O(\ga^3)$ as~$\ga \to 0$ (remember that we have fixed $\bbE[\go_1^2] = 1$).

Notice that the main difference with the periodic case, cf.~Theorem~\ref{th:per2},
is given by the behavior of~$h_c(\gl)$ as~$\gl \to 0$: we could say that
when $\gl$ is small for the copolymer it is easier to localize in the
random case than in the periodic one. This is easily understood by considering that
for~$\gl$ small a major role is played by the {\sl long excursions} of the walk, and
observing that the energetic contribution to an excursion of length~$L$ is~$O(\sqrt{L})$
in the random case by the~CLT, while in the periodic case it is of course~$O(1)$.

We point out that in~\cite{cf:GT}, using Concentration Inequalities techniques,
it has been proven that the limit of $h_c(\gl)/\gl$
as $\gl \to  0$ actually exists and it is {\sl independent of the distribution of $\go_1$},
at least when $\go_1$ is a bounded symmetric random variable or when $\go_1 $
is a standard Gaussian variable. Moreover we will see that the slope at the origin
is also closely related to the phase diagram of a Brownian copolymer model
which emerges as a scaling limit of our copolymer model as~$\gl, h \to 0$,
see~\S~\ref{sec:Brownian} below.
This {\sl universal character} of the slope at the origin
makes this quantity very interesting.

\smallskip

Theorem~\ref{th:sumup} is a mild generalization of the
results proven in \cite{cf:BdH} and \cite{cf:BG2}: the extension
lies in the fact that $\go_1$ is not necessarily symmetric and it requires
minimal changes. Despite of the fact that the {\sl lower bound} $\underline h(\cdot)$
and the {\sl upper bound} $\overline h (\cdot)$ differ only by a scale factor,
their origin is actually quite different, as we are going to see.
We also point out that in the physical literature both the conjectures that
$h_c(\cdot) = \underline h(\cdot)$~\cite{cf:Monthus,cf:SSE} and
that $h_c(\cdot) = \overline h(\cdot)$~\cite{cf:GHLO,cf:TM} have been set forth.

\smallskip
\subsubsection{The upper bound}
\label{sec:ub-ann}

For completeness we report the proof given in~\cite{cf:BdH}
of the upper bound $h_c(\cdot) \leq \overline h(\cdot)$. As a matter of fact, it
is completely elementary: using the fact that the limit in~\eqref{eq:felim2}
holds also in~$\bbL_1(\bbP)$ and applying Jensen's inequality we can write
\begin{equation*}
    \tf (\gl,h) \;=\; \lim_{N\to\infty} \frac 1N \bbE \log Z_{N,\go}^{\gl, h} \;\le\;
    \limsup_{N\to\infty} \frac 1N \log \bbE Z_{N,\go}^{\gl, h}\,,
\end{equation*}
and from~\eqref{eq:Znew} we have
\begin{equation} \label{eq:annealed_first}
    \bbE Z_{N,\go}^{\gl, h} = \bE\left[ \exp\left( \sum_{n=1}^N
    \big( \log \M(-2\gl) -2\gl h \,\big)\, \Delta_n \right) \right]\,.
\end{equation}
Then the upper bound follows immediately, because for~$h \ge \overline h(\gl)$
the argument of the exponential is nonpositive and hence~$\bbE Z_{N,\go}^{\gl, h}\le 1$.
Moreover if~$h < \overline h(\gl)$ it is easy to check that
$\lim_{N\to\infty}\big( \log \bbE Z_{N,\go}^{\gl, h} \big) / N > 0$,
hence~$\overline h(\gl)$ is indeed the best upper bound one can derive
from~\eqref{eq:annealed_first}.

This approach to get an upper bound
by performing the integration~$\bbE$ over the disorder {\sl before} taking the logarithm
is a standard tool in the Statistical Mechanics of disordered systems and it is known
as {\sl annealed bound}. We stress however that in our case this approach is not as trivial as it
may appear: for instance it is easy to check that by making the same steps with the old
partition function~$\tilde Z_{N,\go}^{\gl, h}$ one would end up with an useless bound.
The reason is that $Z_{N,\go}^{\gl, h}$ has been obtained by adding to the Hamiltonian the term
$-\gl \sum_{n=1}^N (\go_n + h)$, that does not depend on~$S$ (and therefore
it leaves the copolymer measure invariant) but that has a strong
dependence on~$\go$, which is able to change in a drastic way the annealed bound.

At this point it is clear that one could go further, searching for other $\go$--depending
terms to add to the Hamiltonian in order to improve the upper bound.
Unfortunately the standard application of this technique, known as {\sl constrained annealing},
to our copolymer model cannot improve the basic annealed bound~$\overline h(\cdot)$ on the
critical line: this point is the object of Chapter~\ref{ch:cg}, where this technique is
explained in more detail.

Of course the difficulties in improving the upper bound~$\overline h(\cdot)$
could be due to the fact that~$\overline h(\cdot)$ is indeed the true critical line. However,
the numerical analysis performed in Chapter~\ref{ch:cgg} suggests that this is not the case.

\smallskip
\subsubsection{The lower bound}

The proof given in~\cite{cf:BG2} of the lower bound
$h_c(\cdot) \ge \underline h(\cdot)$ is obtained by computing explicitly the energy--entropy
contribution to the partition function given by a suitable {\sl strategy} of the
copolymer. Roughly speaking, the strategy chosen is to force the copolymer to spend most of
his time in the upper half plane, making it descend in the lower half plane only in
correspondence of {\sl long stretches} of the sequence of charges~$\go=\{\go_n\}_n$
that have an {\sl atypically negative sample mean}. The statistics of such stretches is governed
by the so--called Large Deviations functional \cite{cf:DZ} for sums of IID random variables
distributed like~$\go_1$, which is nothing but the Legendre transform of~$\log \M(\cdot)$: this is
the reason why also the lower bound~$\underline h(\cdot)$ is of this form.

We do not report here the details of the proof because we will give
an alternative proof of the lower bound in Chapter~\ref{ch:cgg}: see Section~\ref{sec:lb}
for an outline and \S~\ref{app:prooflb} for the details.
We stress that the idea behind the above strategy
(and also behind our proof) takes inspiration from a (non rigorous)
{\sl renormalization scheme} for one--dimensional disordered systems
applied to the copolymer model by C.~Monthus~\cite{cf:Monthus}.
This approach was first proposed by D.~S.~Fisher in the context of quantum Ising model
with transverse random magnetic field~\cite{cf:Fisher} and then applied to random walk in
random environment~\cite{cf:LDMF99} with remarkable success.

We point out that the lower bound~$\underline h(\cdot)$ on the critical line
appears to be a very robust one: several attempts have been performed to enrich
the above strategy (that is to keep many more random walk trajectories)
in order to get a better lower bound, but all of them have failed.
There could be of course the possibility that~$h_c(\cdot) = \underline h(\cdot)$,
but we anticipate that in Chapter~\ref{ch:cgg} we present several numerical
observations and a rigorous statistical test which strongly suggest that
indeed~$h_c(\cdot) > \underline h(\cdot)$.

\smallskip
\subsection{The path behavior}

\label{sec:pathwise}
The question of whether splitting the phase diagram
into the regions $\cL$ and $\cD$, which are defined in terms of the free energy,
does correspond to really different behaviors of the typical paths of the
copolymer measure has a positive answer,
at least if we do not consider the critical case, that is
if we consider the path behavior for $(\gl, h) \in\cL$ and
for $(\gl, h)$ in the interior of $\cD$ (that will be called
{\sl strictly delocalized region}).
However, while the localized regime is rather well understood,
the delocalized one remains somewhat elusive. We first consider
the periodic setting.

\smallskip
\subsubsection{The periodic case}

Strong path localization statement can be obtained applying the technique
used in~\cite{cf:Sinai} by Sinai to study the random case. More precisely,
if $(\gl,h) \in \cL$ then for every~$\gep > 0$ there exist positive
constants $N_0>0$, $L_0>0$ such that
for all~$N\ge N_0$
\begin{equation} \label{eq:per_loc}
    \sup_{n=1,\ldots, N}\bP_{N,\go}^{\gl, h} \big( |S_n| > L \big)
    \;\le\;  \exp \big( - (\tf_\go(\gl, h) - \gep) \,L\, \big)
    \qquad \forall L \ge L_0\,.
\end{equation}
Furthermore, the Large Deviations approach taken in~\cite{cf:BG} gives
detailed information on the returns to zero under the copolymer measure,
that form a set with positive density, see~\cite[\S~1.7]{cf:BG}.

On the other hand, for the delocalized phase the available results are less precise:
the only result known in complete generality is that in the strictly delocalized
regime the polymer spends almost all the time above any prefixed level, that is for any $L>0$
\begin{equation} \label{eq:zero_density}
    \lim_{N\to\infty} \bE_{N,\go}^{\gl,h} \Bigg[ \frac 1N \sum_{n=1}^N
    \ind_{(S_n \ge L)} \Bigg] = 1\,.
\end{equation}
Much stronger results are known to hold in more specific instances:
for example in~\cite{cf:MGO} the case of $\go_n = (-1)^n$ is considered, for
a copolymer model which differs from the ours in the definition
of~$\sign(0)$ (we refer to~\cite[\S~1.5]{cf:BG} for more details on the
implications of this change).
The authors compute the law of the returns to zero under the polymer measure,
from which using the ideas in~\cite{cf:IY} or the general and more robust
approach we take in Chapter~\ref{ch:cgz} one can extract
the Brownian scaling limits of the model. More precisely, one can prove that
for~$(\gl,h)$ in the interior of~$\cD$
the law of the process $\{S_{\lfloor tN \rfloor}/\sqrt{N}\}_{t \in [0,1]}$
under the polymer measure~$\bP_{N,\go}^{\gl,h}$
converges weakly to the law of the Brownian meander process (that is the law of
a standard Brownian motion conditioned not to enter the lower half plane,
cf.~\cite{cf:RevYor}). The analysis can be performed in the critical case
too, that is when $h=h_c(\gl)$, showing that this
time the scaling limit process is the absolute value of a Brownian motion.

We point out that the proof of these results has been obtained
essentially by explicit computations, because by taking the
marginals of a period--2 copolymer over the even sites one gets to
an homogeneous pinning model, which is known to be exactly
solvable~\cite{cf:IY} (see also \S~\ref{sec:pinning} below). However
it is widely believed that these results should hold for any
periodic~$\go$.

In Chapter~\ref{ch:cgz} we are going to show that this is indeed the case,
proving that both the strictly delocalized and the critical Brownian scaling
limit holds in complete generality for a wide class of periodic inhomogeneous polymer
models, including the copolymer near a selective interface (and also pinning/wetting
models, that will be described in~\S~\ref{sec:pinning} below).
Furthermore we will also give a precise description of the {\sl local path
properties} of the copolymer measure (thermodynamic limit) in all regimes,
including the localized case.

\smallskip
\subsubsection{The random case}

Also in the random case it is known that for $(\gl, h) \in\cL$ the copolymer
paths are localized in a strong sense. The random analogue
of~\eqref{eq:per_loc} has been proved by Sinai in~\cite{cf:Sinai}
for the case $\gl>0, h=0$ (but the
method can be extended to the whole localized region, cf.~\cite{cf:G}),
and it requires some care to state it properly.
It is convenient to work with two--sided sequence of charges, that is
for the sake of this section we assume that~$\go=\{\go_n\}_{n\in\Z}$
is an element of the space~$\Omega := \R^\Z$ and~$\bbP$ is of course the product
probability measure on~$\Omega$. We also define for~$n\in\N$ the translation
$\theta^n$ on~$\Omega$ by $(\theta^n \go)_k := \go_{n+k}$.
Then Sinai's result reads as follows: for every~$\gep > 0$ there exist
positive random variables $N_0(\go),\, L_0(\go): \Omega \to \N$ such that
for~$\bbP$--almost every~$\go$
and for all $N\ge N_0(\go)$ the following relation holds:
\begin{equation} \label{eq:ran_loc}
\begin{split}
    & \forall n \in \big\{\log^\gamma N, \ldots, N - \log^\gamma N\big\}
    \qquad \forall L \ge L_0\big(\theta^n(\go)\big) \\
    &\rule{0pt}{13pt} \quad \ \bP_{N,\go}^{\gl, h} \big( |S_n| > L \big)
    \;\le\; \exp \big( - \big(\tf(\gl, h) - \gep\big) \,L\, \big) \,,
\end{split}
\end{equation}
where $\gamma>0$ is an absolute constant (depending neither on~$\gep$ nor on~$\go$).

Some observations are in order. The restriction on the values of~$n$ is made only
for convenience: an analogous statement holding for all~$n\le N$
is possible but the notations become more involved. The key point is rather
the condition $L \ge L_0\big(\theta^n(\go)\big)$, which is saying that the
``radius of localization'' depends on~$n$, and can actually be arbitrarily large
because the the random variable~$L_0$ is essentially unbounded. This fact may make
the estimate~\eqref{eq:ran_loc} appear unsatisfactory,
but in fact it is unavoidable, and the reason for this
is to be sought in the presence of {\sl arbitrary long atypical stretches}
in the sequence~$\go$: in fact if a site~$n$ is surrounded by a stretch
$\{\go_{n-k}, \ldots, \go_{n+k}\}$ with an atypically positive sample mean,
this will have a repulsive effect pushing~$S_n$
to height of order~$\approx \sqrt{k}$
(in Chapter~\ref{ch:cgg} we will see that one can also take advantage of
atypical stretches).

However the situation is not so bad. On the one hand, the random variable
$L_0$ can be chosen such that~$\bbE \big[ \exp(\ga L_0) \big] < \infty$ for some
$\ga > 0$ and therefore for~$\bbP$--a.e.~$\go$ we have that eventually
$L_0(\theta^n \go) \le \ga^{-1} \log n$: thus the radius of localization
is in any case much smaller than the polymer size. On the other hand,
by Birkhoff's Ergodic Theorem we have that, for every~$K>0$, $\bbP(\dd\go)$--a.s.
\begin{equation*}
    \lim_{N \to \infty} \; \sum_{n=1}^N \ind_{\{L_0(\theta^n \go) > K\}}
    \; = \; \bbP \big( L_0  > K \big)\,,
\end{equation*}
hence by choosing~$K$ large we have that the localization radius is smaller
than~$K$ {\sl most of the time}.

We observe that strong localization results are available also for
the thermodynamic limit of the copolymer measure: we do not report them here and
we refer to~\cite{cf:BisdH} and~\cite{cf:AZ} for details.

Turning to the delocalized regime, we point out that almost no result is at
present available for the critical case. In the strictly delocalized case
the situation is somewhat better: for instance it is known that \eqref{eq:zero_density}
holds for~$\bbP$--almost every~$\go$.
However this is quite a weak information on the paths,
above all if compared to what is available for the periodic case (and more generally
for related non disordered models, see e.g. \cite{cf:DGZ} and references therein),
namely Brownian scaling.

The standard way to prove this scaling limit for the strictly
delocalized regime is to show that under the copolymer measure~$\bP
_{N, \go}^{\gl,h}$ the epoch of the last visit to the lower half
plane is~$o(N)$. For non disordered models actually much more is
true: in fact in the limit $N\to\infty$ the polymer becomes
transient and it visits the lower half plane, or any point below a
fixed level, only a {\sl finite number} of times. The situation
appears to be different for the random copolymer: in fact
in~\cite{cf:GT} it has been shown that for $h <\overline{h}(\gl)$
the number of visit to the lower half plane for the \textsl{quenched
averaged measure}~$\bbE \bE_{N,\go}^{\gl,h} [\,\cdot\, ]$ is~$O(\log
N)$. This fact alone does not suffice to yield the scaling limit,
because besides showing that there are $o(N)$ visit to the lower
half plane, one should prove that they all happen close to the
origin: we refer to \cite{cf:GT} for more details and for a
discussion on what is still missing.

We stress that in answering this kind of questions an important role is
played by the asymptotic behavior as~$N\to\infty$
of the partition function~$Z_{N,\go}$ in the interior of the delocalized phase.
In the non--disordered case it is known that~$Z_{N,\go} \approx N^{-1/2}$, see for instance
Theorem~\ref{th:as_Z} of Chapter~\ref{ch:cgz}. On the other hand, this asymptotic
behavior is known not to hold anymore in the random case: more precisely, for every~$(\gl,h)$
in the interior of~$\cD$ there exists~$\gep > 0$ and a subsequence~$\{\tau_N(\go)\}_{N}$
such that $N^{1/2 - \gep} Z_{\tau_N(\go),\go} \to \infty$ as $N\to\infty$, see Proposition~4.1 in~\cite{cf:GT}.

\smallskip

The issue of the delocalized path behavior in the random case
is taken up again in \S~\ref{sec:path1} of Chapter~\ref{ch:cgg}.


\smallskip
\section{Other related polymer models}

\smallskip
\subsection{Pinning at an interface and wetting models}
\label{sec:pinning}

Another problem that has received much attention is the situation in which a polymer
chain is attracted (or repelled) by an interface, which may be penetrable or impenetrable.
We can model this situation by giving a reward (or a penalization)
to each monomer lying on the interface, and this reward/penalization may
vary from one monomer to another if the polymer chain
is heterogeneous. As in the copolymer model analyzed so far, this modification may alter the
paths of the walk inducing a localization/delocalization transition.

Let us define a probabilistic model for these situations when the interface is flat
(for us it will be the~$x$--axis). We start with the case when the interface is
penetrable ({\sl pinning models}): as in the preceding section,
we take a simple random walk $\{S_n\}_n$ with law~$\bP$, and for~$N\in 2\N$, $\gb \in \R$
and~$\go=\{\go_n\}_{n\in\N} \in \R^\N$ we define the new law~$\bP_{N,\go}^\gb$ by
\begin{equation} \label{eq:def_pinning}
    \frac{\dd \bP_{N,\go}^\gb}{\dd \bP} (S) \;\propto\; \exp \Bigg( \gb
    \sum_{n=1}^N \go_n \ind_{(S_n=0)} \Bigg)\,.
\end{equation}
The case of an impenetrable interface is obtained by restricting to paths that
stay nonnegative up to epoch~$N$, that is multiplying the r.h.s. above
by~$\ind_{(S_1\ge 0, \ldots, S_N \ge 0)}$. This second case will be called a
{\sl wetting model}, as it can be also interpreted as the model of an interface
interacting with an impenetrable wall.

Again we will stick to the case when the sequence charges~$\go$ is either
(deterministic and) periodic or (quenched) random.
Of course the main interest is in understanding the behavior of the above measure when~$N$
is large. Localization/delocalization can be defined in terms of the corresponding free
energy, exactly as in the previous section.

A particularly simple case is the homogeneous one, that is when the sequence~$\go$
is constant: $\go_n = \go_1$ for all~$n\in\N$, and up to a redefinition of the
parameter~$\gl$ we may assume that~$\go_1=1$. We point out that in this case
both the pinning and the wetting models
are completely solvable, not only for the purpose of finding
the phase diagram (see e.g.~\cite{cf:G} for an elementary derivation)
but also for a very detailed analysis of the polymer
path behavior~\cite{cf:IY,cf:DGZ}. We do not spend much time here on this issue,
because in Chapter~\ref{ch:cgz} we will treat in full detail the case of periodic~$\go$.
Nevertheless we make some observations: using convexity arguments it is easy to check
that the phase diagram in the homogeneous case is encoded in a single number~$\gb_c$ such that
for~$\gb > \gb_c$ (resp. $\gb \le \gb_c$) the polymer is localized (resp. delocalized). Moreover:
\begin{itemize}
\item in the pinning case $\gb_c=0$, that is an arbitrarily small reward is sufficient to
localize the polymer;
\smallskip
\item in the wetting case on the contrary $\gb_c >0$.
\end{itemize}
The reason why the wetting model is more difficult to localize than the pinning model
is that conditioning the walk to stay nonnegative up to step~$N$ induces a repulsion
effect of order~$\sqrt{N}$ on the paths, a phenomenon which goes under the
name of {\sl entropic repulsion}. In our one dimensional setting, a more precise
version of this statement is provided by the following invariance principle~\cite{cf:Bol76}:
the process $\{S_{\lfloor Nt \rfloor}/\sqrt{N}\}_{t\in[0,1]}$ conditionally on the
event~$\{S_1\ge 0, \ldots, S_N \ge 0\}$ converges weakly as~$N\to\infty$ to the Brownian
meander process, that is to a Brownian motion conditioned to stay nonnegative~\cite{cf:RevYor}.
In Chapter~\ref{ch:llt} we will prove a local version of this weak convergence.

As already anticipated, the periodic version of these models will be analyzed
in Chapter~\ref{ch:cgz}. On the other hand, nothing will be said in this thesis about
the random case. We only mention that, as in the copolymer case, in the physical literature
there is no agreement on the phase diagram of the model, especially for small values of the
coupling constants: for more details on this issue and for the available rigorous
results see~\cite{cf:AS,cf:Pet} and references therein.

\medskip

To conclude we would like to point out the relevance that random walk models
have for the modeling of DNA molecules. DNA is normally in a double--stranded
state, however it may happen that the two strands get detached, for example
when the temperature is sufficiently high ({\sl denaturation transition}) or due to
the effect of an external force ({\sl pulling induced unzipping}). Since the interaction
between the two strands may be described (at least at a first level)
by an Hamiltonian of the
form~\eqref{eq:def_pinning}, the energy/entropy competition that gives origin
to such phase transitions may be understood in terms of suitable modifications
of the pinning/wetting models just described.

\smallskip
\subsection{A Brownian motion model: the coarse graining issue}
\label{sec:Brownian}

One of the main results in the paper of Bolthausen and Den Hollander~\cite{cf:BdH}
is that in the limit of weak coupling the copolymer model described in Section~\ref{sec:main_model}
can be approximated by a continuous model built with Brownian motions instead of random walks.
This continuous model is defined in complete analogy with the discrete one: we take two
Brownian motions $B=\{B_t\}_{t\ge 0}$ (the polymer) and~$\gb=\{\gb_t\}_{t\ge 0}$ (the charges),
with respective laws~$\tilde \bP$ and~$\tilde \bbP$, and for~$t>0$, $\gl, h \ge 0$ and a
$\tilde \bbP$--typical path~$\{\gb_s\}_s$ we introduce
the polymer measure~$\tilde \bP_{t,\gb}^{\gl, h}$ on paths of length~$t$ defined by
\begin{equation*}
    \frac{\dd \tilde \bP_{t,\gb}^{\gl, h}}{\dd \tilde \bP} (B) \;:=\; \frac{1}{Z_{t,\gb}^{\gl, h}} \;
    \exp \bigg( \gl \int_0^t \sign(B_s) \, \big( \dd \gb_s + h \dd s \big) \bigg)\,,
\end{equation*}
where the integral with respect to~$\gb_s$ is an Ito integral. The partition function of the
model is of course
\begin{equation*}
    Z_{t,\gb}^{\gl, h} \;=\; \tilde \bE \exp \bigg( \gl \int_0^t \sign(\gb_s)
    \, \big( \dd \gb_s + h \dd s \big) \bigg)\,,
\end{equation*}
and the free energy~$\tilde f(\gl,h)$ is defined as
\begin{equation*}
    \tilde f(\gl,h) \;:=\; \lim_{t\to \infty} \frac 1t \, \log Z_{t,\gb}^{\gl, h}\,,
\end{equation*}
where the limit holds both $\tilde \bbP$--a.s. and in~$\bbL_1(\tilde \bbP)$ and~$\tilde f(\gl,h)$
is nonrandom (see~\cite{cf:G} for a detailed proof of the existence of such a limit).

As in the discrete case, we have
\begin{equation*}
    \tilde f(\gl, h) \;\ge\; \gl h\,,
\end{equation*}
and consequently we distinguish between a delocalized regime ($\tilde f(\gl,h) = \gl h$)
and a localized regime ($\tilde f(\gl, h) > \gl h$). Notice however that the scaling properties
of Brownian motions entail that for all~$a>0$
\begin{equation*}
    \tilde f(\gl,h) \;= \; \frac{1}{a^2} \tilde f (a\gl, ah)\,,
\end{equation*}
from which it follows immediately that the critical curve of this model is a {\sl straight
line}, that is
\begin{equation*}
    \exists K_c > 0\,: \quad \tilde f(\gl, h)
\begin{cases}
= \gl h & \ \text{if } h \ge K_c \gl \\
> \gl h & \ \text{if } h < K_c \gl
\end{cases}
\,.
\end{equation*}
Despite the apparent simplicity of the phase diagram, we could say that this continuous model
retains the full complexity of the discrete model, which is hidden in the constant~$K_c$. This statement
is made precise by the following fundamental theorem (cf. \cite[Th.~5~and~6]{cf:BdH}), which
also clarifies in which sense the continuous model is an approximation of the discrete one.

\smallskip
\begin{theorem}
\label{th:coarse}
Let us consider the free energy $f(\gl,h)$ of the discrete model (see eq.~\eqref{eq:free_energy_intro})
in the case when $\bbP(\go_1=+1) = \bbP (\go_1 = -1) = 1/2$, and the corresponding critical line
$h=h_c(\gl)$ (see Prop.~\eqref{prop:prel}). Then the following relations hold:
\begin{gather}
    \label{eq:scaling}
    \lim_{a\to 0} \, \frac{1}{a^2} \, f(a \gl,a h) \;=\; \tilde f (\gl, h) \quad \forall \gl, h \ge 0 \\
    \label{eq:slope}
    h'_c(0) \;:=\; \lim_{\gl \to 0} \frac{h_c(\gl)}{\gl} \;=\; K_c\,.
\end{gather}
\end{theorem}
\smallskip
\noindent
In particular by~\eqref{eq:sumupq} it follows that~$2/3 \le K_c \le 1$.

\smallskip
We point out that \eqref{eq:slope} does not follow directly from~\eqref{eq:scaling}: in fact
the scaling limit of the free energy expressed by~\eqref{eq:scaling} yields only the lower bound
$h'_c(0) \ge K_c$. The proof of~\eqref{eq:slope} is achieved through sharp comparison
inequalities between~$f$ and~$\tilde f$ and requires very delicate estimates.

We stress that Theorem~\ref{th:coarse} has been proven
for the case when the charges have
a symmetric Bernoulli law, but its validity should be very general.
The intuitive idea is that
as~$\gl,h \to 0$ what really matters are the {\sl long excursions} of the walk, and consequently
the microscopic details of the model should become irrelevant.

For instance, as we already mentioned, in~\cite{cf:GT} it has been proven
with Concentration Inequalities techniques that actually
\eqref{eq:slope} holds whenever $\go_1$ is bounded and symmetric (and such that~$\bbE[\go_1]=0$
and~$\bbE[\go_1{}^2] = 1$) or if~$\go_1$ is a standard Gaussian.
Alternatively, the original proof of Theorem~\ref{th:coarse} given in~\cite{cf:BdH}
can be adapted to show that indeed both~\eqref{eq:scaling}
and~\eqref{eq:slope} hold for any choice of the law of~$\go_1$ satisfying~\eqref{eq:hypM}
and such that $\bbE[\go_1]=0$, $\bbE[\go_1{}^2] = 1$.

In Chapter~\ref{ch:continuous} we introduce another kind of variation on the discrete model,
namely we will change the law~$\bP$ of the underlying walk, taking into account general random
walks on~$\R$ whose increments are bounded and have a continuous law. This change too is supposed not to
have any influence on the conclusions of Theorem~\ref{th:coarse}, however giving a complete
proof of this fact appears to be very challenging. Some partial steps have been done in the
direction of proving~\eqref{eq:scaling} alone: we refer to Chapter~\ref{ch:continuous} for more details
on this issue.


\smallskip
\section{An overview of the literature}

The {\sl copolymer in the proximity of an interface} problem has a long literature, mostly in
the area of chemistry and physics, but possibly the first article that attracted the
attention of mathematicians is~\cite{cf:GHLO}. The first mathematical study on the subject
has been performed by Sinai in~\cite{cf:Sinai}, where he shows that for~$h=0$ and~$\gl > 0$
(we are referring to the parameter of the model introduced in Section~\ref{sec:main_model})
the copolymer with random charges is localized in a strong pathwise sense (see~\S~\ref{sec:pathwise}
above).
Further path investigations and a detailed analysis of the free energy (always for the random case
and in the symmetric setting $h=0$) have been performed by Albeverio and Zhou in~\cite{cf:AZ}.

As already mentioned, our attention on the random copolymer has been mainly focused on the
issue of investigating the phase diagram, which entails studying the copolymer for~$h>0$.
In this direction the fundamental paper is
the one by Bolthausen and den Hollander~\cite{cf:BdH}, where the existence and some basic properties
of the critical line~$h_c(\cdot)$ (including the upper bound in~\eqref{eq:sumupq}) have been proven.
However the main result of~\cite{cf:BdH} is the coarse graining of the free energy, expressed by
Theorem~\ref{th:coarse} of~\S~\ref{sec:Brownian} below. The other fundamental result on the phase diagram
in the random case, namely the lower bound in~\eqref{eq:sumupq}, has been proven by Bodineau and
Giacomin in~\cite{cf:BG2}.

The strategy used in~\cite{cf:BG2} takes inspiration from the physical
paper by Monthus~\cite{cf:Monthus}, where the lower bound curve~$\underline h (\cdot)$ has been
introduced for the first time, as a conjecture for the true critical line. Again from the physical
literature, we point out that the conjecture~$h_c(\cdot) = \underline h (\cdot)$ has been set
forth also in~\cite{cf:SSE} on the ground of replica computations, while the complementary
conjecture that~$h_c(\cdot) = \overline h(\cdot)$ has been formulated in~\cite{cf:TM}
and in~\cite{cf:GHLO}.

Coming back to mathematical papers, a path analysis for the whole localized region~$\cL$ in the random
case has been performed by Biskup and den Hollander in~\cite{cf:BisdH}: the keywords of their approach
are {\sl thermodynamic limit} and {\sl Gibbs measures}. On the other hand, path results
for the delocalized region appear to be much more challenging: recent progresses in this
direction have been obtained by Giacomin and Toninelli in~\cite{cf:GT}.

Turning to the case of periodic charges, the issue of determining the phase diagram
has received a complete solution in the paper~\cite{cf:BG} by Bolthausen and Giacomin
(see~\S~\ref{sec:per_diagram}). We refer to this
paper also for references to previous works on periodic copolymers.

In the literature one finds also a large number of numerical works on copolymers,
see for example~\cite{cf:CW,cf:SW} and references therein:
with respect to the numerical approach we take in Chapter~\ref{ch:cgg},
the attention is often shifted toward different aspects,
notably the issue of critical exponents and the more
complex model in which the polymer is not directed but rather self--avoiding.

Finally, about other polymer models dealing with pinning/adsorption phenomena we mention
\cite{cf:IY}, \cite{cf:DGZ}, \cite{cf:AS}, \cite{cf:Whi1,cf:Whi2} and references therein.


\smallskip
\section{Outline of the thesis}
\label{sec:out}

\smallskip
The exposition is organized as follows.
\begin{itemize}
\medskip
\item In Chapter~\ref{ch:cgg} we combine numerical computations with rigorous arguments
to study the phase diagram and the path behavior of the copolymer
near a selective interface model defined in
Section~\ref{sec:main_model}. We consider the case of {\sl random
charges}. We provide several evidences for the fact that the
critical line lies strictly in between the two known bounds given
in~\eqref{eq:sumupq} and for the fact that the scaling limit towards
the Brownian meander process holds in the strictly delocalized
region. In particular the conjecture that~$h_c(\cdot) = \underline
h(\cdot)$ can be excluded with a high level of confidence, thanks to
a rigorous statistical test with explicit error bounds. We also give
an alternative self--contained proof of the lower bound~$h_c(\cdot)
\ge \underline h(\cdot)$.

\medskip \noindent
The article~\cite{cf:CGG} has been taken from the content of this chapter.
\bigskip
\item In Chapter~\ref{ch:cg} we address the issue of improving the annealed upper bound
for disordered systems (see \S~\ref{sec:ub-ann})
by adding to the Hamiltonian disorder--dependent terms, a technique known as
{\sl constrained annealing}. We show that for a number of disordered linear chain models
(including the copolymer near selective interfaces and the pinning/wetting model
described in the preceding sections) the standard application
of this technique using empirical averages of local functions cannot improve the basic
annealed bound on the critical curve.

\medskip \noindent
The article~\cite{cf:CG} has been taken from the content of this chapter.
\bigskip
\item In Chapter~\ref{ch:cgz} we consider a general model of a heterogeneous polymer
in the proximity of an interface (including as special cases the copolymer near a selective
interface and the pinning/wetting model) in the case of {\sl periodic
charges}. We propose an approach based on Renewal Theory that yields sharp estimates on the
partition function of the model in all the regimes, including the critical one. From
these results we obtain a very precise description of both the thermodynamic limit
and the scaling limits of the polymer measure.

\medskip \noindent
The preprint~\cite{cf:CGZ} has been taken from the content of this chapter.
\bigskip
\item In Chapter~\ref{ch:continuous} we consider a modification of the copolymer near
a selective interface model in which the reference measure~$\bP$
is not any more the law of the simple symmetric random walk on~$\Z$. More precisely,
we allow~$\bP$ to be the law of a general real random walk whose typical step is
centered, bounded and has an absolutely continuous law. We focus on the case of
{\sl random charges}. Besides giving a proof of the existence of the free energy,
we study the phase diagram of the model, pointing out the close analogies
with the simple random walk case. We finally consider the issue of extending to
this model the coarse graining of the free energy expressed by Theorem~\ref{th:coarse}
(work in progress), giving some partial result in this direction
and discussing what is missing.
\bigskip
\item In Chapter~\ref{ch:llt} we prove a local limit theorem for random walks
conditioned to stay positive which is valid in great generality (whenever the walk
is attracted to the Gaussian law). This theorem provides a local refinement of the well-known
weak convergence of random walks conditioned to stay positive towards the Brownian
meander process. Besides being an interesting result in itself,
it is an important tool for the purpose of dealing with polymer models built over
general random walks, like the one considered in Chapter~\ref{ch:continuous}.

\medskip \noindent
The article~\cite{cf:C} has been taken from the content of this chapter.
\end{itemize}

%
%
\chapter[A numerical study of the random copolymer]
{A numerical study of the phase diagram and path behavior\\
of the copolymer model with random charges}
\label{ch:cgg}

In this chapter we study the copolymer near a selective interface
model, defined in Section~\ref{sec:main_model}
of Chapter~\ref{ch:first}, in the {\sl random} case. We combine numerical computations
with rigorous arguments to get to a better understanding of the phase diagram
and of the path behavior. Our main aim is to provide evidences of the fact that
the critical line of the model~$h_c(\cdot)$ lies \emph{strictly} in between
the two bounds~$\underline h(\cdot)$ and~$\overline h(\cdot)$, defined in
\eqref{eq:sumupq} of Chapter~\ref{ch:first}, and to numerically analyze
the delocalization issues raised in \S~\ref{sec:pathwise} of Chapter~\ref{ch:first}.
A detailed outline of the results is given in~\S~\ref{sec:outline}.

\smallskip
The article~\cite{cf:CGG} has been taken from the content of this chapter.

\smallskip
\section{Introduction and results}

\smallskip
\subsection{Prelimiaries}
\label{sec:notation}

The notations we will use are those introduced in Section~\ref{sec:main_model}
of Chapter~\ref{ch:first}. We recall
in particular the definitions~\eqref{eq:Znew} and~\eqref{eq:pinned} of the
partition functions~$Z_{N,\go}^{\gl,h}$ and~$Z_{N,\go}^{\gl,h}(x)$ (we
will be mainly interested in the case~$x=0$). Also remember that for the
critical line~$h_c(\cdot)$ of our model we have the bounds
\begin{equation} \label{eq:sumup_alt}
    \underline h (\gl) \;\le\; h_c (\gl) \;\le\; \overline{h} (\gl),
\end{equation}
see Theorem~\ref{th:sumup} of Chapter~\ref{ch:first}. For what follows we set
\begin{equation}
\label{eq:hm}
h^{(m)}(\gl) \, = \, \frac 1{2m\gl } \log \M\left( -2m\gl \right),
\end{equation}
for $m >0$, where we recall that $\M (\ga) := \bbE[\exp (\ga \go_1)]$. Observe that the curves $\underline h(\cdot)$ and $\overline h(\cdot)$
correspond respectively to $m=2/3$ and $m=1$, and that for every~$m$ we have
$\frac{\dd}{\dd\gl} h^{(m)}(\gl) |_{\gl=0} = m$.

\medskip

Before proceeding, we present a different viewpoint on the process:
this turns out to be useful for the intuition and it will be used
in some technical steps. We call  $\eta$  the first
return time of the walk $S$ to $0$, that is $\eta:=\inf\left\{  n\ge
1:S_{n}=0\right\} $, and  set $K(2n):=\bP \left(\eta=2n\right)$ for $n\in\N$.
It is well known that $K(\cdot)$ is decreasing on the even natural numbers and
\begin{equation}
\label{eq:asympt}
\lim_{x\in 2\N, x\to \infty} x^{3/2} K(x)=
\sqrt{2/\pi },
\end{equation}
see e.g. \cite[Ch. 3]{cf:Feller}.
Let the IID sequence $\left\{ \eta_j \right\}_{j=1,2, \ldots}$
denote the inter--arrival times at $0$ for $S$, and we set $\tau_k := \eta_0 + \ldots + \eta_k$.
 If we introduce also
 $\ell _N =\max\{j \in \N \cup \{0\} :
\tau_j \le N\}$, then by  exploiting
the up--down symmetry of the excursions of $S$  we directly obtain
\begin{equation}
\begin{split}
\label{eq:reducttoexc}
& Z_{N,\go} (0) \, =\, \bE \left[ \prod_{j=1}^{\ell_N} \varphi \Bigg(\gl \sum_{n=\tau_{j-1}+1}^{\tau_j} \go_n + \gl h  \eta_j \Bigg) ; \tau_{\ell_N} = N \right] \\
    & \quad =\, \sum_{l=0}^N \sumtwo{x_0, \ldots, x_l \in 2\N}{0=:x_0 < \ldots < x_l:= N} \prod_{i=1}^l \; \varphi \Bigg(\gl \sum_{n=x_{i-1}+1}^{x_i} \go_n + \gl h  (x_i - x_{i-1}) \Bigg) \; K(x_i - x_{i-1}) \,,
\end{split}
\end{equation}
with $\varphi (t) : =   \left( 1+\exp(-2t)\right) / 2$.
Of course the formula for $Z_{N,\go}$
is just  slightly different.

Formula \eqref{eq:reducttoexc} reflects the fact that what really
matters for the copolymer are the return times to the interface.

\smallskip
\subsection{Outline of the results}
\label{sec:outline}

Formula \eqref{eq:sumup_alt} leaves an important gap, that
hides the only partial understanding of the nature of this delocalization/localization
transition.
Our purpose is to go toward filling this gap: our results are both of theoretical and numerical nature.
At the same time we address the delocalization issues raised
in \S~\ref{sec:pathwise} of Chapter~\ref{ch:first}, which are intimately related with the precise
asymptotic behavior of $Z_{N,\go}$ and of  $Z_{N,\go}(0)$.
More precisely:

\smallskip
\begin{enumerate}
\item In Section~\ref{sec:testing} we present a statistical test with explicit error bounds, see Proposition~\ref{th:stat}, based on
super--additivity and concentration inequalities, to state that a point $(\gl, h)$ is localized. We apply this test to show that, with a very low level of error, the lower bound $h=\underline h(\gl)$ does not coincide with the critical line.
\item \rule{0pt}{14pt}In Section~\ref{sec:lb} we give the outline of a new proof of the lower bound
$h_c(\cdot) \ge \underline h (\cdot)$. The details of the proof are in \S~\ref{app:prooflb} and we point out
in particular Proposition~\ref{prop:stopping1}, that gives a necessary and sufficient
condition for localization. This viewpoint on the transition, derived from
\cite[\S~4]{cf:GT}, helps substantially in interpreting the {\sl irregularities} in the
behavior  of $\left\{ Z_{N, \go}\right\}_N$ as $N \nearrow \infty$.
\item \rule{0pt}{14pt}In Section~\ref{sec:path} we pick up the conjecture of Brownian scaling
in the delocalized regime both in the intent of testing it and
in trying to asses with reasonable confidence that $(\gl, h)$ is in the interior
of $\cD$. In particular, we present quantitative evidences in favor of the fact that the upper bound $h= \overline h (\gl)$ is strictly greater than the critical line. We stress that this is a very delicate issue, since delocalization, unlike
localization, does not appear to be reducible to a finite volume issue.
\item \rule{0pt}{14pt}Finally, in Section~\ref{sec:guess}, we report the results of a numerical attempt
to determine the critical curve. While this issue has to be treated with care,
mostly for the reasons raised in point 4 above, we observe
a surprising phenomenon: the critical curve appears to be  very close
to $h^{(m)}(\cdot)$ for a suitable value of $m$. By the universality
result proven in \cite{cf:GT}, building on the free energy Brownian scaling result proven
in \cite{cf:BdH}, the slope at the origin of $h_c(\cdot)$ does not depend on the law of $\go$.
Therefore if really $h^{(m)}(\cdot)= h_c (\cdot)$, since the slope at the origin
of $h^{(m)}(\cdot)$ is $m$,
$m$ is the  universal constant we are looking for.
We do not believe that the numerical evidence
allows to make a clear cut  statement, but what we observe is compatible
with such a possibility.
\end{enumerate}
\smallskip

We point out that our numerical results are based on a numerical computation
of the partition function $Z_{N,\go}$, exploiting the standard transfer--matrix
approach (this item is discussed in more details in \S~\ref{app:algo}).


\smallskip
\section{A statistical test for the localized phase}
\label{sec:testing}

\smallskip
\subsection{Checking localization at finite volume}
\label{sec:superadd}
At an intuitive level one is led to believe that, when the copolymer is localized,
it should be possible to detect it
by looking at the system before the infinite volume limit.
This intuition is due to the fact that in the localized phase
the length of each excursion is finite, therefore for $N$ {\sl much larger }
that the {\sl typical} excursion length one should already observe the
localization phenomenon in a quantitative way.
The system being disordered of course does not help,
because it is more delicate to make sense of what
{\sl typicality} means in a non translation invariant set--up.
 However the translation invariance can be recovered
 by averaging and in fact it turns out to be rather easy to give
 a precise meaning to the intuitive idea we have just mentioned.
 The key word here is super--additivity of the averaged free energy.
 \smallskip

 In fact by considering only the $S$ trajectories such that $S_{2N}=0$ and by applying
 the Markov property of $S$ one directly verifies
  that for any $N, M\in \N$
 \begin{equation}
 \label{eq:superadd}
 Z_{2N+2M, \go} (0)\, \ge \, Z_{2N, \go}(0)\, Z_{2M, \theta^{2N} \go}(0),
 \end{equation}
 $(\theta \go)_n=\go_{n+1}$,
 and therefore
 \begin{equation}
 \left\{ \, \bbE\log Z_{2N, \go} (0)
 \right\}_{N=1,2, \ldots}
  \end{equation}
  is a super--additive sequence, which immediately entails
  the existence of the limit of $ \bbE[\log Z_{2N, \go}(0)]/2N$ and the fact that
  this limit coincides with the supremum of the sequence.
  Therefore from the existence of the quenched free energy we have that
  \begin{equation}
  \tf(\gl, h) \, =\, \sup_N \frac 1{2N} \bbE\log Z_{2N, \go}(0)\,.
  \end{equation}
In a more suggestive way one may say that:
\begin{equation}
\label{eq:loc_char}
(\gl, h) \in \cL \  \  \Longleftrightarrow  \ \
\text{there exists } N\in \N \ \, \text{such that} \,  \  \bbE\log Z_{2N, \go}(0)>0\,.
\end{equation}
The price one pays for working with a disordered system
is precisely in taking the $\bbP$--expectation
and from the numerical viewpoint it is an heavy price:
even with the most positive attitude one cannot expect
to have access to  $\bbE\log Z_{2N, \go}(0)$
by direct numerical computation for $N$ above $10$.
Of course in principle small values of $N$ may suffice
(and they do in some cases, see Remark~\ref{rem:N=2}), but they do not
suffice to tackle the specific issue we are interested in.
We elaborate at length on this interesting issue in \S~\ref{sec:computer_assisted}.
\smallskip

\begin{rem}
\label{rem:N=2}
\rm
An elementary application of the localization criterion \eqref{eq:loc_char} is
obtained for $N=1$: $(\gl, h)\in \cL$ if
\begin{equation}
\label{eq:N=2}
\bbE\left[
\log\left(
\frac 12+ \frac 12 \exp\left(
-2\gl \left(\go_1+ \go_2+2h\right)
\right)
\right)
\right]>0.
\end{equation}
In the case $\bbP (\go_1=\pm 1)=1/2$ from
\eqref{eq:N=2}
we obtain that for $\gl$ sufficiently large
$h_c(\gl) > 1- c/\gl$, with $c= (1/4)\log(2\exp(4) -1)\approx1.17$.
From $\underline{h}(\cdot)$ we obtain the same type of bound,
with $c=(3/4)\log 2\approx 0.52$. This may raise some hope that
for $\gl$ large an explicit, possibly computer assisted, computation
for small values of $N$ of $\bbE \log Z_{2N, \go}(0)$ could
lead to new estimates. This is not the case, as we show in~\S~\ref{sec:computer_assisted}.
\end{rem}

\smallskip
\subsection{Testing by using concentration}

In order to decide whether  $\bbE\log Z_{2N, \go}(0)>0$ we
resort to a Montecarlo evaluation of $\bbE\log Z_{2N, \go}(0)$
that can be cast into a statistical test with explicit error bound
by means of concentration of measure ideas.
This procedure is absolutely general, but we have to choose
a set--up for the computations and we take the simplest:
 $\bbP(\go_1=+1)=\bbP(\go_1=-1)=1/2$.
 The reason for this choice is twofold:
 \begin{itemize}
 \item if $\go_1$ is a bounded random variable, a Gaussian concentration
 inequality holds and if $\go$ is symmetric and it takes only two values
 then one can improve on the explicit constant in such an inequality.
 This speeds up in a non negligible   way the computations;
 \item generating {\sl true randomness}
is out of reach, but playing head and tail is certainly
the most elementary case in such a far reaching task (the random numbers
issue is briefly discussed in \S~\ref{app:algo} too).
 \end{itemize}

\smallskip
A third reason to restrict testing to the Bernoulli case is explained at the end of the
caption of Table \ref{tbl:2}.
\smallskip

We start the testing procedure by stating the null hypothesis:
\begin{equation}
\label{eq:H0}
\text{H}0: \  \ \bbE\log Z_{2N, \go}(0)\le 0.
\end{equation}
$N$ in H0  can be chosen
arbitrarily. We stress that refusing H0 implies $\bbE\log Z_{2N, \go}(0)>0$,
which by \eqref{eq:loc_char} implies localization.

The following concentration inequality for Lipschitz functions holds
for the uniform measure on $\{-1,+1\}^N$: for every
function $G_N:\{-1,+1\}^N \to \R$
such that $\vert G_N(\go)-G_N(\go^\prime )\vert
\le C_{\text{Lip}} \sqrt(\sum_{n=1}^N (\go_n-\go^\prime_n)^2)$, where
$C_{\text{Lip}}$ a positive constant
and $G_N(\go)$ is an abuse of notation for $G_N(\go_1, \ldots, \go_N)$, one has
\begin{equation}
\label{eq:concentration}
\bbE\left[\exp\left(\ga \left(G_N(\go)- \bbE[G_N(\go)]\right)\right)\right]\, \le\,
\exp\left(\ga^2C_{\text{Lip}}^2\right),
\end{equation}
for every $\ga$.
Inequality \eqref{eq:concentration} with an extra factor $4$ at the exponent
can be extracted from the proof of Theorem 5.9, page 100 in \cite{cf:Ledoux}.
Such an inequality holds for variables taking values in $[-1,1]$:
the factor $4$ can be removed for the particular case we are considering
(see \cite[p. 110--111]{cf:Ledoux}).
In our case $G_N(\go)= \log Z_{2N,\go}(0)$. By applying the
 Cauchy--Schwarz inequality one obtains that $G_N$ is Lipschitz with
$C_{\text{Lip}}= 2\gl\sqrt{N}$.
Let us now consider an IID sequence  $\{ G^{(i)}_N(\go)\}_i$
with $G^{(1)}_N (\go)=G_N(\go)$: if H0 holds then we have that for every $n\in \N$, $u>0$
and $\ga= un/8\gl^2N$
\begin{equation}
\begin{split}
\bbP\left(
\frac 1n \sum_{i=1}^n G^{(i)}_N(\go) \ge u
\right)
\, &\le  \,
\bbE\left[
\exp\left(
\frac\ga n \left(G_N(\go)- \bbE[G_N(\go)]\right)
\right)
\right]^n
\exp\left(-\ga\left( u- \bbE[G_N(\go)]\right)\right)
\\
&\le \,
\exp\left(\frac{4 \ga^2 \gl^2N}{n}-\ga u\right)
\\
&= \, \exp\left( -\frac{u^2 n}{16 \gl^2 N}\right) .
\end{split}
\end{equation}

Let us sum up what we have obtained:

\medskip
\begin{proposition}
\label{th:stat}
Let us call $\widehat u_n$ the average of a sample of $n$ independent
realizations of $\log Z_{2N,\go}^{\gl,h}(0)$. If $\widehat u_n>0$ then we may
refuse H0, and therefore  $(\gl,h)\in \cL$, with a level
of error not larger than $\exp\left( -{\widehat u_n ^2 n}/{16 \gl^2 N}\right)$.
\end{proposition}
\medskip

\smallskip
\subsection{Numerical tests}
We report in Table \ref{tbl:1}
 the most  straightforward application
of Proposition \ref{th:stat}, obtained by a numerical computation of $\log Z_N$ for a sample of~$n$ independent environments~$\go$.
We aim
at seeing how far above $\underline{h}(\cdot)$ one can go
and still claim localization,
keeping a reasonably small probability of error.

\begin{table}[h]
\begin{center}
\begin{tabular}{|c|c|c|c|}
\hline
$\gl $ &$0.3$&$0.6$&$1$\\
\hline
$h$&  0.22 &0.41&0.58\\
\hline
$p$--value& $1.5\times 10^{-6}$&$9.5\times 10^{-3}$&$1.6 \times 10^{-5}$\\
\hline
$\underline{h}(\gl)$& 0.195 &0.363&0.530\\
\hline
$\overline{h}(\gl)$
 & 0.286 &0.495&0.662\\
\hline
$N$& 300000 &500000&160000\\
\hline
$n$& 225000 &330000&970000\\
\hline
C. I.  99\% & $7.179\pm0.050$ &$9.011\pm 0.045$& $7.643 \pm 0.025$\\
\hline
\end{tabular}
\end{center}
\bigskip
\caption{\label{tbl:1}
According to our numerical computations,
the three pairs $(\gl,h)$ are in $\cL$ and this has been tested with the stated
 $p$--values (or probability/level of error).
We report the values of $\overline{h}(\gl)$ and $\underline{h}(\gl)$ for reference.
Of course in these tests there is quite a bit of freedom in the choice of $n$
and $N$: notice that $N$ enters in the evaluation of the $p$--value also
because a larger value of $N$ yields a larger value of $\bbE \log Z_{2N,\go}^{\gl,h}(0)$.
In the last line we report  standard Gaussian $99\%$ confidence intervals for
$\bbE \log Z_{2N,\go}^{\gl,h}(0)$. Of course the $p$--value under the Gaussian assumption
turns out to be totally negligible.}
\end{table}

\begin{rem}
\rm
One might be tempted to interpolate between the values in Table  \ref{tbl:1},
or possibly to get results for small  values of $\gl$ in order to
extend the result of the test to the slope of the critical curve in the origin.
However the fact that $h_c(\gl)$ is strictly increasing does not help much in this direction
and the same is true for
 the finer result, proven in \cite{cf:BG}, that $h_c(\gl)$ can be written as $U(\gl)/\gl$,
$U(\cdot)$ a convex function.
\end{rem}

\smallskip
\subsection{Improving on $\underline{h}(\cdot)$ is uniformly hard}
\label{sec:computer_assisted}
One can get much smaller $p$--values at little computational cost  by choosing
$h$ {\sl just above} $\underline{h}(\gl)$. As a matter of fact
a natural choice is for example  $h=h ^{(0.67)}(\gl)>\underline{h}(\gl)$, recall
\eqref{eq:hm}, for a set of values of $\gl$, and this is part of the content of
Table \ref{tbl:2}: in particular $\bbE \log Z_{2N_+,\go}^{\gl,h^{(0.67)}(\gl)}(0)  >0$
with a probability of error smaller than $10^{-5}$ for the values of $\gl$
between $0.1$ and $1$. However we stress that
for some of these $\gl$'s we have a much smaller $p$--value, see the caption
of Table \ref{tbl:2}, and that the content of this table is much richer
and it approaches also the question of whether or not
a symbolic computation  or some other form of
computer assisted argument could lead to $h_c(\gl)>\underline{h}(\gl)$ for
some $\gl$, and therefore for $\gl$ in an interval. Since such an argument would require
$N$ to be {\sl small}, intuitively the hope resides in large values of $\gl$, recall also
  Remark~\ref{rem:N=2}. It turns out that one needs in any case
  $N$ larger than $700$ in order to observe a localization phenomenon
  at $h^{(0.67)}(\gl)$.
We now give some details on the procedure that leads to Table  \ref{tbl:2}.

\begin{table}[h]
\begin{center}
\begin{tabular}{|c|c|c|c|c|c|c|c|c|c|}
\hline
$\gl $ &$0.05 (\star)$&$0.1$&$0.2$&$0.4$&$0.6$&$1$&$2(\star)$&$4(\star\star )$&$8(\star\star )$\\
\hline
$N_+$&750000&190000&40000&9500&4250&1800&900&800&800\\
\hline
$N_-$&600000&130000&33000&7500&3650&1550&750&700&700\\
\hline
\end{tabular}
\end{center}
\bigskip
\caption{\label{tbl:2}
For a given $\gl$, both $\bbE \log Z_{2N_+,\go}^{\gl,h^{(0.67)}(\gl)}(0)  >0$ and
$\bbE \log Z_{2N_-,\go}^{\gl,h^{(0.67)}(\gl)}(0)  <0$ with a probability of error
smaller than $10^{-5}$ (and in some cases much smaller than that).
Instead for the two cases marked by a $(\star)$ the level of error
is rather between $10^{-2}$ and  $10^{-3}$. For large values of $\gl$, the two cases marked with
 $(\star \star)$, it becomes computationally
expensive to reach small $p$--values. However, above $\gl=3$ one observes
that the values of $Z_{2N,\go}(0)$
essentially do not depend anymore on the value of $\gl$. This can be interpreted
in terms of convergence to a limit ($\gl \to \infty$) model, as it is explained in Remark~\ref{rem:starstar}.
If we then make the hypothesis that this limit model sharply describes the
copolymer along the curve $(\gl, h^{(m)}(\gl))$ for $\gl$ sufficiently large and we apply the concentration
inequality, then the given values of $N_+$  and $N_-$ are tested
with a very small probability of error. Since the details of such a procedure are
quite lengthy we do not report them here.
We have constructed (partial) tables also for different laws of $\go$,
notably $\go_1 \sim N(0,1)$,  and they turned out
to yield larger, at times substantially larger, values of $N_\pm (\gl)$. 
}
\end{table}

First and foremost, the concentration argument
that leads to Proposition~\ref{th:stat} is symmetric
and it works for deviations below the mean as well as above.
So we can, in the very same way, test the null hypothesis
$\bbE\log Z_{2N, \go}(0)> 0$ and, possibly,  refuse it
if $\hat u _n <0$, exactly with the same $p$--value
as in Proposition~\ref{th:stat}.
Of course an important part of Proposition~\ref{th:stat} was coming
from the finite volume localization condition \eqref{eq:loc_char}:
 we do not have an analogous statement for delocalization
(and we do not expect that there exists one).
But, even if $\bbE\log Z_{2N, \go}(0)\le 0$ does not imply delocalization,
it says at least that it is pointless to try to prove localization
by looking at a system of that size.

In Table~\ref{tbl:2} we show
two values of the system size $N$,  $N_+$ and $N_-$,
for which, at a given $\gl$, one has that
$\bbE\log Z_{2N_+, \go}(0)> 0$ and  $\bbE\log Z_{2N_-, \go}(0)< 0$
with a fixed probability of error (specified in the caption of the Table).
It is then reasonable to guess that the transition from
negative to positive values of $\bbE\log Z_{\cdot , \go}(0)$ happens
for $N\in (N_-,N_+)$. There is no reason whastoever to expect that
 $\bbE\log Z_{N, \go}(0)$ should be monotonic in $N$ but according to our numerical result it is not unreasonable to expect that monotonicity should set in for $N$ large or, at least, that for $N< N_-$ (respectively $N>N_+$)
 $\bbE\log Z_{2N, \go}(0)$ is definitely negative
(respectively positive).

\begin{figure}[h]
\begin{center}
\leavevmode
\epsfysize =8 cm
\psfragscanon
\psfrag{N}[c][l]{$N$}
\psfrag{lambda}[c][l]{ $\gl$}
\epsfbox{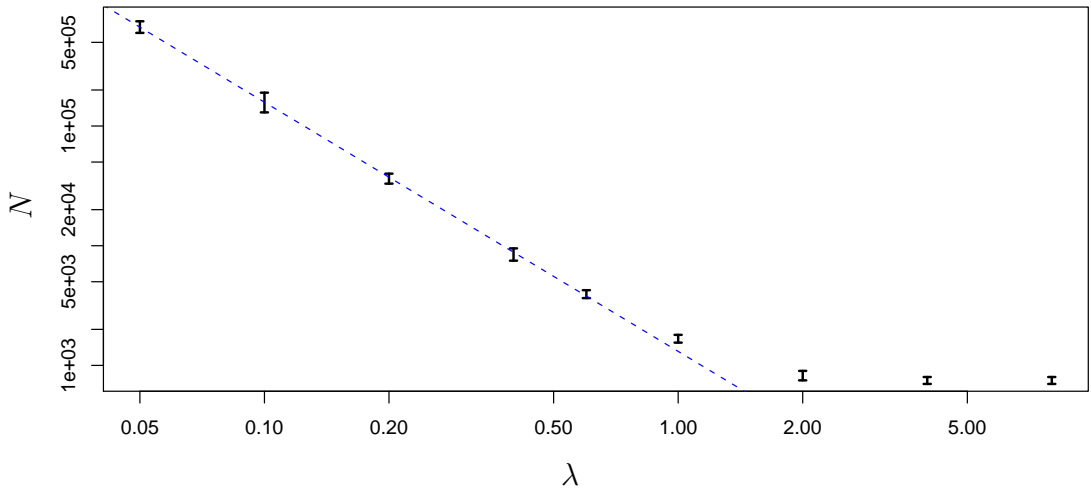}
\end{center}
\caption{\label{fig:NCA} A graphical representation  of
Table~\ref{tbl:2}. The plot is log--log, and a $\gl^{-c}$ behavior is rather
evident, $c$ is about $2.08$. This can be nicely interpreted in terms
of the coarse graining technique in the proof of the weak interaction scaling
limit of the free energy in \cite{cf:BdH}: from that argument one extracts that
if $\gl $ is small the excursions that give a contribution to the free
energy have {\sl typical} length $\gl^{-2}$ and that
in the limit the polymer is just made up by this type of excursions.
One therefore  expects
that it suffices a system of size $N(\gl)$, with $\lim_{\gl \searrow 0} \gl^2 N(\gl)=+\infty$,
to observe localization if $m<h_c^\prime( 0)$, $h=h^{(m)}(\gl)= m \gl (1+o(1))$ and $\gl$ is small.}
\end{figure}

\begin{rem}
\label{rem:starstar}
\rm
As pointed out in the caption of Table  \ref{tbl:2}, from numerics one observes
a very sharp convergence to a $\gl$ independent behavior as $\gl$
becomes large, along the line $h=h^{(m)}(\gl)$.
This is easily interpreted if one observes that $h^{(m)}(\gl)= 1-((\log 2) /2m\gl)+O(\exp(-4m\gl))$
so that
\begin{equation}
\label{eq:lim_mod}
\lim_{\gl \to \infty} \exp\left(-2 \gl \sum_{n=1}^N
\left( \go_n +h\right) \Delta_n \right)\, =\, \exp\left(\frac{\log 2}{m}\sum_{n=1}^N \Delta_n\right)
\ind_{\left\{\sum_{n=1}^N \Delta_n (1+\go_n)=0\right\} } (S).
\end{equation}
This corresponds to the model where a positive charge never enters the lower half-plane and where the energy of a configuration is proportional to the number of negative charges in the lower half-plane.
\end{rem}



\smallskip
\section{Lower bound strategies versus the true strategy}
\label{sec:lb}

\subsection{An approach to lower bounds on the critical curve}
\label{sec:lb_outline}
In this section we give an outline of a new derivation of the lower bound
\begin{equation}
\label{eq:mainBG}
\underline h (\gl)\le h_c(\gl),
\end{equation}
with $\underline h (\gl)$ defined in \eqref{eq:sumupq} of Chapter~\ref{ch:first}.
The complete proof may be found in \S~\ref{app:prooflb}.
The argument takes inspiration from the ideas used in the proof of Proposition 3.1 in \cite{cf:GT} and, even if it is essentially the proof of \cite{cf:BG} in disguise,
in the sense that the selection of the random walk trajectories that are kept
and whose energy contribution is evaluated
does not differ too much (in a word: the {\sl strategy} of the polymer
is similar), it is however conceptually somewhat different and it
will naturally lead to some considerations on the precise
asymptotic behavior of $Z_{N, \go}$ in the delocalized phase and
even in the localized phase close to criticality.

\bigskip

The first step in our proof of (\ref{eq:mainBG}) is a different way of looking at localization. For any fixed positive number $C$ we introduce the stopping time (with respect to the natural filtration of the sequence $\{\go_n\}$) $T^C = T^{C,\gl,h}(\go)$ defined by
\begin{equation}
\label{eq:T^C}
T^{C,\gl,h}(\go) := \inf \{N\in 2\N:\ Z_{N,\go}^{\gl,h}(0) \ge C\}\,.
\end{equation}

The key observation is that if $\bbE[T^C] < \infty$ for some $C>1$, then the polymer is localized. Let us sketch a proof of this fact (for the details, see Proposition~\ref{prop:stopping1} of \S~\ref{app:prooflb}): notice that by the very definition of $T^C$ we have $Z_{T^C(\go),\go}(0) \ge C$. Now the polymer that is in zero at $T^C(\go)$ is equivalent to the original polymer, with a translated environment $\go'=\theta^{T^C} \go$, and setting $T_2(\go) := T^C (\go')$ we easily get $Z_{T_1(\go)+T_2(\go),\go}(0) \ge C^2$ (we have put $T_1(\go):=T^C(\go)$). Notice that the new environment $\go'$ is still typical, since $T^C$ is a stopping time, so that $T_2$ is independent of $T_1$ and has the same law. This procedure can be clearly iterated, yielding an IID sequence $\{T_i(\go)\}_{i=1,2,\ldots}$ that gives the following lower bound on the partition function:
\begin{equation} \label{eq:lowbound}
Z_{T_1(\go)+\ldots +T_n(\go),\go}(0) \ge  C^n\,.
\end{equation}
From this bound one easily obtains that
\begin{equation}
\label{eq:forappB}
\tf(\gl,h) \stackrel{\text{a.s.}}{= } \lim_{n\to\infty} \frac{\log Z_{T_1(\go)+\ldots +T_n(\go),\go}(0)}{T_1(\go)+\ldots +T_n(\go)} \ge \frac{\log C}{\bbE[T^C]}\,,
\end{equation}
where we have applied the strong law of large numbers, and localization follows since by hypothesis $C>1$ and $\bbE[T^C]<\infty$.

\medskip

\begin{rem}
\label{rem:reciprocal}
\rm
It turns out that also the reciprocal of the claim just proved holds true, that is \textsl{the polymer is localized if and only if $\bbE[T^C]<\infty$}, with an arbitrary choice of $C>1$, see Proposition~\ref{prop:stopping1}. In fact the case $\bbE[T^C]=\infty$ may arise in two different ways:
\begin{enumerate}
\item the variable $T^C$ is defective, $\bbP[T^C=\infty]>0$: in this case with positive probability $\{Z_{N,\go}(0)\}_N$ is a bounded sequence, and delocalization follows immediately;
\item\label{en:scenario} the variable $T^C$ is proper with infinite mean, $\bbP[T^C=\infty]=0,\ \bbE[T^C]=\infty$: in this case we can still build a sequence $\{T_i(\go)\}_{i=1,2,\ldots}$ defined as above and this time the lower bound \eqref{eq:lowbound} has \textsl{subexponential} growth. Moreover it can be shown that in this case the lower bound \eqref{eq:lowbound} gives the true free energy, cf. Lemma~\ref{lem:chop}, which therefore is zero, so that delocalization follows also in this case.
\end{enumerate}
As a matter of fact, it is highly probable that in the interior of the delocalized phase
$Z_{N,\go}(0)$ vanishes $\bbP (\dd \go)$--a.s.
when $N \to \infty$ and this would rule out the scenario (\ref{en:scenario}) above, saying that for $C>1$ the random variable $T^C$ must be either integrable or defective. We take up again this point in Sections~\ref{sec:path} and~\ref{sec:guess}: we feel that this issue
is quite crucial in order to fully understand the delocalized phase
of disordered models.
\end{rem}


\begin{rem}
\label{rem:analogy}
\rm
Dealing directly with $T^C$ may be difficult.
Notice however that if one finds
 a random time (by this we mean simply an integer--valued random variable)
 $T=T(\go)$ such that
 \begin{equation} \label{eq:but}
Z_{T(\go), \go}(0)\ge C>1\,, \qquad \text{with \ } \bbE[T]<\infty\,,
\end{equation}
then localization follows. This is simply because
this implies $T^C \le T$ and hence $\bbE[T^C]<\infty$. Therefore localization is equivalent to the condition $\log Z_{T(\go), \go}(0)>0$ for an \textsl{integrable} random time $T(\go)$: we would like to stress the analogy between this and the criterion for localization given in \S~\ref{sec:superadd}, see \eqref{eq:loc_char}.
\end{rem}

\medskip

Now we can turn to the core of our proof: we are going to show that for every $(\gl,h)$ with $h<\underline h (\gl)$ we can build a random time $T=T(\go)$ that satisfies \eqref{eq:but}.
The construction of $T$
 is based on the  idea that for $h>0$ if localization prevails
 is because of rare $\go$--stretches that invite the polymer
 to spend time in the lower half--plane in spite of the action of $h$.

\begin{figure}[h]
\begin{center}
\leavevmode
\epsfysize =5 cm
\psfragscanon
\psfrag{0}[c][l]{$0$}
\psfrag{l}[c][l]{ $\ell$}
\psfrag{n}[c][l]{ $n$}
\psfrag{h}[l][l]{ $L$}
\psfrag{S}[c][l]{$S$}
\psfrag{T}[c][l]{$T(\go)$}
\epsfbox{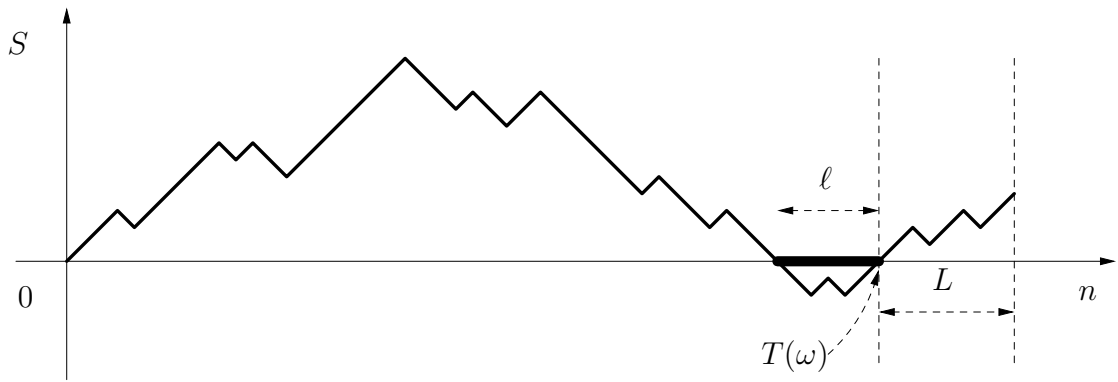}
\end{center}
\caption{\label{fig:ZT}
Inequality \eqref{eq:beyondT}
comes simply from restricting the evaluation of $Z_{T(\go) +L, \go} $ to the trajectories
visiting the {\sl $q$--atypical} stretch of length $\ell$  and by staying away from
the unfavorable solvent after that.
}
\end{figure}

The strategy we use consists in looking
for $q$--atypical stretches of length at least $M\in 2\N$, where $q<-h$ is
the average charge of the stretch. Rephrased a bit more precisely,
we are looking for the smallest $n\in 2\N$ such that
$\sum_{i=n-k+1}^n \go_i /k <q$ for some even integer $k \ge M$.
It is well known that such a random variable grows, in the sense
of Laplace, as $\exp(\Sigma(q) M)$ for $M\to \infty$, where $\gS(q)$ is the Cramer functional
\begin{equation}
\label{eq:cramer}
\gS(q) := \sup_{\ga\in\R} \{\ga q-\log \M(\ga)\}\,.
\end{equation}
One can also show without much effort that the length of such a stretch
cannot be much longer than $M$.
Otherwise stated, this is the familiar statement that the longest $q$--atypical
sub--stretch of $\go_1, \ldots , \go_N$ is of typical length $\sim \log N/\Sigma (q)$.
 So $T(\go)$ is for us the end--point
of a $q$--atypical stretch of length approximately $(\log T(\go))/\Sigma(q)$:
by looking for sufficiently long $q$--atypical stretches
we have always the freedom to choose $T(\go) \gg 1$, in such a
way that also $\log T(\go) \ll T(\go)$ and this is helpful for the estimates.
So let us bound $Z_{T(\go), \go}$ from
below by considering only the trajectories of the walk that
stay in the upper half--plane up to the beginning of the $q$--atypical stretch and
that are negative in the stretch, coming back to zero at step $T(\go)$
(see Fig.~\ref{fig:ZT}: the polymer is cut at the first dashed vertical line).
The contribution of these trajectories is easily evaluated: it is approximately
\begin{equation}
\label{eq:int1_cgg}
\left(
\frac {1} {T(\go)^{3/2}}
\right)
 \exp\left( -2\gl (q+h)  \frac{\log T(\go)}{\Sigma(q)}\right).
\end{equation}
For such an estimate we have used
\eqref{eq:asympt} and $\log T(\go) \ll T(\go)$
both in writing the probability that the first return to zero
of the walk is at the beginning of the $q$--atypical stretch and in neglecting the probability
that the walk is negative inside the stretch.
It is straightforward to see that if
\begin{equation}
\label{eq:int2_cgg}
\frac {4\gl} 3 h < -\frac {4\gl} 3 q -  \Sigma(q),
 \end{equation}
and if $T(\go)$ is large, then also  the quantity in
 \eqref{eq:int1_cgg} is large. We can still optimize this procedure by choosing $q$ (which must be sufficently negative, i.e. $q < -h$).
 By playing with
\eqref{eq:cramer}
one sees that one can choose $q_0\in \R$ such that for $q=q_0$ the right--hand side
in \eqref{eq:int2_cgg} equals $\log \M (-4\gl/3)$ and
if $h<\log \M (-4\gl/3)/(4\gl/3) = \underline{h} (\gl)$
 then $q_0 <-h$. This argument therefore is saying that there
 exists $C>1$ such that
 \begin{equation}
 \label{eq:C1}
 Z_{T(\go), \go}(0) \, \ge \, C,
\end{equation}
for every $\go$. It only remains to show that $\bbE[T]<\infty$: this fact, together with a detailed proof of the argument just presented can be found in \S~\ref{app:prooflb}.

\medskip

\begin{figure}[h]
\begin{center}
\leavevmode
\epsfysize =14 cm
\psfragscanon
\psfrag{x}[c][l]{\small $N$}
\psfrag{y}[c][l]{\small $\log  Z_{2N, \go}^{\gl, h}$}
\psfrag{a}[l][l]{A: $h=0.42$}
\psfrag{b}[l][l]{B:  $h=0.44$}
\psfrag{c}[l][l]{C:  $h=0.43$}
\psfrag{d}[l][l]{D:  $h=0.43$ (Zoom)}
\epsfbox{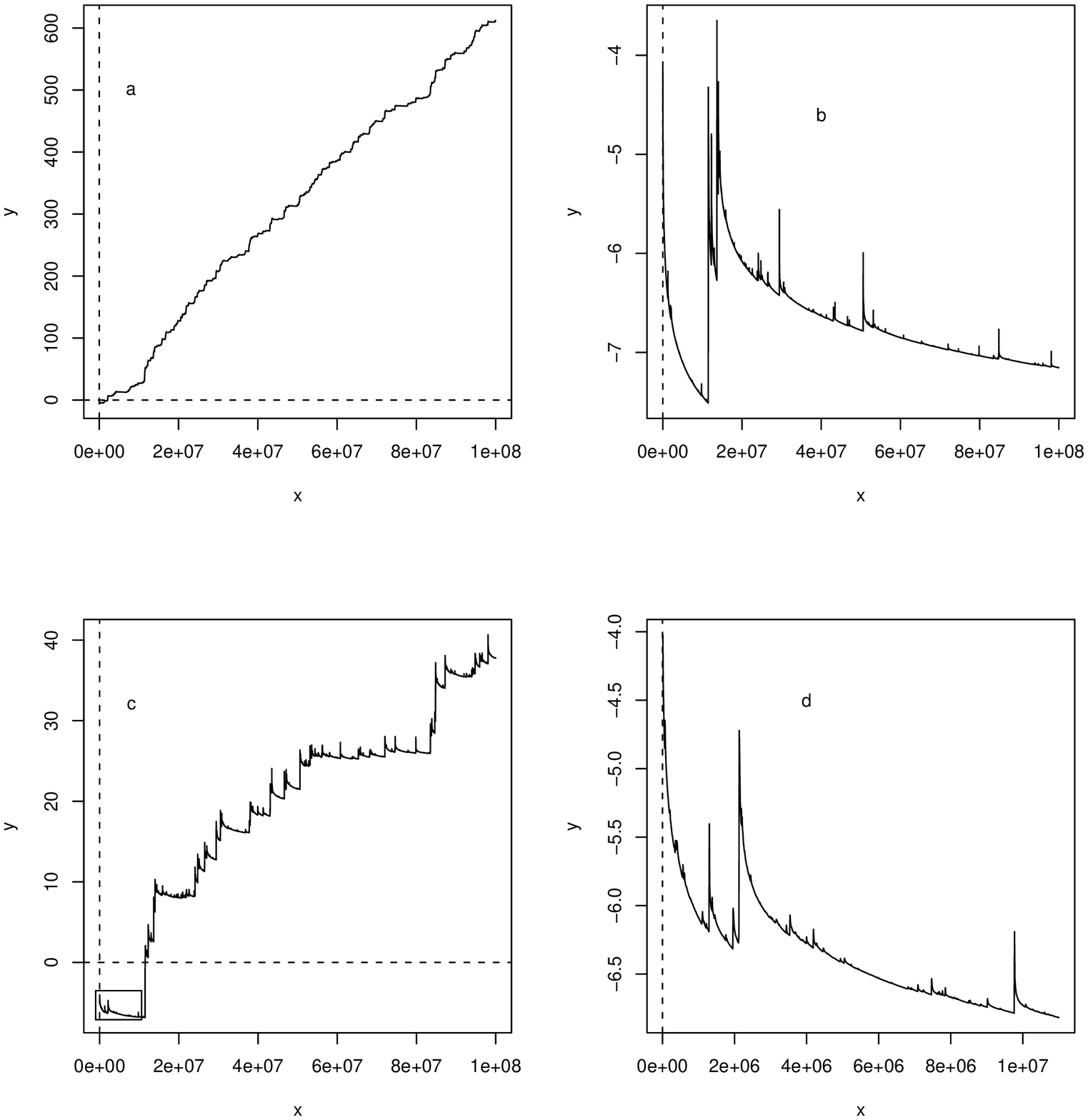}
\end{center}
\caption{
\label{fig:tmp}
For $\gl=0.6$ ($\underline h (0.6)\simeq 0.36$ and $\overline h (0.6)\simeq 0.49$) , the behavior of
$\log Z_{2N, \go}$ for $h=0.42$ (A), $0.43$ (C,D) and $0.44$ (B).
In case A, the polymer is localized with free energy approximately
$3\cdot 10^{-6}$:
the linear growth is quite clear, but a closer look shows sudden jumps,
which correspond to atypically negative stretches of charges.
Getting closer to the critical point, case C,
the linear growth is still evident, but
it is clearly the result of sudden growths followed by slow decays
(approximately polynomial with exponent $-1/2$).
Case B suggests delocalization: a closer analysis reveals a decay
of the type $N ^{-1/2}$, but sharp deviations are clearly visible.
Case D is the zoom of the rectangle in the left corner of C.
The similarity between B and D
make clear that claiming delocalization looking at the behaviour of the partition function is difficult.
}
\end{figure}

\smallskip
\subsection{Persistence of the effect of rare stretches}
As pointed out in the previous section, there is strong evidence
that $h_c(\gl)> \underline{h} (\gl)$.
At this stage Fig.~\ref{fig:tmp}  is of particular interest.
Notice first of all that in spite of being substantially above $\underline{h}(\cdot)$
the copolymer appears to be still localized, see in particular case A.

The rigorous lower bounds that we are able
to prove cannot establish localization in the region we are considering.
All the same, notice that if one does  not cut the polymer at $T(\go)$,
as in the argument above, but at $T(\go)+L$,
a  lower bound of the following type
\begin{equation}
\label{eq:beyondT}
Z_{T(\go)+L, \go}\, \stackrel{\text{roughly}}{\ge}\,
\text{const.} \frac1{T(\go)^{3/2}}\,
\exp\left( -2\gl (q+h)\frac{\log T(\go)}{\Sigma (q)}\right) \, \frac 1 {L^{1/2}},
\end{equation}
is easily established. Of course we are being imprecise, but we just want to convey
the idea, see also Fig.~\ref{fig:ZT}, that after passing through an atypically {\sl negative}
stretch of environment ($q>0$), the effect of this stretch decays at most like $L^{-1/2}$,
that is the probability that a walk stays  positive for a time $L$.

At this point we stress that the argument outlined in \S~\ref{sec:lb_outline}
and re--used for \eqref{eq:beyondT}
may be very well applied to $h> \underline{h}(\gl)$, except that
this time it does not suffice for $\eqref{eq:C1}$.
But it yields nevertheless  that for $h \in \left(\underline{h}(\gl), \overline{h}(\gl)\right)$
the statement $Z_{N,\go}\sim N^{-1/2}$, something a priori expected (for
example \cite{cf:BH})
in the delocalized regime and true for non disordered systems, is
violated. More precisely, one can find a sequence of random
times $\{\tau_j\}_j$, $\lim_j \tau _j= \infty$ such that
$Z_{\tau_j,\go}\ge {\tau_j}^{-1/2+a}$, $a=a(\gl,h)>0$ (see  Proposition~4.1 in \cite{cf:GT}).
These random times are constructed exactly by looking
for $q$--atypical stretches as above and one can appreciate
such an irregular decay for example in case B of Fig.~\ref{fig:tmp},
and this in spite of the fact that the data have been  strongly coarse grained.

\smallskip

Therefore the lower bound
\eqref{eq:beyondT}, both in the localized and in the delocalized regime,
yields the following  picture:
the lower bound we found on  $Z_{N, \go}$ grows suddenly in correspondence of atypical stretches
and after that it decays with an exponent $1/2$, up to another atypical stretch.
This matches Fig.~\ref{fig:tmp}, at least at a qualitative level, see the caption
of the figure.


\smallskip
\section{The delocalized phase: a path analysis}
\label{sec:path}

Let us start with a qualitative observation: if we set the parameters $(\gl,h)$ of the copolymer to $(\gl, h^{(m)}(\gl))$ with $m = 0.9$, then the observed behavior of $\{Z_{N,\go}^{\gl,h}(0)\}_N$ --suitably averaged over blocks in order to eliminate local fluctuations-- is somewhat close to $\text{(const)}/N^{3/2}$. This is true for all the numerically accessible values of $N$ (up to $N\sim 10^8$), at once for a number of values of~$\gl$ and for a great number of typical environments $\go$. Of course this is suggesting that for $m=0.9$ the curve $h^{(m)}(\gl)$ lies in the delocalized region, but it is not easy to convert this qualitative observation into a precise statement, because we do not have a rigorous finite--volume criterion to state that a point $(\gl,h)$ belongs to the delocalized phase (the contrast with the localized phase, see~\eqref{eq:loc_char}, is evident). In other words, we cannot exclude the possibility that the system is still localized but with a characteristic size much larger than the one we are observing.

Nevertheless, the aim of this section is to give an empirical criterion, based on an analysis of the path behavior of the copolymer, that will allow us to provide some more quantitative argument in favor of the fact that the curve $h^{(m)}(\gl)$ lies in the delocalized region even for values of~$m<1$. This of course would entail that the upper bound $\overline h(\gl)$ defined in \eqref{eq:sumup_alt} is not strict.

\smallskip
\subsection{Known and expected path behavior}
\label{sec:path1}
We want to look at the whole \textsl{profile} $\{Z_{N,\go^r}^{\gl,h}(x)\}_{x\in\Z}$ rather than only at $Z_{N,\go^r}^{\gl,h}(0)$, where by $\go^r$ we mean the environment $\go$ in the {\sl
 backward direction}, that is $(\go^r)_n := \go_{N+1-n}$ (the reason for this choice is explained in Remark~\ref{rem:backward} below). The link with the path behavior of the copolymer, namely the law of $S_N$ under the polymer measure $\bP_{N,\go^r}^{\gl,h}$, is given by
\begin{equation}
    \frac{Z_{N,\go^r}^{\gl,h}(x)}{Z_{N,\go^r}^{\gl,h}} = \bP_{N,\go^r}^{\gl,h} (S_N = x)\,.
\end{equation}

We have already remarked in \S~\ref{sec:pathwise} of Chapter~\ref{ch:first} that, although the localized and delocalized phases have been defined in terms of free energy, they do correspond to sharply different path behaviors. In the localized phase it is known \cite{cf:Sinai,cf:BisdH} that the laws of $S_N$ under $\bP_{N,\go^r}^{\gl,h}$ are \textsl{tight}, which means that the polymer is essentially at $O(1)$ distance from the $x$--axis. The situation is completely different in the (interior of the) delocalized phase, where one expects that $S_N = O(\sqrt{N})$: in fact the conjectured path behavior (motivated by the analogy with the known results for non disordered models, see in particular \cite{cf:MGO},  \cite{cf:DGZ} and \cite{cf:CGZ}) should be weak convergence under diffusive scaling to the \textsl{Brownian meander process} (that is Brownian motion conditioned to stay positive on the interval $[0,1]$, see \cite{cf:RevYor}). Therefore in the (interior of the) delocalized phase the law of $S_N/\sqrt{N}$ under $\bP_{N,\go^r}^{\gl,h}$ should converge weakly to the corresponding marginal of the Brownian meander, whose law has density $x \exp(-x^2/2) \ind_{(x\ge 0)}$.

In spite of the lack of precise rigorous results,
the analysis we are going to describe is carried out under the hypothesis that, in the interior of the delocalized phase, the scaling limit towards Brownian meander holds true (as it will be seen, the numerical results provide a sort of \textsl{a posteriori} confirmation of this hypothesis).

\begin{rem} \label{rem:backward} \rm
From a certain point of view attaching the environment backwards does not change too much the model: for example it is easy to check that if one replaces $\go$ by $\go^r$ in \eqref{eq:felim2}, the limit still exists $\bbP(\dd\go)$--a.s. and in $\bbL_1(\dd\bbP)$. Therefore the free energy is the same, because $\{\go^r_n\}_{1\le n \le N}$ has the same law as~$\{\go_n\}_{1\le n \le N}$, for any fixed~$N$.

However, if one focuses  on the law of $S_N$ as a function of~$N$ \textsl{for a fixed environment} $\go$, the behavior reveals to be much smoother under $\bP_{N,\go^r}^{\gl,h}$ than under $\bP_{N,\go}^{\gl,h}$. For instance, under the original polymer measure $\bP_{N,\go}^{\gl,h}$ it is no more true that in the localized region the laws of $S_N$ are tight (it is true only {\sl most of the time}, see \cite{cf:G} for details). The reason for this fact is to be sought in the presence of long {\sl atypical} stretches in every typical~$\go$ (this fact has been somewhat quantified in \cite[Section 4]{cf:GT} and it is at the heart of the approach in Section~\ref{sec:lb}) that are encountered along the copolymer as $N$ becomes larger. Of course the effect of these stretches is very much damped with the backward environment.

A similar and opposite phenomenon takes place also in the delocalized phase. In fancier words, we could say that for fixed $\go$ and as $N$ increases, the way $S_N$ approaches its {\sl limiting behavior} is faster when the environment is attached backwards: it is for this reason that we have
chosen to work with $\bP_{N,\go^r}^{\gl,h}$.
\end{rem}

\smallskip
\subsection{Observed path behavior: a numerical analysis}
In view of the above considerations, we choose as a measure of the delocalization of the polymer the
$\ell_1$ distance $\bigtriangleup _N^{\gl,h}(\go)$ between the numerically computed profile for a polymer of size $2N$ under $\bP_{2N,\go^r}^{\gl,h}$, and the conjectured asymptotic delocalized profile:
\begin{equation}
    \bigtriangleup _N^{\gl, h}(\go) := \sum_{x \in 2\Z} \Bigg| \frac{Z_{N,\go^r}^{\gl,h}(x)}{Z_{N,\go^r}^{\gl,h}} \;-\; \frac{1}{\sqrt{2N}} \, \varphi^+ \bigg( \frac{x}{\sqrt{2N}} \bigg) \Bigg|\,,  \qquad \gp^+(x) := x \, e^{-x^2/2} \ind_{(x\ge 0)}\,.
\end{equation}
Loosely speaking, when the parameters $(\gl,h)$ are in the interior of the the delocalized region we expect $\bigtriangleup _N$ to decrease to~$0$ as $N$ increases, while this certainly will not happen if we are in the localized phase.

\smallskip

The analysis has been carried out at $\gl=0.6$: we recall that the lower and upper bound of \eqref{eq:sumup_alt} give respectively $\underline h(0.6) \simeq 0.36$ and $\overline h(0.6) \simeq 0.49$, while the lower bound we derived with our  test for localization is $h=0.41$, see Table~\ref{tbl:1}. However,
as observed in Section \ref{sec:lb}, Fig.~\ref{fig:tmp}, there is numerical evidence that $h=0.43$ is still localized, and for this reason  we have analyzed the values of $h=0.44, 0.45, 0.46, 0.47$ (see below for an analysis on smaller values of~$h$).

\begin{figure}[h]
\begin{center}
\leavevmode
\epsfysize =7.9 cm
\epsfbox{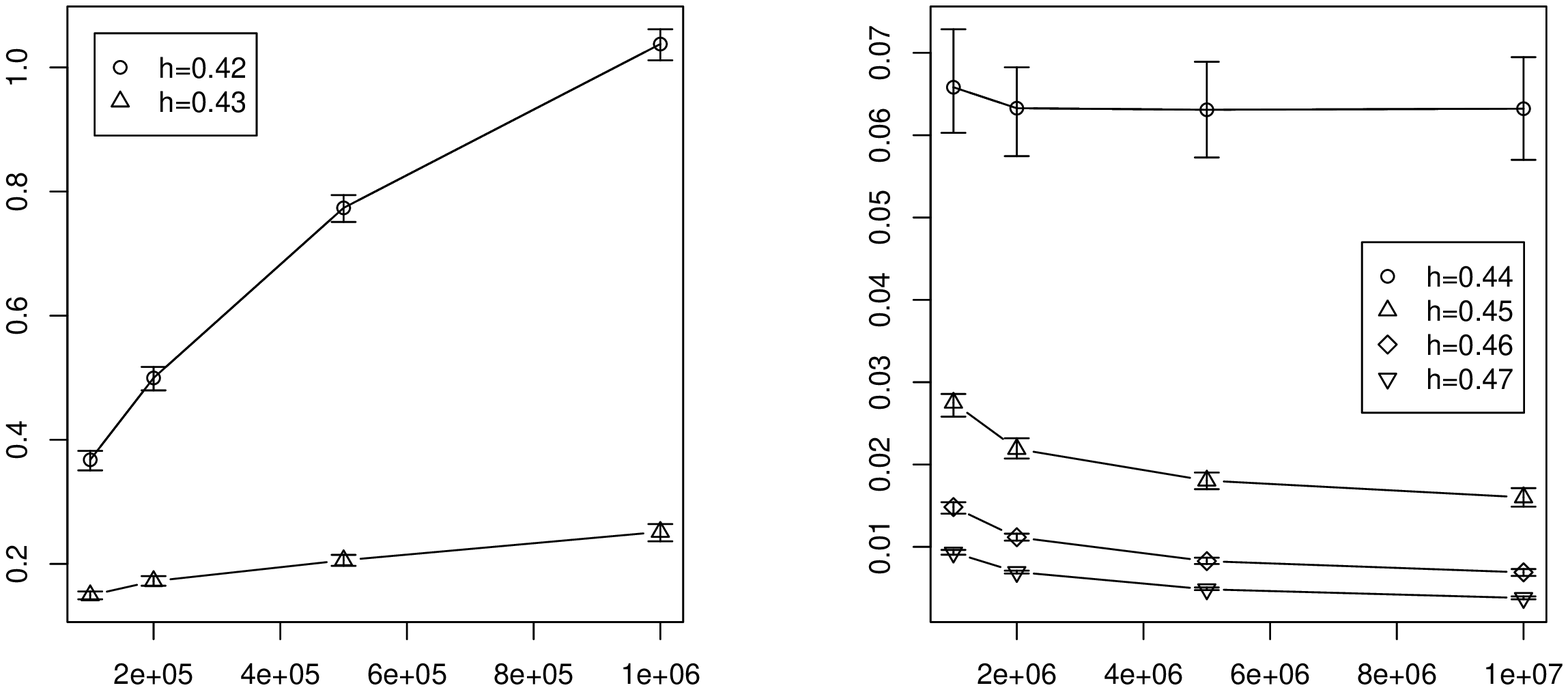}
\end{center}
\caption{
\label{fig:F}
Graphical representation of the data of Tables~\ref{tbl:F} (on the right) and~\ref{tbl:F2} (on the left). The plotted points are the sample medians against the sample size, the error bars correspond to the confidence intervals given in Tables~\ref{tbl:F} and~\ref{tbl:F2}.
}
\end{figure}

For each couple $(\gl, h)$ we have computed $\bigtriangleup _{N}^{\gl,h}(\go)$ for the sizes $N=a\times 10^6$ with $a=1,2,5,10$ and for $500$ independent environments. Of course some type of statistical analysis must be performed on the data in order to decide whether there is a decay of~$\bigtriangleup $ or not. The most direct strategy would be to look at the sample mean of
a family of IID variables distributed like $\bigtriangleup _N(\go)$, but it turns out that the fluctuations are too big to get reasonable confidence intervals for this quantity (in other words, the sample variance does not decrease fast enough), at least for the numerically accessible sample sizes. A more careful analysis shows that the variance is essentially due to a \textsl{very small} fraction of data that have \textsl{large} deviations from the mean, while the most of the data mass is quite concentrated.

For this reason we have chosen to focus on the \textsl{sample median} rather than on the sample mean. Table~\ref{tbl:F} contains the results of the analysis (see also Fig.~\ref{fig:F} for a graphical representation): for each value of $h$ we have reported the standard $95\%$ confidence interval for the sample median (see Remark~\ref{rem:conf_int} below for details) for the four different values of~$N$ analyzed. While for $h=0.44$ the situation is not clear, we see that for the values of $h$ greater than $0.45$ there are quantitative evidences for a decrease in~$\bigtriangleup_{N}$: this leads us to the conjecture that the points $(\gl,h)$ with $\gl=0.6$ and $h\ge 0.45$ (equivalently, the points $(\gl,h^{(m)}(\gl))$ with $m \gtrsim 0.876$) lie in the delocalized region.

\begin{table}[h]
\begin{center}
\begin{tabular}{|c||c|c|c|c|}
\hline
$h \backslash N(\times 10^6) $ & 1 & 2 & 5 & 10\\
\hline
\hline
0.44 & [.0603, .0729] & [.0574, .0682] & [.0572, .0689] & [.0570, .0695] \\
\hline
0.45 & [.0258, .0286] & [.0207, .0232] & [.0170, .0190] & [.0149, .0171] \\
\hline
0.46 & [.0140, .0154] & [.0108, .0116] & [.00792, .00869] & [.00647, .00731] \\
\hline
0.47 & [.00905, .00963] & [.00676, .00711] & [.00475, .00508] & [.00364, .00398] \\
\hline
\end{tabular}
\end{center}
\bigskip
\caption{ \label{tbl:F}
The table contains the standard $95\%$ confidence interval for the median of a sample $\{\bigtriangleup _N^{\gl,h}(\go)\}_{\go}$ of size 500, where $\gl=0.6$ and $h,N$ take the different values reported in the table. For the values of $h \ge 0.45$ the decreasing behavior of $\bigtriangleup_{N}$ is quite evident (the confidence intervals do not overlap), see also Fig.~\ref{fig:F}.}
\end{table}

As already remarked, these numerical observations cannot rule out the possibility that the system is indeed localized, but the system size is too small to see it. For instance, we have seen that there are evidences for $h=0.43$ to be localized (see case~C of Fig.~\ref{fig:tmp}).
In any case, the exponential increasing of $Z_N(0)$ is detectable only at sizes of order$\sim 10^8$, while for smaller system sizes (up to$\sim 10^7$) the qualitative observed behavior of $Z_N(0)$ is rather closer to $(const)/N^{3/2}$, thus apparently suggesting delocalization (see case~D of Fig.~\ref{fig:tmp}).

For this reason it is interesting to look at $\bigtriangleup_N^{0.6,\,h}$ for $h=0.42, 0.43$ and for $N \ll 10^8$. For definiteness we have chosen $N=a\times 10^6$ with $a=1,2,5,10$, performing the computations for $3000$ independent environments: the results are reported in Table~\ref{tbl:F2} (see also Fig.~\ref{fig:F}). As one can see, this time there are clear evidences for an \textsl{increasing} behavior of $\bigtriangleup_N$. On the one hand this fact gives some more confidence on the data of Table~\ref{tbl:F}, on the other hand it suggests that looking at $\{\bigtriangleup_N\}_N$ is a more reliable criterion for detecting (de)localization than looking at $\{Z_N(0)\}_N$.

\begin{table}[h]
\begin{center}
\begin{tabular}{|c||c|c|c|c|}
\hline
$h \backslash N(\times 10^5) $ & 1 & 2 & 5 & 10\\
\hline
\hline
0.42 & [.351, 0.382] & [.480, 0.517] & [.751, 0.794] & [1.01, 1.06] \\
\hline
0.43 & [.143, 0.155] & [.165, 0.180] & [.197, 0.215] & [.236, 0.264] \\
\hline
\end{tabular}
\end{center}
\bigskip
\caption{ \label{tbl:F2}
The table contains the standard $95\%$ confidence interval for the median of a sample $\{\bigtriangleup _N^{\gl,h}(\go)\}_{\go}$ of size 3000, where $\gl=0.6$ and $h,N$ take the values reported in the table. For both values of~$h$ an increasing behavior of $\bigtriangleup_{N}$ clearly emerges, see also Fig.~\ref{fig:F} for a graphical representation.}
\end{table}

\begin{rem} \label{rem:conf_int} \rm
A confidence interval for the sample median can be obtained in the following general way
(the steps below are performed under the assumption that the median
is unique, which is, strictly speaking, not true in our case, but it will be clear that
a finer analysis would not change the outcome). Let $\{Y_k\}_{1\le k \le n}$ denote a sample of size~$n$, that is the variables $\{Y_k\}_k$ are independent with a common distribution, whose median we denote by~$\xi_{1/2}$:
$\bP \left(Y_1 \le \xi_{1/2}\right)=1/2$. Then the variable
\begin{equation}
    \cN_n := \# \{i \le n:\ Y_i \le \xi_{1/2}\}
\end{equation}
has a binomial distribution $\cN_n \sim B(n,1/2)$ and when $n$ is large (for us it will be at least~500) we can approximate $\cN_n/n \approx 1/2 + Z/(2\sqrt{n})$, where $Z \sim N(0,1)$ is a standard gaussian. Let us denote the sample quantiles by $\Xi_q$, defined for $q \in (0,1)$ by
\begin{equation}
    \# \{i \le n:\ Y_i \le \Xi_q\} = \lfloor qn \rfloor\,.
\end{equation}
If we set $a:= |\Phi^{-1}(0.025)|$ ($\Phi$ being the standard gaussian distribution function) then the random interval
\begin{equation}
    \Big[\Xi_{\frac{1}{2}-\frac{a}{2\sqrt{n}}},\; \Xi_{\frac{1}{2}+\frac{a}{2\sqrt{n}}}\Big]
\end{equation}
is a $95\%$ confidence interval for $\xi_{1/2}$, indeed
\begin{align}
    0.95 &= \bP \big( Z \in [-a,a] \,\big) = \bP \bigg( \frac 12 + \frac{1}{2\sqrt{n}}Z \;\in\; \Big[\frac{1}{2}-\frac{a}{2\sqrt{n}}\;,\; \frac{1}{2}+\frac{a}{2\sqrt{n}} \Big] \bigg) \nonumber \\
    &\approx \bP \bigg ( \frac{\cN_n}{n} \in \Big[\frac{1}{2}-\frac{a}{2\sqrt{n}}\;,\; \frac{1}{2}+\frac{a}{2\sqrt{n}} \Big] \bigg) = \bP \bigg( \Xi_{\frac{1}{2}-\frac{a}{2\sqrt{n}}} \le \xi_{1/2} \le \Xi_{\frac{1}{2}+\frac{a}{2\sqrt{n}}} \bigg)\,.
\end{align}
\end{rem}


\smallskip
\section{An empirical observation on the critical curve}
\label{sec:guess}

The key point of this section is that, from a numerical viewpoint,
$h_c(\cdot)$ seems very close to $h^{(m)}(\cdot)$, for a suitable
value of $m$. Of course any kind of
 statement in this direction
 requires first of all a procedure to estimate $h_c(\cdot)$ and
we explain this first.

Our analysis is based on the following
conjecture:
\begin{equation}
\label{eq:conject_h}
(\gl, h)\in \overset{\circ}{\cD} \, \, \Longrightarrow \ \
\lim_{N \to \infty }Z_{2N, \go}^{\gl, h} (0) \, =\, 0, \ \bbP \left( \dd \go \right)-\text{a.s.}.
\end{equation}
The  arguments in Section~\ref{sec:lb}
suggest the validity of such a conjecture, which  is  comforted  by the numerical observation.
Since, if $(\gl, h)\in \cL$,  $Z_{2N, \go}^{\gl, h} (0)$ diverges (exponentially fast)
$\bbP \left( \dd \go \right)$--almost surely and since $Z_{2N, \go}^{\gl, h} (0)$
is decreasing  in $h$, we define $\hat h _{N, \go} (\gl)$ as the only $h$ that solves
$ Z_{2N, \go}^{\gl, h} (0)=1$.
We expect that $\hat h _{N, \go} (\gl)$ converges to $h_c(\gl)$
as $N$ tends to infinity, for typical $\go$'s.
Of course setting the threshold to the value $1$ is rather arbitrary, but it is
somewhat suggested by \eqref{eq:loc_char} and by  the idea behind
the proof of \eqref{eq:mainBG} (Proposition  \ref{prop:stopping1} and
equation
\eqref{eq:T^C}).

\begin{figure}[h]
\begin{center}
\leavevmode
\epsfysize =8.5 cm
\psfragscanon
\psfrag{x}[c][l]{$\gl$}
\psfrag{xg}[c][l]{$\gl$}
\psfrag{y}[c][l]{$\hat h _{N,\go}(\gl)$}
\psfrag{yg}[c][l]{$\hat h _{N,\go}(\gl)$}
\epsfbox{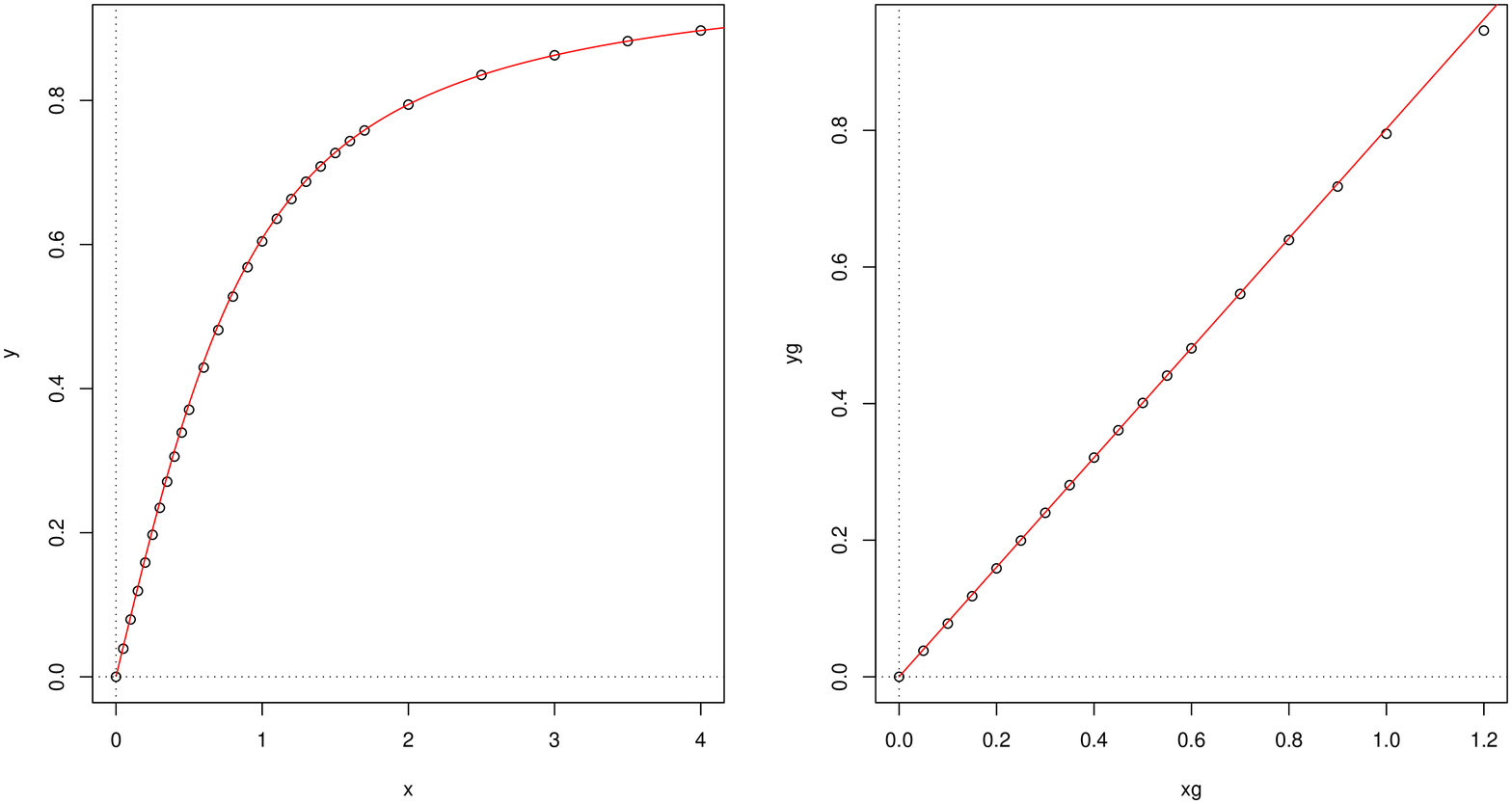}
\end{center}
\caption{
\label{fig:sec5_1}
On the left the case of binary symmetric $\go_1$ and on the right
the case of $\go_1 \sim N(0,1)$, boths for $N= 3.2\cdot 10^7$.
The small circles represent the computed values: the errors on $\hat h _{N,\go}(\gl)$
are
negligible and the plotted points are at the centers of the circles. The continuous
line is instead the curve $h^{(m)}(\cdot)$. In the binary case $m=0.841$ and it has been chosen
by solving $h^{(m)}(4)=\hat h _{N,\go}(4)$. In the Gaussian case $m=0.802$, the maximum
of $ \hat h _{N,\go}(\gl)/\gl$ for the plotted values of $\gl(>0)$.
The rather different values of $\hat m_{N, \go}$ may be somewhat understood
both by considering that these two curves have been obtained for a fixed
realization of $\go$ and by taking into account the remark at the end of the caption
of Table~\ref{tbl:2}: it appears that for  Gaussian charges one needs longer systems in order to get closer
to the values of $m$ observed in the binary case (in particular: for the prolongation,
with the same random
number generator,
of the Gaussian $\go$ sample used here up to $N=5\cdot 10^{7}$ one
obtains $\hat m_{N, \go}=0.812$).
}
\end{figure}

\begin{figure}[h]
\begin{center}
\leavevmode
\epsfysize =7.1 cm
\psfragscanon
\psfrag{Binary}[c][l]{\small $\go_1 =\pm 1, \, \go_1 \sim -\go_1 $}
\psfrag{Gauss}[c][l]{$\go_1 \, \sim \, N(0,1)$}
\psfrag{x}[c][l]{$\gl$}
\psfrag{xg}[c][l]{$\gl$}
\psfrag{y}[c][l]{$\hat h _{N,\go}(\gl)$}
\psfrag{yg}[c][l]{$\hat h _{N,\go}(\gl)$}
\psfrag{tm1}[c][l]{$r_{N,\go}(\gl)$}
\psfrag{tmp1}[c][l]{$r_{N,\go}(\gl)$}
\epsfbox{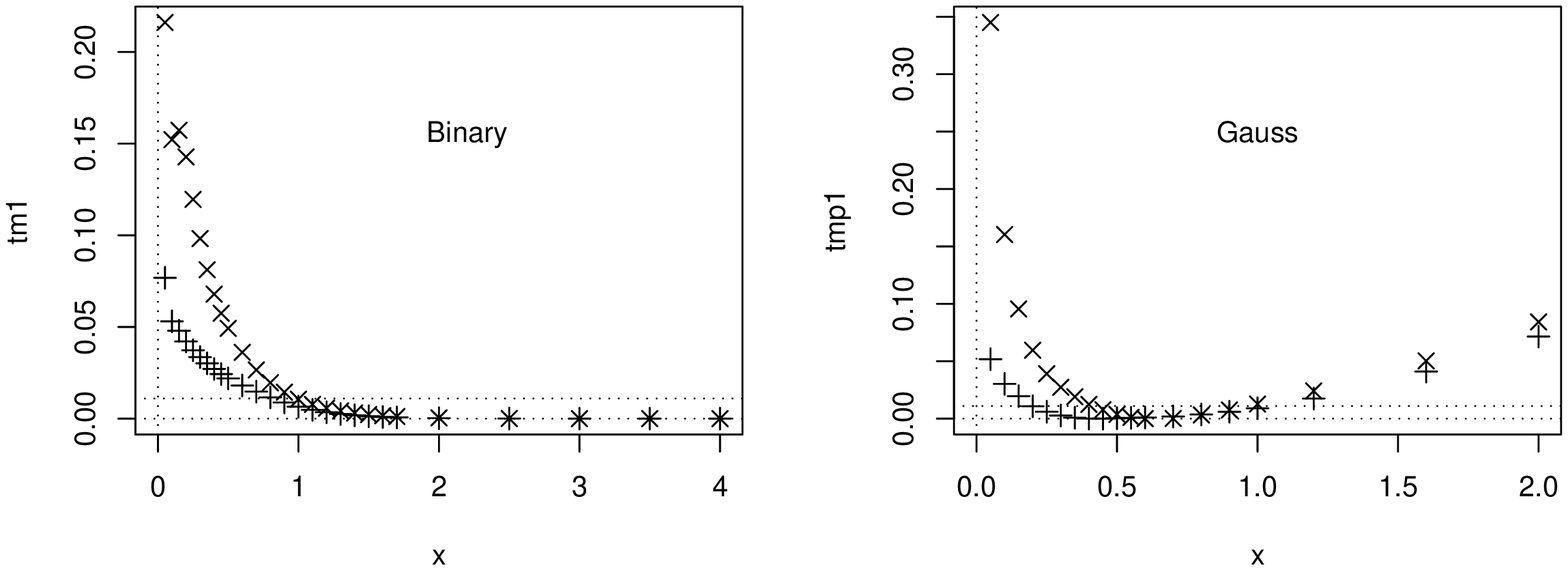}
\end{center}
\caption{
\label{fig:sec5_2}
Relative errors $r_{N,\go}(\gl):= \left( h^{(m)}(\gl) -\hat h _{N, \go} (\gl) \right)/\hat h _{N, \go} (\gl)$,
for the value $m=\hat m_{N,\go}$ explained in the caption of Fig.~\ref{fig:sec5_1} and
for the cases of $N=2.5\cdot 10^5$ ($\times$ dots), and $N=3.2\cdot 10^7$ ($+$ dots).
Notice that in the binary case the error is more important for small values of $\gl$ (recall
Table~\ref{tbl:2} and Fig.~\ref{fig:NCA}). Instead for the Gaussian case there is a deviation
both for small and large values of $\gl$: the deviation for large values is due to
the saturation effect explained in the text. Given the fact that $h_{\text{sat}}$, cf. \eqref{eq:sat},
behaves almost surely and to leading order for $N \to \infty$ as $\sqrt{\log N}$ one
understand  why the slow disappearing of the saturation effect has to be expected.
In both graphs the dotted line above the axis is at level $0.01$.
The fitted values for $\hat m_{N, \go}$, $N=  2.5\cdot 10^5$, are $0.821$ in the binary case
and $0.778$ in the Gaussian case.
}
\end{figure}

What we have observed numerically, see Figures~\ref{fig:sec5_1} and \ref{fig:sec5_2},
may be summed up by
the statement
\begin{equation}
\label{eq:cnj}
\text{there exists } m \text{ such that }
\hat h _{N, \go} (\gl) \, \approx h^{(m)}(\gl).
\end{equation}
Practically this means that $\hat h _{N, \go} (\gl)$, for a set of $\gl$
ranging from $0.05$ to $4$, may be fitted with remarkable precision
by the one parameter family of functions $\left\{ h^{(m)}(\cdot)\right\}_m$.
The fitting value of $m=: \hat m_{N, \go}$ does depend on $N$ and it is essentially increasing.
This is of course expected since localization requires a sufficiently large
system (recall in particular Table~\ref{tbl:2} and Fig.~\ref{fig:NCA} -- see the caption of
Fig.~\ref{fig:sec5_1} for the fitting criterion).
We stress that we are presenting results that have been obtained for one
fixed sequence of $\go$: based on what we have observed for example
in Section~\ref{sec:superadd} for different values of $\gl$ one does expect that
for smaller values of $\gl$ one should use larger values of $N$, but
 changing $N$ corresponds to selecting
a longer, or shorter, stretch of $\go$, that is a different sequence of charges
and this may have a rather strong effect on the value of $\hat m_{N, \go}$.
Moreover there is the problem of deciding which $\gl$-dependence to choose. This may explain
the deviations from \eqref{eq:cnj} that are observed for small values of $\gl$, but these
are in any case rather moderate (see Fig.~\ref{fig:sec5_2}).
\smallskip

A source of stronger (and unavoidable) deviations arises in the cases of unbounded charges:
of course if
\begin{equation}
\label{eq:sat}
h\, \ge \,  h_{\text{sat}} \, :=\, \max_{n\in \{1, \ldots, N \}} \left( -(\go_{2n-1}+ \go_{2n})/2 \right),
\end{equation}
then
$Z^{\gl, h}_{2N, \go} (0) <1$, regardless of the value of $\gl$.
Moreover it is immediate to verify that $\lim_{\gl \to \infty } Z^{\gl, h}_{2N, \go}(0) =+\infty$
for $h<  h_{\text{sat}}$ and therefore $\hat h _{N, \go} (\gl)\nearrow h_{\text{sat}} $
as $\gl \nearrow \infty$. We refer to the captions of Fig.~\ref{fig:sec5_2} for
more on this saturation effect.
 \smallskip

We have tried also alternative definitions of $\hat h _{N, \go} (\gl)$, namely:
\begin{enumerate}
\item the value of $h$ such that $ Z_{2N, \go}^{\gl, h} =1$ (or a different fixed value);
\item the value of $h$ such that the $\ell_1$ distance between
the distribution of the endpoint and the distribution of the meander, cf. Section~\ref{sec:path}, is smaller than a fixed threshold, for example $0.05$.
\end{enumerate}

What we have observed is that \eqref{eq:cnj} still holds.
What is not independent of the criterion is
 $\hat m_{N, \go}$.
 Of course believing deeply in \eqref{eq:cnj} entails
 the expectation  that $\hat m_{N, \go}$ converges to the non random
 quantity $h ^\prime_c (0)$.
  The results reported in this section suggest a value of $h ^\prime_c (0) $ larger than $0.83$
 and the cases presented in Section~\ref{sec:path} suggest that it should be smaller than $0.86$.


\smallskip
\section{Appendix}

\smallskip
\subsection{The algorithm for computing $Z_{N,\go}$}
\label{app:algo}

We are going to briefly illustrate the algorithm we used in the numerical computation of the partition function $Z_{N}=Z_{N,\go}^{\gl,h}$. We recall its definition:
\begin{equation} \label{eq:appZ}
Z_N = \bE \Bigg[ \exp \Bigg( -2\gl \sum_{n=1}^N (\go_n + h) \gD_n \Bigg) \Bigg]\,,
\end{equation}
where $\gD_n := (1-\sign(S_n))/2$ and the convention for $\sign(0)$ described in the introduction.

Observe that a direct computation of $Z_N$ from (\ref{eq:appZ}) would require to sum the contributions of $2^N$ random walk trajectories, making the problem numerically intractable.
 However, here we can make profitably use of the \textsl{additivity} of our Hamiltonian:
 loosely speaking, if we join together two (finite) random walk segments, the energy of the resulting path is the sum of the energies of the building segments.

We can exploit this fact to derive a simple recurrence relation for the sequence of functions $\big\{ \cZ _{M}(y):=Z_{2M}(2y),\ y \in \Z\big\}_{M\in\N}$, where $Z_N(x) = Z_{N,\go}^{\gl,h}(x)$, the latter
defined in \eqref{eq:pinned},
and we recall that we work  with even values of~$N$. Conditioning on $S_{2M}$ and using the Markov property one easily finds
\begin{equation} \label{eq:apprec}
\cZ_{M+1} (y) = \begin{cases}
\frac14 \cZ_{M}(y+1) \;+\; \frac12 \cZ_{M}(y) \;+\; \frac14 \cZ_{M}(y-1) & y>0 \\
\frac14 \Big[ \cZ_{M}(1) + \cZ_{M}(0) \Big] \;+\; \frac14 \ga_M \Big[ \cZ_{M}(0) + \cZ_{M}(-1) \Big] & y=0 \\
\ga_M \Big[ \frac14 \cZ_{M}(y+1) \;+\; \frac12 \cZ_{M}(y) \;+\; \frac14 \cZ_{M}(y-1) \Big] & y<0
\end{cases} \;,
\end{equation}
where we have put $\ga_M := \exp\big(-2\gl\,(\go_{2M+1} + \go_{2M+2} + 2h)\big)$.

From equation (\ref{eq:apprec}) and from the trivial observation that $\cZ_{M}(y)=0$ for $|y|>M$, it follows that $\{\cZ_{M+1}(y),\ y \in \Z\}$ can be obtained from $\{\cZ_{M}(y),\ y \in \Z\}$ with $O(M)$ computations. This means that we can compute $Z_N$  in $O(N^2)$ steps.\footnote{The algorithm just described can be implemented in a standard way: the code we used, written in  C,
 is available on the web page:
 {\tt http://www.proba.jussieu.fr/pageperso/giacomin/C/prog.html}. Graphic representations and standard statistical procedures have been performed
with R \cite{cf:R}. }

\medskip

We point out that sometimes one is satisfied with \textsl{lower bounds} on $Z_N$, for instance in the statistical text for localization described in Section~\ref{sec:superadd}. In this case the algorithm can be further speeded up by restricting the computation to a suitable set of random walk trajectories. In fact when the system size is~$N$ the polymer is at most at distance $O(\sqrt{N})$ (we recall the discussion in Section~\ref{sec:path} on the path behavior), hence a natural choice to get a lower bound on $Z_N$ is to only take into account the contribution coming from those random walk paths $\{s_n\}_{n\in\N}$ for which
\begin{equation}
-A\sqrt{n} \le s_n \le B \sqrt{n} \qquad \text{for } n \ge N_0\,,
\end{equation}
where $A,B,N_0$ are positive constants. Observe that this is easily implemented in the algorithm described above: it suffices to apply relation (\ref{eq:apprec}) only for $y\in[-A\sqrt M, B\sqrt M]$, while setting $\cZ_{M+1}(y)=0$ for the other values of~$y$. In this way the number of computations needed to obtain $Z_N$ is reduced to~$O(N^{3/2})$.

The specific values of $A,B,N_0$ we used in our numerical computations are $3,8,1000$, and we would like to stress that the lower bound on $Z_N$ we got coincides up to the $8^{\text{th}}$ decimal digit with the {\sl true value} obtained applying the complete algorithm.


\smallskip
A final important remark is that for the results we have reported we have used
the Mersenne--Twister~\cite{cf:MT} pseudo--random number generator.
However we have also tried other pseudo--random number generators
and {\sl true randomness} from {\tt www.random.org}:
the results appear not to depend on the generator.


\smallskip
\subsection{Proof of the lower bound on $h_c$}
\label{app:prooflb}
We are going to give a detailed proof of the lower bound \eqref{eq:mainBG} on the critical curve, together with some related result. We stress that this appendix can be made substantially lighter
if one is interested only in the {\sl if} part of Proposition~\ref{prop:stopping1}. In this case
the first part of this appendix is already contained in the first part of \S~\ref{sec:lb_outline}, up
to \eqref{eq:forappB}, and it suffices to look at the proof of the lower bound
starting from page~\pageref{sec:appB2}.

We recall that $Z_{N,\go}^{\gl,h}(0)$ is the partition function corresponding to the polymer pinned at its right endpoint, see \eqref{eq:pinned}, and $T^C=T^{C}(\go)$ is the first $N$ for which $Z_{N,\go}(0)\ge C$, see \eqref{eq:T^C}. In particular, for all $\go$ such that $T^C(\go) < \infty$ we have
\begin{equation} \label{eq:stopping_maj}
Z_{T^C(\go),\go}^{\gl,h}(0) \ge C\,.
\end{equation}
We will also denote by $\cF_n := \gs(\go_1,\ldots,\go_n)$ the natural filtration of the sequence $\{\go_n\}_{n\in\N}$.

\medskip

\subsubsection{A different look at (de)localization}
We want to show that (de)localization can be read from $T^C$. We introduce some notation: given an increasing, $2\N$--valued sequence $\{t_i\}_{i\in\N}$, we set $t_0:=0$ and $\zeta_N:=\max\{k: t_k \le N\}$. Then we  define
\begin{align} \label{eq:low_b}
\begin{split}
\widehat{Z}_{N,\go}(0) = \widehat{Z}_{N,\go}^{\{t_i\},\gl,h}(0) & \;:=\; \bE \left[ e^{-2 \gl \sum_{n=1}^N \left( \go_n +h\right) \Delta_n}; \, S_{t_1}=0,\, \ldots,\, S_{t_{\zeta_N}}=0,\, S_N= 0 \right] \\
&\;=\; \prod_{i=0}^{\zeta_N-1} Z_{t_{i+1}-t_i,\theta^{t_i}\go}^{\gl,h}(0) \,\cdot\, Z_{N-t_{\zeta_N}(\go),\theta^{t_{\zeta_N}}\go}^{\gl,h}(0)\,,
\end{split}
\end{align}
and we recall that $\theta$ denotes the translation on the environment.
One sees immediately that $\widehat{Z}_{N,\go}(0)\le Z_{N,\go}(0)$.
We first establish a preliminary result.

\begin{lemma} \label{lem:chop}
If the sequence $\{t_i\}_i$ is such that $\zeta_N/N \to 0$ as $N\to\infty$, then
\begin{equation}
\lim_{N\to\infty} \frac 1N \log \widehat{Z}_{N,\go}^{\{t_i\},\gl,h}(0) \;=\; \tf(\gl,h)\,,
\end{equation}
both $\bbP(\dd\go)$--a.s. and in $\bbL_1(\bbP)$.
\end{lemma}

\proof
By definition we have $Z_{N,\go}(0) \ge \widehat{Z}_{N,\go}(0)$. On the other hand, we are going to show that
\begin{equation} \label{eq:rough_maj}
Z_{N,\go}^{\gl,h}(0) \;\le\; 4^{\zeta_N} \, A^{2\zeta_N} \, \left( \prod_{i=1}^{\zeta_N} (t_i-t_{i-1}) \cdot (N-t_{\zeta_N}) \right)^{3} \,  \widehat{Z}_{N,\go}^{\{t_i\},\gl,h}(0)\,,
\end{equation}
where $A$ is a positive constant. To derive this bound, we resort to the equation \eqref{eq:reducttoexc} that expresses $Z_{N,\go}(0)$ in terms of random walk excursions. We recall that $K(2n)$ is the discrete probability density of the first return time of the walk $S$ to $0$, and that $K(t) \ge 1/(A\,t^{3/2})$, $t\in 2\N$, for some positive constant~$A$: it follows that for $a_1, \ldots, a_k \in 2\N$
\begin{equation} \label{eq:app_bound}
K(a_1 + \ldots + a_k) \; \le \; 1 \; \le \; A^k \, (a_1 \cdot \ldots \cdot a_k)^{3/2} \, K(a_1) \cdot \ldots \cdot K(a_k)\,.
\end{equation}
This gives us an upper bound to the entropic cost needed to split a random walk excursion of length $(a_1 + \ldots + a_k)$ into $k$~excursions of lengths $a_1, \ldots, a_k$.

Now let us come back to the second line of \eqref{eq:reducttoexc}, that can be rewritten as
\begin{equation}
\label{eq:app_sum}
Z_{N,\go}(0) \,=\, \sum_{\{x_i\} \subseteq \{0, \ldots, N\} \cap 2\N} G(\{x_i\})\,.
\end{equation}
A first observation is that if we restrict the above sum  to the $\{x_i\}$ such that
$\{x_i\} \supseteq \{t_i\}$, then we  get  $\widehat{Z}_{N,\go}^{\{t_i\}}(0)$. Now for each $\{x_i\}$ we
aim at finding an  upper bound on the term $G(\{x_i\})$ of the form
 $c \cdot G(\{x_i\} \cup \{t_i\})$ for some $c>0$ not depending on $\{x_i\}$.
 Each term $G(\{x_i\})$, see \eqref{eq:reducttoexc}, is the product of two terms: an entropic part depending on~$K(\cdot)$ and an energetic part depending on~$\varphi(\cdot)$. Replacing the entropic part costs no more than
\begin{equation}
    c_{\text{ent}} \, :=\,  A^{2\zeta_N} \, \left( \prod_{i=1}^{\zeta_N} (t_i-t_{i-1}) \cdot (N-t_{\zeta_N}) \right)^{3}\,,
\end{equation}
thanks to~\eqref{eq:app_bound}. On the other hand, the cost for replacing the energetic part is easily
bounded above by
\begin{equation}
    c_{\text{energy}} \, :=\, 2^{\zeta_N}\,,
\end{equation}
so that the bound $G(\{x_i\}) \le c \cdot G(\{x_i\} \cup \{t_i\})$ holds true with $c:= c_{\text{ent}} \, c_{\text{energy}}$. Replacing in this way each term in the sum in the r.h.s. of \eqref{eq:app_sum}, we are left with a sum of terms $G(\{y_i\})$ corresponding to sets $\{y_i\}$ such that $\{y_i\} \supseteq \{t_i\}$. It remains to count the {\sl multiplicity} of any such~$\{y_i\}$, that is how many original sets $\{x_i\}$ are such that $\{x_i\} \cup \{t_i\} = \{y_i\}$. Sets $\{x_{i}\}$ satisfying this last condition must differ only for a subset of $\{t_{i}\}$, hence the sought multiplicity is $2^{\zeta_N}$ (the cardinality of the parts of $\{t_{i}\}$) and the bound \eqref{eq:rough_maj} follows.

Therefore we get
\begin{align}
\bigg| \frac{\log \widehat{Z}_{N,\go}^{\{t_i\},\gl,h}(0)}{N}  -  \frac{\log Z_{N,\go}^{\gl,h}(0)}{N} \bigg| & \;\le\; (2\log 2 A) \frac{\zeta_N}{N} \;+\; 3\, \frac{1}{N} \,\log \left( \prod_{i=1}^{\zeta_N} (t_i-t_{i-1}) \cdot (N-t_{\zeta_N}) \right) \\
&\;\le\; (2\log 2A) \frac{\zeta_N}{N} \;+\; 3 \, \frac{\zeta_N+1}{N} \, \log \bigg(\frac{N}{\zeta_N+1}\bigg)\,,
\nonumber
\end{align}
where in the second inequality we have made use of the elementary fact that once the sum of $k$ positive numbers is fixed, their product is maximal when all the numbers coincide (for us $k=\zeta_N+1$). Since by hypothesis $\zeta_N/N \to 0$ as $N\to\infty$, the Lemma is proved.\qed

\bigskip

Now we are ready to prove the characterization of $\cL$ and $\cD$ in terms of $T^C$.
Fix any $C>1$.

\smallskip

\begin{proposition} \label{prop:stopping1}
A point $(\gl,h)$ is localized, that is $h < h_c(\gl)$, if and only if $\bbE[T^C]<\infty$.
\end{proposition}

\proof
We set $\cA:=\{\go: T^C(\go)<\infty\}$. Observe that for $\go\in \cA^{\complement}$ we have $Z_{N,\go}(0) \le C$ for every $N\in2\N$, and consequently $\log Z_{N,\go}^{\gl,h}(0)/N \to 0$ as $N\to\infty$.

Consider first the case when the random variable $T^C$ is defective, that is $\bbP[\cA^{\complement}]>0$ (this is a particular case of $\bbE[T^C]=\infty$). Since we know that $\log Z_{N,\go}^{\gl,h}(0)/N \to \tf(\gl,h)$, $\bbP(\dd\go)$--a.s., from the preceding observation it follows that $\tf(\gl,h)=0$ and the Proposition is proved in this case.

\smallskip

Therefore in the following we can assume that $T^C$ is proper, that is $\bbP(\cA)=1$, so that equation (\ref{eq:stopping_maj}) holds for almost every~$\go$. Setting $\theta^{-1}\cA := \{\go : \theta\go \in \cA\}$, we have $\bbP\left(\theta^{-1}\cA\right)=1$ since $\bbP$ is $\theta$--invariant, and consequently $\bbP\left(\cap_{k=0}^\infty \theta^{-k}\cA
 \right)=1$, which amounts to saying that (\ref{eq:stopping_maj}) can be actually strengthened to
\begin{equation} \label{eq:im_stopping_maj}
Z_{T^C(\theta^k\go),\theta^k\go}^{\gl,h}(0) \ge C \qquad \forall k \ge 0,\ \bbP(\dd\go)\text{--a.s.}\,.
\end{equation}
Observe that the sequence $\{( \theta^{T^C(\go)} \go )_n\}_{n\in\N}$ has the same law as $\{\go_n\}_{n\in\N}$ and it is independent of $\cF_{T^C}$. We can define inductively an increasing sequence of stopping times $\{T_n\}_{n\in\N}$ by setting $T_0:=0$ and $T_{k+1}(\go) - T_k(\go) := T^C(\theta^{T_k(\go)}\go) =: S_k(\go)$. We also set $\zeta_N(\go) := \max \{n:\ T_n(\go) \le N\}$. Since $\{S_k\}_{k\in\N}$ is an IID sequence, by the strong law of large numbers we have that, $\bbP(\dd\go)$--a.s., $T_n(\go)/n \to \bbE[T^C]$ as $n\to\infty$, and consequently $\zeta_N(\go)/N \to 1/\bbE[T^C]$ as $N\to\infty$ (with the convention that $1/\infty=0$).

Now let us consider the lower bound $\widehat{Z}_{N,\go}(0)$ corresponding to the sequence $\{t_i\}=\{T_i(\go)\}$: from \eqref{eq:low_b} and \eqref{eq:im_stopping_maj} we get that $\bbP(\dd\go)$--a.s.
\begin{align} \label{eq:major}
\begin{split}
\widehat{Z}_{N,\go}^{\{T_i(\go)\},\gl,h}(0) &\;=\; \prod_{i=0}^{\zeta_N(\go)-1} Z_{T^C(\theta^{T_i}\go),\theta^{T_i}\go}^{\gl,h}(0) \,\cdot\, Z_{N-T_{\zeta_N(\go)}(\go),\theta^{T_{\zeta_N(\go)}}\go}^{\gl,h}(0) \\
&\;\ge\; C^{\zeta_N(\go)} \cdot \frac{c}{N^{3/2}}\,,
\end{split}
\end{align}
where $c$ is a positive constant (to estimate the last term we have used the  lower bound $Z_k(0) \ge c/k^{3/2}$, cf. \eqref{eq:step_deloc}), and consequently
\begin{equation}
\tf(\gl,h) \;=\; \lim_{N\to\infty} \frac{\log Z_{N,\go}^{\gl,h}(0)}{N} \;\ge\; \liminf_{N\to\infty} \frac{\log \widehat{Z}_{N,\go}^{\{T_i(\go)\},\gl,h}(0)}{N} \;\ge\; \frac{\log C}{\bbE[T^C]}\,.
\end{equation}
It follows that if $\bbE[T^C]<\infty$  then $\tf(\gl,h)>0$, that is $(\gl,h)$ is localized.

\smallskip

It remains to consider the case $\bbE[T^C]=\infty$, and we want to show that this time $\widehat{Z}_{N,\go}(0)$, defined in \eqref{eq:major}, gives a null free energy. In fact, as $T^C(\eta)$ is defined as the \textsl{first} $N$ such that $Z_{N,\eta}(0) \ge C$, it follows that $Z_{T^C(\eta),\eta}(0)$ cannot be much greater than $C$. More precisely, one has that
\begin{equation}
Z_{T^C(\eta),\eta}(0) \le  C \,\exp(2\gl|\eta_{T^C(\eta)-1} + \eta_{T^C(\eta)}|)\,,
\end{equation}
and from the first line of \eqref{eq:major} it follows that
\begin{equation}
\frac 1N \log \widehat{Z}_{N,\go}(0) \;\le\; \frac{\zeta_N(\go)+1}{N} \log C \;+\;  \frac{2\gl}{N} \sum_{i=1}^{\zeta_N(\go)} \Big( |\go_{T_i(\go)}| + |\go_{T_i(\go)-1}| \Big)\,.
\end{equation}
We estimate the second term in the r.h.s. in the following way:
\begin{align}
    & \frac 1N \sum_{i=1}^{\zeta_N(\go)} \Big( |\go_{T_i(\go)}| + |\go_{T_i(\go)-1}| \Big) = \frac 1N \sum_{k=1}^{N} \ind_{\{\exists i:\, T_i(\go) = k\}} \Big( |\go_k| + |\go_{k-1}| \Big)\nonumber  \\
    & \qquad \le \left( \frac 1N \sum_{k=1}^{N} \ind_{\{\exists i:\, T_i(\go) = k\}} \right) ^{1/2} \left( \frac 1N \sum_{k=1}^{N} \Big( |\go_k| + |\go_{k-1}| \Big)^2 \right)^{1/2}\\
    &\qquad \le \sqrt{\frac{\zeta_N(\go)}{N}} \cdot 2 \sqrt{ \frac 1N \sum_{k=1}^{N} |\go_k|^2} \le A \sqrt{\frac{\zeta_N(\go)}{N}} \,, \nonumber
\end{align}
for some positive constant~$A=A(\go)$ and  eventually as $N \to \infty$, having used the Cauchy--Schwartz inequality and the law of large numbers for the sequence $\{|\go_k|^2\}_{k\in\N}$. Therefore
\begin{equation}
    \frac 1N \log \widehat{Z}_{N,\go}(0) \;\le\; \frac{\zeta_N(\go)+1}{N} \log C \;+\; 4\gl A \sqrt{\frac{\zeta_N(\go)}{N}}\,,
\end{equation}
and since $\bbE[T^C]=\infty$ implies $\zeta_N(\go)/N \to 0$, $\bbP(\dd\go)$--a.s., we have $\log \widehat{Z}_{N,\go}(0) /N \to 0$, $\bbP(\dd\go)$--a.s.. Then Lemma~\ref{lem:chop} allows us to conclude that $\tf (\gl,h)=0$, and the proof of the Proposition is completed.\qed


\medskip

\subsubsection{Proof of the lower bound on $h_c$}
\label{sec:appB2}
To prove equation \eqref{eq:mainBG}, we are going to build, for every $(\gl,h)$ such that $h < \underline h(\gl)$, a random time $T$ such that $\bbE[T]<\infty$ and $Z_{T(\go),\go}^{\gl,h}(0) \ge C$, for some $C>1$. It follows that $T^C \le T$, yielding that $\bbE[T^C]<\infty$ and by Proposition~\ref{prop:stopping1} $(\gl,h)$ is localized, that is, $\underline h(\gl) \le h_c(\gl)$.

\smallskip

Given $M\in 2\N$
and $q<-h$, we start defining the stopping time
\begin{equation}
\tau_M(\go) = \tau_{M,q}(\go) := \inf \bigg\{n\in2\N:\ \exists k \in 2\N,\ k \ge M:\ \frac{\sum_{i=n-k+1}^n \go_i}{k} \le q \bigg\}\,.
\end{equation}
This is the first instant at which a $q$--atypical stretch of length at least $M$ appears along the sequence $\go$. The asymptotic behavior of $\tau_M$ is given by Theorem 3.2.1 in \cite[\S~3.2]{cf:DZ} which says that $\bbP(\dd\go)$--a.s.
\begin{equation} \label{eq:as_tau}
\frac{\log \tau_{M}(\go)}{M} \to \gS(q) \qquad \text{as } M\to\infty\,,
\end{equation}
where $\gS(q)$ is Cramer's Large Deviations functional for $\go$, \eqref{eq:cramer}.
We also give a name to the shortest of the terminal stretches in the definition of $\tau_M$:
\begin{equation}
R_{M}(\go) = R_{M,q}(\go) := \inf \bigg\{k \in 2\N,\ k \ge M:\ \frac{\sum_{i=\tau_M-k+1}^{\tau_M} \go_i}{k} \le q \bigg\}\,,
\end{equation}
and it is not difficult to realize that $R_M \le 2M$.

\smallskip

We are ready to give a simple lower bound on the partition function of size $\tau_{M,q}$ (for any $M\in 2\N$ and $q<-h$): it suffices to consider the contribution of the trajectories that are negative in correspondence of the last (favorable) stretch of size $R_M$, and stay positive the rest of the time. Recalling that
we use $K(\cdot)$ for the discrete density of the first return time to the origin and that by
\eqref{eq:asympt} we have $K(2n)\ge c/n^{-3/2}$ for a constant $c>0$,
we estimate
\begin{equation}\label{eq:needs}
\begin{split}
Z_{\tau_M(\go), \go}^{\gl,h}(0) & \ge \frac 1 4 \, K\left({\tau_M - R_M}\right) \, K\left(
{R_M}\right) \, e^{-2\gl (q+h) R_M} \ge \frac{c^2}{4 \tau_M^{3/2} (2M)^{3/2}} e^{-2\gl (q+h) M} \\
& \ge c'\, \exp \bigg\{ \frac 3 2 M \bigg[ (-4\gl /3)q - \frac{\log \tau_M}{M} - ( 4 \gl/3) h - \frac{\log M}{M} \bigg] \bigg\} \,,
\end{split}
\end{equation}
where $c' := c^2/(8\sqrt{2})$.

Having in mind (\ref{eq:as_tau}), we define a random index $\ell=\ell_{A,\gep,q}$ depending on the two parameters $A \in 2\N,\ \gep>0$ and on $q$:
\begin{equation} \label{eq:def_ell}
\ell(\go)  = \ell_{A,\gep,q}(\go) := \inf \bigg\{ k\in 2 \N,\ k \ge A:\ \frac{\log \tau_{k,q}(\go)}{k} \le \gS(q) + \gep \bigg\} \,,
\end{equation}
and we finally set
\begin{equation}\label{eq:def_T_cgg}
    T (\go) = T_{A,\gep,q} (\go) := \tau_{\ell(\go)}(\go)\,.
\end{equation}
Then for the partition function of size $T(\go)$ we get
\begin{equation} \label{eq:almost_maj}
Z_{T(\go), \go}^{\gl,h}(0) \ge c'\, \exp \bigg\{ \frac 3 2 A \bigg[ (-4\gl /3)q - \gS(q) - (4\gl /3)h - \frac{\log A}{A} - \gep \bigg] \bigg\} \,.
\end{equation}

The fact that $\bbE[T_{A,\gep,q}]<\infty$ for any choice of $A,\gep,q$ (with $q<-h$) is proved in Lemma~\ref{lem:stopping2} below. It only remains to show that for every fixed $(\gl,h)$ such that $h<\underline h(\gl)$, or equivalently
\begin{equation}
\label{eq:cond_low}
(4\gl/3)h < \log \M(-4\gl/3)\,,
\end{equation}
the parameters $A,\gep,q$ can be chosen such that the right--hand side
 of equation (\ref{eq:almost_maj}) is greater than~$1$.

\smallskip

The key point is the choice of  $q$. Note that
the generating function $ \M(\cdot)$ is smooth, since  finite on the whole real line.
 Moreover for all $ \gl \in\R$ there exists some $q_0\in\R$ such that
\begin{equation}
\log \M(-4\gl/3) = (-4\gl/3) q_0 - \gS(q_0)\,,
\end{equation}
and from (\ref{eq:cond_low}) it follows that $q_0 < -h$. Therefore we can take $q=q_0$, and equation (\ref{eq:almost_maj}) becomes
\begin{equation}
\label{eq:lbonZTom}
Z_{T(\go), \go}^{\gl,h}(0) \ge c'\, \exp \bigg\{ \frac 3 2 A \bigg[ \log \M(-4\gl/3) - (4\gl/3) h - \frac{\log A}{A} - \gep \bigg] \bigg\} \,.
\end{equation}
It is now clear that for every $(\gl,h)$, such that (\ref{eq:cond_low}) holds, by choosing
$\gep$ sufficiently small and $A$ sufficiently large, the right--hand side of \eqref{eq:lbonZTom} is greater than~1, and the proof of \eqref{eq:mainBG} is complete.

\medskip

\begin{lemma}
\label{lem:stopping2}
For every $A\in 2\N$,  $\gep > 0$ and $ q <-h$
the random variable $T(\go)=T_{A,\gep,q}(\go)$ defined
below \eqref{eq:def_ell} is integrable: $\bbE[T]<\infty$.
\end{lemma}

\medskip

\proof
By the definition (\ref{eq:def_ell}) of $\ell=\ell_{A,\gep,q}$ we have
\begin{equation}
T_{A,\gep,q} \le \exp\big( (\gS(q) + \gep)\, \ell_{A,\gep,q} \big)\,,
\end{equation}
so it suffices to show that for any $\gb>0$ the random variable $\exp(\gb\, \ell_{A,\gep,q})$ is integrable.

For any $l\in2\N$, we introduce the IID sequence of random variables $\{Y_n^{l}\}_{n\in\N}$ defined by
\begin{equation}
Y^{l}_n := \frac 1 {l}{\sum_{i=(n-1)l+1}^{n l} \go_i}\,.
\end{equation}
By Cramer's Theorem \cite{cf:DZ} we have that for any fixed $q<0$ and $ \gep>0$ there exists  $l_0$ such that $\bbP\left(Y_1^{l} \le q\right) \ge e^{-l (\gS(q) + \gep/2)}$ for every $l\ge l_0$.
By (\ref{eq:def_ell}) have that
\begin{equation}
\begin{split}
\{\ell > l\} \subseteq \left\{\tau_l > \exp((\gS(q)+\gep) l)\right\}
\subseteq \inter_{i=1}^{\lfloor M/l \rfloor} \{Y_i^{l} > q\}\,,
\end{split}\end{equation}
with $M:=\exp((\gS(q)+\gep) l)$, so that
\begin{equation}
\begin{split}
\bbP\left(\ell > l\right)
  \le \left( 1- e^{-l (\gS(q) + \gep/2)} \right)^{\lfloor M/l \rfloor}
  &\le
  \exp \left( - \lfloor M/l \rfloor e^{-l (\gS(q) + \gep/2)} \right)
  \\
& \le  \exp \left(- \exp\left(l \gep/4\right)\right),
\end{split}
\end{equation}
where the last step holds if $l$ is sufficiently large   (we have also used
 $1-x \le e^{-x}$). Therefore
\begin{equation}
\bbP\left(\exp(\gb\,\ell) > N\right) = \bbP\left(\ell > (\log N) /\gb \right)
\le  \exp \left( - N ^{\gep / 4\gb }\right),
\end{equation}
when $N$ is large, and the proof is complete. \qed


\chapter[On improving the annealed bound]
{On improving the annealed bound\\ for polymer chains with random charges}
\label{ch:cg}

In this chapter we address the issue of improving the anneal bound
on the critical line~$h_c(\cdot)$ of the random copolymer
via the so--called constrained annealing, that means nothing but
applying the annealing procedure (which is just Jensen's inequality)
after having added to the Hamiltonian a disorder--dependent term
(sometimes interpreted as a Lagrange multiplier) in a way that
the quenched expressions are left unchanged, remember \S~\ref{sec:ph_random}
of Chapter~\ref{ch:first}.

A popular class of multipliers is the one consisting of empirical
averages of local functions of the disorder.
These multipliers are particularly suitable for computations,
and it is often believed that in this
class one can  approximate arbitrarily well the quenched free energy.

We are going to prove that {\sl this is not the case}
for a wide family of polymer models, including
the copolymer near a selective interface and the pinning/wetting models
defined in Chapter~\ref{ch:first}.
More precisely we show that the multipliers in the above class
cannot improve on the basic annealed bound from the viewpoint
of characterizing the phase diagram. For simplicity the proof has
been carried out under the assumption that the random variable~$\go_1$
takes only a finite number of values, however the statement remains true
also the general case, provided one makes some suitable boundedness
assumptions on the multiplier.

\smallskip
The article~\cite{cf:CG} has been taken from the content of this chapter.

\smallskip
\section{The framework and the main result}
\label{sec:intro_CG}

\smallskip
\subsection{The general set--up}
A number of disordered models of linear chains undergoing
localization or pinning effects can be put into the following general framework.
Let  $S:=\left\{ S_n \right\}_{n=0,1,\ldots}$ be a process with $S_n$ taking
values in $\Z^d$, $d \in \N :=\{1,2,\ldots\}$ and
 law $\bP$.

The disorder in the system is given by a sequence $\go :=\left\{\go _n\right\}_n$
of IID random variables taking values in a finite set $\Gamma$ with law $\bbP$, acting on the path of $S$ via an Hamiltonian
that, for a system of size $N$, is a function $H_{N,\go}$
of the trajectory $S$, but depending only on $S_0, S_1, \ldots, S_N$.
One is interested in the properties of the probability measures
$\bP_{N, \go}$ defined by giving the density with respect to $\bP$:
\begin{equation}
\label{eq:1}
\frac{\dd \bP_{N, \go}}{\dd\bP} \left( S \right) \, =\,
\frac1{Z_{N, \go}} \exp\left(H_{N, \go}\left(S\right)\right),
\end{equation}
where $Z_{N, \go}:= \bE \left[\exp\left(H_{N, \go}\left(S\right)\right)\right]$ is the normalization constant.
Our attention focuses on the
asymptotic behavior   of $\log Z_{N, \go}$.

In the sequel we will assume:

\medskip
\begin{assumption} \rm
There exists a sequence $\left\{ D_n \right\}_n$
of subsets of $\Z^d$ such that
$\bP ( S_n \in D_n \text{ for } n=1,2, \ldots,  N)\stackrel{N\to \infty}{\asymp} 1 $, namely
\begin{equation}
\lim_{N\to \infty} \frac 1N \log \bP \left( S_n \in D_n \text{ for } n=1,2, \ldots, N\right)=0,
\end{equation}
and $H_{N, \go}(S)=0$ if $S_n \in D_n$ for $n=1, 2, \ldots, N$.
\end{assumption}

\vskip 0.1 cm

One sees directly that this hypothesis implies
\begin{equation}
\label{eq:first}
\liminf_{N \to \infty}
\frac 1N \log Z_{N, \go} \, \ge
\lim_{N \to \infty}
\frac 1N \log \bP \left( S_n \in D_n \text{ for } n=1,2, \ldots, N\right)\, =\, 0,
\end{equation}
$\bbP (\dd \go)$--a.s.. We will assume that $\left\{(1/N)\log Z_{N, \go}\right\}_N$
is a sequence of integrable random variables that converges
 in the $L^1\left( \bbP (\dd \go)\right)$ sense and
$P(\dd \go)$--almost surely to a constant, the free energy, that we will call $\free$.
These assumptions are verified in the large majority of the interesting situations,
for example whenever super/sub--additivity tools are applicable.

Of course \eqref{eq:first} says that $\free \ge 0$ and one is lead to
the natural question of whether $\free =0$ or $\free >0$.
In the instances that we are going to consider the free energy may be
zero or positive according to some parameters from which the $H_{N, \go}$
depends: $\free =0$ and $\free >0$ are associated to sharply different
behaviors of the system.

\medskip

In order to establish upper bounds on $\free$
one may apply directly Jensen inequality ({\sl annealed bound}) obtaining
\begin{equation}
\label{eq:J}
\free  \, \le \, \liminf_{N \to \infty} \frac 1N \log \bbE \left[ Z_{N,\go}
\right]\, =: \, \freea,
\end{equation}
and,  in our context,
  if $\freea=0$ then $\free =0$.
The annealed bound
 may be improved by adding to
$H_{N, \go}(S)$ an integrable function $A_N: \Gamma^\N \to \R$ such that
$\bbE\left[ A_N (\go) \right]=0$: while the left--hand side is unchanged,
$\freea$ may depend on the choice of $\{A_N\}_N$.
We stress that not only $\free$ is left
unchanged by $H_{N, \go}(S) \to H_{N, \go}(S)+A_N (\go)$, but $\bP_{N,\go}$
itself is left unchanged (for every $N$).
Notice that the choice
$A_N (\go) = -\log Z_{N,\go} + \bbE\left[  \log Z_{N,\go}\right]$
yields the equality in \eqref{eq:J}.

In the sequel when we refer to $\freea$ we mean that  $Z_{N,\go}$ is defined with respect to $H_{N, \go}$
satisfying the Basic Hypothesis (no $A_N$ term added).

\smallskip
\subsection{The result}
What we prove in this note is that
\medskip

\begin{proposition}
\label{th:main}
If $\freea >0$ then
for every local function $F: \Gamma ^\N \longrightarrow \R$ such that $\bbE \left[ F (\go) \right]=0$
one has
\begin{equation}
\label{eq:main_CG}
\liminf_{N \to \infty} \frac 1 N \log
\bbE \bE \left[ \exp \left( H_{N, \go}(S) + \sum_{n=0}^N F(\theta_n \go)\right)\right]\, > \, 0,
\end{equation}
where $(\theta_n \go)_m= \go_{n+m}$.
\end{proposition}

\medskip

We can sum up this result by saying that when  $ \free =0$ but $ \freea >0$
it is of no use  modifying  the Hamiltonian by  adding the
empirical average of a (centered) local function.

On a  mathematical level it is clear that we are playing with
an exchange of limits and that it is not obvious that the
free energy, recall the optimal choice of $A_N$ above,
may be approximated via empirical averages of a
local function of the disorder. But we remark
that in the physical literature the approach of approximating
the free energy via what can be viewed
as a constrained annealed computation, the term $\sum_{n=0}^N F(\theta_n \go)$
being interpreted as a Lagrange multiplier,
is often considered as an effective way of approximating the quenched
free energy. Here we mention in particular \cite{cf:Morita}
and \cite{cf:Kuhn} in which this approach is taken up in a systematic way:
the aim is to approach the quenched free energy by constrained
annealing via local functions $F$ that are more and more complex,
the most natural example being linear combinations of correlations of higher and higher
order.
\smallskip

The proof of
Proposition \ref{th:main} is based on the simple observation
that whenever $A_N$ is centered
\begin{multline}
\label{eq:argproof}
 \frac 1 N \log
\bbE \bE\left[ \exp \left( H_{N, \go}(S) + A_N(\go) \right)\right]\, \ge
\\
 \frac 1 N \log
\bbE \left[ \exp \left( A_N(\go) \right)\right] +
\frac 1 N \log \bP \left( S_n \in D_n \text{ for } n=1, 2, \ldots , N
\right)\, =:\, Q_N+P_N.
\end{multline}
By hypothesis $P_N=o(1)$ so one has to consider the asymptotic behavior
of $Q_N$. If $\liminf_N Q_N >0$ there is nothing to prove.
So let us assume that $\liminf_{N}Q_N=0$:
in this case the inferior limit of the left--hand side of \eqref{eq:argproof}
may be  zero and we want to exclude this possibility
when $\tilde f >0$ and $A_N (\go)= \sum_{n=0}^N F(\theta_n \go)$, $F$ local and centered
(of course in this case $\lim_N Q_N$ does exist).
And in Theorem~\ref{th:main1} below  in fact
we show that
if $\log \bbE\left[ \exp\left(A_N(\go)\right)\right] =o(N)$, then
there exists a local function $G$ such that $F( \go)= G(\theta_1 \go)- G(\go)$
so that $ \{\sum_{n=0}^N F(\theta_n \go)\}_N$ is just a boundary term and
 the corresponding constrained annealing is just the standard annealing.

Notice that having chosen $\Gamma$ finite frees us from
integrability conditions.

\medskip

\begin{rem}
\label{th:rem}
\rm
We stress that our Basic Hypothesis is more general than it may look at first.
As already observed, one has  the freedom of adding to the Hamiltonian $H_{N, \go}(S)$
any term  that does not depend on $S$ (but  possibly does depend on $\go$ and $N$)
without changing the model $ \bP_{N,\go}$.
 It may therefore happen
 that the {\sl natural} formulation of the Hamiltonian does not satisfy
 our Basic Hypothesis, but it does after a suitable additive correction.
This  happens for example for the Copolymer near a selective interface model,
as we have seen in~\S~\ref{rem:Z} of Chapter~\ref{ch:first} (see also
\S~\ref{sec:cop} below):
the additive correction in this case is linear in $\go$ and it corresponds to
  what in \cite{cf:ORW} is called {\sl first order} Morita approximation.
  In these terms, Proposition \ref{th:main} is saying that
  {\sl higher order} Morita approximations cannot improve
  the bound on the critical curve found with the first order computation.
\end{rem}

\smallskip

Let us now look at applications of Proposition \ref{th:main}.

\smallskip
\subsection{Random rewards or penalties at the origin}
Let $S$, $S_0=0\in \Z ^d$, be a random walk with centered IID non degenerate increments $\{X_n\}_n$,
$(X_n)_j\in \{-1,0,1\}$ for $j=1, 2, \ldots, d$,
and
\begin{equation}
H_{N, \go} = \gb
  \sum_{n=1}^N \left(1+ \gep \go_n\right) \ind_{\{S_n=0\}}.
\end{equation}
for $\gb \ge 0$ and $\gep \ge 0$. This model is a $d$--dimensional
version (with somewhat different notations) of
the pinning model introduced in~\S~\ref{sec:pinning} of Chapter~\ref{ch:first}.
The random variable $\go_1$ is chosen
such that $\bbE[\exp (\gl\go_1)]<\infty$ for every $\gl\in\R$, and centered.
We write $\free (\gb , \gep)$ for $\free$:
 by super--additive arguments $\free$ exists and it is self--averaging (this
 observation is valid for all the models we consider and will not be repeated).
As we already remarked in Chapter~\ref{ch:first}, for $\gep =0$ the model can be solved, see
e.g. \cite{cf:G}, and in particular $\free (\gb, 0)=0$ if and only if $\gb \le \gb_c (d):= -\log (1-\bP(S$ never
comes back to $0))$. Adding the disorder makes this
model much more complex: the annealed bound yields $\free (\gb ,\gep) =0$
if $ \gb \le \gb_c(d)- \log \bbE \left[\exp (\gep \go_1)\right]
=:\tilde {\gb_c}$.
It is an open question whether $\tilde {\gb_c}$ coincides with
the quenched critical value  or not, that is whether $\free (\gb , \gep) =0$
implies $ \gb \le \tilde {\gb_c}$ or not. For references about this issue
we refer to \cite{cf:AS,cf:Pet},
see however also the next
paragraph:   the model we are considering  can in fact be exactly mapped
to the wetting problem (\cite{cf:AS}, \cite{cf:G}).
Proposition~\ref{th:main} applies to this context with $D_n= \{0\}^\complement$ for every $n$
\cite[Ch. 3]{cf:Feller}
and says that one cannot answer this question via constrained annealed bounds.

\smallskip
\subsection{Wetting models in $1+d$ dimensions}
Let $S$ and $\go$ as in the previous example
and
\begin{equation}
\label{eq:Wetting}
H_{N, \go} =
\begin{cases}
 \gb \sum_{n=1}^N \left( 1 + \gep \go_n\right) \ind_{\{(S_n)_d=0\}}
 &\text{ it } (S_n)_d\ge 0 \text{ for } n=1,2, \ldots , N \\
 -\infty &\text{ otherwise.}
 \end{cases}
\end{equation}
with $\gb \ge 0$ and
$\gep\ge 0$.
If one takes  the directed walk viewpoint,
that is if one considers the walk $\{(n,S_n)\}_n$,
then this is a model of a walk constrained above
the (hyper--)plane $x_d=0$ and rewarded $\gb$, on the average,
when touching this plane. If $d=1$ then this is an effective
model for a (1+1)--dimensional interface above a  wall which mostly attracts it.
As a matter of fact in this case there is no loss of generality
in considering $d=1$, since in the directions parallel to the wall
the model is just the original walk.
Once again if $\gep =0$ the model can be solved in detail,
see e.g. \cite{cf:G}. Computing the critical $\gb$
and deciding whether the annealed bound is sharp, at least
for small $\gep$, is an unresolved and disputed question
in the physical literature, see e.g. \cite{cf:FLNO}, \cite{cf:DHV} and \cite{cf:TC}.
Proposition \ref{th:main} applies with the choice
$D_n =\bbZ ^{d-1}\times \N$.

\smallskip
\subsection{Copolymer and adsorption models}
\label{sec:cop}
Choose $S$ as above and take the directed walk viewpoint.
Imagine that above the axis ($x_d>0$) is filled of a solvent $A$,
while below ($x_d<0$) there is a solvent $B$. At $x_d=0$ there is the
interface. We choose $\go = \{A,B\}$ and for example
\begin{equation} \label{eq:hab}
H^{AB}_{N, \go}(S)\, =\,
\sum_{n=1}^N \left(a \ind_{\{\sign(S_n)=+1,\,  \go_n =A\}} +
b\ind_{\{\sign (S_n) =-1, \, \go_n =B\}}+ c\ind_{\{S_n=0\}}
\right)
\end{equation}
with $a$, $b$ and $c$ real parameters, $\sign (S_n) := \sign \big( (S_n)_d \big)$
and the convention $\sign(S_n)= \sign(S_{n-1})$ if $(S_n)_d=0$
already used in Chapter~\ref{ch:first}.
In order to apply Proposition \ref{th:main}
one has to subtract a disorder dependent term, cf. Remark~\ref{th:rem}: if $a \ge b$
we change the Hamiltonian
\begin{equation}
\label{eq:H1}
H_{N, \go}(S)\, :=\, H^{AB}_{N, \go}(S)-   \sum_{n=1}^N a \ind_{\{ \go_n =A\}}.
\end{equation}
without changing the measure $\bP_{N, \go}$  while the free energy
has the trivial shift
 from $\free$ to $\free - a \bbP \left( \go_1=A\right)$.
One can therefore choose $D_n =\bbZ ^{d-1}\times \N$ and Proposition \ref{th:main}
applies. This model has been considered for example
in \cite{cf:ORW}.

Note that if~$c=0$ and~$d=1$ the model is nothing but the copolymer model introduced
in Chapter~\ref{ch:first}, that is we can cast~\eqref{eq:hab} in the form
\begin{equation}
\label{eq:H2}
H_{N,\go} (S)\, = \, \gl  \sum_{n=1}^N  \left(\go_n +h \right) \sign (S_n) \,,
\end{equation}
with $\go$ taking values in $\R$. Once again the Hamiltonian
has to be corrected by subtracting
the term $\gl \sum_n (\go_n +h)$ (which is exactly what was done in \S~\ref{rem:Z} of
Chapter~\ref{ch:first}) in order to apply Proposition \ref{th:main}.
One readily sees that  \eqref{eq:H1}
and  \eqref{eq:H2} are the same model when in the second case
$\go$ takes only the values $\pm 1$, $A=+1$ and $B=-1$,
and $h= (a-b)/(a+b)$, $\gl= (a+b)/4$.

Proposition \ref{th:main} acquires some interest in this context:
in fact we have already remarked that
the physical literature is rather split on the
precise value of the critical curve and on whether the annealed
bound is sharp or not.
We recall that the numerical analysis performed in Chapter~\ref{ch:cgg}
is suggesting that the annealed curve does not coincide with
the quenched one, and in view of Proposition~\ref{th:main} this would mean that
constrained annealing via local functions cannot capture
the phase diagram of the quenched system.

\smallskip
\subsection{Further models and  observations}
In spite of substantial numerical evidence
that in several  instances $\free =0$ but $\freea>0$,
we are unaware of an {\sl interesting} model for which
this situation is rigorously known to happen.
Consider however the case $\bbP (\go_n =\pm 1)=1/2$ and
\begin{equation}
H_{N, \go} (S)=
 \gb \sum_{n=1}^N \left( 1 + \gep \go_n\right) \ind_{\{S_n=n\}},
 \end{equation}
 with $\gb$ and $\gep$ real numbers and $S$ the simple random walk on $\Z$.
 We observe that
 Proposition \ref{th:main} applies to this case with $D_n =\{ n\}^\complement$ and
 that the model is solvable in detail.
 In particular $\free (\gb, \gep) = (\gb - \log 2)\vee 0$, regardless of the value of $\gep$.
 The annealed computation instead yields $\freea (\gb , \gep)=
 (\gb +\log \cosh(\gep)-\log 2)\vee 0$. Notice in particular
 that the critical values of $\gb$, respectively $\log 2$ and
 $\log 2 -\log \cosh(\gep)$, differ as long as there is disorder
 in the system ($\gep \neq 0$).
 It is interesting to see in this toy model how
 $A_N$ has to be chosen {\sl very non local} in order to
 improve on the annealed bound.

\smallskip

\begin{rem}\rm
We point out that we restricted our examples
only to cases in which $S$ is a simple random walk, but
 in principle our approach goes through for much more general models,
 like walks with correlated  increments or self--interacting walks, see
 \cite{cf:OTW}  for an example.
 And of course $S_n$ takes  values in $\Z^d$ only for ease of exposition
 and can be easily generalized.
 It is however unclear whether our argument applies
 to the disordered wetting problem in $d+1$ dimensions, $d>1$.
 In this case $S$ is a random interface, the Hamiltonian is like
 in
\eqref {eq:Wetting}, but $n \in \{0,1, 2, \ldots\} ^d$, $S_n\in \Z$ or
$\R$. We set for example $S_n=0$ when one of the coordinates
of $n$ is zero. The missing ingredient is an analog of Theorem
\ref{th:main1} in higher dimensions.
\end{rem}

\smallskip
\section{On cocycles with null free energy}
\label{sec:cocycles}

Let $\{\bs{\go}_n\}_{n\in\N}$ be an IID sequence of random variables under the probability measure $\bbP$, taking values in a finite space $\Gamma$ (we have switched the notation $\go \to \bs{\go}$ for clarity). The law of $\bs{\go}_1$ on $\Gamma$ is denoted by $\nu$: we will assume that $\nu(\ga)>0$ for all $\ga\in\Gamma$.

We are interested in families $A = \{A_N\}_{N\in\N}$ of random variables of the form of empirical averages of a centered local function $F$, that is
\begin{equation} \label{eq:hamiltonian}
A_N = \sum_{n = 1}^N F(\bs{\go}_n,\ldots,\bs{\go}_{n+k})\,,
\end{equation}
where $k\in\{0\} \cup \N$ and $F$ is a real function defined on $\Gamma ^{k+1}$ such that $\int F \dd\nu^{*(k+1)}=0$. We will call $A = \{A_N\}_{N\in\N}$ a centered \textsl{cocycle}, and with some abuse of notation we will speak of {\sl the cocycle~$F$} to mean the cocycle~$\{A_N\}_{N\in\N}$ defined by ~(\ref{eq:hamiltonian}).

A cocycle $F:\Gamma^{k+1} \to \R$ is said to be a \textsl{coboundary} if (when $k\geq 1$) there exists a function $G:\Gamma ^k \to \R$ such that
\begin{equation} \label{eq:grad}
F(\ga_1,\ldots ,\ga_{k+1}) = G(\ga_2,\ldots \ga_{k+1}) - G(\ga_1,\ldots ,\ga_{k})
\end{equation}
for all $\ga_1,\ldots, \ga_{k+1} \in \Gamma$. When $k=0$, we say that $F$ is a \textsl{coboundary} if it is identically zero: $F(\ga) = 0$ for every $\ga \in \Gamma$.

\smallskip

For $\gb \in\R$ we define the free energy $L^F(\gb)$ of a cocycle $F$ as
\begin{equation} \label{eq:free_energy}
L^F(\gb) := \lim_{N\to\infty} \frac 1N \log \bbE \Big[e^{\gb A_N}\Big] \,.
\end{equation}
The limit above is easily seen to exist by a standard superadditive argument, and Jensen's inequality yields immediately $L^F(\gb)\geq0$. Of course, if $F$ is a coboundary then the corresponding free energy vanishes for all $\gb\in\R$. That also the converse is true is the object of the following theorem.

\begin{theorem} \label{th:main1}
Let $F$ be a centered cocycle, and let $L^F(\gb)$ be the corresponding free energy, defined by (\ref{eq:free_energy}). The following conditions are equivalent:
\begin{enumerate}
\item $F$ is a coboundary;
\item $L^F(\gb)=0$ for all $\gb\in\R$;
\item $L^F(\gb_0)=0$ for some $\gb_0\in\R\setminus\{0\}$.
\end{enumerate}
\end{theorem}

The proof is obtained combining convexity ideas with the following combinatorial reformulation of the condition that a function be a coboundary.

\begin{lemma} \label{lem:main}
A function  $F:\Gamma ^{k+1} \to \R$ is a coboundary if and only if for every $N\in\N$ and for every $(\eta_1,\ldots , \eta_{N}) \in \Gamma^N$ the following relation holds:
\begin{align} \label{eq:grad1}
\sum_{i=1}^{N} F(\eta_i,\eta_{i\oplus_N 1},\ldots, \eta_{i\oplus_N k}) = 0 \,,
\end{align}
where for $a,b \in \N$ we have set $a\oplus_N b := (a + b)\mod N$.
\end{lemma}

\proof
The {\sl if} part trivially follows from the definition of a coboundary (see (\ref{eq:grad})), so we can focus on the
{\sl only if} part. As a matter of fact, we will use the hypothesis of the Lemma only for two values of $N$, namely $N=2k$ and $N=2k+1$.

Let us take $k$ elements  $\gamma_1,\ldots,\gamma_k\in \Gamma$,
arbitrarily chosen, that will be kept \textsl{fixed} throughout the proof; moreover, let $\ga_1,\ldots , \ga_{k+1}$ denote generic elements of~$\Gamma$. We start rewriting equation~(\ref{eq:grad1}) for $N=2k+1$, with $(\eta_1,\ldots , \eta_{N})=(\ga_1,\ldots , \ga_{k+1},\gamma_1,\ldots,\gamma_k)$,  as
\begin{align} \label{eq:first_step}
F(\ga_1,\ldots , \ga_{k+1}) = -\sum_{i=1}^{k} F(\ga_{i+1},\ldots, \ga_{k+1}, \gamma_1,\ldots, \gamma_i) - \sum_{i=1}^{k} F(\gamma_i,\ldots, \gamma_k,\ga_1,\ldots, \ga_i) \,.
\end{align}
In order to determine an alternative expression for the second sum in the r.h.s., we use again equation~(\ref{eq:grad1}), this time with~$N=2k$ and~$(\eta_1,\ldots , \eta_{N})=(\ga_1,\ldots , \ga_{k},\gamma_1,\ldots,\gamma_k)$, getting
\begin{equation} \label{eq:second_step}
\sum_{i=1}^{k} F(\gamma_i,\ldots, \gamma_k,\ga_1,\ldots, \ga_i) = - \sum_{i=1}^{k} F(\ga_i,\ldots, \ga_k,\gamma_1,\ldots, \gamma_i)\,.
\end{equation}

If now we introduce a function $G:\Gamma^k \to \R$, defined by
\begin{align*}
G(\zeta_1,\ldots, \zeta_k) := - \sum_{i=1}^{k} F(\zeta_i,\ldots, \zeta_k,\gamma_1,\ldots, \gamma_i) \,,
\end{align*}
we can combine equations (\ref{eq:first_step}) and (\ref{eq:second_step}) to get
\begin{align*}
F(\ga_1,\ldots , \ga_{k+1})
= G(\ga_2,\ldots \ga_{k+1}) - G(\ga_1,\ldots \ga_{k}) \,,
\end{align*}
so that the proof is completed.
\qed

\bigskip

\noindent\textbf{Proof of Theorem \ref{th:main1}.} It has already been remarked that $(1) \Rightarrow (2)$, and of course $(2) \Rightarrow (3)$ holds trivially. In the following we are going to prove that $(3) \Rightarrow (2) \Rightarrow (1)$.

We start determining an explicit expression for the free energy. For this, we define a slight modification of the cocycle $A$ defined by~(\ref{eq:hamiltonian}), by setting
\begin{equation} \label{eq:mod_ham}
\tilde{A}_N := \sum_{n = 1}^N F(\bs{\go}_n,\bs{\go}_{n\oplus_N 1},\ldots,\bs{\go}_{n\oplus_N k})\,,
\end{equation}
where by $\oplus_N$ we mean addition modulo $N$. Of course, only the last $k$ addends in the sum are really changed: as $F$ is a bounded function (the space $\Gamma$ is finite), it easily follows that the free energies of $A$ and $\tilde{A}$ are the same, so that we can write
\begin{equation} \label{eq:zeta}
L^F(\gb) = \lim_{N\to\infty} \frac 1 N \log Z_N(\gb) \qquad \text{where} \qquad Z_N(\gb) = Z^F_N(\gb) = \bbE \Big[ e^{\gb \tilde{A}_N} \Big] \,.
\end{equation}

Now we introduce the $\Gamma^{k+1} \times \Gamma^{k+1}$ matrix $A_\gb$, defined for $\ga_i,\gamma_i \in \Gamma,\ i=1,\ldots,k+1$ by
\begin{equation} \label{eq:matrix}
A_\gb \big[ (\ga_1, \ldots, \ga_{k+1}), (\gamma_1, \ldots, \gamma_{k+1}) \big] := \gd_{\gamma_1,\ga_2} \cdots \gd_{\gamma_{k}, \ga_{k+1}} \cdot e^{\gb F(\gamma_1, \ldots, \gamma_{k+1})} \cdot \nu(\gamma_{k+1})\,.
\end{equation}
Developing the expectation defining $Z_N(\gb)$ we get
\begin{align}
Z_N(\gb) &\;=\; \sum_{\zeta_1,\ldots,\zeta_N \in \Gamma} e^{\gb \sum_{i=1}^N F(\zeta_i,\zeta_{i\oplus_N 1}, \ldots, \zeta_{i\oplus_N k})} \cdot \nu(\zeta_1) \cdots \nu(\zeta_N) \nonumber \\
\label{eq:trace}
& \;=\; \text{Tr} \big[ A_\gb ^N \big] \;=\; \sum_{i=1}^{|\Gamma|^{2(k+1)}} e_i(\gb)^N \,,
\end{align}
where $\{e_i(\gb),\ i=1,\ldots,|\Gamma|^{2(k+1)}\}$ are the (possibly complexes) eigenvalues of the matrix $A_\gb$ (counted repeatedly according to their algebraic multiplicity). It's immediate to check that $A_\gb$ is an irreducible, aperiodic matrix, and since its entries are nonnegative we can apply Perron--Frobenius theory \cite{cf:Asm}: there exists a real positive simple eigenvalue, say $e_1(\gb)$, such that $|e_i(\gb)| < e_1(\gb)$ for every $i\geq 2$. To lighten the notation, from now on we will let $e(\gb) := e_1(\gb)$. Combining (\ref{eq:zeta}) with (\ref{eq:trace}) we get
\begin{equation} \label{eq:zeta_as}
Z_N(\gb) = e(\gb)^N \cdot \Bigg( 1 + \sum_{i=2}^{|\Gamma|^{2(k+1)}} \bigg( \frac{e_i(\gb)}{e(\gb)} \bigg) ^N \Bigg) \,,
\end{equation}
so that
\[
Z_N(\gb) \cdot e(\gb)^{-N} \to 1 \quad \text{as }N\to\infty \,.
\]
From this \textsl{sharp asymptotics} for $Z_N(\gb)$ we obtain in particular the explicit expression of $L^F(\gb)$ we were looking for:
\begin{equation} \label{eq:log_e}
L^F(\gb) = \log e(\gb) \,.
\end{equation}

This equation shows that $L^F(\gb)$ is a \textsl{real analytic} function of $\gb\in\R$, since $e(\gb)$ is so: this is because the Perron--Frobenius eigenvalue is a simple root of the characteristic polynomial and the entries of $A_\gb$ are real--analytic functions of $\gb\in\R$.

From (\ref{eq:zeta}) it is clear that $\log Z_N(\gb)$ is a convex function of $\gb\in\R$, for every $N\in\N$. Moreover, we have $Z_N(\gb) \geq 1$ for every $\gb\in\R$ by Jensen's inequality, and trivially $Z_N(0)=1$. It follows immediately that $L^F(\gb)$ is a convex function too, being the pointwise limit of $\log Z_N(\gb)/N$, that $L^F(\gb) \geq 0$ for every $\gb\in\R$, and $L^F(0)=0$.

\smallskip

Let's assume that condition $(3)$ in the statement of the theorem holds, that is $L^F(\gb_0)=0$ for some $\gb_0 > 0$ (the case $\gb_0<0$ is completely analogous): the preceding observations yield $L^F(\gb)=0$ for every $\gb \in [0,\gb_0]$, and by analyticity we conclude that indeed $L^F(\gb)=0$ for every $\gb \in \R$. We have thus shown that $(3) \Rightarrow (2)$.

\smallskip

Now we assume that condition $(2)$ holds: by (\ref{eq:log_e}) this means $e(\gb)=1$ for every $\gb\in\R$, and (\ref{eq:zeta_as}) we have that
\begin{equation} \label{eq:ineq}
|Z_N(\gb)| \leq   e(\gb)^N \cdot |\Gamma|^{2(k+1)}  =  |\Gamma|^{2(k+1)} \qquad \forall N\in\N\,,\ \forall \gb\in\R\,.
\end{equation}
Since $\log Z_N(\gb)$ is a convex function, $Z_N(\gb)$ is convex too; furthermore, we have already remarked that $Z_N(\gb) \geq 1$ for every $\gb\in\R$ and that $Z_N(0)=1$. Since (\ref{eq:ineq}) shows that $|Z_N(\gb)|$ is bounded, by elementary convex analysis it follows that $Z_N$ must be constant, therefore $Z_N(\gb)=1$ for all $\gb\in\R$ and $N\in\N$. This means that for every $\gb\in\R$ Jensen's inequality for $Z_N(\gb)$ it's not strict: since for any $\gb>0$ the function $\{x \mapsto e^{\gb x}\}$ is a strictly convex function, this can happen if and only if $\tilde{A}_N$ is $\bbP$--a.s. constant, for every $N\in\N$. Recalling (\ref{eq:mod_ham}) and the fact that by hypothesis $\nu(\ga)>0$ for every $\ga\in\Gamma$, this amounts to saying that
\[
\sum_{i=1}^{N} F(\eta_i,\eta_{i\oplus_N 1},\ldots, \eta_{i\oplus_N k}) = 0 \,,
\]
for every $N\in\N$ and for every $\eta_1, \ldots, \eta_N \in \Gamma$: applying Lemma~\ref{lem:main} we conclude that $F$ is a coboundary, and the proof is complete.\qed

\chapter[A renewal theory approach to periodic polymers]
{A renewal theory approach to polymers\\ with periodic
distribution of charges}
\label{ch:cgz}

In this chapter we consider a general model
of an heterogeneous polymer chain in the proximity of
an interface between two selective solvents, which
includes as special cases the copolymer near a selective interface
and the pinning model introduced in Chapter~\ref{ch:first}.
The heterogeneous character of the model comes from the fact that
the interaction of each {\sl monomer unit} is governed by a
{\it charge} that it carries. We consider the model in the {\sl periodic setting},
that is when the charges repeat themselves along the chain in a periodic fashion.
The main question is of course whether the polymer remains tightly close
to the interface ({\sl localization})
or there is a marked preference for one solvent ({\sl delocalization}).

We propose an approach based on renewal theory
that yields sharp estimates on the partition function of the model
in all the regimes (localized, delocalized and critical).
This in turn allows to get a very precise
description of the polymer measure, both in a local sense
({\sl thermodynamic limit}) and in a global sense ({\sl scaling limits}):
see~\S~\ref{sec:spirit} for an outline of our results
and~\S~\ref{sec:results} for a detailed exposition.
A key point, but also a byproduct, of our analysis is the closeness of
the polymer measure to a suitable {\sl Markov Renewal Process}.

\smallskip
The preprint~\cite{cf:CGZ} has been taken from the content of this chapter.

\smallskip
\section{Introduction and main results}
\label{sec:intro_cgz}

\smallskip

\subsection{Two motivating models}
\label{sec:2examples}
We slightly enlarge our setting with respect to Chapter~\ref{ch:first}, namely
we work with a random walk $S:=\left\{ S_n\right\}_{n=0,1, \ldots}$
with IID symmetric increments $\{X_j\}_{j \ge 1}$ taking values in $\{-1, 0,+1\}$.
Hence the law of the walk is identified by
$p:=\bP \left( X_1 =1\right)(=\bP \left( X_1 =-1\right))$, and we assume that $p\in (0,1/2)$.
Note that we have excluded the case $p=1/2$ and this has been done in order to
lighten the exposition: all the results we  present have a close analog in the case $p=1/2$,
however the statements
require a minimum of notational care because of the periodicity of the walk. We also consider
a sequence $\go:=\left\{ \go_n\right\}_{n\in \N=\{1,2,\ldots\} }$ of real numbers
with the property that $\go_n=\go_{n+T}$ for some $T\in \N$ and for every $n$:
we denote by  $T(\go)$ the minimal value of $T$.

\smallskip

Before defining the general model that will be the object of our analysis,
we recall the two motivating models that were introduced in Chapter~\ref{ch:first}.

 \smallskip
 \begin{enumerate}
 \item {\it Pinning and wetting models.} For $\gl \ge 0$ consider the probability measure $\bP_{N,\go}$
  defined by
 \begin{equation}
 \label{eq:pinning}
 \frac{\dd \bP_{N, \go}}{\dd \bP} (S) \, \propto \, \exp \left( \gl \sum _{n=1}^ {N}
 \go_n \ind_{\left\{S_n=0 \right\}} \right).
 \end{equation}
The walk receives a {\sl pinning reward}, which may be negative or positive, each time
it visits the origin. By considering the directed walk viewpoint, that is $\left\{ (n, S_n)\right\}_n$,
one may interpret this model in terms of a directed linear chain
receiving an energetic contribution when it touches an interface. In this context it is natural
to introduce the asymmetry parameter $h:= \sum_{n=1}^T \go_n /T$, so that one isolates
a constant drift term from the {\sl fluctuating} behavior of $\go$.
The question is whether for large $N$
the measure $\bP_{N, \go}$ is rather attracted or repelled by the interface (there is in principle
the possibility for the walk to be essentially indifferent of such a change of measure, but we anticipate
that this happens only in trivially degenerate cases while in {\sl critical} situations a more subtle
scenario shows up).

By multiplying the right--hand side of  \eqref{eq:pinning} by $\ind_{\left\{S_n\ge 0:\,  n=1, \ldots, N\right\}}$
one gets to a so called {\sl wetting model}, that is the model of an interface interacting with an impenetrable
wall. The {\sl hard--wall} condition induces a repulsion effect of purely entropic origin which
is in competition with attractive energy effects: one expects that in this case $h$ needs to
be positive for the energy term to overcome the entropic repulsion effect, but quantitative estimates
are not a priori obvious.

There is an extensive literature on periodic pinning and wetting models, the majority of which is
 restricted
to the $T=2$ case, we mention for example \cite{cf:GG,cf:NZ}.
\medskip
 \item {\it Copolymer near a selective interface.} Much in the same way we introduce
 \begin{equation}
 \label{eq:copolymer}
 \frac{\dd \bP_{N, \go}}{\dd \bP} (S) \, \propto \, \exp \left( \gl \sum _{n=1}^ {N}
 \go_n \sign \left( S_n\right) \right),
 \end{equation}
 where if $S_n=0$ we set $\sign (S_n) := \sign(S_{n-1}) \, \ind_{\{S_{n-1} \neq 0\}}$. This
 convention for defining $\sign(0)$, that will be kept throughout the chapter,
 has the following simple interpretation: $\sign(S_{n}) = +1,0,-1$ according to whether
 the bond joining $S_{n-1}$ and $S_n$ lies above, on,  or below the $x$--axis.

 Also in this case we take a directed walk viewpoint and then
 $\bP_{N,\go}$ may be interpreted as a polymeric chain in which the
 monomer units, the bonds of the walk, are charged. An interface, the $x$--axis,
 separates two solvents, say oil above and water below:
 positively charged monomers are hydrophobic and negatively charged ones
 are instead hydrophilic.
In this case one expects a competition between three possible scenarios: polymer
preferring water, preferring oil or undecided between the two and choosing to
fluctuate in the proximity of the interface. We will therefore talk of delocalization in water (or oil)
or of localization at the interface. Critical cases are of course of particular  interest.

We select \cite{cf:MGO,cf:SD} from the physical literature on periodic copolymers,
keeping however in mind that periodic copolymer modeling  has a central
role in applied chemistry and material science.
 \end{enumerate}


\smallskip

\subsection{A general model}
We point out that the models presented in \S~\ref{sec:2examples}
are particular examples of the polymer measure with Hamiltonian
\begin{equation}
 \label{eq:genH}
\cH _N (S)\, =\,
 \sum_{i =\pm 1} \sum _{n=1}^ {N}
 \go_n^{(i)} \ind_{\left\{\sign \left( S_n\right)=i\right\}}+
 \sum _{n=1}^ {N}
 \go^{(0)}_n \ind_{\left\{S_n=0 \right\}}+
  \sum _{n=1}^ {N}
 \widetilde\go^{(0)}_n \ind_{\left\{\sign\left(S_n\right)=0 \right\}},
 \end{equation}
where $\go^{(\pm1)}$, $ \go^{(0)}$ and  $\widetilde\go^{(0)}$
are periodic sequences of real numbers. Observe that, by our conventions on~$\sign(0)$, the last term gives an energetic contribution (of pinning/depinning type) to the bonds lying on the interface.

Besides being a natural model, generalizing and interpolating between
pinning and copolymer models, the general model we consider is the
one considered at several instances, see e.g. \cite{cf:SW} and references therein.
\smallskip

\begin{rem} \label{rem:1.1}\rm
The \textsl{copolymer} case corresponds to $\go^{(+1)}=-\go^{(-1)}= \gl \, \go$ and $ \go^{(0)}=\widetilde\go^{(0)}=0$,
while the \textsl{pinning} case corresponds to $\go^{(0)}=\gl \, \go$ and
$\go^{(+1)}=\go^{(-1)}=\widetilde\go^{(0)}=0$. We stress
 that the \textsl{wetting} case can be included
too, with the choice $\go^{(0)}=\gl \go$, $\go^{(-1)}_n=-\infty$ for every $n$ and
$\go^{(+1)}=\widetilde\go^{(0)}=0$. Of course plugging $\go^{(-1)}_n=-\infty$ into the Hamiltonian~\eqref{eq:genH}
 is a bit formal, but it
 simply corresponds to a constraint on $S$ in the polymer measure associated to~$\cH_N$, see \eqref{eq:genP} below. For ease of exposition we will restrict to finite values of the charges~$\go$, but
 the generalization is straightforward.
\end{rem}
\smallskip

\begin{rem} \rm
\label{rem:intper}
We take this occasion for stressing that, from an applied viewpoint, the interest  in periodic
models of the type we consider appears to be  at least two--fold.
On one hand periodic models are often chosen as caricatures of the
{\sl quenched disordered} models, like the ones in which the charges are
a typical realization of a sequence of independent random variables
(e.g. \cite{cf:AS,cf:BdH,cf:G,cf:SW} and references therein).
In this respect and taking a mathematical standpoint, the relevance
of periodic models, which may be viewed as {\sl weakly inhomogeneous},
for understanding the strongly inhomogeneous quenched set--up is
at least questionable and the approximation of quenched models
with periodic ones, in the limit of large period, poses very interesting and
challenging questions. In any case, the precise description
of the periodic case that we have obtained in this work
highlights limitations and perspectives of periodic modeling
for strongly inhomogeneous systems.
One the  other hand, as already mentioned above, periodic models
are absolutely natural and of direct relevance for application,
for example when dealing with {\sl molecularly engineered} polymers~\cite{cf:NN,cf:SD}.
\end{rem}

\smallskip

Starting from the Hamiltonian \eqref{eq:genH}, for $a= \rc$ ({\sl constrained}) or
$a=\rf$ ({\sl free})
we introduce the \textsl{polymer measure}~$\bP^a_{N, \go}$ on $\Z^\N$, defined by
\begin{equation}
\label{eq:genP}
\frac{\dd \bP^a_{N, \go}}{\dd \bP} (S)\, =\,
\frac{\exp\left(\cH_N (S)\right)}{\tilde Z_{N, \go}^a} \left(\ind_{\left\{a=\rf\right\} } +
\ind_{\left\{a=\rc\right\} } \ind_{\left\{ S_N=0 \right\}}
\right),
\end{equation}
where $\tilde Z_{N, \go}^a := \bE [\exp (\cH_N) \, (\ind_{\{a=\rf\} } +
\ind_{\{a=\rc\} } \ind_{\{ S_N=0 \}}) ]$ is the {\sl partition function}, that is the normalization constant.
Here $\go$ is a shorthand for the four periodic sequences appearing in the definition~\eqref{eq:genH} of~$\cH_N$, and we will use $T=T(\go)$ to denote the smallest common period of the sequences.

\smallskip

The Laplace asymptotic behavior of $\tilde Z_{N,\go}$
plays an important role
and the quantity
\begin{equation}
\label{eq:fe}
f_\go \, := \, \lim_{N \to \infty } \frac 1N \log \tilde Z_{N, \go} ^{\rc},
\end{equation}
is usually called {\sl free energy}. The existence of the limit above follows from
a direct super--additivity argument, and it is easy to check that
$\tilde Z_{N, \go} ^{\rc}$ can be replaced by $\tilde Z_{N, \go} ^{\rf}$
without changing the value of $f_\go$,
see e.g.~\cite{cf:G}. The standard free energy approach to this type of models
starts from the observation that
\begin{equation}
\label{eq:fdeloc0}
\begin{split}
f_\go \, & \;\ge\;
\lim_{N\to \infty} \; \frac{1}{N} \, \log \bE
\Big[\exp\big(\cH_N(S)\big)\, ; \, S_n>0 \text{ for } n=1, \ldots, N\, \Big]\\
& \;=\; \frac 1{T(\go)}\sum_{n=1}^{T(\go)} \go_n^{(+1)} \;+\;
\lim_{N\to \infty} \;
\frac{1}{N} \, \log \bP\big(S_n>0 \text{ for } n=1, \ldots, N\, \big) \,.
\end{split}
\end{equation}
It is a classical result~\cite[Ch. XII.7]{cf:Fel2} that $\bP (S_n>0 \text{ for } n=1, \ldots, N) \sim c N^{-1/2}$, as~$N\to\infty$,
for some $c \in (0,\infty)$ (by
 $a_N \sim b_N$ we mean $a_N/b_N \to 1$). Hence the  limit of the last term of \eqref{eq:fdeloc0} is zero and one
easily concludes that
\begin{equation}
\label{eq:fdeloc}
f_\go \, \ge\, f_{\go}^{\cD} \, :=\, \max _{i=\pm 1} h_\go(i),
\qquad h_\go (i):=\frac 1{T(\go)}\sum_{n=1}^{T(\go)} \go_n^{(i)}
.
\end{equation}
Having in mind the steps in \eqref{eq:fdeloc0}, one is led to the following basic
\begin{definition} \label{def:main}
The polymer chain defined by \eqref{eq:genP} is said to be:
\begin{itemize}
\item \textsl{localized (at the interface)} if $f_\go>f_{\go}^{\cD}$;
\smallskip
\item \textsl{delocalized above the interface} if $f_\go= h_\go (+1)$;
\smallskip
\item \textsl{delocalized below the interface} if $f_\go= h_\go (-1)$.
\end{itemize}
\end{definition}
\noindent
Notice that, with this definition, if $h_\go(+1) = h_\go(-1)$ and the polymer is delocalized, it is delocalized both above and below the interface.

\smallskip

\begin{rem}\rm
Observe that the polymer measure~$\bP^a_{N,\go}$ is invariant under the joint transformation $S \to -S$, \ $\go^{(+1)} \to \go^{(-1)}$, hence by symmetry we may (and will) assume that
\begin{equation}\label{eq:ass_h}
h_\go \, := \, h_\go(+1) \, - \, h_\go(-1) \, \geq \, 0\,.
\end{equation}
It is also clear that we can add to the Hamiltonian~$\cH_N$ a constant term (with respect to~$S$) without changing the polymer measure. Then we set
\begin{align*}
    \cH'_N (S) \, &:=\, \cH _N (S) \;-\; \sum _{n=1}^ {N}\go_n^{(+1)}\,,
\end{align*}
which amounts to redefining $\go_n^{(+1)} \to 0$, $\go_n^{(-1)} \to (\go_n^{(-1)}-\go_n^{(+1)})$ and $\widetilde\go_n^{(0)} \to (\widetilde\go_n^{(0)}-\go_n^{(+1)})$, and we can write
\begin{equation}
\label{eq:newP}
\frac{\dd \bP^a_{N, \go}}{\dd \bP} (S)\, =\,
\frac{\exp\left(\cH_N' (S)\right)}{Z_{N, \go}^a} \left(\ind_{\left\{a=\rf\right\} } +
\ind_{\left\{a=\rc\right\} } \ind_{\left\{ S_N=0 \right\}}
\right),
\end{equation}
where $Z_{N, \go}^a$ is a new partition function which coincides with
$\tilde Z_{N, \go}^a \exp(- \sum _{n=1}^ {N}\go_n^{(+1)})$. The corresponding free energy~$\tf_\go$ is given by
\begin{equation}
\label{eq:F}
    \tf_\go \;:=\; \lim_{N\to\infty} \frac 1N \log Z_{N,\go}^a \;=\; f_\go - f_\go^\cD\,,
\end{equation}
and notice that in terms of~$\tf_\go$ the condition for localization (resp. delocalization) becomes $\tf_\go > 0$ (resp.~$\tf_\go = 0$). From now on, speaking of partition function and free energy we will always mean~$Z_{N,\go}^a$ and~$\tf_\go$.
\end{rem}


\smallskip

\subsection{From free energy to path behavior}
\label{sec:spirit}

In order to understand the spirit of our approach, let us briefly outline our results (complete results are given in~\S~\ref{sec:results} below).

Our first goal is to give necessary and sufficient \textsl{explicit conditions} in terms of the charges~$\go$ for the (de)localization of the polymer chain, see Theorem~\ref{th:as_Z}. We point out that the content of this theorem is in fact much richer, as it gives the \textsl{sharp asymptotic behavior} (and not only the Laplace one \cite{cf:BG}) as~$N\to\infty$ of the constrained partition function~$Z^\rc_{N,\go}$. In particular we show that when the polymer is delocalized ($\tf_\go = 0$) the constrained partition function $Z^\rc_{N,\go}$ is actually vanishing as~$N\to\infty$. Moreover the rate of the decay induces a further distinction in the delocalized regime between a \textsl{strictly delocalized regime} ($Z_{N,\go}^\rc \sim c_1 N^{-3/2}$, $c_1 \in (0, \infty)$) and a \textsl{critical regime} ($Z_{N,\go}^\rc \sim c_2 N^{-1/2}$, $c_2 \in (0, \infty)$).

These asymptotic results are important because they allow to address further interesting issues. For example, it has to be admitted that defining (de)localization in terms of the free energy is not completely satisfactory, because one would like to characterize the
 polymer path properties. In different terms, given a polymer measure which is (de)localized according to Definition~\ref{def:main}, to what extent are its typical paths really (de)localized? Some partial answers to this question are known, at least in some particular instances: we mention
 here the case of $T(\go)=2$ copolymers \cite{cf:MGO}
 and the case of homogeneous pinning and wetting models \cite{cf:DGZ,cf:IY,cf:Upton}.

Our main aim is to show that, for the whole class of models we are considering, free energy (de)localization does correspond to a \textsl{strong form} of path (de)localization. More precisely, we look at path behavior from two different viewpoints.
\smallskip
\begin{itemize}
\item \emph{Thermodynamic limit.} We show that the measure $\bP_{N,\go}^a$ converges weakly as~$N\to\infty$ toward a measure~$\bP_\go$ on~$\Z^\N$, of which we give an explicit construction, see Section~\ref{sec:infvol}. It turns out that the properties of~$\bP_\go$ are radically different in the three regimes (localized, strictly delocalized and critical), see Theorem~\ref{th:infvol}.
It is natural to look at these results as those characterizing the {\sl local} structure of the polymer chain.
\item \rule{0pt}{15pt}\emph{Brownian scaling limits.} We prove that the diffusive rescaling of the polymer measure $\bP_{N,\go}^a$ converges weakly in $C([0,1])$ as~$N\to\infty$. Again the properties of the limit process, explicitly described in Theorem~\ref{th:scaling}, differ considerably in the three regimes.
Moreover we stress that scaling limits describe {\sl global} properties of the chain.
\end{itemize}
\smallskip
We insist on the fact  that the path analysis just outlined has been obtained exploiting heavily the sharp asymptotic behavior of~$Z_{N,\go}^\rc$ as~$N\to\infty$.
In this sense our results are the direct sharpening of the Large Deviations
approach taken in \cite{cf:BG}, where a formula  for $\tf_\go$
was obtained for periodic copolymers (but the method of course directly extends to the general
case considered here).
Such a formula (see \S~\ref{sec:as_Z_loc}), that reduces the problem of computing the free energy
to a finite dimensional problem connected to a suitable
Perron--Frobenius matrix, in itself suggests the new approach taken here
since it makes rather apparent the link between periodic
copolymers and the class of {\sl Markov renewal processes}~\cite{cf:Asm}.
On the other hand, with respect to \cite{cf:BG}, we leave aside any issue
concerning the phase diagram (except for \S~\ref{sec:P} below).


\smallskip

\subsection{The order parameter~$\gd^\go$}
\label{sec:order_par}

It is a remarkable fact that the dependence of our results on the charges~$\go$ is essentially encoded in one single parameter~$\gd^\go$, that can be regarded as the \textsl{order parameter} of our models. For the definition of this parameter, we need some preliminary notation. We start with the law of the first return to zero of the original walk:
\begin{equation} \label{eq:first_ret}
    \tau_1 := \inf \{n>0:\ S_n = 0\} \qquad \quad K(n) \;:=\; \bP \big( \tau_1 = n \big)\,.
\end{equation}
It is a classical result \cite[Ch. XII.7]{cf:Fel2} that
\begin{equation} \label{eq:as_K}
    \exists \lim_{n\to\infty} n^{3/2} \, K(n) \;=:\; c_K \in (0,\infty)\,.
\end{equation}
Then we introduce the Abelian group $\bbS:= \Z/ (T\Z)$ and to indicate that an integer~$n$ is in the equivalence class $\gb\in\bbS$ we write equivalently $[n]=\gb$ or $n\in\gb$. Notice that the charges~$\go_n$ are functions of~$[n]$, and with some abuse of notation we can write~$\go_{[n]}:=\go_n$. The key observation is that, by the $T$--periodicity of the charges~$\go$ and by the definition~\eqref{eq:ass_h} of~$h_\go$, we can write
\begin{equation*}
    \sum_{n=n_1+1}^{n_2} (\go_n ^{(-1)}-\go_n ^{(+1)}) \;=\; -(n_2-n_1) \, h_\go \;+\; \Sigma_{[n_1],[n_2]}\,.
\end{equation*}
Thus we have decomposed the above sum into a drift term and a more fluctuating term, where the latter has the remarkable property of depending on~$n_1$ and~$n_2$ only through their equivalence classes~$[n_1]$ and~$[n_2]$. Now we can define three basic objects:
\begin{itemize}
\smallskip
\item for $\ga, \gb\in \bbS$ and $\ell \in \N $ we set
\begin{equation} \label{eq:def_Phi}
\Phi^\go_{\ga, \gb}(\ell)\, :=\,
\begin{cases}
    \go ^{(0)}_{\gb} \;+\; \Big(\tilde  \go ^{(0)}_{\gb} - \go^{(+1)}_\gb \Big) & \text{if }
    \ell=1,\ \ell \in \gb-\ga \\
     \go ^{(0)}_{\gb} +
    \rule{0pt}{22pt} \displaystyle  \log  \bigg[ \frac 12 \Big(1+ \exp\big( -\ell \, h_\go  + \Sigma_{\ga,\gb}\big)\Big) \bigg] &\text{if } \ell > 1,\ \ell \in \gb-\ga\\
    \rule{0pt}{16pt} 0 & \text{otherwise}
\end{cases}
\;,
\end{equation}
which is a sort of integrated version of our Hamiltonian;

\medskip
\item for $x\in\N$ we introduce the $\bbS \times \bbS$ matrix $M^\go_{\ga,\gb}(x)$ defined by
\begin{equation} \label{eq:matrix_cgz}
M^\go_{\ga,\gb}(x) \;:=\; e^{\Phi^\go_{\ga,\gb}(x)} \, K(x) \, \ind_{(x \in \gb -\ga)}\,;
\end{equation}

\medskip
\item summing the entries of~$M^\go$ over~$x$ we get a~$\bbS\times\bbS$ matrix that we call~$B^\go$:
\begin{equation} \label{eq:def_B}
    B^\go_{\ga,\gb} := \sum_{x\in\N} M^\go_{\ga,\gb}(x)\,.
\end{equation}
\end{itemize}
\smallskip
The meaning and motivation of these definitions, that at this point might appear artificial,
are explained in detail in \S~\ref{sec:RWexcurs}. For the moment we only stress that the
above quantities are \textsl{explicit functions} of the charges~$\go$ and of the law of the
underlying random walk (to lighten the notation, the $\go$--dependence of these quantities will
be often dropped in the following).

\smallskip
We can now define our order parameter~$\gd^\go$. Observe that $B_{\ga,\gb}$ is a finite dimensional matrix with nonnegative entries, hence the Perron--Frobenius (P--F) Theorem (see e.g. \cite{cf:Asm}) entails that $B_{\ga,\gb}$ has a unique real positive eigenvalue, called the Perron--Frobenius eigenvalue, with the property that it is a simple root of the characteristic polynomial and that it coincides with the spectral radius of the matrix. This is exactly our parameter:
\begin{equation} \label{eq:def_delta}
\gd^\go \ := \ \text{Perron--Frobenius eigenvalue of $B^\go$}\,.
\end{equation}

\smallskip


\subsection{The main results}
\label{sec:results}

Now we are ready to state our results. We start characterizing the (de)localization of the polymer chain in terms of~$\gd^\go$.


\medskip

\begin{theorem}[\bf Sharp asymptotics] \label{th:as_Z}
The polymer chain is localized if and only if $\gd^\go > 1$.
More precisely, the asymptotic behavior of $Z_{N,\go}^\rc$ as
$N\to\infty$, $[N]=\eta$ is given by
\smallskip
\begin{enumerate}
\item for $\gd^\go > 1$ (\textsl{localized regime}) \ $Z_{N,\go}^\rc \;\sim\; C^>_{\go,\eta} \, \exp \big(\tf_\go N \big)$ \,;
\medskip
\item for $\gd^\go < 1$ (\textsl{strictly delocalized regime}) \ $Z_{N,\go}^\rc \;\sim\; C^<_{\go,\eta} \, / \, N^{3/2}$ \,;
\medskip
\item for $\gd^\go = 1$ (\textsl{critical regime}) \ $Z_{N,\go}^\rc \;\sim\; C^=_{\go,\eta} \, / \, \sqrt{N} $ \,,
\end{enumerate}
\smallskip
where $\tf_\go > 0$ is the free energy and its explicit definition
in terms of~$\go$ is given in~\S~\ref{sec:as_Z_loc}, while~$C^>_{\go,\eta}$,
$C^<_{\go,\eta}$ and $C^=_{\go,\eta}$ are
explicit positive constants, depending on $\go$ and $\eta$,
whose value is given in Section~\ref{sec:as_Z}.
\end{theorem}

\medskip
\begin{rem}
\label{sec:compdis}
\rm
Theorem~\ref{th:as_Z}
is the building block of all the path analysis that follows. It is therefore
important to stress that, in the quenched disordered case,
cf. Remark~\ref{rem:intper},  such a strong statement in general
does not hold, see \cite[Section 4]{cf:GT}.
\end{rem}
\medskip

Next we investigate the thermodynamic limit, that is the
weak limit as~$N\to\infty$ of the sequence of
measures~$\bP^{a}_{N,\go}$ on~$\Z^\N$ (endowed with the
standard product topology). The next theorem provides a
first connection between free energy (de)localization
and the corresponding path properties.

Before stating the result, we need a notation: we denote
by $\cP$ the set of $\omega$ such that:
\begin{equation}\label{sit1}
\cP \, := \, \left\{ \go \, : \ \gd^\go \leq 1, \quad
h_\go=0, \quad \exists \ \ga,\gb: \ \Sigma_{\ga,\gb} \ne 0
\right\},
\end{equation}
\[
\cP^< \, := \, \cP \cap\{\gd^\go < 1\},
\qquad \cP^= \, := \, \cP \cap\{\gd^\go = 1\}.
\]
Here $\cP$ stands for {\sl problematic}, or {\sl pathologic}.
Indeed, we shall see that for $\go\in\cP$ the results are
weaker and more involved than for $\go\notin\cP$. We stress however
that these restrictions do not concern localized regime,
because $\cP\subset\{\go:\delta^\go\leq 1\}$. We also notice that
for the two motivating models of \S~\ref{sec:2examples}, the pinning
and the copolymer models, $\go$ {\it never} belongs to
$\cP$. This is clear for the pinning case,
where by definition $\Sigma\equiv 0$. On the other hand, in the copolymer case
it is known that if $h_\go=0 $ and $\exists \ \ga,\gb: \ \Sigma_{\ga,\gb} \ne 0$
then $\delta^\go>1$: see \S~\ref{app:loc} or~\cite{cf:BG}.
In reality
 the {\sl pathological} aspects
observed for $\go \in \cP$
may be understood in statistical mechanics terms
and we sketch an interpretation  in   \S~\ref{sec:P} below:
this goes rather far from the point of view adopted here,
since it is an issue
tightly entangled with the analysis of the free energy.
It will therefore be taken up in a further work.

\medskip

\begin{theorem}[\bf Thermodynamic limit]\label{th:infvol}
If $\go\notin\cP^<$, then
both the polymer measures $\bP^\rf_{N, \go}$ and
$\bP^\rc_{N, \go}$ converge as $N\to\infty$ to the same limit
$\bP_\go$, law of an irreducible Markov process on $\Z$ which is:
\smallskip
\begin{enumerate}
\item positive recurrent if $\gd^\go>1$ (\textsl{localized regime})\,;
\medskip
\item transient if $\gd^\go<1$
(\textsl{strictly delocalized regime})\,;
\medskip
\item null recurrent if $\gd^\go=1$
(\textsl{critical regime})\,.
\end{enumerate}
\smallskip
If $\go\in\cP^<$ (in particular $\gd^\go < 1$), for all $\eta\in\bbS$ and $a=\rf,\rc$ the measure
$\bP^a_{N, \go}$ converges as $N\to\infty$, $[N]=\eta$ to
$\bP^{a,\eta}_{\go}$, law of an irreducible transient Markov
chain on $\Z$.
\end{theorem}

\medskip\noindent
We stress that in all regimes the limit law~$\bP_\go$ or
$\bP^{a,\eta}_{\go}$ has an explicit
construction in terms of~$M^\go_{\ga,\gb}(x)$, see
Section~\ref{sec:infvol} for details.

%
%
%
%
%

\bigskip

We finally turn to the analysis of the \textsl{diffusive rescaling}
of the polymer measure~$\bP_{N,\go}^a$. More precisely, let us
define the map $X^N: {\mathbb R}^N \mapsto C([0,1])$:
\[
X^N_t(x) \, = \, \frac{ x_{\lfloor Nt\rfloor}}{\sigma N^{1/2}} + (Nt-\lfloor Nt\rfloor) \,
\frac {x_{\lfloor Nt\rfloor+1}-x_{\lfloor Nt\rfloor}}{\sigma N^{1/2}}, \qquad t\in[0,1],
\]
where $\lfloor \,\cdot\, \rfloor$ denotes the integer part and $\gs^2 := 2p$ is the variance
of~$X_1$ under the original random walk measure~$\bP$.
Notice that $X^N_t(x)$ is nothing but the linear interpolation of~$\{x_{\lfloor Nt \rfloor}/(\gs\sqrt{N})\}_{t \in \frac{\N}{N} \cap [0,1]}$. For $a=\rf,\rc$ we set:
\[
Q^{a}_{N,\go} \, := \, {\bf P}^{a}_{N,\go} \circ (X^N)^{-1},
\]
Then $Q_{N,\go}^a$ is a measure on $C([0,1])$, the space of real
continuous functions defined on the interval~$[0,1]$, and we want
to study the behavior as~$N\to\infty$ of this sequence of measures.

\smallskip

We start fixing a notation for the following standard processes:
\smallskip
\begin{itemize}
\item the Brownian motion $\left\{B_\tau\right\}_{\tau\in [0,1]}$;
\smallskip
\item the Brownian bridge $\left\{\beta_\tau\right\}_{\tau\in
[0,1]}$ between $0$ and $0$;
\smallskip
\item the Brownian motion {\sl conditioned to stay
non-negative on $[0,1]$} or, more precisely, the Brownian meander
$\{m_{\tau}\}_{\tau\in[0,1]}$, see \cite{cf:RevYor};
\smallskip
\item the Brownian bridge {\sl conditioned to stay
non-negative on $[0,1]$} or, more precisely, the normalized
Brownian excursion $\{e_{\tau}\}_{\tau\in[0,1]}$, also known as the
Bessel bridge of dimension 3 between $0$ and $0$, see
\cite{cf:RevYor} .
\end{itemize}
\smallskip
Then we introduce a modification of the above processes labeled by a
parameter $p \in [0,1]$:
\begin{itemize}
\item \rule{0pt}{12pt}the process $\{B^{(p)}_\tau\}_{\tau \in [0,1]}$
is the so--called {\sl skew Brownian motion of parameter~$p$},
cf.~\cite{cf:RevYor}. More explicitly, $B^{(p)}$ is a process such
that~$|B^{(p)}| = |B|$ in distribution, but in which the sign of each
excursion is chosen to be $+1$ (resp.~$-1$) with probability~$p$
(resp.~$1-p$) instead of~$1/2$. In the same way, the process
$\{\gb^{(p)}_\tau\}_{\tau \in [0,1]}$ is the skew Brownian bridge
of parameter~$p$. Notice that for $p=1$ we have $B^{(1)}=|B|$ and
$\beta^{(1)}=|\beta|$ in distribution.
\smallskip
\item the process $\{m^{(p)}_\tau\}_{\tau \in [0,1]}$ is defined by
\[
\bbP(m^{(p)}\in dw) \, := \, p \, \bbP(m\in dw) \, + \,
(1-p) \, \bbP(-m\in dw),
\]
i.e. $m^{(p)}=\sigma  m$, where $\bbP(\sigma=1)=1-\bbP(\sigma=-1)=p\,$
and $(m,\sigma)$ are independent. The process
$\{e^{(p)}_\tau\}_{\tau \in [0,1]}$ is defined in exactly the same manner.
For $p=1$ we have $m^{(1)}=m$ and $e^{(1)}=e$.
\end{itemize}
\smallskip
Finally, we introduce a last process, labeled by two
parameters $p,q\in[0,1]$:
\begin{itemize}
\item \rule{0pt}{12pt}consider a r.v. $U\mapsto [0,1]$ with
the arcsin law: $\bbP(U\leq t)=\frac 2\pi\, \arcsin\sqrt t$,
and processes $\gb^{(p)}$, $m^{(q)}$ as defined above,
with $(U,\gb^{(p)},m^{(q)})$ independent triple.
Then we denote by $\{B^{(p,q)}_\tau\}_{\tau \in [0,1]}$
the process defined by:
\[
B^{(p,q)}_\tau \, := \,
\begin{cases}
\sqrt U \, \gb^{(p)}_{\frac\tau U}
& \text{if } \tau \leq U
\\
\rule{0pt}{24pt}
\sqrt{1-U} \, m^{(q)}_{\frac{\tau-U}{1-U}}
& \text{if } \tau > U
\end{cases} \;.
\]
Notice that the process $B^{(p,q)}$ differs from the $p$--skew Brownian motion
$B^{(p)}$ only for the last excursion in~$[0,1]$, whose sign is~$+1$ with
probability~$q$ instead of~$p$.
\end{itemize}

\medskip

We are going to show that the sequence $\{Q_{N,\go}^a\}$ has a
weak limit as~$N\to\infty$ (with a weaker statement if $\go\in\cP$).
Again the properties of the limit process differ considerably in
the three regimes $\gd^\go>1$, $\gd^\go<1$ and $\gd^\go=1$.
However for the precise description of the limit processes, for the regimes
$\gd^\go = 1$ and~$\gd^\go < 1$ we need to distinguish between
$a\in\{\rf,\rc\}$ and to introduce
further parameters~$\tt p_\go, \tt q_\go$, defined as follows:
\begin{itemize}
\item case $\gd^\go=1$\, :
\begin{itemize}
\item $\tt p_\go := \tt p_\go^=$,  defined in \eqref{eq:p_om3}. We point out two special cases:
if~$h_\go>0$ then~$\tt p^=_\go =1$, while if~$h_\go=0$
and~$\Sigma \equiv 0$ then~$\tt p^=_\go = 1/2$;
\item for each $\eta\in\bbS$,
$\tt q_\go := \tt q_{\go,\eta}^=$, defined by \eqref{eq:p_om4}.
\end{itemize}
%
%
\item case $\gd^\go<1$\, :
\begin{itemize}
\item
 $\go \notin \cP^<$:\: if~$h_\go>0$ we set
$\tt p_\go:=\tt p_\go^<:=1$
while if~$h_\go = 0\,$ we set $\tt p_\go:=\tt p_\go^< := 1/2$;
\item  $\go \in \cP^<$:\: for each $\eta\in\bbS$
and $a=\rf,\rc$,
${\tt p}_\go:={\tt p}_{\go, \eta}^{<, a}$ is defined in \eqref{eq:p_om}
and \eqref{eq:p_om2}.
\end{itemize}
\end{itemize}

\medskip
\begin{theorem}[\bf Scaling limits]
\label{th:scaling}
If $\go\notin\cP$, then
the sequence of measures~$\{Q_{N,\go}^a\}$ on $C([0,1])$ converges weakly
as~$N\to\infty$. More precisely:
\smallskip
\begin{enumerate}
\item for $\delta^\go >1$ (\textsl{localized regime}) $Q^a_{N,\go}$ converges to the measure
concentrated on the constant function taking the value zero\,;
\medskip
\item for $\delta^\go<1$ (\textsl{strictly delocalized regime}):
\begin{itemize}
\item
$Q^\rf_{N,\go}$ converges to the law of $m^{(\rm p_\go^<)}$\,;
\item
$Q^\rc_{N,\go}$ converges to the law of $e^{(\rm p_\go^<)}$\,;
\end{itemize}
\medskip
\item for $\delta^\go=1$ (\textsl{critical regime}):
\begin{itemize}
\item
$Q^\rf_{N,\go}$ converges to the law of $B^{{(\rm p_\go^=)}}$\,;
\item $Q^\rc_{N,\go}$ converges to the law of $\beta^{(\rm p_\go^=)}$\,.
\end{itemize}
\smallskip
\end{enumerate}
If $\go\in\cP$, then for all $\eta\in\bbS$ the measures
$Q^\rc_{N, \go}$ and $Q^\rf_{N, \go}$ converge as $N\to\infty$,
$[N]=\eta$ to, respectively:
\begin{enumerate}
\item for $\delta^\go<1$,
the law of $e^{{(\tt p_{\go,\eta}^{<,\rc})}}$ and $m^{{(\tt p_{\go,\eta}^{<,\rf})}}$.
\item for $\delta^\go=1$,
the law of $\beta^{{(\tt p_{\go}^{=})}}$ and $B^{{(\tt p_\go^=,\tt q_{\go,\eta}^=)}}$.
\end{enumerate}
\end{theorem}

\medskip

Results on thermodynamic limits in the direction of Theorem \ref{th:infvol}
have been obtained in the physical literature by exact computations either for homogeneous
polymers or for $T=2$ pinning models and copolymers, see e.g. \cite{cf:MGO},
while Brownian scaling limits have been heuristically derived
at several instances, see e.g. \cite{cf:Upton}.
Rigorous results corresponding to our three main theorems have been obtained
for homogeneous pinning/wetting models
in \cite{cf:DGZ,cf:IY}. We would like to stress the very much richer variety
of limit processes that we have obtained in our general context.

\smallskip

\subsection{About the regime $\cP$}
\label{sec:P}
We have seen, cf. Theorem~\ref{th:infvol},
that if $\go \in \cP^<$ the infinite volume limit (in particular
 the probability that the walk
escapes either to $+\infty$ or to $-\infty$) depends on
$a=\rc$ or $\rf$ and on the subsequence $[N]=\eta \in \bbS$.
 This reflects  directly into
Theorem~\ref{th:scaling} and  in this case also the $\cP^=$ regime is affected,
but only for $a=\rf$ and the change is restricted to the
sign of the very last excursion of the process.
It is helpful to keep in mind that $\go \in \cP$ if and only if
there is a non trivial unbiased copolymer part, that is $h_\go=0$ but
the matrix $\gS$ is non trivial, and at the same time the polymer is delocalized.
It is known (\S~\ref{app:loc} and \cite{cf:BG}) that in absence of
 pinning terms, that is $\go^{(0)}_n= \tilde \go^{(0)}_n=0$ for every $n$,
  the polymer is localized. However if the pinning rewards are
 sufficiently large and negative, one easily sees that (de)pinning takes over
 and the polymer delocalizes. This is the phenomenon that characterizes
 the regime $\cP$ and its lack of uniqueness of limit measures.
\smallskip

\begin{figure}[h]
\begin{center}
\leavevmode
\epsfysize =7 cm
\psfragscanon
\psfrag{D}[c][c]{\Large $\cD$}
\psfrag{L}[c][c]{ \Large $\cL$}
\psfrag{r1}[c][c]{$\varrho=+1$}
\psfrag{r2}[c][c]{ $\varrho=-1$}
\psfrag{0}[c][c]{$0$}
\psfrag{h}[c][c]{ $h$}
\psfrag{bc}[c][c]{$\gb _c$}
\psfrag{b}[c][c]{$\gb$}
\epsfbox{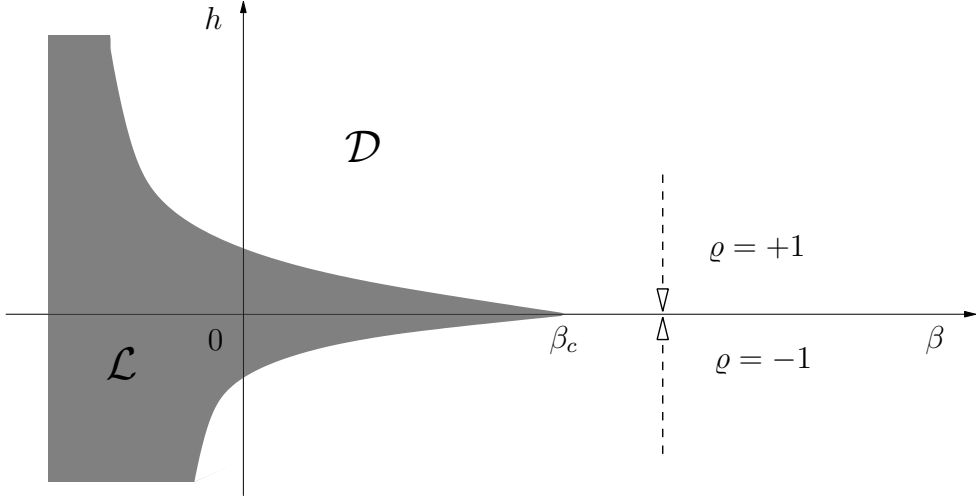}
\end{center}
\caption{\label{fig:P}
A sketch of the phase diagram for  the model \eqref{eq:Pfig}.
In this case, with abuse of notation, $\cP=\{ (h, \gb): \, h=0, \, \gb \ge \gb_c\}$.
Approaching $\cP$ in the sense of the dashed arrowed lines one
observes the two sharply different behaviors of paths completely delocalized
above ($\varrho=+1$) or below ($\varrho=-1$) the interface.
}
\end{figure}

Lack of uniqueness of infinite volume measures and   dependence
on  boundary conditions do not come as a surprise if one takes
a statistical mechanics viewpoint and if one notices that the system undergoes
a {\sl first order} phase transition exactly at $\cP $. In order to be more precise
let us consider the particular case of
\begin{equation}
\label{eq:Pfig}
 \frac{\dd \bP_{N, \go}}{\dd \bP} (S) \, \propto \, \exp \left( \sum _{n=1}^ {N}\left(
 \go_n +h \right)\sign \left( S_n\right) -\gb \sum _{n=1}^ {N}
  \ind_{\{ S_n=0\}}\right),
\end{equation}
with $h$ and $\gb$ two real parameters and $\go$ a fixed non trivial centered
($\sum_{n=1}^T\go_n=0$) periodic configuration of charges.
The phase diagram of such a model is sketched in Figure~\ref{fig:P}.
In particular it is easy to show that for $h=0$ and  for $\gb$ large and positive
the polymer is delocalized and, recalling that for $\gb=0$
the polymer is localized, by monotonicity of the free energy in $\gb$ one immediately
infers that there exists $\gb_c>0$ such that localization prevails for $\gb<\gb_c$,
while the polymer is delocalized (both above and below the interface)
if $\gb\ge \gb_c$. However the two regimes of delocalization above
or below the interface, appearing for example
as soon as $h$ is either positive or negative and $\gb \ge \gb_c$,
are characterized by opposite values ($\pm 1$) of $\varrho=\varrho (h,\gb):=
\lim_{N\to \infty}\bE_{N, \go}
\left[ N^{-1}\sum _{n=1}^ {N}
 \sign \left( S_n\right) \right]$ and of course $\varrho$ is the derivative of the
 free energy with respect to $h$. Therefore the free energy is not
 differentiable at $h=0$ and we say that there is a {\sl first order phase transition}.
First order phase transitions are usually associated to multiple infinite
volume limits ({\sl phase coexistence}).
A detailed analysis of this interesting phenomenon
will be given elsewhere.

\smallskip

\subsection{Outline of the exposition}

In Section~\ref{sec:as_Z} we study the asymptotic behavior of $Z^\rc_{N,\go}$,
proving Theorem~\ref{th:as_Z}. In Section 3 we compute the thermodynamic
limits of $\bP_{N,\go}^a$, proving Theorem~\ref{th:infvol}. In Section
4 we compute the scaling limits of $\bP_{N,\go}^a$, proving Theorem~\ref{th:scaling}.
Finally, in Section~\ref{app_cgz} we give the proof of some technical
results and some additional material.


\smallskip
\section{Sharp asymptotic behavior of the partition function}
\label{sec:as_Z}

In this section we are going to derive the precise asymptotic behavior of~$Z_{N,\go}^\rc$, in particular proving Theorem~\ref{th:as_Z}. The key observation is that the study of the partition function for the models we are considering can be set into the framework of the theory of \textsl{Markov renewal processes}, see~\cite[Ch.~VII.4]{cf:Asm}. We start recalling the basic notions of this theory and setting the relative notation.

\smallskip
\subsection{Markov Renewal Theory}

Given a finite set $\bbS$ (for us it will always be $\Z/(T\Z)$), by a \textsl{kernel} we mean a family of nonnegative $\bbS\times\bbS$ matrices $F_{\ga,\gb}(x)$ depending on a parameter $x\in\N$. We say that the kernel $F_{\ga,\gb}(x)$ is \textsl{semi--Markov} if $F_{\ga,\cdot}(\cdot)$ is a probability mass function on $\bbS\times\N$ for every $\ga\in\bbS$, that is if $\sum_{\gb,x}F_{\ga,\gb}(x)=1$.

\smallskip

A semi--Markov kernel~$F_{\ga,\gb}(x)$ has a simple probabilistic interpretation: it defines a Markov chain $\{(J_k,T_k)\}$ on $\bbS\times\N$ through the transition kernel given by
\begin{equation} \label{eq:trans_kernel}
\bbP\big[(J_{k+1},T_{k+1}) = (\gb,x) \:\big|\: (J_{k},T_{k}) = (\ga,y) \big] = F_{\ga,\gb}(x)\,.
\end{equation}
In this case we say that the process $\{J_k,T_k\}$ is a (discrete) \textsl{Markov--renewal process}, the~$\{T_k\}$ being thought of as interarrival times. This provides a generalization of classical renewal processes, since the $\{T_k\}$ are no longer IID but their laws are rather \textsl{modulated} by the process~$\{J_k\}$. Since the r.h.s. of \eqref{eq:trans_kernel} does not depend on~$y$, it follows that $\{J_k\}$ is a Markov chain, and it is called the \textsl{modulating chain} of the Markov renewal process (observe that in general the process $\{T_k\}$ is \textsl{not} a Markov chain). The transition kernel of~$\{J_k\}$ is given by $\sum_{x\in\N} F_{\ga,\gb}(x)$. We will assume that this chain is irreducible (therefore positive recurrent, since $\bbS$ is finite) and we denote by $\{\nu_\ga\}_{\ga\in\bbS}$ its invariant measure.

\smallskip

Given two kernels $F$ and $G$, their convolution $F*G$ is the kernel defined by
\begin{equation} \label{eq:convolution}
(F*G)_{\ga,\gb}(x) := \sum_{y\in\N} \sum_{\gamma\in\bbS} F_{\ga,\gamma}(y) G_{\gamma,\gb}(x-y) = \sum_{y\in\N} \big[ F(y) \cdot G(x-y) \big]_{\ga,\gb} \,,
\end{equation}
where $\cdot$ denotes matrix product. Observe that if $F$ and $G$ are semi--Markov kernels, then $F*G$ is semi--Markov too. With standard notation, the $n$--fold convolution of a kernel~$F$ with itself will be denoted by~$F^{*n}$, the $n=0$ case being by definition the identity kernel $[F^{*0}]_{\ga,\gb}(x):=\ind_{(\gb=\ga)} \ind_{(x=0)}$.

\smallskip

A fundamental object associated to a semi--Markov kernel $F$  the so--called \textsl{Markov--Green function} (or Markov--renewal kernel), which is the kernel $\bs U$ defined by
\begin{equation} \label{eq:renewal_kernel}
\bs{U}_{\ga,\gb}(x) := \sum_{k=0}^\infty \big[ F^{*k} \big]_{\ga,\gb}(x)\,.
\end{equation}
Of course the kernel $\bs U$ is the analog of the Green function of a classical renewal process, and it has a similar probabilistic interpretation in terms of the associated Markov renewal process~$\{(J_k,T_k)\}$:
\begin{equation}\label{eq:prob_interpret}
    \bs U_{\ga,\gb} (x) = \bbP_\ga \big[ \exists k \ge 0 :\; T_0+\ldots +T_k=x\;,\; J_{k} = \gb \big]\,,
\end{equation}
where $\bbP_\ga$ is the law of $\{(T_k,J_k)\}$ conditioned on $\{J_0=\ga, T_0 = 0\}$.

\smallskip

We need some notation to treat our periodic setting: we say that a kernel~$F_{\ga,\gb}(x)$ has period $T\in\N$ if the set $\{x: \bs{U}_{\ga,\ga}(x) \neq 0\}$ is contained in $T\Z$, for the least such~$T$ (this definition does not depend on $\ga$ because the chain $\{J_k\}$ is supposed to be irreducible, see the discussion at p.~208 of \cite{cf:Asm}). It follows that the set $\{x: \bs{U}_{\ga,\gb}(x) \neq 0\}$ is contained in the translated lattice $\gamma(\ga,\gb) + T\N$, where $\gamma(\ga,\gb)\in\{0,\ldots,T-1\}$ (for us it will be $\gamma(\ga,\gb) = [\gb-\ga]$).

\smallskip

In analogy to the classical case, the asymptotic behavior of $\bs U_{\ga,\gb}(x)$
as $x\to\infty$ is of particular interest. Let us define the (possibly infinite)
\textsl{mean} $\mu$ of a semi--Markov kernel $F_{\ga,\gb}(x)$ as
\begin{equation} \label{eq:mu}
\mu := \sum_{\ga,\gb \in \bbS} \sum_{x\in\N} x\, \nu_\ga\, F_{\ga,\gb}(x) \,.
\end{equation}
Then we have an analog of Blackwell's Renewal Theorem, that in our periodic setting reads as
\begin{equation} \label{eq:per_renewal_theorem}
\exists \limtwo{x\to\infty}{[x] = \gb-\ga} \bs{U}_{\ga,\gb}(x) \;=\; T \frac{\nu_\gb}{\mu}\,,
\end{equation}
cf. Corollary~2.3 p.~10 of \cite{cf:Asm} for the classical case.

\smallskip

We will see that determining the asymptotic behavior of $\bs U_{\ga,\gb}(x)$ when the kernel $F_{\ga,\gb}(x)$ is no more semi--Markov is the key to get the asymptotic behavior of the partition function~$Z_{N,\go}^\rc$.


\smallskip
\subsection{A random walk excursion viewpoint}
\label{sec:RWexcurs}

Now we are ready to make explicit the link between the partition function for our model and the Theory of Markov Renewal Processes. Let us look back to our Hamiltonian \eqref{eq:genH}: its specificity comes from
the fact that it can be decomposed in an efficient way by considering the
return times to the origin of $S$. More precisely we set for~$j \in \N$
\begin{equation*}
    \tau_0=0 \qquad \quad \tau_{j+1}=\inf\{ n>\tau_j:\, S_n=0 \}\,,
\end{equation*}
and for $\bP$--typical trajectories of $S$ one has an infinite sequence $\tau:=\{\tau_j\}_j$ of stopping
times. We set $T_j= \tau_j-\tau_{j-1}$ and of course
$\{ T_j\}_{j=1, 2, \ldots}$ is, under~$\bP$,  an IID sequence.
By conditioning on $\tau$ and integrating on the up--down symmetry
of the random walk excursions one easily obtains the following expression for the constrained partition function:
\begin{equation} \label{eq:step1}
Z_{N,\go}^{\rc}\, = \,
\bE \left[
\prod_{j=1}^{\iota_N} \exp\big( \Psi^\go(\tau_{j-1}, \tau_j)\big);
\, \tau_{\iota_N}=N
\right],
\end{equation}
where $\iota_N = \sup\{k:\, \tau_k\le N\}$ and we have introduced
the \textsl{integrated Hamiltonian} $\Psi^\go (n_1,n_2)$, which gives the energetic contribution of an excursion from~$n_1$ to~$n_2$:
\begin{equation}
\label{eq:def_Psi}
\Psi^\go (n_1,n_2)=\begin{cases}
    \go ^{(0)}_{n_2} \,+\, \big( \tilde  \go ^{(0)}_{n_2} - \go^{(+1)}_{n_2} \big) & \text{if } n_2=n_1+1 \\
    \rule{0pt}{24pt} \go ^{(0)}_{n_2} \,+\, \displaystyle \log \bigg[ \frac 12 \Big( 1 + \exp \sum_{n=n_1+1}^{n_2} \big(\go_n ^{(-1)}-\go_n ^{(+1)}\big)
\Big) \bigg]  & \text{if } n_2>n_1 +1 \\
    \rule{0pt}{11pt} 0 & \text{otherwise}\,.
\end{cases}
\end{equation}

\smallskip

Now we are going to use in an essential way the fact that our charges are~$T$--periodic. In fact a look at \eqref{eq:def_Psi} shows that the energy $\Psi^\go (n_1,n_2)$ of an excursion from~$n_1$ to~$n_2$ is
a function only of $(n_2-n_1)$, $[n_1]$ and~$[n_2]$, where by~$[\,\cdot\,]$ we mean the equivalence class modulo~$T$, see \S~\ref{sec:order_par}. More precisely for $n_1 \in \ga$, $n_2 \in \gb$ and $\ell = n_2-n_1$ we have $\Psi^\go(n_1,n_2) = \Phi^\go_{\ga,\gb} (\ell)$, where~$\Phi^\go$ was defined in~\eqref{eq:def_Phi}. Then recalling the law~$K(n)$ of the first return, introduced in~\eqref{eq:first_ret}, we can rewrite \eqref{eq:step1} as
\begin{equation} \label{eq:zero_set_dec}
Z_{N,\go}^{\rc}\, = \, \sum_{k=1}^N \sumtwo{t_0,\ldots,t_k \in \N}{0=:t_0 < t_1 < \ldots < t_k:=N} \, \prod_{j=1}^k K\left( t_j - t_{j-1}\right) \; \exp\big(\Phi ^\go_{[t_{j-1}], [t_j]} (t_j - t_{j-1}) \big)\,.
\end{equation}
This decomposition of $Z_{N,\go}^{\rc}$ according to the random walk excursions makes explicit the link with Markov Renewal Theory. In fact using the kernel $M_{\ga,\gb}(x)$ introduced in \eqref{eq:matrix_cgz} we can rewrite it as
\begin{equation} \label{eq:conv1}
\begin{split}
    Z_{N,\go}^{\rc} &= \sum_{k=1}^N \sumtwo{t_0,\ldots,t_k \in \N}{0=:t_0 < t_1 < \ldots < t_k:=N} \, \prod_{j=1}^k M_{[t_{j-1}],[t_j]}(t_j-t_{j-1}) \\
    &= \sum_{k=1}^N \sumtwo{t_0,\ldots,t_k \in \N}{0=:t_0 < t_1 < \ldots < t_k:=N} \big[ M(t_1) \cdot M(t_2-t_{1}) \cdot \ldots \cdot M(N-t_{k-1}) \big]_{0,[N]} \\
    &= \sum_{k=0}^\infty \big[ M^{*k} \big] _{[0],[N]}(N)\,.
\end{split}
\end{equation}
Therefore it is natural to introduce the kernel $\cZ_{\ga,\gb}(x)$ defined by
\begin{align} \label{eq:expr_zeta1}
    \cZ_{\ga,\gb}(x) = \sum_{k=0}^\infty \big[ M^{*k} \big] _{\ga,\gb}(x)\,,
\end{align}
so that $Z_{N,\go}^{\rc} = \cZ_{[0],[N]}(N)$. More generally $\cZ_{\ga,\gb}(x)$ for $[x] = \gb - \ga$ can be interpreted as the partition function of a directed polymer of size~$x$ that starts at a site $(M,0)$, with $[M]=\ga$, and which is pinned at the site~$(M+x,0)$.

\smallskip

Our purpose is to get the precise asymptotic behavior of $\cZ_{\ga,\gb}(x)$ as~$x\to\infty$, from which we will obtain the asymptotic behavior of $Z_{N,\go}^{\rc}$ and hence the proof of Theorem~\ref{th:as_Z}.
It is clear that equation \eqref{eq:expr_zeta1} is the same as equation \eqref{eq:renewal_kernel}, except for the fact that in general the kernel~$M$ has no reason to be semi--Markov. Nevertheless we will see that with some transformations one can reduce the problem to a semi--Markov setting.

\smallskip

It turns out that for the derivation of the asymptotic behavior of~$\cZ_{\ga,\gb}(x)$ it is not necessary to use the specific form~\eqref{eq:matrix_cgz} of the kernel~$M_{\ga,\gb}(x)$, the computations being more transparent if carried out in a general setting. For these reasons, in the following we will assume that $M_{\ga,\gb}(x)$ is a generic $T$--periodic kernel such that the matrix~$B_{\ga,\gb}$ defined by~\eqref{eq:def_B} is finite. While these assumption are sufficient to yield the asymptotic behavior of~$\cZ_{\ga,\gb}(x)$ when $\gd^\go > 1$, for the cases $\gd^\go < 1$ and $\gd^\go = 1$ it is necessary to know the asymptotic behavior as~$x\to\infty$ of~$M_{\ga,\gb}(x)$ itself. Notice that our setting is an \textsl{heavy--tailed} one: more precisely we will assume that
for every $\ga,\gb\in\bbS$:
\begin{equation} \label{eq:asympt_cgz}
\exists \limtwo{x\to\infty}{[x]=\gb-\ga} x^{3/2}\, M_{\ga,\gb}(x) \;=:\; L_{\ga,\gb} \in (0,\infty) \,.
\end{equation}
From equation \eqref{eq:def_Phi} it is easy to check that the kernel $M_{\ga,\gb}(x)$ defined by \eqref{eq:matrix_cgz} does satisfy \eqref{eq:asympt_cgz} (see Section~\ref{sec:infvol} for more details on this issue).

\smallskip

For ease of exposition, we will treat separately the three cases $\gd^\go>1$, $\gd^\go <1$ and $\gd^\go = 1$.

\smallskip
\subsection{The localized regime $(\gd^\go > 1)$}
\label{sec:as_Z_loc}

The key idea is to introduce the following exponential perturbation of the kernel $M$ (cf. \cite[Theorem 4.6]{cf:Asm}), depending on the positive real parameter~$b$:
\begin{equation*}
    A^{b}_{\ga,\gb}(x) := M_{\ga,\gb}(x)\, e^{-bx}\,.
\end{equation*}
Let us denote by $\Delta(b)$ the Perron--Frobenius eigenvalue of the matrix $\sum_x A_{\ga,\gb}^{b}(x)$. As the entries of this matrix are analytic and nonincreasing functions of~$b$, $\Delta(b)$ is analytic and nonincreasing too, hence strictly decreasing because $\Delta(0) = \gd^\go > 1$ and $\Delta(\infty) = 0$. Therefore there exists a single value $\tf_\go>0$ such that $\Delta\big(\,\tf_\go\,\big)=1$, and we denote by $\{\zeta_\ga\}_\ga$, $\{\xi_\ga\}_\ga$ the Perron--Frobenius left and right eigenvectors of $\sum_x A_{\ga,\gb}^{\tf_\go}(x)$, chosen to have (strictly) positive components and normalized in such a way that $\sum_\ga \zeta_\ga\, \xi_\ga = 1$ (of course there is still a degree of freedom in the normalization, which however is immaterial).

\smallskip

Now we set
\begin{equation} \label{eq:def_Gamma>}
    \Gamma^>_{\ga,\gb} (x) := A^{\tf_\go}_{\ga,\gb} (x) \, \frac{\xi_\gb}{\xi_\ga} = M_{\ga,\gb} (x)\, e^{-\tf_\go x}\, \frac{\xi_\gb}{\xi_\ga}\,,
\end{equation}
and it is immediate to check that $\Gamma^>$ is a semi--Markov kernel. Furthermore, we can rewrite~\eqref{eq:expr_zeta1} as
\begin{equation} \label{eq:Z_loc}
    \cZ_{\ga,\gb}(x) := e^{\tf_\go x}\, \frac{\xi_\ga}{\xi_\gb}\, \sum_{k=0}^\infty \big[ (\Gamma^>)^{*k} \big] _{\ga,\gb}(x) = e^{\tf_\go x}\, \frac{\xi_\ga}{\xi_\gb}\, \cU_{\ga,\gb}(x) \,,
\end{equation}
where $\cU_{\ga,\gb}(x)$ is nothing but the Markov--Green function associated to the semi--Markov kernel~$\Gamma^>_{\ga,\gb}(x)$. Therefore the asymptotic behavior of $\cZ_{\ga,\gb}(x)$ is easily obtained applying Blackwell's Renewal Theorem~\eqref{eq:per_renewal_theorem}. To this end, let us compute the mean $\mu$ of the semi--Markov kernel $\Gamma^>$: it is easily seen that the invariant measure of the associated modulating chain is given by $\{\zeta_\ga \xi_\ga\}_\ga$, therefore
\begin{align*}
    \mu & = \sum_{\ga,\gb \in \bbS} \sum_{x\in\N} x \, \zeta_\ga \, \xi_\ga \, \Gamma_{\ga,\gb}^>(x) = \sum_{\ga,\gb \in \bbS} \sum_{x\in\N} x e^{-\tf_\go x} \, \zeta_\ga\, M_{\ga,\gb}(x)\, \xi_\gb\\
    & = - \bigg( \frac{\partial}{\partial b} \Delta(b) \bigg) \bigg|_{b = \tf_\go} \ \ \in \ (0,\infty)\,,
\end{align*}
(for the last equality see for example \cite[Lemma 2.1]{cf:BG}). Coming back to \eqref{eq:Z_loc}, we can now apply Blackwell's Renewal Theorem \eqref{eq:per_renewal_theorem} obtaining the desired asymptotic behavior:
\begin{equation}\label{eq:as_Z_loc}
    \cZ_{\ga,\gb}(x) \; \sim \; \xi_\ga \, \zeta_\gb \,  \frac{T}{\mu} \,
\exp \left( \tf_\go\, x \right) \qquad \quad x\to\infty,\quad [x]=\gb-\ga\,.
\end{equation}
In particular, for $\ga=[0]$ and $\gb=\eta$ we have part~(1) of Theorem~\ref{th:as_Z},
where $C^>_{\go,\eta} = \xi_0 \zeta_\eta T /\mu$.

\smallskip
\subsection{The strictly delocalized case $(\gd^\go < 1)$}
\label{sec:as_Z_del}

We prove that the asymptotic behavior of~$\cZ_{\ga,\gb}(x)$ when $\gd^\go < 1$ is given by
\begin{equation}\label{eq:as_Z_del}
    \cZ_{\ga,\gb}(x) \; \sim \; \Big( \big[(1-B)^{-1} L \, (1-B)^{-1} \big]_{\ga,\gb} \Big) \: \frac{1}{x^{3/2}} \qquad \quad x\to\infty,\quad [x]=\gb-\ga \,,
\end{equation}
where the matrixes $L$ and $B$ have been defined in \eqref{eq:asympt_cgz} and \eqref{eq:def_B}.
In particular, taking $\ga=[0]$ and $\gb=\eta$, \eqref{eq:as_Z_del} proves part~(2) of
Theorem~\ref{th:as_Z} with
\begin{equation*}
    C^<_{\go,\eta} := \big[(1-B)^{-1} L \, (1-B)^{-1} \big]_{0,\eta} \,.
\end{equation*}

\smallskip

To start with, we prove by induction that for every $n\in\N$
\begin{equation} \label{eq:sum}
\sum_{x\in\N} [M^{*n}]_{\ga,\gb} (x) = [B^n]_{\ga,\gb}\,.
\end{equation}
The $n=1$ case is the definition of~$B$, while for $n\geq 1$
\begin{align*}
\sum_{x\in\N} M^{*(n+1)}(x) & \;=\; \sum_{x\in\N} \sum_{z\le x} M^{*n}(z) \cdot M(x-z) \;=\; \sum_{z\in\N}\, M^{*n}(z) \cdot \sum_{x\ge z}\, M(x-z) \\
& \;=\; \sum_{z\in\N}\, M^{*n}(z) \cdot B \;=\; B^n \cdot B \;=\; B^{n+1} \,.
\end{align*}

\smallskip

Next we claim that, if \eqref{eq:asympt_cgz} holds, then for every $\ga,\gb \in \bbS$
\begin{equation} \label{eq:as_M_k}
\exists \limtwo{x\to\infty}{[x]=\gb-\ga} \, x^{3/2} \big[M^{*k}\big]_{\ga,\gb}(x) \;=\; \sum_{i=0}^{k-1} \big[B^i\cdot L \cdot B^{(k-1)-i}\big]_{\ga,\gb}\,.
\end{equation}
We proceed by induction on~$k$. The $k=1$ case is given by~\eqref{eq:asympt_cgz}, and we have that
\[
M^{*(n+1)}(x) = \sum_{y=1}^{x/2} \bigg( M(y)\cdot M^{*n}(x-y) \;+\; M(x-y) \cdot M^{*n}(y) \bigg)
\]
(strictly speaking this formula is true only when $x$ is even, however the odd~$x$ case is analogous). By the inductive hypothesis equation \eqref{eq:as_M_k} holds for every $k\leq n$, and in particular this implies that $\{x^{3/2} [M^{*k}]_{\ga,\gb}(x)\}_{x\in\N}$ is a bounded sequence. Therefore we can apply Dominated Convergence and \eqref{eq:sum}, getting
\begin{align*}
    \exists \limtwo{x\to\infty}{[x]=\gb-\ga} & \, x^{3/2} \big[M^{*(n+1)}\big]_{\ga,\gb}(x) \\
    & \;=\; \sum_\gamma \sum_{y=1}^{\infty} \bigg( M_{\ga,\gamma}(y)\, \sum_{i=0}^{n-1} \big[B^i\cdot L \cdot B^{(n-1)-i}\big]_{\gamma,\gb} \;+\; L_{\ga,\gamma} \big[M^{*n}\big]_{\gamma,\gb}(y) \bigg)\\
& \;=\; \sum_\gamma  \bigg( B_{\ga,\gamma}\, \sum_{i=0}^{n-1} \big[B^i\cdot L \cdot B^{(n-1)-i}\big]_{\gamma,\gb} \;+\; L_{\ga,\gamma} \big[ B^{*n} \big]_{\gamma,\gb} \bigg)\\
& \;=\; \sum_{i=0}^{n} \big[B^i\cdot L \cdot B^{n-i}\big]_{\ga,\gb} \,.
\end{align*}

\smallskip

Our purpose is to apply the asymptotic result \eqref{eq:as_M_k} to the terms of \eqref{eq:expr_zeta1}, hence we need a bound to apply Dominated Convergence. What we are going to show is that
\begin{equation}\label{eq:dom_con}
    x^{3/2} \, \big[M^{*k}\big]_{\ga,\gb}(x) \;\le\; C\, k^3\, \big[ B^k \big]_{\ga,\gb}
\end{equation}
for some positive constant $C$ and for all $\ga,\gb\in\bbS$ and $x,k\in\N$. Observe that the r.h.s. above, as a function of~$k$, is a summable sequence because the matrix~$B$ has spectral radius $\gd^\go<1$. We proceed again by induction: for the $k=1$ case, thanks to \eqref{eq:asympt_cgz}, it is possible to find~$C$ such that \eqref{eq:dom_con} holds true (this fixes $C$ once for all). Now assuming that \eqref{eq:dom_con} holds for all $k < n$ we show that it does also for~$k=n$ (we suppose for simplicity that $n=2m$ is even, the odd~$n$ case being analogous). Then we have (assuming that also~$x$ is even for simplicity)
\begin{align*}
    x^{3/2} \, \big[M^{*2m}\big]_{\ga,\gb}(x) &\;=\; 2 \sum_{y=1}^{x/2} \sum_{\gamma \in \bbS} \big[M^{*m}\big]_{\ga,\gamma}(y)\; x^{3/2} \big[M^{*m}\big]_{\gamma,\gb}(x-y)\\
    &\;\le\; 2 \cdot 2^{3/2} \, C \, m^{3} \sum_{y=1}^{x/2} \sum_{\gamma \in \bbS} \big[M^{*m}\big]_{\ga,\gamma}(y)\; \big[B^{m}\big]_{\gamma,\gb} \\
    &\;\le\; C\, (2m)^3 \big[ B^{2m} \big]_{\ga,\gb}\,,
\end{align*}
where we have applied \eqref{eq:sum}, and \eqref{eq:dom_con} is proven.

\smallskip

We can finally obtain the asymptotic behavior of~$\cZ_{\ga,\gb}(x)$ applying the bound \eqref{eq:as_M_k} to \eqref{eq:expr_zeta1}, using Dominated Convergence thanks to \eqref{eq:dom_con}. In this way we get
\begin{align*}
    \exists \limtwo{x\to\infty}{[x]=\gb-\ga} & \; x^{3/2} \cZ_{\ga,\gb}(x) \;=\; \sum_{k=1}^\infty \sum_{i=0}^{k-1} \big[B^i\cdot L \cdot B^{(k-1)-i}\big]_{\ga,\gb} \\
    &\;=\; \sum_{i=0}^\infty \sum_{k=i+1}^\infty \big[B^i\cdot L \cdot B^{(k-1)-i}\big]_{\ga,\gb} \;=\; \sum_{i=0}^\infty \big[B^i\cdot L \cdot (1-B)^{-1}\big]_{\ga,\gb} \\
    &\;=\; \big[(1-B)^{-1}\cdot L \cdot (1-B)^{-1}\big]_{\ga,\gb}\,,
\end{align*}
and equation \eqref{eq:as_Z_del} is proven.

\smallskip

\subsection{The critical case $(\gd^\go = 1)$}
\label{sec:as_Z_cri}

In the critical case the matrix $B$ defined in~\eqref{eq:def_B} has Perron--Frobenius eigenvalue equal to~1. Let $\{\zeta_\ga\}_\ga$, $\{\xi_\ga\}_\ga$ denote its corresponding left and right eigenvectors, always chosen to have positive components and normalized so that $\sum_\ga \zeta_\ga\, \xi_\ga = 1$. Then it is immediate to check that the kernel
\begin{equation}\label{eq:def_Gamma=}
    \Gamma^=_{\ga,\gb}(x) \;:=\; M_{\ga,\gb} (x) \, \frac{\xi_\gb}{\xi_\ga}
\end{equation}
is semi--Markov, and the corresponding Markov--Green function $U_{\ga,\gb}(x)$ is given by
\begin{equation} \label{eq:def_U}
    U_{\ga,\gb}(x) \;:=\; \sum_{k=0}^\infty \big[ (\Gamma^=)^{*k} \big]_{\ga,\gb} (x) \;=\; \frac{\xi_\gb}{\xi_\ga} \, \cZ_{\ga,\gb}(x)\,,
\end{equation}
where the last equality follows easily from~\eqref{eq:expr_zeta1}. We are going to derive the asymptotic behavior of~$U_{\ga,\gb}(x)$, and from the above relation we will get the analogous result for~$\cZ_{\ga,\gb}(x)$.

\smallskip

Denoting by $\{(T_k,J_k)\}$ under~$\bbP$ the Markov--renewal process generated by the semi--Markov kernel $\Gamma^=_{\ga,\gb}(x)$, for $U_{\ga,\gb}(x)$ we have the probabilistic interpretation \eqref{eq:prob_interpret}, that we rewrite for convenience
\begin{equation} \label{eq:prob_interpret2}
    U_{\ga,\gb} (x) = \bbP_\ga \big[ \exists k \ge 0 :\; T_0+\ldots +T_k=x\;,\; J_{k} = \gb \big]\,.
\end{equation}
For $\gb \in \bbS$ we introduce the sequence of stopping times~$\{\kappa^{(\gb)}_n\}_{n \geq 0}$
corresponding to the visit of the chain~$\{J_k\}$ to the state~$\gb$:
\begin{equation}\label{defkappa}
    \kappa_0^{(\gb)} := \inf \{k \geq 0:\ J_k = \gb\} \qquad \quad \kappa^{(\gb)}_{n+1}
:= \inf \{k > \kappa^{(\gb)}_n:\ J_k = \gb\}\,,
\end{equation}
and we define the process $\{T^{(\gb)}_n\}_{n \geq 0}$ by setting
\begin{equation}\label{defTkappa}
    T_0^{(\gb)} := T_0 + \ldots + T_{\kappa^{(\gb)}_0} \qquad \quad T_{n}^{(\gb)} :T_{\kappa^{(\gb)}_{n-1} + 1} + \ldots + T_{\kappa^{(\gb)}_{n}} \,.
\end{equation}
The key point is that under~$\bbP_\ga$ the random variables $\{T^{(\gb)}_n\}$ are the interarrival times of a (possibly delayed) \textsl{classical renewal process}, equivalently the sequence $\{T^{(\gb)}_n\}_{n \geq 1}$ is IID and independent of~$T_0^{(\gb)}$. We denote for $x\in\N$ by $q^{(\gb)}(x)$ the (mass function of the) law of~$T_n^{(\gb)}$ for $n\geq 1$, while the law of~$T_0^{(\gb)}$ under~$\bbP_\ga$ is denoted by~$q^{(\ga;\gb)}(x)$. Since clearly
\begin{equation*}
    \big\{ \exists k \ge 0 :\; T_0+\ldots +T_k=x\;,\; J_{k} = \gb \big\} \quad \iff \quad \big\{ \exists n \ge 0 :\; T^{(\gb)}_0+\ldots +T^{(\gb)}_n=x \big\}\,,
\end{equation*}
from \eqref{eq:prob_interpret2} we get
\begin{equation} \label{eq:cl_ren_pr}
    U_{\ga,\gb} (x) = \bbP_\ga \big[ \exists n \ge 0 :\; T^{(\gb)}_0+\ldots +T^{(\gb)}_n=x  \big] = \bigg( q^{(\ga;\gb)} * \sum_{n=0}^\infty \big( q^{(\gb)} \big) ^{*n} \bigg) (x) \,,
\end{equation}
which shows that $U_{\ga,\gb}(x)$ is indeed the Green function of the classical renewal process whose interarrival times are the $\{T_n^{(\gb)}\}_{n\ge 0}$.

\smallskip

Now we claim that the asymptotic behavior of $q^{(\gb)}(x)$ as $x\to\infty$, $x\in\gb$, is given by
\begin{equation} \label{eq:cb1}
    q^{(\gb)}(x) \; \sim \; \frac{c_\gb}{x^{3/2}} \qquad \quad
c_\gb \;:=\; \frac1{\zeta_\gb\, \xi_\gb}
\, \sum_{\ga,\gamma} \zeta_\ga \, L_{\ga,\gamma}\, \xi_\gamma  \;>\;0\,,
\end{equation}
see \S~\ref{app:qb} for a proof of this relation. Then the asymptotic behavior of \eqref{eq:cl_ren_pr} is given by
\begin{equation}\label{eq:doney}
    U_{\ga,\gb} (x) \;\sim\; \frac{T^2}{2\pi\,c_\gb} \, \frac{1}{\sqrt{x}} \qquad \quad x\to\infty,\quad [x] = \gb-\ga\,,
\end{equation}
as it follows by~\cite[Th.~B]{cf:Don97} (the factor~$T^2$ is due to our periodic setting). Combining equations \eqref{eq:def_U}, \eqref{eq:cb1} and \eqref{eq:doney} we finally get the asymptotic behavior of~$\cZ_{\ga,\gb}(x)$:
\begin{equation}\label{eq:as_Z_cri}
    \cZ_{\ga,\gb}(x) \;\sim\; \frac{T^2}{2\pi} \, \frac{\xi_\ga\, \zeta_\gb}
{ \sum_{\gamma,\gamma'} \zeta_\gamma \, L_{\gamma,\gamma'}\, \xi_{\gamma'} }
\: \frac{1}{\sqrt{x}} \qquad \quad x\to\infty,\quad [x]=\gb-\ga \,.
\end{equation}
Taking $\ga=[0]$ and $\gb=\eta$, we have the proof of part~(3) of Theorem~\ref{th:as_Z}.


\smallskip
\section{Thermodynamic limits}
\label{sec:infvol}

In this section we study the limit as $N\to\infty$ of the polymer
measure $\bP^a_{N, \go}$, using the sharp asymptotics for the partition
function obtained in the previous section.
We recall that  $\bP^\rc_{N,\go}$ is a probability
measure on $\Z^\N$, which we endow with the product topology.
In particular, weak convergence on $\Z^\N$ means convergence
of all finite dimensional marginals.

\medskip

We start giving a very useful decomposition of $\bP^a_{N,\go}$.
The intuitive idea is that a path $(S_n)_{n\le N}$ can be split
into two main ingredients:
\smallskip
\begin{itemize}
\item the family $(\tau_k)_{k=0,1,\ldots}$ of \textsl{returns to zero}
of $S$ (defined in \S~\ref{sec:RWexcurs});
\item \rule{0pt}{13pt}the family of \textsl{excursions from zero}
$(S_{i+\tau_{k-1}}:0\le i\le \tau_{k}-\tau_{k-1})_{k=1, 2, \ldots}$
\end{itemize}
\smallskip
Moreover, since each excursion can be
either positive or negative, it is also useful to consider separately
the signs of the excursions $\sigma_k := \sign (S_{\tau_{k-1} + 1})$ and the
absolute values $(e_k(i):=|S_{i+\tau_{k-1}}|:\, i=1,\ldots,\tau_{k}-\tau_{k-1})$.
Observe that these are trivial for an
excursion with length~$1$: in fact if $\tau_{k}=\tau_{k-1} + 1$
then $\gs_k =0$ and~$e_k(0)=e_k(1)=0$.

\smallskip

Let us first consider the returns $(\tau_k)_{k\le\iota_N}$ under~$\bP_{N,\go}^a$,
where $\iota_N = \sup\{k:\, \tau_k\le N\}$.
The law of this process can be viewed as a probability measure $p^a_{N,\go}$
on the class $\cA_N$ of subsets of $\{1,\ldots, N\}$: indeed
for $A\in\cA_N$, writing
\begin{equation}\label{notA}
A  = \{t_1,\ldots,t_{|A|}\}, \qquad
0 \, =: \, t_0<t_1<\cdots<t_{|A|} \, \leq \, N,
\end{equation}
we can set
\begin{equation}\label{defp^a}
p^a_{N,\go}(A) \, := \, \bP^a_{N,\go}(\tau_i= t_i, \ i\leq\iota_N).
\end{equation}
The measure $p^a_{N,\go}$ describes the zero set of the polymer of size~$N$,
and it is analyzed in detail below. From the inclusion of $\cA_N$ into $\{0,1\}^\N$,
the family of all subsets of~$\N$, $p_{N,\go}^a$ can be viewed as a measure
on~$\{0,1\}^\N$ (this observation will be useful in the following).

\smallskip

Now we pass to the signs: we can see that, given $(\tau_j)_{j\leq\iota_N}$,
under $\bP^a_{N,\go}$ the signs $(\sigma_k)_{k\leq\iota_N}$ form
an independent family. Conditionally on~$(\tau_j)_{j\leq\iota_N}$, the law of~$\gs_k$ is specified by:
\begin{itemize}
\item[-] if~$\tau_{k} = 1+\tau_{k-1} $, then $\gs_k = 0$;
\item[-] \rule{0pt}{12pt}if~$\tau_{k} > 1+\tau_{k-1} $, then $\gs_k$ can take the two values~$\pm 1$ with
\begin{equation}\label{defsigma}
\bP^a_{N,\go}\Big(\gs_k=+1 \, \Big| \ (\tau_j)_{j\leq\iota_N}\Big)
\, = \, \frac{1}
{1 + \exp\left\{ -(\tau_{k} - \tau_{k-1}) \, h_\go + \Sigma_{[\tau_{k-1}],[\tau_{k}]} \right\}}\,.
\end{equation}
\end{itemize}
Observe that when~$\tau_{\iota_N} < N$ (which can happen only for~$a=\rf$)
there is a last (incomplete) excursion in the interval
$\{0, \ldots, N\}$, and the sign of this excursion is also expressed by~\eqref{defsigma} for~$k=\iota_{N+1}$,
provided we set~$\tau_{\iota_{N} + 1}:= N$.

\smallskip

Finally we have the moduli: again, once $(\tau_{k-1},\sigma_k)_{1\leq k\leq\iota_N+1}$ are given,
the excursions $(e_k)_{k=1,\dots,\iota_N+1}$ form an independent family. The conditional law of~$e_k(\cdot)$
on the event $\{\tau_{k-1}=\ell_0, \ \tau_{k}=\ell_1\}$ and for $f=(f_i)_{i=1, \ldots, \ell_1 - \ell_0}$
is, for $k\leq\iota_N$, given by
\begin{align}\label{defexc}
\begin{split}
& \bP^a_{N,\go}\Big( e_k(\cdot)=f \ \Big| \ (\tau_{j-1},\sigma_j)_{1\leq j\leq\iota_N+1} \Big) \\
& = \ \bP\Big(S_i=f_i: \ i=1,\ldots,\ell_1-\ell_0 \ \Big| \
S_i > 0: \ i=1, \ldots, \ell_1-\ell_0-1, \ S_{\ell_1-\ell_0}=0\Big)\,.
\end{split}
\end{align}
In the case~$\tau_{\iota_N} < N$ we have a last excursion
$e_{\iota_N + 1}(\cdot)$: its conditional law, on the event $\{\tau_{\iota_N}=\ell<N\}$
and for $f=(f_i)_{i=1, \ldots, N - \ell}$, is given by
\begin{align}\label{defexcl}
\begin{split}
&\bP^a_{N,\go}\Big( e_{\iota_N + 1}(\cdot)=f \ \Big| \ (\tau_{j-1},\sigma_j)_{1\leq j\leq\iota_N+1}
\Big)\\
& = \ \bP\Big(S_i=f_i: \ i=1,\ldots,N-\ell \ \Big| \
S_i> 0: \ i = 1, \ldots, N-\ell\Big),
\end{split}
\end{align}

\smallskip

We would like to stress that the above relations fully characterize
the polymer measure $\bP^a_{N,\go}$. A remarkable fact is that,
conditionally on $(\tau_k)_{k\in \N}$, the joint distribution
of $(\gs_j,e_j)_{j\leq \iota_N}$ \textsl{does not depend on $N$}: in this sense,
all the $N$--dependence is contained in the measure $p^a_{N,\go}$.

For this reason, this section is mainly devoted to the study of the
asymptotic behavior of the zero set measures~$p^a_{N,\go}$ as $N\to\infty$.
The main result is that $p^\rc_{N,\go}$ and $p^\rf_{N,\go}$  have the
same weak limit $p_\go$ on~$\{0,1\}^\N$ as~$N\to\infty$ (with some
restrictions when~$\go \in \cP^<$).
Once this is proven, it follows easily that also the polymer
measure~$\bP_{N,\go}^a$ converges
to a limit measure $\bP_\go$ on~$\Z^\N$, constructed by pasting
the excursion over the limit zero set. More precisely, $\bP_\go$
is the measure under which the processes $(\tau_j)$, $(\sigma_j)$
and~$(e_j)$ have the following laws:
\begin{itemize}
\smallskip
\item the law of the~$(\tau_j)_{j\in\N}$ is determined in an obvious way by the
limiting zero set measure~$p_\go$;
\smallskip
\item then, conditionally on the~$(\tau_j)_{j\in\N}$, the process $(\sigma_j)_{j\in\N}$ is an independent
one with marginal laws given by~\eqref{defsigma};
\smallskip
\item finally, conditionally on~$(\tau_j,\gs_j)_{j\in\N}$, on the event $\{\tau_{k-1}=\ell_0, \tau_{k}=\ell_1\}$
with $\ell_0<\ell_1<\infty$ the law of $e_k$ is given by the r.h.s. of (\ref{defexc}).
We have to consider also the case $\ell_0 < \infty$, $\ell_1 = \infty$, because
in the regime $\delta^\go<1$ it turns out that $\bP_\go(\tau_{k}=\infty) > 0$ (see
below and \S~\ref{app:ther}): in this case
the law of $e_k$ is given for any~$n\in\N$ and for $f=(f_i)_{i=1, \ldots, n}$ by:
\begin{align}\label{defexclul}
\begin{split}
&\bP_{\go}\Big(e_{k}(i)=f_i\,:\ i = 1, \ldots, n \ \Big| \ (\tau_j,\gs_j)_{j\in\N} \Big) \, = \,
\bP^+\Big(S_i=f_i\,:\ i = 1, \ldots, n\Big)
\\ & \qquad := \ \lim_{N\to\infty} \bP\Big(S_i=f_i\,: \ i=1,\ldots,n \ \Big| \
S_i> 0: \ i = 1, \ldots, N\Big),
\end{split}
\end{align}
where the existence of such limit is well known: see e.g. \cite{cf:G}.
\end{itemize}

\smallskip
\subsection{Law of the zero level set in the free and constrained cases}
\label{sec:zero_set}
Let us describe more explicitly $p^a_{N,\go}(A)$, using the (strong)
Markov property of $\bP^a_{N,\go}$. We use throughout the chapter the notation
(\ref{notA}). Recalling the definition \eqref{eq:matrix_cgz}
of $M_{\ga,\gb}(t)$, we have:
\begin{itemize}
\item for $a= \rc$ and $A\in\cA_N$: \ $p^\rc_{N,\go}(A)\ne 0$ \ if and
only if \  $t_{|A|}=N$, \ and in this case:
\[
 p^\rc_{N,\go}(A) \, = \,
\frac 1{Z^\rc_{N,\go}} \, \prod_{i=1}^{|A|}
 M_{[t_{i-1}],[t_i]}(t_i-t_{i-1})
\]
\item \rule{0pt}{13pt}for $a= \rf$ and $A\in\cA_N$:
\[
 p^\rf_{N,\go}(A) \, = \,
\frac 1{Z^\rf_{N,\go}} \, \left[\prod_{i=1}^{|A|}
 M_{[t_{i-1}],[t_i]}(t_i-t_{i-1})\right] P(N-t_{|A|}) \,
 \exp \Big( \tPhi_{[t_{|A|}],[N]}(N-t_{|A|}) \Big).
\]
where $P(n) := \sum_{k=n+1}^\infty K(k) = \sum_{k=n+1}^\infty \bP (\tau_1 = k)$ and we have introduced
\begin{equation} \label{eq:def_tPhi}
    \tPhi_{\ga,\gb}(\ell) \;:=\; \log \bigg[ \frac 12 \Big( 1 + \exp \big( -\ell h_\go + \Sigma_{\ga,\gb} \big) \Big) \bigg] \; \ind_{(\ell\, > 1)} \; \ind_{(\ell \, \in \gb-\ga)}\,,
\end{equation}
which differs from~$\Phi$ in not having the terms of interaction with the interface, cf.~\eqref{eq:def_Phi}.
\end{itemize}

\smallskip

We are going to show that, for any value of~$\gd^\go$, the
measure~$p^a_{N,\go}$ on~$\{0,1\}^\N$ converges as~$N\to\infty$
(with some restrictions if~$\go \in \cP^<$) to a limit measure
under which the process $([\tau_k],\tau_k-\tau_{k-1})_{k\in\bbN}$
is a \textsl{Markov renewal process}. Moreover, we will compute
explicitly the corresponding semi--Markov kernel, showing that
the returns to zero are
\begin{enumerate}
\item integrable if $\delta^\go>1$ (localized regime);
\item \rule{0pt}{13pt}defective if $\delta^\go<1$ (strictly delocalized regime);
\item \rule{0pt}{13pt}non integrable if $\delta^\go=1$ (critical regime).
\end{enumerate}
Thanks to the preceding observations, this will complete the proof of~Theorem~\ref{th:infvol}. We stress that the key result in our derivation is given by the sharp asymptotics of the partition function $Z^\rc_{N,\go}$ obtained in the previous section.

\smallskip

Before going into the proof, we give some preliminary material which is useful for all values of~$\delta^\go$.
For $k\in\N$ we define the shift operator:
\[
\theta_k:\R^\N\mapsto\R^\N, \qquad
\theta_k\zeta := \zeta_{[k+\cdot]},
\]
and it is easy to check that the following relations hold true:
\begin{equation}\label{eq:ZcZ}
Z_{N-k,\theta_k\go}^{\rc} = \cZ_{[k],[N]}(N-k), \qquad k\leq N.
\end{equation}
\begin{equation}\label{eq:relfc}
Z_{N,\go}^{\rf} \, = \, \sum_{t=0}^N Z_{t,\go}^{\rc} \, P(N-t) \,
\exp\left(\tPhi_{[t],[N]}(N-t) \right),
\end{equation}
\begin{equation}\label{eq:relpmz}
\bP^a_{N, \go} \left( \tau_1 =k \right)
\, = \, M_{0,[k]}(k) \ \frac{Z^a_{N-k,\theta_k\go}}{Z^a_{N,\go}},
\qquad 1\leq k\leq N, \quad a=\rc,\rf.
\end{equation}
Finally, using \eqref{eq:as_K}, \eqref{eq:matrix_cgz} and~\eqref{eq:def_Phi} it is easy to see that
\eqref{eq:asympt_cgz} holds true, namely
\begin{equation}\label{eq:as_L_div}
\exists \limtwo{x\to\infty}{[x]=\gb-\ga} \, x^{3/2} \, M_{\ga,\gb}(x)
\;=\;  L_{\ga, \gb},
\end{equation}
where:
\begin{equation}\label{defL}
L_{\ga, \gb} \, = \,
\begin{cases}
\displaystyle  c_K \, \frac 12 \, \Big(1+
\exp\big( \Sigma_{\ga,\gb}\big)\Big) \exp(\go ^{(0)}_{\gb})
& \text{if\ \ } h_\go = 0\\
\displaystyle\rule{0pt}{22pt} c_K \, \frac 12 \, \exp(\go ^{(0)}_{\gb})
& \text{if\ \ } h_\go > 0
\end{cases}
\; .
\end{equation}
Since also the asymptotic behavior of $P(\ell) \exp(\tPhi_{\ga,\gb}(\ell))$
will be needed, we set
\begin{equation} \label{eq:as_tPhi}
    \tL_{\ga,\gb} \;:=\; \lim_{\ell \to \infty,\ \ell \in \gb-\ga} \sqrt{\ell}
\, P(\ell) \, e^{\tPhi_{\ga,\gb}(\ell)} \;=\;
\begin{cases}
\displaystyle  c_K \big( 1 + \exp(\Sigma_{\ga,\gb}) \big) & \text{if\ \ } h_\go = 0\\
\displaystyle\rule{0pt}{22pt} c_K & \text{if\ \ } h_\go > 0
\end{cases}
\;,
\end{equation}
as it follows easily from~\eqref{eq:def_tPhi} and from the fact
that~$P(\ell) \sim 2 \, c_K / \sqrt{\ell}$ as~$\ell\to\infty$.

%
%

\smallskip
\subsection{The localized regime $(\delta^\go>1)$}

We prove point (1) of Theorem \ref{th:infvol}. More precisely, we
prove the following:
\begin{proposition}\label{pr:infvol1}
If $\delta^\go>1$ then
the polymer measures $\bP^\rf_{N, \go}$ and $\bP^\rc_{N, \go}$
converge as $N\to\infty$ to the same limit
$\bP_\go$, under which $([\tau_k],\tau_k-\tau_{k-1})_{k\in\bbN}$
is a Markov renewal process with semi-Markov kernel
$(\Gamma^>_{\ga,\gb}(x): \ga,\gb\in\bbS, x\in\N)$.
\end{proposition}
\noindent
For the definition of $\Gamma^>$ see (\ref{eq:def_Gamma>}).

\smallskip
\subsubsection{Proof of Proposition \ref{pr:infvol1}.}
We prove first the case $a=\rc$.
By \eqref{eq:ZcZ}, (\ref{eq:relpmz})
and by the asymptotics of $\cZ$ in (\ref{eq:as_Z_loc}) above,
we have for all $\ga,\gb,\gga\in\bbS$ and $\ell\in\ga$, $m\in\gb$
\[
\exists \limtwo{N\to\infty}{N\in\gga} \
\frac{Z^\rc_{N-m,\theta_m\go}}{Z^\rc_{N-\ell,\theta_\ell\go}}
\, = \, \limtwo{N\to\infty}{N\in\gga} \
\frac{\cZ_{\gb,\gga}(N-m)}{\cZ_{\ga,\gga}(N-\ell)}
\, = \, e^{-\tf_\go k} \, \frac{\xi_\gb}{\xi_\ga},
\]
and since the right hand side does not depend on $\gga$, then the
limit exists as $N\to\infty$. It follows that for $\ell\in\ga$, $k+\ell\in\gb$:
\[
\lim_{N\to \infty}
\bP^\rc_{N-\ell,\theta_\ell\go} \left( \tau_1 =k \right)\, = \, M_{\ga,\gb}(k)
\, e^{-\tf_\go k} \, \frac{\xi_{\gb}}{\xi_\ga} \, = \,
\Gamma^>_{\ga,\gb}(k).
\]
By the Markov property of $\bP^\rc_{N, \go}$ this yields
\[
\lim_{N\to \infty}
\bP^\rc_{N, \go} \left( \tau_1 =k_1, \ldots, \tau_j=k_j \right)\, =\,
\prod_{i=1}^j \Gamma^>_{[k_{i-1}],[k_i]}(k_i-k_{i-1}),
\qquad k_0:=0.
\]

\smallskip
\noindent
The argument for $\bP_{N,\go}^\rf$ goes along the very same line:
by (\ref{eq:relfc}),
\begin{align*}
& e^{-\tf_\go N} \, Z_{N-k,\theta_k\go}^\rf =
e^{-\tf_\go N} \sum_{t=0}^{N-k} \cZ_{ [k],[N-t]}(N-k-t) \: P(t) \,
\exp\left(\tPhi_{[N-t],[N]}(t) \right) \\
& \qquad = e^{-\tf_\go k}\; \sum_{\eta \in \bbS}
\sum_{t=0}^{N-k} e^{-\tf_\go \, t} \, P(t)
\left[\exp\left(\tPhi_{\eta,[N]}(t) \right)
\, e^{-\tf_\go \, (N-k-t)} \cZ_{ [k],\eta}(N-k-t)\right].
\end{align*}
Since by~\eqref{eq:as_Z_loc} the expression in brackets converges
as $N\rightarrow \infty$ and $N\in[t]+\eta$, we obtain
\begin{align*}
\exists \limtwo{N\to\infty}{N\in\gamma}
e^{-\tf_\go N} \, Z_{N-k,\theta_k\go}^\rf \, = \,
\xi_{[k]} \, e^{-\tf_\go k}
\Bigg( \frac{T}{\mu} \, \sum_{\eta \in \bbS}
\sumtwo{t \in \N}{[t]=\gamma - \eta} e^{-\tf_\go \, t} \, P(t) \,
\exp\left(\tPhi_{\eta,\, \gamma}(t) \right)
\zeta_\eta \Bigg) \,.
\end{align*}
Observe that the term in parenthesis is just a function of~$\gamma$. Having found
the precise asymptotics of $Z_{N,\go}^\rf$, we can argue as for $\bP_{N,\go}^\rc$
to conclude the proof.
\qed

\smallskip
\subsection{The critical regime $(\delta^\go=1)$}

We prove point (3) of Theorem \ref{th:infvol}. More precisely, we
prove the following:
\begin{proposition}\label{pr:infvol3}
If $\delta^\go=1$ then
the polymer measures $\bP^\rf_{N, \go}$ and $\bP^\rc_{N, \go}$
converge as $N\to\infty$ to the same limit
$\bP_\go$, under which $([\tau_k],\tau_k-\tau_{k-1})_{k\in\bbN}$
is a Markov renewal process with semi-Markov kernel
$(\Gamma^=_{\ga,\gb}(x): \ga,\gb\in\bbS, x\in\N)$.
\end{proposition}
\noindent
For the definition of $\Gamma^=$ see (\ref{eq:def_Gamma=}).

\smallskip
\subsubsection{Proof of Proposition \ref{pr:infvol3}.}
We prove first the case $a=\rc$.
By (\ref{eq:as_L_div}) and
and by the asymptotics of $\cZ$ in (\ref{eq:as_Z_cri}) above, we obtain
for all $k\in\ga$:
\[
\exists \ \limtwo{N\to\infty}{N\in\gb} \ N^{1/2}
\ \cZ_{\ga,\gb}(N-k) \, = \,
\frac{T^2}{2\pi} \, \frac{\xi_\ga\, \zeta_\gb}
{\sum_{\gamma,\gga'} \zeta_\gamma \,
L_{\gga,\gamma'}\, \xi_{\gamma'} }.
\]
It follows for all $\ga,\gb,\gga\in\bbS$ and $\ell\in\ga$, $m\in\gb$
\[
\exists \ \limtwo{N\to\infty}{N\in\gga} \
\frac{Z^\rc_{N-m,\theta_m\go}}{Z^\rc_{N-\ell,\theta_\ell\go}}
\, = \, \limtwo{N\to\infty}{N\in\gga} \
\frac{\cZ_{\gb,\gga}(N-m)}{\cZ_{\ga,\gga}(N)}
\, = \, \frac{\xi_\gb}{\xi_\ga},
\]
and since the right hand side does not depend on $\gga$, then the
limit exists as $N\to\infty$. It follows for $\ell\in\ga$, $k+\ell\in\gb$:
\[
\lim_{N\to \infty}
\bP^\rc_{N-\ell,\theta_\ell\go} \left( \tau_1 =k \right)\, = \, M_{\ga,\gb}(k)
\, \frac{\xi_{\gb}}{\xi_\ga} \, = \,
\Gamma^=_{\ga,\gb}(k).
\]
By the Markov property of $\bP^\rc_{N, \go}$ this yields
\[
\lim_{N\to \infty}
\bP^\rc_{N, \go} \left( \tau_1 =k_1, \ldots, \tau_j=k_j \right)\, =\,
\prod_{i=1}^j \Gamma^=_{[k_{i-1}],[k_i]}(k_i-k_{i-1}),
\qquad k_0:=0.
\]

\smallskip

For $\bP_{N,\go}^\rf$,
by (\ref{eq:relfc}) we have for $N\in\gb$ and $k\leq N$:
\[
Z_{N-k,\theta_k\go}^\rf \, = \, \sum_\gga
\sum_{t=0}^{N-k} \cZ_{ [k],\gga}(t) \, P(N-k-t) \,
\exp\Big(\tPhi_{\gga, \, \gb}(N-k-t) \Big).
\]
By the previous results and using~\eqref{eq:as_tPhi} we obtain that for every~$k\in\N$
\begin{align}\label{asZf}
\begin{split}
\exists \; \lim_{N\to\infty,\; N\in\gb} \, Z_{N-k,\theta_k\go}^\rf \ & = \ \xi_{[k]}\, \frac{T}{2\pi} \,
\frac{\sum_\eta \zeta_\eta \, \tL_{\eta,\gb}}{\sum_{\eta,\eta'} \zeta_\eta \,
L_{\eta,\eta'}\, \xi_{\eta'}}
\int_0^1 \frac{dt}{t^{\frac12}(1-t)^{\frac 12}}\\
& = \ \xi_{[k]}\, \Bigg( \frac{T}{2} \, \frac{\sum_\eta \zeta_\eta \, \tL_{\eta,\gb}}{\sum_{\eta,\eta'} \zeta_\eta \,
L_{\eta,\eta'}\, \xi_{\eta'}} \Bigg) \, .
\end{split}
\end{align}
To conclude it suffices to argue as in the constrained case. \qed

\smallskip
\subsection{The strictly delocalized regime $(\delta^\go<1)$}
\label{sec:tl_del}

We prove point (2) and the last assertion of Theorem \ref{th:infvol}.
In this case the result is different
according to whether $\go\in\cP^<$ or $\go\notin\cP^<$ (recall the definition
\eqref{sit1}). To be more precise, there is
first a weak formulation for {\sl all} $\omega$ 
which gives a thermodynamic
limit of $\bP^a_{N,\go}$ depending on the sequence
$\{N:[N]=\eta\}$ and on $a=\rf,\rc$; secondly, there
is a stronger formulation only for $\go\notin\cP^<$,
which says that such limits coincide for all $\eta\in\bbS$
and $a=\rf,\rc$.

\smallskip

It will turn out that in the strictly delocalized regime
there exists a.s. a last return to zero, i.e. the process  $(\tau_k)_{k\in\bbN}$ is defective.
In order to express this with the language of Markov renewal processes,
we introduce the sets
$\obbS:=\bbS\cup\{\infty\}$ and $\obbN:=\bbN\cup\{\infty\}$,
extending the equivalence relation to $\obbN$ by $[\infty] =\infty$.
Finally we set for all $\ga,\eta\in\bbS$:
\[
\Lambda_{\ga,\eta}^\rc \, := \, \big[(1-B)^{-1}L\, (1-B)^{-1}\big]_{\ga,\eta}, \qquad
\mu_{\ga,\eta}^\rc \, := \, \big[L\, (1-B)^{-1}\big]_{\ga,\eta},
\]
\[
\Lambda_{\ga,\eta}^\rf \, := \, \big[(1-B)^{-1}\tL\big]_{\ga,\eta}, \qquad
\mu_{\ga,\eta}^\rf \, := \, \tL_{\ga,\eta},
\]
and for all $\eta\in\bbS$ and $a=\rf,\rc$
we introduce the semi-Markov kernel on $\obbS\times\obbN$:
\[
    \Gamma_{\ga,\gb}^{\eta,a}(x) \, := \,
\begin{cases}
    {\displaystyle M_{\alpha,\gb}(k) \, \Lambda_{\gb,\eta}^a/\Lambda_{\ga,\eta}^a}
\quad & \ga\in\bbS, \ x\in\bbN, \ \beta=[x]\in\bbS\\
    \rule{0pt}{15pt}{\displaystyle \mu_{\ga,\eta}^a/\Lambda_{\ga,\eta}^a} \quad
& \ga\in\bbS, \ x=\infty, \ \beta=[\infty] \\
    \rule{0pt}{15pt}1 & \ga=\gb=[\infty], \ x=0 \\
    \rule{0pt}{15pt}0 & {\rm otherwise}.
\end{cases}
\]
Notice that $\Gamma^{\eta,a}$ is really a semi-Markov kernel, since for $\ga\in\bbS$:
\begin{eqnarray*}
\sum_{\gb\in\obbS} \sum_{x\in\obbN} \Gamma^{\eta,a}_{\ga,\gb}(x)
& = & \frac{\mu_{\ga,\eta}^a}{\Lambda_{\ga,\eta}^a} +
\sum_{\gb\in\bbS} \sum_{x\in\bbN}
\frac{M_{\alpha,\gb}(x) \, \Lambda_{\gb,\eta}^a}{\Lambda_{\ga,\eta}^a}
= \frac{\mu_{\ga,\eta}^a}{\Lambda_{\ga,\eta}^a} + \frac 1{\Lambda_{\ga,\eta}^a}
[B\cdot \Lambda^a]_{\ga,\eta}
\\ \\ & = & \frac{\mu_{\ga,\eta}^a}{\Lambda_{\ga,\eta}^a} +
\frac 1{\Lambda_{\ga,\eta}^a} (\Lambda_{\ga,\eta}^a-{\mu_{\ga,\eta}^a}) \, = \, 1.
\end{eqnarray*}

\smallskip
We are going to prove the following:

\smallskip
\begin{proposition}\label{pr:infvol2}
Let $\delta^\go<1$. Then:
\begin{enumerate}
\item for $a=\rf,\rc$,
$\bP^a_{N, \go}$
converges as $N\to\infty$, $[N]=\eta$ to a measure
$\bP_{\go}^{a,\eta}$, under which $([\tau_k],\tau_k-\tau_{k-1})_{k\in\bbN}$
is a Markov renewal process with semi-Markov kernel
$(\Gamma^{\eta,a}_{\ga,\gb}(x): \ga,\gb\in\obbS, x\in\obbN)$.
\smallskip
\item
if $\go\notin\cP^<$, then $\bP_{\go}^{a,\eta}=:\bP_\go$ and $\Gamma^{\eta,a}=:\Gamma^<$
depend neither on $\eta$ nor on $a$, and
both $\bP^\rf_{N, \go}$ and $\bP^\rc_{N, \go}$
converge as $N\to\infty$ to $\bP_\go$, under which $([\tau_k],\tau_k-\tau_{k-1})_{k\in\bbN}$
is a Markov renewal process with semi-Markov kernel~$\Gamma^<$.
\end{enumerate}
\end{proposition}
\smallskip

\begin{rem}\label{sit2}{\rm
Part (2) of Proposition~\ref{pr:infvol2} is an easy consequence of part~(1). In fact
from equations \eqref{defL} and \eqref{eq:as_tPhi} it follows immediately that
when $\go\notin\cP^<$ then both matrices $(L_{\ga,\gb})$ and $(\tL_{\ga,\gb})$ are constant in $\ga$, and
therefore $\Lambda^a$ factorizes into a tensor product, i.e.
\[
\Lambda_{\ga,\eta}^a \, = \, \lambda_\ga^a \, \nu_\eta^a, \qquad
\ga,\eta\in\bbS,
\]
where $(\lambda_\ga^a)_{\ga\in\bbS}$ and $(\nu_\ga^a)_{\ga\in\bbS}$
are easily computed.
But then it is immediate to check that the semi--Markov kernel~$\Gamma^{\eta,a}=:\Gamma^<$ depends
neither on $\eta$ nor on $a$.}
\end{rem}
\smallskip
\subsubsection{Proof of Proposition \ref{pr:infvol2}}
By the preceding Remark it suffices to prove part~(1). For all $k\in\ga$, by \eqref{eq:as_Z_del} we have
\begin{equation}\label{eq:as_Z_div}
\exists \limtwo{N\to\infty}{[N]=\beta} \; N^{3/2}
\ \cZ_{\ga,\gb}(N-k)
\, = \, \big[(1-B)^{-1} L \, (1-B)^{-1} \big]_{\ga,\gb}
\, = \,  \Lambda_{\ga,\gb}^\rc.
\end{equation}
In particular, we have for all $\ga,\gb,\eta\in\bbS$ and $\ell\in\ga$, $m\in\gb$:
\[
\exists \limtwo{N\to\infty}{N\in\eta} \
\frac{Z^\rc_{N-m,\theta_m\go}}{Z^\rc_{N-\ell,\theta_\ell\go}}
\, = \, \limtwo{N\to\infty}{N\in\eta} \
\frac{\cZ_{\gb,\eta}(N-m)}{\cZ_{\ga,\eta}(N)}
\, = \, \frac{\Lambda_{\gb,\eta}^\rc}{\Lambda_{\ga,\eta}^\rc},
\]
Then by (\ref{eq:relpmz}) we get
\[
\limtwo{N\to\infty}{N\in\eta} \ \bP^\rc_{N, \go}(\tau_1=k)
\, = \,
\frac{M_{0,[k]}(k) \, \Lambda_{[k],\eta}^\rc}{\Lambda_{0,\eta}^\rc} =
\Gamma^{\eta,\rc}_{0,[k]}(k).
\]
By the Markov property of $\bP^\rc_{N, \go}$ this generalizes to
\[
\limtwo{N\to\infty}{N\in\eta} \
\bP^\rc_{N, \go} \left( \tau_1 =k_1, \ldots, \tau_j=k_j \right)\, =\,
\prod_{i=1}^j \Gamma^{\eta,\rc}_{[k_{i-1}],[k_i]}(k_i-k_{i-1}),
\qquad k_0:=0.
\]

\smallskip
We prove now the case $a=\rf$.
Recalling (\ref{eq:relfc}) above, we see here that
\[
N^{1/2} \, Z_{N-k,\theta_k\go}^\rf \, = \,
\sum_{t=0}^{N-k} \cZ_{ [k],[t+k]}(t) \, N^{1/2} \, P(N-k-t) \,
\exp\Big(\tPhi_{[t+k],[N]}(N-k-t) \Big).
\]
Then by~\eqref{eq:as_tPhi} we obtain
\begin{equation}\label{as_Z_f}
\exists \ \limtwo{N\to\infty}{N\in\eta}
N^{1/2} \, Z_{N-k,\theta_k\go}^\rf \, = \, \sum_{t=0}^\infty \cZ_{ [k],[t+k]}(t)\,
\tL_{[t+k],\eta} \, = \, \big[(1-B)^{-1}\tL\big]_{[k],\eta}
\, = \,  \Lambda_{[k],\eta}^\rf\,,
\end{equation}
since
\begin{equation}\label{mmu_gb}
\sum_{t=0}^\infty \cZ_{\ga,\gga}(t) \, = \,
\sum_{t=0}^\infty \sum_{k=0}^\infty M^{*k}_{\ga,\gga}(t)
\, = \, \sum_{k=0}^\infty B^{*k}_{\ga,\gga} \, = \,
\left[(I-B)^{-1}\right]_{\ga,\gga}.
\end{equation}
Arguing as for $\bP_{N,\go}^\rc$, we conclude the proof.
\qed


\smallskip
\section{Scaling limits}
\label{sec:scaling_limits}

In this section we prove that the measures ${\bf P}^{a}_{N,\go}$
converge under Brownian rescaling. The results and proof follow
very closely those of \cite{cf:DGZ} and we shall refer to this
paper for several technical lemmas.

The first step is tightness of $(Q^a_{N,\go})_{N\in\N}$ in $C([0,1])$.
\begin{lemma}\label{lemma5}
For any $\go$ and $a=\rc,\rf$ the sequence $(Q^a_{N,\go})_{N\in\N}$
is tight in $C([0,1])$.
\end{lemma}
\noindent
For the standard proof we refer to Lemma 4 in \cite{cf:DGZ}.

\smallskip\noindent
In the rest of the section we prove Theorem \ref{th:scaling}.

\smallskip
\subsection{The localized regime $(\delta^\go>1)$}

We prove point (1) of Theorem \ref{th:scaling}.
By Lemma~\ref{lemma5}
it is enough to prove that ${\bf P}^{a}_{N,\go}(|X^N_t|>\gep)\to 0$
for all $\gep>0$ and $t\in[0,1]$ and  one can obtain this
estimate explicitly. We point out however that in this
regime one can avoid using the compactness lemma
and one can obtain a stronger result by elementary means:
observe that for any $k, n \in \N$ such that  $n>1$ and $k+n\le N$, we have
\begin{multline}
\label{eq:insert}
\bP^a_{N,\go} \big(
S_k=S_{k+n}=0, \, S_{k+i}\neq 0 \text{ for } i=1, \ldots, n-1
\big) \\
\le\;\frac{\frac 12 \left( 1+\exp\left( \sum_{i=1}^n \left(
\go_{k+i}^{(-1)}-\go_{k+i}^{(+1)}
\right)
\right)\right)}
{Z_{n, \theta_k \go}^{\rc}}
\, =:\, \widehat K_k  (n),
\end{multline}
and this holds both for $a=\rc$ and $a=\rf$.
Inequality \eqref{eq:insert} is obtained by using the Markov
property of $S$ both in the numerator and the denominator of the expression
\eqref{eq:newP}
defining $\bP^a_{N,\go} \left(\cdot \right)$ after having bounded
$Z_{N, \go} ^a$ from below  by inserting the event $S_k=S_{k+n}=0$.
Of course $\lim_{n \to \infty} (1/n)\log  \widehat K_k  (n)= -\tf _\go$
uniformly in $k$ (notice that $\widehat K_{k+T}  (n) = \widehat K_k  (n)$).
Therefore if we fix $\gep >0$ by the union bound we obtain (we recall that
$\{\tau_j\}_j$ and $\iota_N$ were defined in Section~\ref{sec:infvol})
\begin{equation*}
\begin{split}
\bP^a_{N,\go} \bigg(
\max_{j=1,2, \ldots, \iota_N } &
\tau_j-\tau_{j-1} > \left( 1+\gep \right)
\log N/{\tf_\go}
\bigg)\\
&\qquad  \le \,
\sum_{k \le N-(1+\gep) \log N/\tf_\go} \;
\sum_{n> (1+\gep) \log N /\tf_\go} \widehat K_k (n) \\
&\qquad \le \; N \sum_{n> (1+\gep) \log N /\tf_\go} \max_{k=0, \ldots , T-1}\widehat K_k (n)
\;\le\; \frac{c}{N^{\gep}},
\end{split}
\end{equation*}
for some $c>0$.

Let us start with the constrained case: notice that $\bP_{N, \go} ^\rc(\dd S)$--a.s.
we have $\tau_{\iota_N}=N$ and hence $\max_{j \le  \iota_N }
\tau_j-\tau_{j-1} \ge \max_{n=1, \ldots, N} \vert S_n \vert$, since $|S_{n+1}-S_{n}| \le 1$.
Then we immediately obtain that for any $C>1/\tf_\go$
\begin{equation}
\label{eq:longexc}
\lim_{N \to \infty}
\, \bP_{N, \go} ^\rc
\left( \max_{n=1, \ldots, N} \vert S_n \vert > C \log N\right) \, =\, 0,
\end{equation}
which is of course a much stronger statement than the
scaling limit of point (1) of Theorem~\ref{th:scaling}.
If we consider instead the measure $\bP_{N, \go} ^\rf$,
the length of the last excursion has to be taken into account too:
however, an argument very close to the one used in  \eqref{eq:insert}
yields also that the last excursion is exponentially bounded
(with the same exponent) and the proof of point (1) of Theorem~\ref{th:scaling} is complete.

\smallskip
\subsection{The strictly delocalized regime $(\delta^\go<1)$}
\label{4.2}

We prove point (2) of Theorem \ref{th:scaling}. We set for $t\in\{1,\ldots,N\}$:
\[
D_t \, := \, \inf\{ k=1,\ldots,N: \, k> t, \ S_k=0 \},
\qquad G_t \, := \, \sup\{ k=1,\ldots,N: \, k\leq t, \ S_k=0 \}.
\]
The following result shows that in the strictly delocalized regime,
as $N\to\infty$, the visits to zero under $\bP^a_{N,\go}$ tend to
be very few and concentrated at a finite distance from the origin
if $a=\rf$ and from $0$ or $N$ if $a=\rc$.
\begin{lemma}\label{pro2}
If $\delta^\go<1$ there exists a constant $C>0$ such that for all $L>0$:
\begin{equation}
\label{eq:zeta_cgz}
\limsup_{N\to\infty} \,
\bP^\rf_{N,\go}\left(G_N \ge L\right)  \, \leq \, C \, L^{-1/2},
\end{equation}
\begin{equation}
\label{eq:zeta_cgz2}
\limsup_{N\to\infty} \
\bP^\rc_{N,\go}\left(G_{N/2} \ge L\right) \, \leq \, C \, L^{-1/2},
\end{equation}
\begin{equation}
\label{eq:zeta_cgz3}
\limsup_{N\to\infty} \
\bP^\rc_{N,\go}\left(D_{N/2} \le N-L\right) \, \leq \, C \, L^{-1/2}.
\end{equation}
\end{lemma}
\medskip\noindent
Lemma~\ref{pro2} is a quantitative version of point (2) of Theorem~\ref{th:infvol}
and it is a rather straightforward complement: the proof is sketched in
\S~\ref{app:ther}, in particular \eqref{eq:last0}.

\smallskip
\subsubsection{The signs}

In order to prove point (2) of Theorem \ref{th:scaling}, it is now enough
to argue as in the proof of Theorem 9 in \cite{cf:DGZ}, with the difference
that now the excursions are not necessarily in the upper half plane, i.e.
the signs are not necessarily positive.
So the proof is complete if we can show that there exists the
limit (as~$N\to\infty$ along~$[N] = \eta$) of the probability that the
process (away from $\{0,1\}$) lives in the upper half plane.
In analogy with Section~\ref{sec:tl_del}, in the general case we
have different limits depending on the sequence $[N] = \eta$
and on $a=\rf,\rc$, while if $\go\notin\cP^<$ all such limits
coincide.

We start with the constrained case: given Lemma~\ref{pro2},
it is  sufficient to  show that
\begin{equation}
\label{eq:pom}
\exists   \ \limtwo{N\to\infty}{N\in\eta} \
\bP_{N,\go}^\rc (S_{N/2} > 0) =:\texttt p^{<,\rc}_{\go,\eta}.
\end{equation}
Formula \eqref{eq:pom} follows from the fact that
\begin{align*}
    \bP_{N,\go}^\rc (S_{N/2} > 0) =
\sum_{\ga,\gb}
\sum_{x < N/2} \sum_{y>N/2}  \frac{\cZ_{0,\ga}(x) \, \rho^+_{\ga,\gb} (y-x)
\, M_{\ga,\gb}(y-x) \, \cZ_{\gb,[N]}(N-y)}{\cZ_{0,[N]}(N)} \, ,
\end{align*}
where for all $z\in\N$ and $\ga,\gb\in\bbS$:
\begin{equation}\label{defrho}
\rho^+_{\ga,\gb} (z) \, := \, \frac{1}
{1+\exp\left(- z \, h_\go +\Sigma_{\ga,\gb}\right)},
\end{equation}
cf. \eqref{defsigma}.
By Dominated Convergence and by (\ref{defL}) and \eqref{mmu_gb}:
\[
\begin{split}
\exists \ \limtwo{N\to\infty}{N\in\eta} \ N^{3/2}
\sum_{x < N/2} \sum_{y>N/2}  &  \cZ_{0,\ga}(x) \, \rho^+_{\ga,\gb} (y-x) \,
M_{\ga,\gb}(y-x) \, \cZ_{\gb,\eta}(N-y)\\
& \, = \,
\big[(1-B)^{-1}\big]_{0,\ga} \, c_K \, \frac 12 \,
\exp(\go ^{(0)}_{\gb}) \, \big[(1-B)^{-1}\big]_{\gb,\eta} \, .
\end{split}
\]
By \eqref{eq:as_Z_del} we obtain \eqref{eq:pom} with
\begin{equation}
\label{eq:p_om}
 \texttt p^{<,\rc}_{\go,\eta} \;:=\;
\frac{\sum_{\ga,\gb}\big[(1-B)^{-1}\big]_{0,\ga}\, c_K \, \frac 12 \, \exp(\go ^{(0)}_{\gb}) \,
\big[(1-B)^{-1}\big]_{\gb,\eta}}{\big[(1-B)^{-1} L \, (1-B)^{-1} \big]_{0,\eta}} .
\end{equation}
Observe that by \eqref{defL}:
\begin{itemize}
\item if $h_\go > 0$ then
in \eqref{eq:p_om} the denominator is equal to the numerator, so that
$\texttt p^{<,\rc}_{\go,\eta} =1$ for all~$\eta$.
\item if $h_\go = 0$ and $\Sigma\equiv 0$ then
in \eqref{eq:p_om} the denominator is equal to {\it twice}
the numerator, so that
$\texttt p^{<,\rc}_{\go,\eta} =1/2$ for all~$\eta$.
\item in the remaining case, i.e. if $\go\in\cP^<$,
in general~$\texttt p^{<,\rc}_{\go,\eta}$ depends on~$\eta$.
\end{itemize}

\medskip
Now let us consider the free case. This time it is  sufficient to show that
\begin{equation}
\label{eq:pom2}
\exists  \ \limtwo{N\to\infty}{N\in\eta} \  \bP_{N,\go}^\rf (S_{N} > 0) =:
\, \texttt p^{<,\rf}_{\go,\eta}.
\end{equation}
Formula \eqref{eq:pom2} follows from the fact that
\begin{equation*}
    \bP_{N,\go}^\rf (S_{N} > 0) =
\sum_{\ga}\sum_{x < N}  \frac{\cZ_{0,\ga}(x) \, \cdot \frac 1 2 P(N-k)
}{Z^\rf_{N,\go}}\,,
\end{equation*}
and using \eqref{eq:relfc}, \eqref{mmu_gb} and~\eqref{eq:as_tPhi} we obtain that \eqref{eq:pom2} holds with
\begin{equation}
\label{eq:p_om2}
 \texttt p^{<,\rf}_{\go,\eta} \, = \,
\frac{\sum_{\ga}\big[(1-B)^{-1}\big]_{0,\ga}\, c_K}
{\big[(1-B)^{-1} \tL\big]_{0,\eta}} .
\end{equation}
Again, observe that by \eqref{eq:as_tPhi}:
\begin{itemize}
\item if $h_\go > 0$ then
in \eqref{eq:p_om2} the denominator is equal to the numerator and
$\texttt p^{<,\rf}_{\go,\eta} =1$ for all~$\eta$.
\item if $h_\go = 0$ and $\Sigma\equiv 0$ then
in \eqref{eq:p_om2} the denominator is equal to {\it twice}
the numerator, so that
$\texttt p^{<,\rf}_{\go,\eta} =1/2$ for all~$\eta$.
\item in the remaining case, i.e. if $\go\in\cP^<$,
in general~$\texttt p^{<,\rf}_{\go,\eta}$ depends on~$\eta$
and is different from $\texttt p^{<,\rc}_{\go,\eta}$.
\end{itemize}

\smallskip
\subsection{The critical regime $(\delta^\go=1)$}

In this section we prove point~(3) of Theorem~\ref{th:scaling}.
As in the previous section, we first determine the the asymptotic behavior of the
zero level set of the copolymer and then we pass to the study of the signs of the excursions.

\medskip

We introduce the random closed subset $\cA^a_N$ of $[0,1]$, describing the zero
set of the polymer of size~$N$ rescaled by a factor~$1/N$:
\[
\bbP(\cA^a_N = A/N) \, = \, p^a_{N,\go}(A), \qquad
A\subseteq\{0,\ldots,N\},
\]
where we recall that $p^a_{N,\go}(\cdot)$ has been defined in~\S~\ref{sec:zero_set}.
Let us denote by ${\cF}$ the class of \textsl{all closed subsets} of $\R^+:=[0,+\infty)$.
We are going to put on~$\cF$ a topological and measurable structure, so that we can view
the law of~$\cA^a_N$
as a probability measure on (a suitable $\gs$--field of)~$\cF$ and we can
study the weak convergence of~$\cA^a_N$.

\smallskip

We endow $\cF$ with the topology of~Matheron, cf.~\cite{cf:Matheron}
and \cite[\S~3]{cf:FFM}, which is a metrizable topology. To define it,
to a closed subset $F \subseteq \R^+$ we associate the closed
\textsl{nonempty} subset~$\tilde F$ of the compact interval~$[0,\pi/2]$
defined by $\tilde F := \arctan\!\big( F \cup \{+\infty \} \big)$.
Then the metric~$\rho(\cdot, \cdot)$ we take on~$\cF$ is
\begin{equation} \label{eq:Haus}
\rho(F,F') \, := \, \max\bigg\{
\sup_{t\in \tilde F}\, d(t,\tilde F')\,, \ \sup_{t'\in \tilde F'} \,d(t',\tilde F)\, \bigg\} \qquad \quad F,\,F' \in \cF\,,
\end{equation}
where $d(s,A):=\inf\{|t-s|, t\in A\}$ is the standard distance between a point and a set.
We point out that the r.h.s. of~\eqref{eq:Haus} is the so--called Hausdorff metric between
the compact sets~$\tilde F,\, \tilde F'$.
Thus given a sequence $\{ F_n\}_n\subset{\cF}$ and $F\in{\cF}$, we say that~$F_n \to F$ in~$\cF$
if and only if $\rho(F_n , F) \to 0$. We observe that this is equivalent to requiring
that for each open set $G$ and each compact $K$
\begin{equation}\label{hau1}
\begin{split}
& \rule{0pt}{12pt} F\cap G\ne \emptyset \quad \Longrightarrow \quad F_n\cap G\ne \emptyset
\ \ {\rm eventually} \\
& \rule{0pt}{12pt} F\cap K = \emptyset \quad  \Longrightarrow \quad
F_n \cap K = \emptyset \ \ {\rm eventually}
\end{split}
\;.
\end{equation}
Another necessary and sufficient condition for~$F_n \to F$ is that
$d(t,F_n) \to d(t,  F)$ for every~$t\in\R^+$.

This topology makes $\cF$ a separable and compact metric space~\cite[Th.~1-2-1]{cf:Matheron},
in particular a Polish space. We endow $\cF$ with the Borel $\gs$--field,
and by standard theorems on weak convergence
we have that also the space~$\cM_1(\cF)$ of probability measures on~$\cF$ is compact.

\smallskip

The main result of this section is to show that the law of the
random set~$\cA_N^a \in \cM_1(\cF)$ converges as~$N\to\infty$ to the
law of the zero set of a Brownian motion $\{B(t)\}_{t\in[0,1]}$ for~$a = \rf$
or of a Brownian bridge~$\{\beta(t)\}_{t\in[0,1]}$ for~$a=\rc$.
\smallskip
\begin{proposition}\label{c0bm}
If~$\gd^\go = 1$ then as $N\to \infty$
\begin{equation}\label{an1}
{\cA}_N^\rf \ \Longrightarrow \
\left\{t\in [0,1]: B(t)=0 \right\},
\end{equation}
\begin{equation}\label{an2}
{\cA}_N^\rc \  \Longrightarrow \
\left\{t\in [0,1]: \beta(t)=0 \right\}.
\end{equation}
\end{proposition}
\medskip

The proof of Proposition \ref{c0bm} is achieved comparing the law
of ${\cA}_N^\rf$ and ${\cA}_N^\rc$ with the law of a random
set ${\cR}_N$ defined as follows: recalling that
$\{\tau_k\}_{k\in\N}$ denotes the sequence of return times of~$S$ to
zero, we set
\[
\cR_N \, := \, \text{range\;}\{\tau_i/N, \ i\geq 0\}
\]
and we look at the law~$\cR_N$ under the critical infinite volume measure~$\bP_\go$
of Proposition~\ref{pr:infvol3}. Observe that under~$\bP_\go$ the process
$([\tau_k],\tau_k-\tau_{k-1})_{k\in\N}$ is a Markov renewal process, whose semi--Markov kernel
is given by $\Gamma^=$. The key point of the proof is given by the following result:

\smallskip
\begin{lemma}\label{th:JB}
The law of $\{\cR_N\}_N$ under $\bP_\go$ converges weakly to the law of
the random set $\{ t \geq 0: B(t) =0\}$.
\end{lemma}
\smallskip

The core of the proof (see Step~1 below) uses the theory of {\sl regenerative sets}
and their connection with the concept of {\sl subordinator}, see
\cite{cf:FFM}. However we point out that it is also possible to give a more standard proof,
using tightness and checking ``convergence of the finite dimensional distributions'':
this approach is outlined in \S~\ref{app:weak_conv}.

\smallskip

\noindent
{\it Proof of Lemma~\ref{th:JB}} We introduce the random set
\begin{equation*}
    \cR_N^{(\gb)} := \text{range}\{\tau_k/N: k\geq 0, \, [\tau_k]=\beta\} \qquad \beta\in\bbS \,.
\end{equation*}
Notice that $\cR_N=\cup_\beta \cR_N^{(\gb)}$.
Let us also recall the definitions \eqref{defkappa} and~\eqref{defTkappa}:
\[
\gk_0^{(\gb)} \, := \, \inf\{k\geq 0: [\tau_k] = \beta\}, \qquad
\gk_{i+1}^{(\gb)} \, := \, \inf\{k> \gk_i^{(\gb)}: [\tau_k] = \beta\},
\]
\[
T_0^{(\gb)} \, := \, \tau_{\gk_0^{(\gb)}}, \qquad
T_i^{(\gb)} \, := \, \tau_{\gk_i^{(\gb)}}-\tau_{\gk_{i-1}^{(\gb)}},
\qquad i\geq 1.
\]
Then $(T_i^{(\gb)})_{i\geq 1}$ is under~$\bP_\go$ an IID sequence, independent of
$T_0^{(\gb)}$: see the discussion before~\eqref{eq:cl_ren_pr}.
We divide the rest of the proof in two steps.

\medskip\noindent{\bf Step 1.}
This is the main step: we prove that the law of~$\cR_N^{(\gb)}$ under~$\bP_\go$ converges
to the law of $\{ t \geq 0: B(t) =0\}$. For this we follow the proof of
Lemma 5 in \cite{cf:DGZ}.

Let $\{P (t)\}_{t\ge 0}$ be a Poisson process with
rate $\gamma>0$, independent of $(T_i^{(\gb)})_{i\geq 0}$. Then
$\sigma_t =[T_1^{(\gb)}+\cdots+T_{P(t)}^{(\gb)}]/N$
forms a non decreasing CAD process with independent
stationary increments and $\sigma_0=0$: in other words
$\sigma=(\sigma_t)_{t\ge 0}$ is a subordinator. Notice that
\[
\cR_N^{(\gb)} \, := \, T_0^{(\gb)}/N \, + \, \widehat\cR_N^{(\gb)},
\qquad \widehat\cR_N^{(\gb)} \,:= \,\text{range\,} \{\sigma_t: t\ge 0\}.
\]
Thus ${\widehat\cR}_N^{(\gb)}$ is the (closed) range of the
subordinator $\sigma$, i.e. by~\cite{cf:FFM} a regenerative
set. As for any Levy process, the law of
$\sigma$ is characterized by the Laplace transform of the one-time
distributions:
\[
\mathbb E\left[ \exp\left(- \lambda\sigma_t \right)\right] \, = \,
\exp\left( -t \phi_N(\lambda)\right), \qquad \lambda\ge 0, \ t\ge 0,
\]
for a suitable function $\phi_N:[0,\infty)\mapsto[0,\infty)$, called
L\'evy exponent, which has a canonical representation, the
L\'evy--Khintchin formula (see e.g. (1.15) in \cite{cf:FFM}):
\begin{eqnarray*}
\phi_N (\lambda) & = & \int_{(0,\infty)} \left( 1-
e^{-\lambda s}\right) \, \gamma \,
{\mathbb P}(T_1^{(\gb)}/N\in \, ds)
\\ \\ & = & \gamma \sum_{n=1}^\infty
\left( 1-\exp(-\lambda n/N)\right) \, q^{(\gb)}(n) \;.
\end{eqnarray*}
Notice now that the law of
the regenerative set ${\widehat\cR}_N^{(\gb)}$ is invariant under
the change of time scale $\sigma _t \longrightarrow \sigma_{ct}$,
for $c>0$, and in particular independent of $\gamma>0$.
Since $\phi_N \longrightarrow c\,  \phi_N$ under this change of scale, we
can fix $\gamma=\gamma_N$ such that $\phi_N(1)=1$ and
this will be implicitly assumed from now on.
By Proposition (1.14) of \cite{cf:FFM}, the law of
${\widehat\cR}_N^{(\gb)}$ is uniquely determined by $\phi_N$.

By the asymptotics of $q^{(\gb)}$ given in (\ref{eq:cb1}),
one directly obtains that
$\phi_N (\lambda) \to \lambda^{1/2}=:
\Phi_{BM}(\lambda)$ as $N\to\infty$.
It is now a matter of
applying the result in \cite[\S 3]{cf:FFM} to obtain that
${\widehat\cR}_N^{(\gb)}$ converges in law to the regenerative set
corresponding to~$\Phi_{BM}$. However by direct computation
one obtains that the latter is nothing but the zero level set of a Brownian motion,
hence~${\widehat\cR}_N^{(\gb)} \Rightarrow \{t \in [0,1]: \ B(t) = 0\}$.
From the fact that $T_0^{(\gb)}/N$ tends to~$0$ a.s., the same weak convergence for
$\cR_N^{(\gb)}$ follows immediately.

\medskip
\noindent{\bf Step 2.}
We notice now that $\cR_N=\cup_\beta \cR_N^{(\gb)}$ is the union
of non independent sets. Therefore, although we know that each $\cR_N^{(\gb)}$
converges in law to $\{ t \geq 0: B(t) =0\}$,
it is not trivial that $\cR_N$ converges to the same limit. We start
showing that for every positive $t\geq 0$, the distance
between the first point in $\cR_N^{(\ga)}$ after $t$ and the first point
in $\cR_N^{(\gb)}$ after $t$ converges to zero in probability.
More precisely, for any closed set $F\subset[0,\infty)$ we set:
\begin{equation}\label{d_t(F)}
d_t(F)  \, := \, \inf (F\cap(t,\infty)).
\end{equation}
and we claim that for all $\ga,\beta\in\bbS$ and $t\geq 0$,
$|d_t(\cR_N^{(\ga)})- d_t(\cR_N^{(\gb)})| \to 0$ in probability.

Recalling \eqref{eq:cl_ren_pr} and
setting $q^{(\ga;\gb)}(t)=\bP_{\theta_\ga\go}(T_0^{(\gb)} =t)$, for
all $\epsilon>0$:
\begin{eqnarray*} & &
\bP_\go \left(d_t(\cR_N^{(\ga)})\geq d_t(\cR_N^{(\gb)})+\epsilon\right) \\ \\ & &
= \sum_k \sum_{\gga} \sum_{y=0}^{\lfloor Nt\rfloor}
\bP_\go (\tau_k=y, [\tau_k]=\gga) \,
\sum_{z=\lfloor Nt\rfloor -y+1}^\infty \bP_{\theta_\gga \go} (T_0^{(\gb)} =z) \
\bP_{\theta_\gb \go} (T_0^{(\ga)}\geq \lfloor N\epsilon\rfloor)
\\ \\ & & = \sum_{\gga} \sum_{y=0}^{\lfloor Nt\rfloor}
U_{0,\gga}(y) \,
\sum_{z=\lfloor Nt\rfloor-y+1}^\infty q^{(\gga;\gb)} (z) \,
\sum_{w=\lfloor N\epsilon\rfloor}^\infty q^{(\gb;\ga)}(w).
\end{eqnarray*}
Arguing as in the proof of \eqref{eq:cb1}, it is easy to obtain
the bound: $q^{(\gb;\ga)}(w)\leq C_1 \, w^{-3/2}$, and
by \eqref{eq:doney}: $U_{0,\gga}(y)\leq C_2 \, y^{-1/2}$, where $C_1,\,C_2$ are positive constants.
Then asymptotically
\[
\bP_\go \left(d_t(\cR_N^{(\ga)})\geq d_t(\cR_N^{(\gb)})+\epsilon\right)
\;\leq\; \frac{C_3}{N^{1/2}} \,
\bigg( \int_0^{t/T} dy \int_{(t-y)/T}^\infty dz \int_{\epsilon/T}^\infty dw
\, \frac  1{y^{1/2} \, z^{3/2}\, w^{3/2}} \bigg)
\]
for some positive constant~$C_3$, having used the convergence of the
Riemann sums to the corresponding integral. The very same computations
can be performed exchanging~$\ga$ with~$\gb$, hence the claim is proven.

Now notice that $d_t(\cR_N) = \min_{\ga \in \bbS} d_t(\cR_N^{(\ga)})$, and
since $\bbS$ is a finite set we have that also
$|d_t(\cR_N) - d_t(\cR_N^{(\gb)})| \to 0$ in probability for any fixed $\gb \in \bbS$.
Since we already know that~$\cR_N^{(\gb)}$ converges weakly to the law
of~$\{t \ge 0: B(t) = 0\}$, the analogous statement for~$\cR_N$ follows
by standard arguments. More precisely,
let us look at $(\cR_N, \cR_N^{(\gb)})$ as a random element of the space~$\cF \times \cF$:
by the compactness of~$\cF$ it suffices to take any convergent subsequence
$(\cR_{k_n}, \cR_{k_n}^{(\gb)}) \Rightarrow (\mathfrak B, \mathfrak C)$
and to show that~$\bbP(\mathfrak B \neq \mathfrak C) = 0$.
By the Portmanteau Theorem it is sufficient to prove that
$\lim_{N\to\infty} \bP_\go (\cR_{N} \neq \cR_{N}^{(\gb)}) = 0$,
and this is an immediate consequence of the decomposition
\begin{equation*}
    \big\{ \cR_{N} \neq \cR_{N}^{(\gb)} \big\} \;=\; \union_{t \in \Q^+} \union_{n \in \N}
    \big\{ |d_t(\cR_N) - d_t(\cR_N^{(\gb)})| > 1/n \big\}\,,
\end{equation*}
which holds by the right--continuity of~$t\mapsto d_t$.\qed

\bigskip

\noindent {\it Proof  of  (\ref{an2}).} First, we compute the
Radon-Nykodim density of the law of $\cA^\rc_N\cap[0,1/2]$ with respect to
the law of $\cR^{1/2}_N:=\cR_N\cap[0,1/2]$: for
$F=\{t_1/N,\ldots,t_k/N\}\subset[0,1/2]$ with $0=:t_0<t_1<\cdots<t_k$
integer numbers, the Radon--Nykodim derivative of the law of $\cA_N^\rc\cap[0,1/2]$
with respect to the law of $\cR^{1/2}_N$ for $\cR^{1/2}_N=F$ is:
\[
f_N^\rc(g_{1/2}(F)) \, = \, f_N^\rc(t_k/N) \, = \,
\frac{\sum_{n=N/2}^N  M_{[t_k],[n]}(n-t_k)
\, \cZ_{[n],[N]}(N-n)} {\cZ_{0,[N]}(N) \ Q_{[t_k]}(N/2-t_k)}
\ \frac{\xi_0}{\xi_{[t_k]} },
\]
where $Q_\ga(t) := \sum_\gb \sum_{s=t+1}^\infty \Gamma^=_{\ga,\gb}(s)$
and for any closed set $F\subset[0,\infty)$ we set:
\begin{equation}\label{g_t(F)}
g_t(F) \, := \, \sup(F\cap[0,t]).
\end{equation}

By \eqref{eq:as_Z_cri}, for all $\gep>0$ and uniformly in $g\in[0,1/2-\gep]$:
\begin{eqnarray*}
f_N^\rc(g) & \sim & \frac {\sum_\gga L_{[Ng],\gga}
\ \frac{T^2}{2\pi} \, \frac{\xi_\gga\, \zeta_{[N]}}
{\sum_{\gamma,\gga'} \zeta_\gamma \,
L_{\gga,\gamma'}\, \xi_{\gamma'}} \ T^{-1} \int_0^{1/2} y^{-1/2}\, (1-y-g)^{-3/2} \, dy}
{\frac{T^2}{2\pi} \, \frac{\xi_0\, \zeta_{[N]}}
{\sum_{\gamma,\gga'} \zeta_\gamma \, L_{\gga,\gamma'}\, \xi_{\gamma'}}
 \, T^{-1} \, \sum_\gga L_{[Ng],\gga} \, \xi_\gga/\xi_{[Ng]}\ 2 \, (1/2-g)^{-1/2}}
\ \frac{\xi_0}{\xi_{[Ng]} }
\\ \\ & = & \frac{\sqrt{1/2}}{1-g} \, =: \, r(g).
\end{eqnarray*}
If $\Psi$ is a bounded continuous functional on $\cF$
such that $\Psi(F)=\Psi(F\cap[0,1/2])$ for all $F\in\cF$, then,
setting $Z_B:=\left\{t\in [0,1]: B(t)=0 \right\}$ and
$Z_\beta:=\left\{t\in [0,1]: \beta(t)=0 \right\}$, we get:
\[
{\mathbb E}[\Psi(Z_\beta)] \, = \, {\mathbb
E}\left[\Psi(Z_B) \, r(g_{1/2}(Z_B))\right],
\]
see formula (49) in \cite{cf:DGZ}.
By the asymptotics of $f_N^\rc$ we obtain that
\[
\bbE\left[\Psi(\cA^\rc_N)\right] \, = \,
\bbE\left[\Psi(\cR^{1/2}_N) \, f_N^\rc(g_{1/2}(\cR^{1/2}_N))\right]
\, \to \,
\bbE\left[\Psi(Z_B) \, r(g_{1/2}(Z_B))\right] \, = \,
\bbE\left[\Psi(Z_\beta)\right]
\]
i.e. $\cA^\rc_N\cap[0,1/2]$ converges to $Z_\beta\cap[0,1/2]$.
Notice now that the distribution of
the random set $\{1-t: \, t\in\cA^\rc_N\cap[1/2,1]\}$ under
$\bP^\rc_{N,\go}$ is the same as the distribution of $\cA^\rc_N\cap[0,1/2]$
under $\bP^\rc_{N,\overline\go}$, where $\overline\go_{[i]}:=\go_{[N-i]}$.
Therefore we obtain that $\cA^\rc_N\cap[1/2,1]$
converges to $Z_\beta\cap[0,1/2]$ and the proof is complete.

\bigskip
\noindent
{\it Proof of (\ref{an1}).} By conditioning on the last zero, we see that
if $\Psi$ is a bounded continuous functional on $\cF$ then:
\[
\bbE\left[\Psi(\cA^\rf_N)\right] \, = \, \sum_{t=0}^N
\bbE\left[\Psi(\cA^\rc_t)\right] \, \frac{Z^\rc_{t,\go}}
{Z^\rf_{N,\go}}
\, P \left(N -t\right) \,
\exp\left(\tPhi_{[t],[N]}(N -t)\right).
\]
We denote by $\beta^{t}$ a Brownian bridge over the interval
$[0,t]$, i.e. a Brownian motion over $[0,t]$ conditioned to
be 0 at time $t$, and we set
$Z_{\beta^{t}}:=\left\{s\in [0,t]: \beta^{t}(s)=0 \right\}$.
By (\ref{an2}), \eqref{eq:as_Z_cri} and (\ref{asZf}) we obtain
as $N\to\infty$:
\begin{eqnarray*} & &
\bbE\left[\Psi(\cA^\rf_N)\right] \, = \, \sum_{t=0}^N \sum_\gga
1_{(t\in\gga)} \,
\bbE\left[\Psi(\cA^\rc_t)\right] \, \frac{Z^\rc_{t,\go}}
{Z^\rf_{N,\go}}
\, P \left(N -t\right) \,
\exp\left(\tPhi_{\gga,[N]}(N -t)\right)
\\ \\ & & \sim \int_0^1 \bbE[\Psi(Z_{\beta^{t}})] \,
\frac 1{t^{\frac 12}(1-t)^{\frac 12}} \, dt \cdot \sum_\gga \frac 1T \,
\frac{T^2}{2\pi} \frac{\xi_0 \, \zeta_\gga}{\sum_{\eta,\eta'} \zeta_\eta \,
L_{\eta,\eta'}\, \xi_{\eta'}} \, \frac{\tL_{\gga,[N]}}{\xi_0 \, \frac T2 \,
\frac{\sum_\eta \zeta_\eta \, \tL_{\eta,[N]}}{\sum_{\eta,\eta'} \zeta_\eta \,
L_{\eta,\eta'}\, \xi_{\eta'}}}
\\ \\ & & = \int_0^1 \bbE[\Psi(Z_{\beta^{t}})] \,
\frac 1{\pi \, t^{\frac 12} (1-t)^{\frac 12} }\, dt \, = \, \bbE[\Psi(Z_B)].
 \qed
\end{eqnarray*}

\smallskip
\subsubsection{The signs}

\medskip\noindent
To complete the proof of point~(3) of Theorem~\ref{th:scaling} in the critical
case $(\gd^\go = 1)$ we follow closely the proof given in Section~8 of~\cite{cf:DGZ}.
We have already proven the convergence of the set of zeros and we have to
``paste'' the excursions. From Section~\ref{sec:infvol} we know that, {\sl conditionally
on the zeros}:
\begin{itemize}
\item the signs $\{\gs_k\}_k$ and the absolute values $\{e_k(\cdot)\}_k$ of the excursions are independent;
\item the (conditional) law of $e_k(\cdot)$ is the same as under the original random walk measure~$\bP$.
\end{itemize}
Furthermore, the weak convergence under diffusive rescaling on~$e_k(\cdot)$ towards
the Brownian excursion~$e(\cdot)$ follows by the arguments described in~\cite{cf:DGZ}.
Then it only remains to concentrate on the signs.

\smallskip

We start with the constrained case: we are going to show that for all $t\in(0,1)$
\begin{equation}
\label{eq:pom3}
\exists   \ \lim_{N\to\infty} \
\bP_{N,\go}^\rc (S_{\lfloor tN\rfloor} > 0)
\, =: \, \texttt p^=_\go,
\end{equation}
and the limit is independent of $t$. We point out that actually we should fix the extremities
of the excursion embracing~$t$, that is we should rather prove that
\begin{equation} \label{eq:fixing}
    \lim_{N\to\infty} \bP_{N,\go}^\rc \big(S_{\lfloor tN\rfloor} > 0\,,\; G_{\lfloor tN \rfloor}/N \in (a,a+\gep)\,,\; D_{\lfloor tN \rfloor}/N \in (b,b+\gep) \big) \;=\; \texttt p^=_\go\,,
\end{equation}
for $a < t < b$ and~$\gep>0$ (recall the definition of~$G_t$ and~$D_t$ in~\S~\ref{4.2}).
However to lighten the exposition we will stick to~\eqref{eq:pom3}, since proving~\eqref{eq:fixing}
requires only minor changes.

We have, recalling
\eqref{defrho}:
\begin{align*}
    \bP_{N,\go}^\rc (S_{\lfloor tN\rfloor} > 0) =
\sum_{\ga,\gb}
\sum_{x < \lfloor tN\rfloor} \sum_{y>\lfloor tN\rfloor}
\frac{\cZ_{0,\ga}(x) \, \rho^+_{\ga,\gb} (y-x)
\, M_{\ga,\gb}(y-x) \, \cZ_{\gb,[N]}(N-y)}{\cZ_{0,[N]}(N)} \, .
\end{align*}
By Dominated Convergence and by \eqref{eq:as_Z_cri}:
\[
\begin{split} &
\exists \ \limtwo{N\to\infty}{N\in\eta} \ N^{1/2}
\sum_{x < \lfloor tN\rfloor} \sum_{y>\lfloor tN\rfloor}
\cZ_{0,\ga}(x) \, \rho^+_{\ga,\gb} (y-x) \,
M_{\ga,\gb}(y-x) \, \cZ_{\gb,\eta}(N-y)\\
& \, = \,  \frac 1{T^2} \int_0^t dx \int_t^1 dy \, [x(y-x)^3(1-y)]^{-\frac 12}
\ \left(\frac{T^2}{2\pi}\right)^2 \, \frac{\xi_0\, \zeta_\ga
\, \xi_\gb \, \zeta_\eta}
{(\sum_{\gamma,\gga'} \zeta_\gamma \,
L_{\gga,\gamma'}\, \xi_{\gamma'})^2}
\, c_K \, \frac 12 \, \exp(\go ^{(0)}_{\gb})
\end{split}
\]
see (\ref{defL}). We obtain \eqref{eq:pom3} with
\begin{equation}
\label{eq:p_om3}
 \texttt p^=_\go \, \;:=\;
\frac{\sum_{\ga,\gb} \zeta_\ga \, c_K \, \frac 12 \,
\exp(\go ^{(0)}_{\gb}) \, \xi_{\gb}}
{\sum_{\ga,\gb} \zeta_\ga \,
L_{\ga,\gb}\, \xi_{\gb}}
\end{equation}
Observe the following: by (\ref{defL}),
\begin{itemize}
\item if $h_\go > 0$ then
in \eqref{eq:p_om3} the denominator is equal to the numerator, so that
$\texttt p^=_\go =1$.
\item if $h_\go = 0$ and $\Sigma\equiv 0$ then
in \eqref{eq:p_om3} the denominator is equal to {\it twice}
the numerator, so that $\texttt p^=_\go =1/2$.
\end{itemize}

\medskip
Now let us consider the free case. We are going to show that
for all $t\in(0,1]$:
\begin{equation}
\label{eq:pom4}
\exists  \ \limtwo{N\to\infty}{[N]=\eta}
\bP_{N,\go}^\rf (S_{\lfloor tN\rfloor} > 0) \, = \,
\left( 1 -\frac{2\, \arcsin\sqrt t}\pi\right) \texttt p^=_\go
\, + \, \frac{2\, \arcsin\sqrt t}\pi \,
\texttt q^=_{\go,\eta}
\, =: \, \texttt p^{=,\rf}_{\go,\eta}(t)\,,
\end{equation}
where $\texttt p^=_\go$ is the same as above, see~\eqref{eq:p_om3},
while $\texttt q^=_{\go,\eta}$ is defined in \eqref{eq:p_om4} below.
We stress again that we should actually fix the values of~$G_{\lfloor tN \rfloor}$
and~$D_{\lfloor tN \rfloor}$ like in~\eqref{eq:fixing}, proving that
the limiting probability is either~$\texttt p^=_\go$ or~$\texttt q^=_{\go,\eta}$
according to whether~$D_{\lfloor tN \rfloor} \le N$ or~$D_{\lfloor tN \rfloor} > N$,
but this will be clear from the steps below.
Formula \eqref{eq:pom4} follows from the fact that
\[
\begin{split}
    \bP_{N,\go}^\rf (S_{\lfloor tN\rfloor} > 0) = &
\sum_{\ga,\gb}\sum_{x < \lfloor tN\rfloor}
\sum_{y>\lfloor tN\rfloor}  \frac{\cZ_{0,\ga}(x) \,
\rho^+_{\ga,\gb} (y-x) \, M_{\ga,\gb}(y-x) \,
Z^\rf_{N-y,\theta_{[y]}\go}}{Z^\rf_{N,\go}}
\\ & +  \sum_{\ga}\sum_{x < \lfloor tN\rfloor} \frac{\cZ_{0,\ga}(x) \,
\rho^+_{\ga,[N]} (N-x) \, P(N-x) \, \exp\left(\tPhi_{[x],[N]}(N-x)\right)}
{Z^\rf_{N,\go}} \; .
\end{split}
\]
By \eqref{asZf}, letting $N\to\infty$,
the first term in the r.h.s. converges to:
\[
\begin{split} &
\int_0^t \frac{dx}{x^{\frac 12}} \int_t^1 \frac{dy}{(y-x)^{\frac 32}} \ \cdot
\\ & \ \cdot \ \sum_{\ga,\gb}
\frac 1{T^2} \, \frac{T^2\, \xi_0\, \zeta_\ga}
{2\pi \sum_{\gamma,\gga'} \zeta_\gamma \, L_{\gga,\gamma'}\, \xi_{\gamma'}}
\ c_K \, \frac 12 \, \exp(\go ^{(0)}_{\gb})
\ \frac{\xi_\gb \, \frac{T}{2} \, \sum_\gga \zeta_\gga \, \tL_{\gga,\eta}}
{\sum_{\gamma,\gga'} \zeta_\gamma \, L_{\gga,\gamma'}\, \xi_{\gamma'}}
\cdot \frac{\sum_{\gamma,\gga'} \zeta_\gamma \, L_{\gga,\gamma'}\, \xi_{\gamma'}}
{\xi_0\, \frac{T}{2} \, \sum_\gga \zeta_\gga \, \tL_{\gga,\eta}}
\\ & = \left( 1 -\frac{2\, \arcsin\sqrt t}\pi\right) \cdot \texttt p^=_\go
\end{split}
\]
while letting $N\to\infty$ with $[N]=\eta$, the second term converges to
\[
\begin{split} &
\int_0^t \frac{dx}{x^{\frac 12}(1-x)^{\frac 12}}  \
\frac 1T \, \sum_\ga \frac{T^2\, \xi_0\, \zeta_\ga}
{2\pi \sum_{\gamma,\gga'} \zeta_\gamma \, L_{\gga,\gamma'}\, \xi_{\gamma'}}
\ c_K
\cdot \frac{\sum_{\gamma,\gga'} \zeta_\gamma \, L_{\gga,\gamma'}\, \xi_{\gamma'}}
{\xi_0 \, \frac{T}{2} \, \sum_\gga \zeta_\gga \, \tL_{\gga,\eta}}
\\ & = \, \frac{2\, \arcsin\sqrt t}\pi \cdot \frac{c_K \sum_\gga \zeta_\gga}
{\sum_\gga \zeta_\gga \, \tL_{\gga,\eta}} \; .
\end{split}
\]
Therefore we obtain \eqref{eq:pom4} with:
\begin{equation}
\label{eq:p_om4}
\texttt q^=_{\go,\eta} \, = \,
\frac{c_K \sum_\gga \zeta_\gga}
{\sum_\gga \zeta_\gga \, \tL_{\gga,\eta}}  \; .
\end{equation}
We observe that, by (\ref{eq:as_tPhi}):
\begin{itemize}
\item if $h_\go > 0$ or if $h_\go = 0$ and $\Sigma\equiv 0$,
then $\texttt p^{=,\rf}_{\go,\eta}(t)
=\texttt q^=_{\go,\eta} = \texttt p^=_\go$ for all $t$ and $\eta$
\item in the remaining case, i.e. if $\go\in\cP^=$, in general
$\texttt p^{=,\rf}_{\go,\eta}(t)$ depends on $t$ and $\eta$.
\end{itemize}

\smallskip

Now that we have proven the convergence of the probabilities of the signs
of the excursion, in order to conclude the proof of point (3) of Theorem \ref{th:scaling}
it is enough to argue as in the proof of Theorem 11 in \cite{cf:DGZ}.


\smallskip
\section{Appendix}

\label{app_cgz}

\smallskip
\subsection{An asymptotic result}
\label{app:qb}

We are going to prove that equation \eqref{eq:cb1} holds true.
Before starting, let us recall an elementary fact about Markov chains.
Let $Q_{\ga,\gb}$ denote the transition matrix of an irreducible, positive
recurrent Markov chain, and let us introduce the matrix $Q^{(\gamma)}$ and
the (column) vector~$|\gamma\rangle$ defined by
\begin{equation*}
    \big[ Q^{(\gamma)} \big] _{\ga,\gb} \;:=\;
Q_{\ga,\gb} \, \ind_{(\gb \neq \gamma)} \qquad \quad
\big[ |\gamma\rangle \big]_\ga := \ind_{(\ga = \gamma)} \,.
\end{equation*}
Since for any~$\gamma$ the matrix $Q^{(\gamma)}$ has spectral radius
strictly smaller than~$1$, we can define the geometric series
\begin{equation*}
    (1-Q^{(\gamma)})^{-1} \;:=\; \sum_{k=0}^\infty \big( Q^{(\gamma)} \big) ^k \,.
\end{equation*}
The interesting point is that, for every~$\gamma$, the row vector
$\langle \gamma |\cdot (1-Q^{(\gamma)})^{-1}$ is (a multiple of) the left
Perron--Frobenius eigenvector of the matrix~$Q$ (by $\langle \gamma |$ we
denote the transposed of $| \gamma \rangle$). Similarly the column vector
$(1-Q^{(\gamma)})^{-1} \cdot Q \cdot | \gamma \rangle$ is (a multiple of)
the right Perron--Frobenius eigenvector of~$Q$. More precisely we have
\begin{equation} \label{eq:markov_chains_fact}
    \big[ \langle \gamma | \cdot (1-Q^{(\gamma)})^{-1} \big]_{\ga}=
     \frac{\nu_\ga}{\nu_\gamma} \qquad \quad \big[ (1-Q^{(\gamma)})^{-1}
\cdot Q \cdot | \gamma \rangle \big]_{\ga} = 1\,,
\end{equation}
where $\{\nu_\ga\}_\ga$ is the invariant measure of the chain, that is
$\sum_\ga \nu_\ga Q_{\ga,\gb} = \nu_\gb$ and $\sum_\ga \nu_\ga = 1$.
Equation \eqref{eq:markov_chains_fact} can be proved by exploiting its
probabilistic interpretation in terms of expected number of visits to
state~$\ga$ before the first return to site~$\gamma$, see~\cite[\S~I.3]{cf:Asm}.

\medskip

Next we turn to the asymptotic behavior of $q^{(\gb)}(x)$, giving the
law of $T^{(\gb)}_0$ under~$\bbP_\gb$ (recall the notations introduced in \S~\ref{sec:as_Z_cri}).
With a standard renewal
argument, we can express it as
\begin{equation} \label{eq:conv_qb}
    q^{(\gb)}(x) \;=\; \sum_{y=0}^{x-1}
\sum_{\gamma \in \bbS} \; V^{(\gb)}_{\gb,\gamma}(y) \;
\Gamma^=_{\gamma,\gb}(x-y) \;=\; \big( V^{(\gb)} * \Gamma^= \big)_{\gb,\gb}(x) \,,
\end{equation}
where the kernel $V^{(\gb)}$ is defined by
\begin{equation*}
    V^{(\gb)}_{\ga,\gamma}(x) \;=\; \sum_{k=0}^\infty
\big[ \big( \Gamma^{(\gb)} \big) ^{*k} \big]_{\ga,\gamma}(x)\,,
\end{equation*}
and we have set $\Gamma^{(\gb)}_{\ga,\gamma}(x) :=
\Gamma^=_{\ga,\gamma}(x) \ind_{(\gamma \neq \gb)}$. Let us look more closely
at both terms in the r.h.s. of~\eqref{eq:conv_qb}.
\begin{itemize}
\item For the semi--Markov kernel $\Gamma^=$, recall its
definition~\eqref{eq:def_Gamma=}, the asymptotic behavior as
$x\to\infty$, $[x] = \gb-\gamma$ is given by
\begin{equation} \label{eq:as_Gamma=}
    \Gamma^=_{\gamma,\gb}(x) \;\sim\; \frac{\tL_{\gamma,\gb}}{x^{3/2}}
\qquad \quad \tL_{\gamma,\gb} := L_{\gamma,\gb} \frac{\xi_\gb}{\xi_\gamma}\,.
\end{equation}
Moreover, we have that
\begin{equation} \label{eq:sum_Gamma=}
    \sum_{x\in\N} \Gamma^=_{\gamma,\gb}(x) \;=\; B_{\gamma,\gb}
\frac{\xi_\gb}{\xi_\gamma} \;=:\; \tB_{\gamma,\gb} \,.
\end{equation}
\item On the other hand, for the kernel $V^{(\gb)}$ we can apply the
theory developed in~\S~\ref{sec:as_Z_del} for the case~$\gd^\go < 1$, because the matrix
\begin{equation*}
    \sum_{x\in\N} \Gamma^{(\gb)}_{\ga,\gamma}(x) \;=\; \big[ \tB^{(\gb)} \big]_{\ga,\gamma}
\end{equation*}
has Perron--Frobenius eigenvalue strictly smaller than~$1$ (we recall the
convention $[ Q^{(\gb)} ]_{\ga,\gamma} := Q_{\ga,\gamma} \ind_{(\gamma \neq \gb)}$
for any matrix~$Q$). Since
\begin{equation*}
    \Gamma^{(\gb)}_{\ga,\gamma}(x) \;\sim\;
\frac{\big[ \tL^{(\gb)} \big]_{\ga,\gamma}}{x^{3/2}} \qquad \quad x\to\infty\,,
\quad [x] = \gamma-\ga\,,
\end{equation*}
we can apply \eqref{eq:as_Z_del} to get that as $x\to\infty$, $[x] = \ga-\gamma$
\begin{equation} \label{eq:as_V}
    V^{(\gb)}_{\ga,\gamma}(x) \;\sim\;
\Big( \big[(1-\tB^{(\gb)})^{-1} \tL^{(\gb)} (1-\tB^{(\gb)})^{-1} \big]_{\ga,\gamma} \Big)
\,\frac{1}{x^{3/2}}\,.
\end{equation}
Moreover applying an analog of \eqref{eq:sum} we get that
\begin{equation} \label{eq:sum_V}
    \sum_{y \in \N} V^{(\gb)}_{\ga,\gamma}(y)=
    \sum_{k=0}^\infty \big[ \big( \tB^{(\gb)} \big)^k\big]_{\ga,\gamma}
    =
    \big[ (1-\tB^{(\gb)})^{-1} \big]_{\ga,\gamma} \,.
\end{equation}
\end{itemize}

\smallskip

We are finally ready to get the asymptotic behavior of~$q^{(\gb)}$.
As both $V^{(\gb)}$ and $\Gamma^=$ have a $x^{-3/2}$--like tail, it is easy
to check from~\eqref{eq:conv_qb} that as $x\to\infty$, $x \in T\N$
\begin{equation*}
    q^{(\gb)}(x) \;\sim\;
\sum_{\gamma \in \bbS} \Bigg\{ \bigg( \sum_{y\in\N} V^{(\gb)}_{\gb,\gamma}(y) \bigg)
\Gamma^=_{\gamma,\gb} (x)  \;+\;  V^{(\gb)}_{\gb,\gamma}(x) \bigg( \sum_{y\in\N}
\Gamma^=_{\gamma,\gb} (y) \bigg) \Bigg\}\,,
\end{equation*}
and applying \eqref{eq:sum_V}, \eqref{eq:as_Gamma=}, \eqref{eq:as_V} and \eqref{eq:sum_Gamma=}
we get that $q^{(\gb)}(x) \sim c_\gb / x^{3/2}$ as $x\to\infty$, $x \in T\N$, with
\begin{align*}
    c_\gb &\;=\; \Big[ (1-\tB^{(\gb)})^{-1} \cdot \tL \Big]_{\gb,\gb} \;+\;
\Big[ (1-\tB^{(\gb)})^{-1} \cdot \tL^{(\gb)} \cdot (1-\tB^{(\gb)})^{-1} \cdot \tB \Big]_{\gb,\gb} \\
    & \;=\; \Big[ (1-\tB^{(\gb)})^{-1} \cdot \tL \cdot (1-\tB^{(\gb)})^{-1} \cdot \tB \Big]_{\gb,\gb}\\
    &\;=\; \langle \gb | \cdot  (1-\tB^{(\gb)})^{-1} \cdot \tL \cdot
(1-\tB^{(\gb)})^{-1} \cdot \tB   \cdot | \gb \rangle\,.
\end{align*}
To obtain the second equality we have used the fact that
\begin{equation*}
    \Big[ (1-\tB^{(\gb)})^{-1} \cdot \tB \Big]_{\gb,\gb} \;=\;
\Big[ \langle \gb | \cdot (1-\tB^{(\gb)})^{-1} \cdot \tB \Big]_\gb \;=\; 1\,,
\end{equation*}
which follows from \eqref{eq:markov_chains_fact} applied to the matrix~$Q=\tB$.
Again from~\eqref{eq:markov_chains_fact} we get
\begin{equation*}
    c_\gb \;=\; \frac{1}{\nu_\gb} \sum_{\ga,\gamma \in \bbS} \nu_\ga \tL_{\ga,\gamma}\,,
\end{equation*}
where $\{\nu_\ga\}_\ga$ is the invariant measure (that is the
normalized left Perron--Frobenius eigenvector) of the matrix~$\tB$.
However from the definition~\eqref{eq:sum_Gamma=} of~$\tB$ it is
immediate to see that $\{\nu_\ga\} = \{\zeta_\ga\, \xi_\ga\}$, and
recalling the definition~\eqref{eq:as_Gamma=} of~$\tL$ we obtain the
expression for~$c_\gb$ we were looking for:
\begin{equation}\label{eq:cb}
     c_\gb \;=\; \frac 1{\zeta_\gb\, \xi_\gb}
\sum_{\ga,\gamma} \zeta_\ga \, L_{\ga,\gamma}\, \xi_\gamma.
\end{equation}


\smallskip
\subsection{Some computations on the thermodynamic limit measure}

\label{app:ther}

We want now to give a description of the typical
paths under $\bP^{a,\eta}_\go$ in the delocalization regime, i.e.
when $\delta^\go<1$. We are going
to compute the distribution of two interesting
random variables under $\bP^{a,\eta}_\go$ in this case: the last return
to zero and the total number of returns to zero. Other
analogous computations are possible using the same procedure.

\smallskip
\subsubsection{The last return to zero}
We want to study the law under $\bP^{a,\eta}_\go$ of the last zero
$\ell:=\sup\{i\in\N: S_i=0\}$ in the strictly delocalized regime.
For simplicity we consider the case $a=\rc$, the case
$a=\rf$ being completely analogous.
We compute first the law of $\ell_k:=\sup\{i\leq k: S_i=0\}$ with
$k\in\N$: for $x\leq k< N$ and $N\in\eta$:
\begin{equation}\label{bbb}
\bP^\rc_{N, \go}(\ell_k \geq x) \, = \,
\sum_{y=x}^k {\cZ_{0,[y]}(y)} \sum_{z=k+1}^N
\frac{ M_{[y],[z]}(z-y) \, {\mathcal Z_{[z],\eta}(N-z)}}{\cZ_{0,\eta}(N)}
\end{equation}
By (\ref{eq:as_L_div}) and (\ref{eq:as_Z_div}) we obtain:
\[
\limtwo{N\to\infty}{N\in\eta} \, \bP^\rc_{N, \go}(\ell_k \geq x)
=
\sum_{y=x}^k  \cZ_{0,[y]}(y) \left[ \sum_{z=0}^\infty
\frac{ L_{[y],\eta-[z]}}{ \Lambda^\rc_{0,\eta} } \,
\cZ_{\eta-[z],\eta}(z) + \sum_{z=k+1}^\infty M_{[y],[z]}(z-y) \,
\frac{\Lambda^\rc_{[z],\eta}}{\Lambda^\rc_{0,\eta}}\right]
\]
Notice now that, by \eqref{mmu_gb}:
\begin{equation}\label{mu_gb}
\sum_{z=0}^\infty L_{[y],\eta-[z]} \, \cZ_{\eta-[z],\eta}(z)
\, = \, \sum_\gga L_{[y],\gga} \sum_{z=0}^\infty \cZ_{\gga,\eta}(z) \, = \,
\left[ L \cdot (I-B)^{-1}\right]_{[y],\eta} \, = \, \mu^\rc_{[y],\eta}.
\end{equation}
Therefore, we have proven that:
\[
\bP^{\rc,\eta}_\go(\ell_k\geq x) \, = \, \limtwo{N\to\infty}{N\in\eta}
\ \bP^\rc_{N, \go}(\ell_k \geq x) \, =
\sum_{y=x}^k \cZ_{0,[y]}(y) \left[ \frac{\mu^\rc_{[y],\eta}}{\Lambda^\rc_{0,\eta}} \,
+ \sum_{z=k+1}^\infty M_{[y],[z]}(z-y) \, \frac{\Lambda^\rc_{[z],\eta}}{\Lambda^\rc_{0,\eta}}\right]
\]
and letting $k\to\infty$ we obtain:
\[
\bP^{\rc,\eta}_\go(\ell\geq x) \, = \, \sum_{y=x}^\infty \cZ_{0,[y]}(y) \,
\frac{\mu^\rc_{[y],\eta}}{\Lambda^\rc_{0,\eta}}.
\]
For the proof of Lemma \ref{pro2} above, notice for instance that by (\ref{bbb}):
\begin{eqnarray}\label{eq:last0}
& &
\bP^\rc_{N, \go}(G_{N/2} \geq L) \, = \,
\bP^\rc_{N, \go}(\ell_{N/2} \geq L)
\\ \nonumber \\ \nonumber
& & \leq \,
C_1 \, N^{3/2} \sum_{t=L}^{\lfloor N/2\rfloor} t^{-3/2}
\sum_{k= \lfloor N/2\rfloor +1}^{N+1} (k-t)^{-3/2} \,
(N+2-k)^{-3/2} \, \le \, C_2 \, L^{-1/2},
\end{eqnarray}
where $C_1,C_2$ are positive constants.

\smallskip

\subsubsection{The number of returns to zero}
Analogously, we want to study the law of the total number of returns to
zero $\cN:=\#\{i\in\N: S_i=0\}$ under $\bP^{\rc,\eta}_{\go}$. We study first
Let $\cN_K:=\#\{i: 1\leq i\leq K: S_i=0\}$ for $k\in\N$.
For $k\leq K$ and $N\in\eta$:
\[
\bP^\rc_{N, \go}(\cN_K = k) \, = \,
\sum_{x=1}^K M^{*k}_{0,[x]}(x) \sum_{y=K+1}^N
\frac{M_{[x],[y]}(y-x) \, {\mathcal Z_{[y],\eta}(N-y)}}{\cZ_{0,\eta}(N)}
\]
Then by (\ref{eq:as_L_div}) and (\ref{eq:as_Z_div}):
\begin{eqnarray*}& &
\limtwo{N\to\infty}{N\in\eta} \ \bP^\rc_{N, \go}(\cN_K = k)
\\ \\
& & = \sum_{x=0}^K M^{*k}_{0,[x]}(x) \left[\sum_{y=0}^\infty
\frac{ L_{[x],\eta-[y]}}{\Lambda^\rc_{0,\eta}} \cZ_{\eta-[y],\eta}(y)
+ \sum_{z=K+1}^\infty M_{[x],[y]}(y-x)
\frac{\Lambda^\rc_{[y],\eta}}{\Lambda^\rc_{0,\eta}}\right].
\end{eqnarray*}
By (\ref{mu_gb}), letting $K\to\infty$ we obtain:
\[
\bP^{\rc,\eta}_\go(\cN = k) \, = \, \frac 1{\Lambda^\rc_{0,\eta}} \left[B^k
\cdot \mu^\rc\right]_{0,\eta}.
\]


\smallskip

\subsection{On the weak convergence of the critical zero set}
\label{app:weak_conv}

We are going to outline an alternative proof of Lemma~\ref{th:JB},
that is we are going to show that when~$\gd^\go = 1$ as~$N\to\infty$
\begin{equation} \label{eq:app_goal}
    \cR_{N} \text{ under } \bP_\go \;\Longrightarrow\; \{ t \geq 0: B(t) =0\}.
\end{equation}
To keep the notation transparent, it is convenient to denote by~$\cG_N \in \cM_1(\cF)$
the image law of~$\cR_N$ under~$\bP_\go$. That is $\cG_N$~is a probability law on~$\cF$
(the class of all closed subsets of~$\R^+$) defined for a measurable subset~$A \subseteq \cF$ by
\begin{equation*}
    \cG_N ( A) \;:=\; \bP_\go \big( \cR_N \in A \big) \,.
\end{equation*}
In the same way the law of~$\{ t \geq 0: B(t) =0\}$ will be denoted by~$\cG^{(BM)}$.
Then we can reexpress our goal~\eqref{eq:app_goal} as
\begin{equation} \label{eq:app_goal1}
    \cG_N \;\Longrightarrow \; \cG^{(BM)} \;.
\end{equation}

\smallskip

Remember the definition \eqref{d_t(F)} of the mapping $d_t : \cF \mapsto \R^+ \cup \{+\infty\}$.
We claim that to prove \eqref{eq:app_goal1} it suffices to show that, for every~$n\in\N$ and for all
$t_1, \ldots, t_n \in \R$, the law of the vector~$(d_{t_1}, \ldots, d_{t_n})$ under~$\cG_N$
converges to the law of the same vector under~$\cG^{(BM)}$:
\begin{equation}\label{eq:finite_dim}
    (d_{t_1}, \ldots, d_{t_n}) \circ \big( \cG_N \big) ^{-1} \; \Longrightarrow \;
(d_{t_1}, \ldots, d_{t_n}) \circ \big( \cG^{(BM)} \big) ^{-1} \,.
\end{equation}

\smallskip

The intuitive explanation of why~\eqref{eq:finite_dim} should imply~\eqref{eq:app_goal1}
is that an element $\xi \in \cF$ can be {\sl identified} with the process~$\{d_t(\xi)\}_{t\in\R^+}$,
since $\xi = \{t\in\R^+:\; d_{t-}(\xi) = t\}$. Hence the convergence in~$\cM_1(\cF)$ can be read in
terms of the random process~$\{d_t(\cdot)\}_{t\in\R^+}$, and using the compactness of~$\cM_1(\cF)$
it turns out that~\eqref{eq:finite_dim} is indeed sufficient to ensure~\eqref{eq:app_goal1}. Let
us sketch more in detail these arguments.
\begin{enumerate}
\item \label{point:sigma}The Borel $\gs$--field of~$\cF$ coincides with $\gs(\{d_t\}_{t\in\R^+})$,
i.e. with the $\gs$-field generated by~$\{d_t\}_{t\in\R^+}$, and also with~$\gs(\{d_t\}_{t\in I})$ where~$I$
is any dense subset of~$\R^+$.
\item \label{point:caglad} Suppose that we are given $\{\nu_k\},\,\nu \in \cM_1(\cF)$ such that
$\nu_k \Rightarrow \nu$: this fact does not entail the convergence of all the finite
dimensional marginals of~$\{d_t\}$, that is \textsl{it is not true} that the law of the vector
$(d_{t_1}, \ldots, d_{t_n})$ under~$\nu_k$ converges to the law of the same vector under~$\nu$,
because the mappings $d_t(\cdot)$ are not continuous on~$\cF$. Nevertheless one can show that
this convergence does hold for \textsl{almost all} choices of the indexes~$t_1, \ldots, t_n$.
More precisely, given any measure $\nu \in \cM_1(\cF)$ there exists a subset~$I_\nu \subseteq \R^+$
with $Leb(I_\nu{}^\rc)=0$ with the following property: for any sequence $\{\nu_k\}$ with
$\nu_k \Rightarrow \nu$, for any~$n\in\N$ and for all $t_1, \ldots, t_n \in I_\nu$, the law of
the vector~$(d_{t_1}, \ldots, d_{t_n})$ under~$\nu_k$ converges as~$k\to\infty$ to the law of
the same vector under~$\nu$. This is a well-known feature of processes whose discontinuity points
form a negligible set, in particular CADLAG processes: in fact the set~$I_\nu$ can be chosen as
the set of~$t\in\R^+$ such that $\nu\big\{\xi:\, d_{t-}(\xi) = d_{t}(\xi)\big\} = 1$, because
$d_{t-}(\xi) = d_t(\xi)$ implies that $d_t(\cdot)$ is continuous at~$\xi$.
\item Since~$\cM_1(\cF)$ is compact, to prove~\eqref{eq:app_goal1} it suffices to show that any
convergent subsequence of~$\{\cG_N\}_N$ converges to~$\cG^{(BM)}$. Thus we take a convergent
subsequence $\cG_{k_n} \Rightarrow \nu$ for some~$\nu\in \cM_1(\cF)$ and we want to prove
that~$\nu = \cG^{(BM)}$. By point~\eqref{point:caglad} there exists a dense subset~$I_\nu \subseteq \R^+$
such that for $t_1, \ldots, t_n \in I_\nu$ the law of the vector~$(d_{t_1}, \ldots, d_{t_n})$
under~$\cG_{k_n}$ converges to the law of the same vector under~$\nu$, and since we are assuming
that~\eqref{eq:finite_dim} holds this means that the vector~$(d_{t_1}, \ldots, d_{t_n})$ has the
same law under~$\nu$ and under~$\cG^{(BM)}$. This is equivalent to say that $\nu$ and $\cG^{(BM)}$
coincide on the $\gs$--field $\gs(\{d_t\}_{t\in I_\nu})$, and by point~\eqref{point:sigma} it
follows that indeed~$\nu = \cG^{(BM)}$.
\end{enumerate}

\smallskip

Thus it only remains to show that~\eqref{eq:finite_dim} holds, and this can be done by direct
computation. For simplicity we consider only the case~$n=1$ of the one--time marginals, but
everything can be extended to the case~$n>1$.

For any~$t>0$ the law of~$d_t$ under~$\cG^{(BM)}$ is given by
\begin{equation*}
    \cG^{(BM)} \big( d_t \in \dd y \big) \;=\; \frac{t^{1/2}}{\pi \, y(y-t)^{1/2}} \,
\ind_{(y > t)} \, dy \;=:\; \rho_t(y) \, dy \,,
\end{equation*}
see \cite{cf:RevYor}. Hence we have to show that for every~$x\in\R^+$
\[
\lim_{N\to\infty} \, \bP_\go \big(d_t(\cR_N) > x\big)
\, = \,\int_x^\infty \rho_t(y) \, dy \,.
\]
We recall that $\cR_N = \text{range}\{\tau_n/N:\; n \ge 0\}$ is the range of the
process~$\{\tau_n\}_{n\in\N}$ rescaled by a factor~$1/N$, and that under~$\bP_\go$
the process~$\{\tau_n\}_{n\in\N}$ is a Markov--renewal process with semi--Markov
kernel~$\Gamma^=_{\ga,\gb}(x)$ defined by~\eqref{eq:def_Gamma=}. We also use the
notation~$U_{\ga,\gb}(x)$ for the corresponding Markov--Green function, defined
by~\eqref{eq:def_U}. Then using the Markov property we get
\begin{align*}
& \bP_\go \big(d_t(\cR_N) > x\big) = \sum_{k\in\N} \bP_\go \big(\tau_k \leq Nt\,,\;
\tau_{k+1}> Nx\big) \\
&\qquad = \sum_{\ga,\gb \in \bbS} \sum_{y=1}^{Nt} \sum_{w=Nx}^\infty \sum_{k\in\N}
\bP_\go \big( \tau_k=y\,,\; [\tau_k] = \ga \big) \, \bP_{\theta_y \go} \big( \tau_1
= w-y \,,\; [\tau_1] = \gb - \ga \big) \\
&\qquad = \sum_{\alpha,\beta\in{\mathbb S}}
\sum_{y=1}^{Nt} U_{0,\alpha}(y) \,
\sum_{w=Nx}^\infty \Gamma^=_{\alpha,\beta}(w-y)
\end{align*}
The asymptotic behavior of the terms appearing in the expression can be extracted
from~\eqref{eq:doney} and~\eqref{eq:asympt_cgz}: the net result is that as~$z\to\infty$
\begin{align*}
    &\sqrt{z} \; U_{0,\alpha}(z) \; \xrightarrow{[z]=\ga} \; \frac{T^2}{2\pi} \,
\frac{\zeta_\ga \xi_\ga}{\sum_{\gamma, \gamma'} \zeta_\gamma L_{\gamma,\gamma'}
\xi_{\gamma '}} \;=:\; c^U_{0,\ga} \\
    &z^{3/2} \; \Gamma^=_{\ga,\gb}(z) \; \xrightarrow{[z]=\gb-\ga} \;
\frac{\xi_\gb}{ \xi_\ga} \, L_{\ga,\gb} \;=:\; c^\Gamma_{\ga,\gb}\,.
\end{align*}
Therefore we have as~$N\to\infty$
\begin{align*}
& \bP_\go \big(d_t(\cR_N) > x\big) \;\sim\; \sum_{\alpha,\beta\in{\mathbb S}}
c^U_{0,\ga} \, c^\Gamma_{\ga,\gb}
\sum_{y=1}^{Nt}  \frac{1}{\sqrt{y}} \, \ind_{([y]=\ga)}
\sum_{w=Nx}^\infty \frac{1}{(w-y)^{3/2}} \, \ind_{([w]=\gb)} \\
& \qquad \sim \; \frac{1}{T^2} \Bigg( \sum_{\alpha,\beta\in{\mathbb S}}
c^U_{0,\ga} \, c^\Gamma_{\ga,\gb} \Bigg)
\frac{1}{N^2} \sum_{s \in (0, \frac{t}{T}) \cap \frac{\Z}{N}}  \frac{1}{\sqrt{s}}
\sum_{u \in (\frac{x}{T}, \infty) \cap \frac{\Z}{N}} \frac{1}{(u-s)^{3/2}}\,,
\end{align*}
and from the explicit expressions for~$c^U_{0,\ga},\, c^\Gamma_{\ga,\gb}$ together with
the convergence of the Riemann sums to the corresponding integral we get
\begin{align*}
    & \exists \lim_{N\to\infty} \bP_\go \big(d_t(\cR_N) > x\big) \;=\; \frac{1}{2\pi}
\int_0^{t/T} \dd s \, \frac{1}{\sqrt{s}} \int_{x/T}^\infty \dd u \, \frac{1}{(u-s)^{3/2}} \\
    &\qquad \;=\; \frac{1}{\pi} \int_0^{t/T} \dd s \, \frac{1}{\sqrt{s}} \, \frac{1}{\sqrt{x/T-s}}
    \;=\; \frac{1}{\pi} \int_0^{t} \dd y \, \frac{1}{\sqrt{y}} \, \frac{1}{\sqrt{x-y}}
    \;=\; \int_x^\infty \dd z \, \rho_t(z) \;,
\end{align*}
that is what was to be proven.

%
%


\smallskip

\subsection{A localization argument}

\label{app:loc}

Let us give a proof that for the \textsl{copolymer near a selective interface}
model, described in~\S~\ref{sec:2examples}, the charge~$\go$ never belongs
to~$\cP$ (see~\eqref{sit1} for the definition of~$\cP$). More precisely, we
are going to show that if~$h_\go=0$ and $\Sigma\not \equiv 0$ then~$\gd^\go > 1$,
that is the periodic copolymer with zero--mean, nontrivial charges is always
localized. As a matter of fact this is an immediate consequence of the estimates
on the critical line obtained in~\cite{cf:BG}. However we want to give here an
explicit proof, both because it is more direct and because the model studied
in~\cite{cf:BG} is built over the simple random walk measure, corresponding
to~$p=1/2$ with the language of~\S~\ref{sec:intro_cgz}, while we consider the
case~$p < 1/2$.

\smallskip

We give some preliminary notation: given an \textsl{irreducible} $T\times T$ matrix $Q_{\ga,\gb}$ with \textsl{nonnegative entries}, its Perron--Frobenius eigenvalue (= spectral radius) will be denoted by ${\tt Z}= {\tt Z}(Q)$ and the corresponding left and right eigenvectors (with any normalization) will be denoted by $\{\zeta_\ga\}, \{\xi_\ga\}$. Being a simple root of the characteristic polynomial, ${\tt Z}(Q)$ is an \textsl{analytic} function of the entries of~$Q$, and one can check that
\begin{equation} \label{eq:der}
    \frac{\partial {\tt Z}}{\partial Q_{\ga,\gb}}  \;=\; \frac{\zeta_\ga\,\xi_\gb}{\big( \sum_{\gamma} \zeta_\gamma \xi_\gamma \big)}\,.
\end{equation}
Hence ${\tt Z}(Q)$ is a strictly increasing function of each of the entries of~$Q$. We also point out a result proved by Kingman~\cite{cf:Kin}: if the matrix is a function of a real parameter $Q=Q(t)$ such that all the entries $Q_{\ga,\gb}(t)$ are \textsl{log--convex} functions of~$t$ (that is $t \mapsto \log Q_{\ga,\gb} (t)$ is convex for all~$\ga,\gb$), then also $t \mapsto {\tt Z}(Q(t))$ is a log--convex function of~$t$.

\smallskip

Next we come to the \textsl{copolymer near a selective interface} model: with reference to the general Hamiltonian~\eqref{eq:genH}, we are assuming that $\go_n^{(0)} = \tilde \go_n^{(0)} = 0$ and~$h_\go = 0$ (where~$h_\go$ was defined in~\eqref{eq:ass_h}). In this case the integrated Hamiltonian~$\Phi_{\ga,\gb} (\ell)$, see~\eqref{eq:def_Phi}, is given by
\begin{equation*}
    \Phi_{\ga,\gb} (\ell) = \begin{cases}
0 & \text{if} \ \ \ell = 1 \ \ \text{or} \ \ \ell \notin \gb-\ga \\
\rule{0pt}{18pt}\log \Big[ \frac 12 \Big( 1 + \exp \big( \Sigma_{\ga,\gb} \big) \Big) \Big] & \text{if} \ \ \ell > 1 \ \ \text{and} \ \ \ell \in \gb-\ga
\end{cases} \;.
\end{equation*}
We recall that the law of the first return to zero of the original walk
is denoted by $K(\cdot)$, see \eqref{eq:first_ret}, and we introduce the
function~$q: \bbS \to \R^+$ defined by
\begin{equation*}
    q(\gamma) := \sum_{x \in \N, \ [x]=\gamma} K(x)
\end{equation*}
(notice that $\sum_\gamma q(\gamma)=1$). Then the matrix~$B_{\ga,\gb}$ defined by~\eqref{eq:def_B} becomes
\begin{equation} \label{eq:B}
    B_{\ga,\gb} \;=\;
\begin{cases}
\frac 12 \Big( 1 + \exp \big( \Sigma_{\ga,\gb} \big) \Big) \, q(\gb-\ga)  & \ \text{if} \ \ \gb-\ga \neq [1]\\
\rule{0pt}{18pt}K(1) \;+\; \frac 12 \Big( 1 + \exp \big( \Sigma_{\ga,\ga + [1]} \big) \Big) \cdot \big( q([1]) - K(1) \big)  & \ \text{if} \ \ \gb-\ga = [1]
\end{cases}
\end{equation}
By \eqref{eq:def_delta}, to prove localization we have to show that the Perron--Frobenius eigenvalue of the matrix~$(B_{\ga,\gb})$ is strictly greater than~$1$, that is ${\tt Z}(B) > 1$.

\smallskip

Applying the elementary convexity inequality $(1+\exp(x))/2 \ge \exp(x/2)$ to \eqref{eq:B} we get
\begin{equation} \label{eq:step_matrix}
    B_{\ga,\gb} \;\ge\; \tilde  B_{\ga,\gb} \;:=\;
\begin{cases}
\exp \big( \Sigma_{\ga,\gb} / 2 \big) \, q(\gb-\ga)  & \ \text{if} \ \ \gb-\ga \neq [1]\\
\rule{0pt}{15pt}K(1) \;+\; \exp \big( \Sigma_{\ga,\ga + [1]} / 2 \big) \cdot \big( q([1]) - K(1) \big)  & \ \text{if} \ \ \gb-\ga = [1]
\end{cases} \;.
\end{equation}
By hypothesis $\Sigma_{\ga_0,\gb_0} \neq 0$ for some $\ga_0,\gb_0$, therefore the inequality above is strict for $\ga=\ga_0$, $\gb=\gb_0$. We have already observed that the P--F eigenvalue is a strictly increasing function of the entries of the matrix, hence ${\tt Z}(B) > {\tt Z}(\tilde B)$. Therefore it only remains to show that ${\tt Z}(\tilde B)\ge 1$, and the proof will be completed.

Again an elementary convexity inequality applied to the second line of~\eqref{eq:step_matrix} yields
\begin{equation} \label{eq:Bhat}
    \tilde B_{\ga,\gb} \;\ge\; \widehat  B_{\ga,\gb} \;:=\; \exp \big(\, c(\gb-\ga) \, \Sigma_{\ga,\gb} / 2 \, \big) \cdot q(\gb-\ga)
\end{equation}
where
\begin{equation*}
    c(\gamma) \;:=\;
\begin{cases}
1 & \ \text{if} \ \ \gamma \neq [1]\\
\rule{0pt}{15pt}  \frac{ q([1]) - K(1)}{q([1])}  & \ \text{if} \ \ \gamma = [1]
\end{cases} \;.
\end{equation*}
We are going to prove that~${\tt Z} (\widehat B) \ge 1$. Observe that setting $v_\ga := \Sigma_{[0],\ga}$ we can write
\begin{equation*}
    \Sigma_{\ga,\gb} = \Sigma_{[0],\gb} - \Sigma_{[0],\ga} = v_\gb - v_\ga\,.
\end{equation*}
Then we make a similarity transformation via the matrix~$L_{\ga,\gb} := \exp( v_\gb /2 ) \, \ind_{(\gb = \ga)}$, getting
\begin{align*}
    C_{\ga,\gb} &\;:=\; \big[ L \cdot \widehat B \cdot L^{-1} \big]_{\ga,\gb} \;=\; \exp \Big( \big( c(\gb - \ga) - 1 \big) \Sigma_{\ga,\gb} / 2 \Big) \cdot q(\gb - \ga)\\
    & \;=\; \exp \Big( d \, \Sigma_{\ga,\ga + [1]} \, \ind_{(\gb - \ga = 1)} \Big) \cdot q(\gb - \ga) \,,
\end{align*}
where we have introduced the constant~$d := -K(1) / (\,2\, q([1])\,)$. Of course~${\tt Z}(\widehat B) = {\tt Z}(C)$. Also notice that by the very definition of~$\Sigma_{\ga,\gb}$ we have $\Sigma_{\ga,\ga + [1]} = \go^{(-1)}_{\ga + [1]} - \go^{(+1)}_{\ga + [1]}$, hence the hypothesis~$h_\go = 0$ yields~$\sum_{\ga \in \bbS} (\Sigma_{\ga,\ga + [1]}) = 0$.

Thus we are finally left with showing that ${\tt Z}(C) \ge 1$ where $C_{\ga,\gb}$ is an~$\bbS \times \bbS$ matrix of the form
\begin{equation*}
    C_{\ga,\gb} \;=\; \exp \big(\, w_\ga \, \ind_{(\gb - \ga = 1)} \,\big) \cdot q(\gb - \ga) \qquad \text{where} \qquad \sum_{\ga} w_{\ga} = 0 \qquad \sum_\gamma q(\gamma) = 1\,.
\end{equation*}
To this end, we introduce an interpolation matrix
\begin{equation*}
    C(t)_{\ga,\gb} \;:=\; \exp \big(\, t \cdot w_\ga \, \ind_{(\gb - \ga = 1)} \,\big) \cdot q(\gb - \ga) \,,
\end{equation*}
defined for $t \in \R$, and notice that $C(1) =C$. Let us denote by~$\eta(t):= {\tt Z}\big( C(t) \big)$ the Perron--Frobenius eigenvalue of~$C(t)$: as the entries of $C(t)$ are log--convex functions of~$t$, it follows that also $\eta(t)$ is log--convex, therefore in particular convex. Moreover $\eta(0)=1$ (the matrix~$C(0)$ is bistochastic) and using~\eqref{eq:der} one easily checks that $\frac{\dd}{\dd t} \eta(t) |_{t=0} = 0$.
Since clearly $\eta(t) \ge 0$ for all~$t\in\R$, by convexity it follows that indeed $\eta(t) \ge 1$ for all~$t\in\R$, and the proof is complete.

\chapter{A general copolymer model with continuous increments}
\label{ch:continuous}


In this chapter we introduce and study a modification of the copolymer
near a selective interface model defined
in Chapter~\ref{ch:first}. The difference is that we change the reference measure~$\bP$
on which the model is built: instead of the law of the simple symmetric random walk
on~$\Z$, we allow~$\bP$ to be the law of a more general real random walk
(see~\S~\ref{sec:motivations} for some motivations for this choice). More precisely
we will consider the case when the typical increment of the walk
is bounded, centered and has an absolutely continuous law (plus a standard regularity
hypothesis on the density in order to apply the Central Local Limit Theorem).
About the charges~$\{\go_n\}_n$, we focus on the {\sl random case}.

Besides giving a proof of the existence of the free energy (which in this setting
is not trivial) we analyze the phase diagram of the model, pointing out the
close analogies with the simple random walk case described in Chapter~\ref{ch:first}.
We also consider briefly
the issue of extending to this model the coarse graining of
the free energy expressed by Theorem~\ref{th:coarse} of Chapter~\ref{ch:first}
(work in progress),
giving some partial results in this direction and discussing what is
still missing.


\smallskip
\section{The model}

\smallskip
\subsection{Motivations}
\label{sec:motivations}

Up to now all the polymer models we have worked on were built as modifications
of the law of a $(1+1)$--dimensional directed walk, the latter being
of the form~$\{(n,S_n)\}_{n}$ where~$\{S_n\}$ is a symmetric nontrivial random walk
on~$\Z$ with increments in~$\{0,\pm 1\}$. One could object that
from the viewpoint of modeling a real polymer chain these restrictions are too severe,
that is we are working with oversimplified models. A possible answer to this objection
is that the phenomena that we want to understand, like localization/delocalization,
{\sl should not} depend too much on the microscopic details of the model,
at least at a qualitative level. Even more, one could maintain that the
essential reasons of the phenomenon we are investigating may be even more visible
in an extremely simplified model. The paradigm in this direction is given by
the {\sl Ising model}, which despite its extreme simplicity is 
able to explain the origin of the ferromagnetic behavior.

Nevertheless, it would be certainly very interesting to be able to study more refined models,
at least for the purpose of understanding to what extent the results
one obtains are indeed independent of the microscopic details of the models.
In our situation, possibly the more direct refinement
that one could think of considering is to work with a
$(1+1)$--dimensional directed walk~$\{(n,S_n)\}_{n}$ in which $\{S_n\}_n$ is allowed
to be a generic real random walk.

\begin{figure}[t]
\bigskip
\psfragscanon \psfrag{Sn}[l][l]{\small $S_n$}
\psfrag{n}[l][l]{\small $n$}
\psfrag{0}[l][l]{\small $0$}
\centerline{\psfig{file=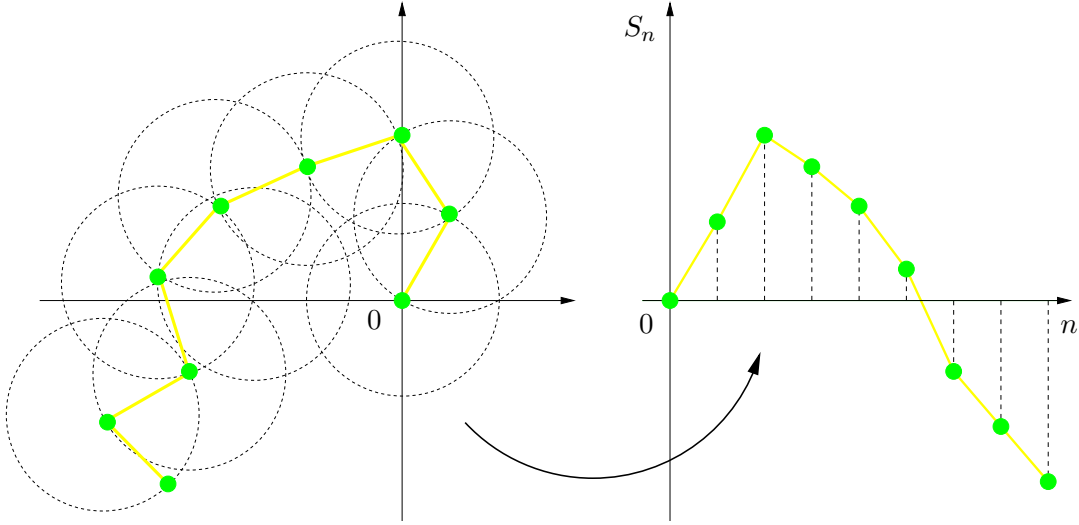,width=14.2 cm}}
\medskip
\caption{ \label{fig:general}
The correspondence between a two--dimensional random walk and a $(1+1)$--dimensional
directed walk, for the purpose of modeling a polymer chain in the proximity of an
interface (the $x$--axis in this case).}
\end{figure}

It may not be a priori evident why this should be a more realistic model: after all
it is always a directed walk model in which the first component is deterministic.
However we claim that, for the purpose of modeling a copolymer in the proximity of
a flat interface, any $d$--dimensional random walk~$\{Y_n\}_n$
is essentially equivalent to a $(1+1)$--dimensional directed walk $\{(n,S_n)\}_{n}$
for a suitable choice of the real random walk $\{S_n\}_n$.
In fact, assuming that the interface is the hyperplane~$\{x_d = 0\}$ and denoting
by~$\bs{Q}$ the law of the $d$--dimensional random walk~$\{Y_n\}_n$,
the analogue of the polymer measure introduced in Chapter~\ref{ch:first},
see equation~\eqref{eq:Boltzmann}, can be written as
\begin{equation*}
    \frac{\dd\bs{Q}_{N,\go}^{\gl, h}}{\dd\bs{Q}} (Y) \;\propto\; \gl \sum_{n=1}^N (\go_n + h)
    \sign\big( (Y_n)_d \big)\,,
\end{equation*}
where by~$(Y_n)_d$ we mean the $d$--th component of the vector $Y_n \in \R^d$.
Now observe that for the purpose of investigating the localization/delocalization
phenomenon it is sufficient to look at the~$d$--th coordinate~$\{(Y_n)_d\}_n$
under the polymer measure, which simply amounts to defining the copolymer model
over the $(1+1)$--dimensional directed
walk $(n, (Y_n)_d )$ (observe that $\{(Y_n)_d\}_n$ is a real random walk).
A graphical representation of this correspondence is given in Fig.~\ref{fig:general}
for the case of a two--dimensional random walk in which the step law is concentrated
on the surface of a sphere (which means that the distance between monomers is fixed).

We take the above considerations as sufficient motivation and we proceed to the
definition and analysis of the model.

\smallskip
\subsection{Definition of the model}
\label{sec:deff}

We take a real random walk $\{S_n\}_{n \ge 0}$, that is $S_0 = 0$ and
$S_{n} - S_{n-1} =: X_n$ with $\{X_n\}_n$ an IID sequence.
The law of the walk will be denoted by $\bP$.
Our assumptions are that:
\begin{itemize}
\item \rule{0pt}{12pt}the typical step of the walk is bounded (to be definite we take~$|X_1| \le 1$)
and centered: $\bE[X_1] = 0$;
\item \rule{0pt}{12pt}the law of~$X_1$ is absolutely continuous w.r.t. Lebesgue measure, with
density~$f$: $\bP(X_1 \in \dd x) = f(x)\,\dd x$;
\item \rule{0pt}{12pt}for some $n_0 \in \N$ the density $f_{n_0}(x) := f^{\ast n_0}(x)$ of~$S_{n_0}$
is essentially bounded: $f_{n_0}(x) \in L^\infty (\R, \dd x)$.
\end{itemize}
We point out that the last hypothesis is made in order to apply the so--called
Local Central Limit Theorem, see~\S~\ref{sec:llt_fluct} below. We denote by
$\gs^2 := \bE \big[ X_1{}^2\big] < \infty$ the variance of the typical step of the walk.

For the {\sl charges} we place ourself in the {\sl random setting}: we take the sequence
$\go = \{\go_n\}_{n \ge 1}$ to be a typical realization of a sequence of
IID random variables, whose global law is denoted by $\bbP$. The assumptions we make on the law
of~$\go_1$ are exactly the same as in Chapter~\ref{ch:first}, namely that it has finite exponential
moments: $\M(\ga) := \bE[\exp (\ga \go_1)] < \infty$ for every $\ga \in \R$ and that it is
centered: $\bbE[\go_1]=0$. We also fix $\bbE[\go_1{}^2]=1$.

For technical reasons it will be convenient to assume that $\go$ is a double--sided sequence,
that is $\go=\{\go_n\}_{n\in\Z}$, though for the definition of the copolymer model we will only
need the~$\go_n$ with~$n\ge 1$. The enlarged $\go$--space will be denoted by~$\Omega$, and of
course we look at~$\bbP$ as a probability measure on~$\Omega$. We also recall for~$k\in\Z$
the notation~$\theta^k$ for the translation on~$\Omega$, defined by $(\theta^k \go)_n := \go_{k+n}$.

\smallskip

Now we are ready to define the copolymer measure in our setting:
for $\gl, h \ge 0$ and $N\in\N$ we define $\bP^{\gl, h}_{N, \go}$ through its Radon--Nikodym
derivative:
\begin{equation}\label{eq:DefCop}
    \frac{\dd \bP^{\gl, h}_{N, \go}}{\dd \bP} (S) = \frac{1}{\tilde Z_{N,\go}}
\exp \left( \gl \sum_{n=1}^{N} (\go_n + h) \sign (S_n) \right) =: \frac{\cG_{N,\,\go}^{\gl, h}(S)}{\tilde Z_{N,\,\go}}\,.
\end{equation}
For definiteness we put $\sign (0) := 0$, but observe that in this new setting this has no role,
because the event that $S_n = 0$ for some $n$ has zero probability. We point out that in this
continuous model the charges are assigned to the points rather than to the bonds of the polymeric
chain (we recall the discussion in the caption of Fig.~\ref{fig:cop_bonds} of~Chapter~\ref{ch:first}
for the discrete setting).

The normalization factor ({\sl partition function}) $\tilde Z_{N,\go} = \tilde Z_{N,\go}^{\gl, h}$
appearing in~\eqref{eq:DefCop} is of course given by
\begin{equation} \label{eq:part_funct_cont}
    \tilde Z_{N,\go} = \bE \left[ \exp \left( \gl \sum_{n=1}^{N} (\go_n + h) \sign (S_n) \right) \right] = \bE \big[ \cG_{N,\,\go} \big] \,,
\end{equation}
and the corresponding {\sl free energy} $f(\gl, h)$ is defined by
\begin{equation}\label{eq:free_energy_cont}
f(\gl, h) := \lim_{n \to \infty} \frac{1}{N} \log \tilde Z_{N, \go}^{\gl, h} \,.
\end{equation}
A proof of the existence of such a limit, both $\bbP(\dd\go)$--a.s. and in~$\bbL_1(\bbP)$,
and of the fact that $f(\gl,h)$ is nonrandom (a phenomenon called {\sl self--averaging})
will be given in full detail in Section~\ref{sec:ex_free}.

\smallskip

Before proceeding with the analysis of the phase diagram of the model,
it is convenient to recall some basic results.

\smallskip
\subsection{Local Limit Theorem and Fluctuation Theory}
\label{sec:llt_fluct}

Since the random walk we consider has a typical step with finite nonzero variance~$\gs^2$,
the Central Limit Theorem (CLT) holds, that is we have the weak convergence as~$N\to\infty$
of the law of~$S_N/(\gs \sqrt{N})$ towards the standard Gaussian Law:
\begin{equation*}
    \forall t \in \R\,:\quad \bP \bigg[ \;\frac{S_N}{\gs \sqrt N} \le t\; \bigg] \ \longrightarrow  \ \int_{-\infty}^t \dd x\, \frac{e^{-x^2/2}}{\sqrt{2\pi}}
    \qquad \quad(N\to\infty)\,.
\end{equation*}
However in the following we will need rather precise estimates on the {\sl density}
of~$S_N$ for large~$N$. This does not follow automatically from the CLT, and some further
assumptions are required. It turns out that the (mild) assumption that for
some $n_0 \in \N$ the density $f_{n_0}(x) := f^{\ast n_0}(x)$ of~$S_{n_0}$ is essentially
bounded is sufficient to guarantee the {\sl uniform convergence} of the density
of~$S_N/(\gs \sqrt{N})$ towards the standard Gaussian density: this is the content of the
so--called Local Central Limit Theorem (LLT), cf.~\cite{cf:GneKol}.

\begin{theorem}[LLT]
Under the above assumptions, the density~$f_N(x)$ of~$S_N$ is bounded and continuous
for large~$N$. Moreover, the (continuous version of the) density of $S_N/(\gs \sqrt N)$
converges uniformly to the standard Normal density:
\begin{equation} \label{eq:unif_llt}
    \sup_{x\in\R} \bigg|\; \gs \sqrt N f_n \big(\gs \sqrt N x \big) \;-\; \frac{e^{-x^2/2}}{\sqrt{2\pi}} \;\bigg| \to 0 \qquad (n\to\infty)\,.
\end{equation}
\end{theorem}

The usefulness of the LLT is that it allows a precise control of the probability of events
like $\{S_N \in I_N\}$ when the area of $I_N$ grows slower than~$\sqrt{N}$. A typical example
in this direction is provided by the following lemma, which is an immediate consequence
of~\eqref{eq:unif_llt}.

\begin{lemma} \label{lem:LLT}
$\forall \ x \in \R$
\begin{align} \label{LLTeq}
    \sqrt{N} \cdot f_N(x) & \;\to\; \frac{1}{\sqrt{2\pi}} \qquad \quad (N \to \infty)\\
    \label{LLTcorolleq} \bP \big[ \, |S_N| \le x \, \big] &\;\sim\; \frac{2x}{\sqrt{2\pi}} \cdot \frac{1}{\sqrt{N}}
    \qquad \quad (N \to \infty)\,,
\end{align}
where in both relations the convergence is uniform for~$x$ in compact sets.
\end{lemma}

\smallskip

We conclude this section by recalling some results from the Fluctuation Theory
for random walks about conditioning a random walk to stay positive (for more details
see Section~\ref{sec:fluct} of Chapter~\ref{ch:llt}). We start with the asymptotic behavior
of the probability of the first entrance of the walk in the negative half--line:
it is a classical result~\cite{cf:Fel2} that, whenever the step of the walk has zero mean and
finite nonzero variance~$\gs^2$, as~$n\to\infty$
\begin{equation} \label{eq:rw_tail}
    \bP \big[ \, S_1 > 0, \ldots, S_n > 0 \, \big] \;=\;
    \sum_{k = n+1}^\infty  \bP \big[\, S_1 > 0, \ldots, S_{k-1} > 0, S_k \le 0 \,\big]
    \;\sim\; \frac{2\,C}{\sqrt{n}}\,,
\end{equation}
for some~$C\in(0,\infty)$. The local version of this relation holds too, namely
\begin{equation}\label{eq:rw_point}
    \bP \big[\, S_1 > 0, \ldots, S_{n-1} > 0, S_n \le 0 \,\big]
    \;\sim\; \frac{C}{n^{3/2}} \qquad \quad (n\to\infty)\,,
\end{equation}
and it has been proven in~\cite{cf:AliDon}. Of course both~\eqref{eq:rw_tail} and~\eqref{eq:rw_point}
hold also for the first entrance in the positive half--line (with possibly a different
constant~$C$).

We will also need to control the probabilities of the first entrances in the case when
the random walk does not necessarily starts from~$0$. More precisely, for~$x\in\R$ let us denote
by~$\bP_x$ the law of~$\{S_n + x\}_n$ under~$\bP$. Then for~$x \ge 0$ we have
\begin{align}
    \label{eq:rw_tail_x}
    \bP_x \big[ \, S_1 > 0, \ldots, S_n > 0 \, \big]
    &\;\sim\; \frac{2\,C_x}{\sqrt n} \qquad \quad (n\to\infty)\\
    \label{eq:rw_point_x}
    \bP_x \big[\, S_1 > 0, \ldots, S_{n-1} > 0, S_n \le 0 \,\big]
    &\;\sim\; \frac{C_x}{n^{3/2}} \qquad \quad (n\to\infty)\,,
\end{align}
where~$C_x$ is a positive nondecreasing function of~$x\ge 0$. We point out that, up to multiplicative
constants, the function~$C_x$ coincides with the renewal function associated to the descending
ladder heights process of the random walk, see~\cite[\S~3]{cf:BerDon} for more on this issue.
A direct proof of~\eqref{eq:rw_point_x} can be also given using the methods of Chapter~\ref{ch:llt}.
Again, an analog of~\eqref{eq:rw_point_x} is valid also for the first entrance in the positive
half--line, when the random walk stars from~$x\le 0$.

\smallskip
\subsection{The phase diagram}

Next we turn to the analysis of the phase diagram of the model we have introduced. As in the
discrete case, the first step is the identification of the
free energy coming from delocalized paths: restricting to trajectories that stay positive up
to epoch~$N$ we have that for~$\bbP$--a.e.~$\go$
\begin{equation}
\label{eq:del_free}
\begin{split}
\frac 1N \log {\tilde Z}_{N,\go}^{\gl, h}  &\;\ge\;
\frac 1N \log \bE
\left[
\exp\left(
\gl \sum_{n=1}^N
\left( \go_n +h\right) \sign \left(S_n\right)
\right)
;\,  S_1 > 0, \ldots, S_N > 0
\right]
\\
&\;=\; \frac {\gl} {N} \sum_{n=1}^N \left( \go_n +h\right) \ + \
\frac 1N
\log \bP \big( S_1 > 0, \ldots, S_N > 0 \big)\, \stackrel{N \to \infty}{\longrightarrow}\, \gl h ,
\end{split}
\end{equation}
where in the last line we have used the strong law of large numbers and the asymptotic
behavior given by~\eqref{eq:rw_tail}.

Arguing as in Chapter~\ref{ch:first}, we partition the $(\gl,h)$--space in two regions:
\smallskip
\begin{itemize}
\item the localized region: $\cL = \left\{ (\gl , h): \, f(\gl, h)>\gl h\right\}$;
\smallskip
\item the delocalized region:  $\cD = \left\{ (\gl , h): \, f(\gl, h) = \gl h\right\}$.
\end{itemize}
\smallskip
For the critical line~$h_c(\cdot)$ separating the two regions we have the following result,
in complete analogy with the discrete case:
\smallskip
\begin{proposition}
\label{prop:prel_cont}
There exists a continuous {\sl increasing} function $h_c: [0,\infty) \to [0,\infty)$
with~$h_c(0)=0$ such that
\begin{equation*}
    \cD \,=\, \big\{(\gl,h):\: h \ge h_c(\gl) \big\} \qquad \cL \,=\,
    \big\{(\gl,h):\: h < h_c(\gl)\big\} \,.
\end{equation*}
\end{proposition}
\smallskip

\noindent
Thus the picture of the phase diagram in this continuous setting looks quite similar to
the discrete case analyzed in Chapter~\ref{ch:first}, at least at a qualitative level.
Now we are going to make this statement quantitative, showing that for the critical
curve~$h_c(\cdot)$ of our continuous model we have {\sl exactly the same upper and lower bound}
that hold in the discrete case, namely
\begin{equation} \label{eq:sumup_cont}
    h^{(2/3)}(\cdot) \; =: \; \underline h (\cdot) \ \le \ h_c(\cdot) \ \le \ \overline h(\cdot)
    \; := \; h^{(1)}(\cdot) \;,
\end{equation}
where we recall the definition of~$h^{(m)}(\cdot)$ for~$m>0$:
\begin{equation*}
    h^{(m)}(\gl) \; := \; \frac{\log \M ( -2 m \gl)}{2 m \gl}\,,
\end{equation*}
and $\M (\ga) := \bbE \big[ \exp (\ga \go_1) \big]$ is the moment generating function of
the environment.

Before proceeding, let us spend some words on Proposition~\ref{prop:prel_cont}:
using convexity arguments as in~\cite[\S~1.2]{cf:BG2} it is not difficult
to prove the existence of the critical line, together with some of its properties.
However showing that $h_c(\cdot)$ is indeed increasing and not only nondecreasing,
that it is continuous also at~$\gl=0$ and that~$h_c(\gl) < \infty$ for every~$\gl \ge 0$
does not follow immediately. A rather cheap (if not elementary) way of proving these properties
is to supply convex analysis with the knowledge of the bounds~\eqref{eq:sumup_cont}
on~$h_c(\cdot)$ (whose proof is independent of Proposition~\ref{prop:prel_cont}).

\smallskip
\subsubsection{Upper bound}

The proof of the upper bound in~\eqref{eq:sumup_cont} is completely analogous to the one given in
Chapter~\ref{ch:first} for the discrete setting, that is it suffices to
apply the annealing procedure. However, in order
not to end up with a useless bound, we have to suitably modify the partition function, as
in~\S~\ref{rem:Z} of Chapter~\ref{ch:first}. More precisely, subtracting to the Hamiltonian
the term~$\gl \sum_{n=1}^N (\go_n + h)$ (that does not depend on~$S$ and that once averaged
on the environment is simply $\gl h N$) and using the fact
that the limit~\eqref{eq:free_energy_cont} holds also in~$\bbL_1(\bbP)$ we can write
\begin{align*}
    f (\gl,h) - \gl h \;= \; \lim_{N\to\infty} \frac 1N \bbE \log \bE
    \left[ \exp\left( -2 \gl \sum_{n=1}^N \left( \go_n +h\right) \ind_{\{ \sign(S_n) = -1 \}}
    \right)  \right]\,.
\end{align*}
However by Jensen's inequality we can bring the expectation~$\bbE$ inside the~$\log$, and
performing the integration over the disorder we get
\begin{align} \label{eq:ann_cont}
    f (\gl,h) - \gl h \;\le\; \lim_{N\to\infty} \frac 1N \log \bE
    \left[ \exp\left( \sum_{n=1}^N \big( \log \M (-2 \gl) -2 \gl h\big) \ind_{\{ \sign(S_n) = -1 \}}
    \right) \right] \,.
\end{align}
For~$h \ge \overline h(\gl)$ the argument of of the exponential is nonpositive: thus
$f(\gl,h)-\gl h \le 0$ and by~\eqref{eq:del_free} we have~$f(\gl,h)=\gl h$, hence
we have proven that~$h_c(\cdot) \le \overline h(\cdot)$.

Observe that for~$h < \overline h(\gl)$
the r.h.s. of~\eqref{eq:ann_cont} equals $\big( \log \M (-2 \gl) -2 \gl h\big) > 0$,
hence~$\overline h(\cdot)$ is indeed the best upper bound on~$h_c(\cdot)$ that
one can extract from~\eqref{eq:ann_cont}.
Also notice that the arguments of Chapter~\ref{ch:cg} can be applied to our continuous
setting with essentially no change: therefore the technique of {\sl constrained
annealing} via empirical averages of local functions cannot improve the upper bound we have found.

\smallskip
\subsubsection{Lower bound}

A proof of the lower bound in~\eqref{eq:sumup_cont} can be obtained by following very
closely the proof in the discrete setting given in \S~\ref{app:prooflb} of Chapter~\ref{ch:cgg}.
For this reason we simply outline the main steps. Let us introduce
a notation for the modified Hamiltonian
\begin{equation*}
    \cH'_{N,\go} \;:=\; -2 \gl \sum_{n=1}^N \left( \go_n +h\right)
    \ind_{\{ \sign(S_n) = -1 \}} \,,
\end{equation*}
so that the reduced free energy $f(\gl,h) - \gl h$ can be expressed for~$\bbP$--a.e.~$\go$ as
\begin{equation} \label{eq:alt_free}
    f(\gl, h) - \gl h \;=\; \lim_{N\to\infty} \frac 1N \log I_{N,\go}^{\gl, h}\,,
\end{equation}
where
\begin{equation*}
    I_{N,\go}^{\gl,h} \;:=\; \inf_{-1 \le x \le 1} \bE_x \Big[ \exp \big( \cH'_{N,\go} \big)
    \, \ind_{\{|S_N| \le 1\}} \Big]\,.
\end{equation*}
The proof of relation~\eqref{eq:alt_free} is the core of Section~\ref{sec:ex_free}.

We stress that $I_{N,\go}$ takes for the continuous setting the role
that the pinned partition function~$Z_{N,\go}(0)$ (see \eqref{eq:pinned}
of Chapter~\ref{ch:first}) has in the discrete setting. In fact, using~\eqref{eq:alt_free}
together with the superadditivity of the process~$\{I_{N,\go}\}_N$
(proved in~\S~\ref{sec:a}), the arguments of the first
part of~\S~\ref{app:prooflb} of Chapter~\ref{ch:cgg} can be easily adapted to the continuous
setting. In particular, in order to prove that a point $(\gl, h)$ is localized, it suffices to find
a number~$C>1$ and a random variable $T:\Omega \to \N$ with the following two properties:
\begin{equation} \label{eq:cond_T}
    (1)\quad I_{T(\go), \go}^{\gl, h} \ge C \quad \bbP(\dd\go)\text{--a.s.}
    \qquad \qquad \quad (2)\quad \bbE\big[ T \big] < \infty \,.
\end{equation}

Thus it only remains to show that for every~$(\gl,h)$ with $h<\underline h(\gl)$ one can
build a random time~$T$ satisfying~\eqref{eq:cond_T}. However, if we define~$T=T_{A,\gep,q}$
as in~\eqref{eq:def_T_cgg} of Chapter~\ref{ch:cgg},
then using the asymptotic relation~\eqref{eq:rw_point_x} we can easily
get a lower bound on~$I_{T(\go),\go}$ like~\eqref{eq:almost_maj}, for a possibly
different value of the constant~$c'$ (see also~\eqref{eq:needs}).
Therefore one can tune the parameters~$A,\gep,q$ exactly as it is done
in the end of~\S~\ref{app:prooflb} of Chapter~\ref{ch:cgg}, see page~\pageref{eq:lbonZTom},
and condition~\eqref{eq:cond_T} will be satisfied.


\smallskip
\section{Existence of the free energy}
\label{sec:ex_free}

\label{existence}

In this section we give a proof of the existence of the free energy,
that is of the limit~\eqref{eq:free_energy_cont}. The standard procedure to get the
existence of such a limit is to modify the partition function of the model
(without changing the Laplace asymptotic behavior) in order to perform
superadditivity arguments. For instance in the discrete case it is sufficient
to restrict the sum defining the partition function~$\tilde Z_{N,\go}$ to the trajectories
such that~$S_N=0$: in our continuous setting this is no longer possible, because
the event~$\{S_N=0\}$ has probability~$0$. This obstacle is easily bypassed and
it is not difficult to find a useful modification of the partition function.
The drawback is that showing that the modified partition function yields
the same Laplace asymptotic behavior as the original one is no longer trivial.

\smallskip

Remember that by hypothesis the steps of our random walk are bounded by~$1$: $|S_n - S_{n-1}| \le 1$.
We also recall the notation~$\cG_{N,\go}:=\cG_{N,\go}^{\gl,h}$ for the Boltzmann
factor appearing in the definition of the copolymer measure~\eqref{eq:DefCop}, and the
expression~\eqref{eq:part_funct_cont} for the partition function~$\tilde Z_{N,\go}$.
The modified partition function to which we will apply superadditivity arguments
will be
\begin{equation} \label{eq:def_I}
    I_{N,\go} \;=\; I_{N,\go}^{\gl,h} \;:=\; \inf_{-1 \le x \le +1} \bE_x \big[ \cG_{N,\go} \, \ind_{\{|S_N| \le 1\}} \big]\,,
\end{equation}
where $\bP_x$ is the law of the random walk starting at~$x\in\R$, introduced in the
preceding section.

The proof is organized in three steps: in \S~\ref{sec:a} we show that the
limit~\eqref{eq:free_energy_cont} exists if we replace~$\tilde Z_{N,\go}$ by~$I_{N,\go}$,
and then in \S~\ref{sec:b} and~\S~\ref{sec:c} we prove some comparison inequalities
showing that~$I_{N,\go}$ and~$\tilde Z_{N,\go}$ are equivalent for the sake of computing
the free energy. To this purpose it will be convenient to consider an intermediate
partition function~$J_{N,\go}$ defined as follows:
\begin{equation} \label{eq:def_J}
    J_{N,\go} \;=\; J_{N,\go}^{\gl,h} \;:=\; \bE \big[ \cG_{N,\go} \, \ind_{\{|S_N| \le 1\}} \big]\,.
\end{equation}

\smallskip

\subsection{Step 1}
\label{sec:a}

We start showing that the sequence of random variables $\{\log I_{N,\,\go}\}_N$
satisfies the hypothesis of \textsl{Kingman's Superadditive Ergodic Theorem}~\cite{cf:Kin2}.

We begin with the upper bound on~$\bbE [\log I_{N,\go}]$: using Jensen's inequality, the
definition~\eqref{eq:def_I} of~$I_{N,\go}$ and a rough bound on~$\cG_{N,\go}$ we get
\begin{align*}
    \bbE \, \big[ \log  I_{N,\,\go} \big] &\;\le\; \log \bE \bbE \bigg[ \exp \bigg( \gl \sum_{n=1}^N
    (|\go_n| + h) \bigg)  \, \ind_{\{|S_N| \le 1\}} \bigg] \;\le\; \big( \log
    \bbE\big[ e^{\gl |\go_1|} \big] + \gl h \big) N \,,
\end{align*}
hence $\sup_N \{ \bbE [ \log  I_{N,\,\go} ] /N \} < + \infty$. Also the superadditivity
is easily obtained: making explicit the functional dependence of $\cG_{N,\,\go}$
on the path $(S_1, \dots, S_N)$ when it is convenient
and using the Markov property, we obtain that $\forall \ x \in [-1,+1]$
\begin{align*}
\bE_x & \left[ \: \cG_{N+M,\,\go} \, \ind_{\{|S_{N+M}|\le 1\}} \right]  \;\ge\; \bE_x \left[ \: \cG_{N+M,\,\go} \, \ind_{\{|S_{N}|\le 1\}} \ind_{\{|S_{N+M}|\le 1\}} \right] \\
& \;=\; \int_{-1}^{1} \dd z \: f_N(z-x) \: \bE_x \big[ \: \cG_{N,\,\go} (S_1, ...\, , S_{N-1}, z) \big] \cdot \bE_z \left[ \: \cG_{M,\,\theta^{N}\go} \, \ind_{\{|S_{M}|\le 1\}} \right] \\
& \;\ge\; \bigg( \int_{-1}^{1} \dd z \: f_N(z-x) \: \bE_x \big[ \: \cG_{N,\,\go} (S_1, ...\, ,z) \big] \bigg) \cdot \inf_{z \in [-1,1]} \bE_z \big[ \: \cG_{M,\,\theta^{N}\go} \, \ind_{\{|S_{M}|\le 1\}} \big] \\
& \;=\; \bE_x \big[ \: \cG_{N,\,\go} \cdot \ind_{\{|S_N| \le 1\}} \big] \cdot I_{M,\,\theta^{N}\go} \  \ge \ I_{N,\,\go} \cdot I_{M,\,\theta^{N}\go}
\end{align*}
Reading only the extremities of this chain of inequalities, we have
\[
\bE_x \left[ \: \cG_{N+M,\,\go} \, \ind_{\{|S_{N+M}|\le 1\}} \right]   \ge  \ I_{N,\,\go} \cdot I_{M,\,\theta^{N}\go} \qquad \Rightarrow \qquad I_{N+M,\,\go} \ge  \ I_{N,\,\go} \cdot I_{M,\,\theta^{N}\go}
\]
so that the superadditivity of the process $\{\log I_{N,\,\go}\}_N$ is proved.
We can thus apply Kingman's Theorem, concluding that the sequence $\{\log I_{N,\,\go}^{\gl,h} / N\}_N$ converges $\bbP(\dd\go)$--a.s. and in $\bbL_1(\bbP)$ to a limit $f_\go(\lambda, h)$ which is $\theta$--invariant. By tail triviality, $f$ is $\bbP(\dd\go)$--a.s. constant and we consequently omit the $\go$ dependence: \mbox{$f = f(\gl, h)$}.

\smallskip

\subsection{Step 2}
\label{sec:b}

Now we show that also the sequence $\{\log J_{N,\,\go}^{\gl, h} / N\}_N$ has, $\bbP(\dd\go)$--a.s. and in $\bbL_1(\bbP)$, the limit $f(\gl, h)$ as $N \to \infty$. We start noting that by definition
\begin{equation} \label{liminf}
J_{N,\,\go}^{\gl, h} \: \ge \: I_{N,\,\go}^{\gl,h} \qquad \Rightarrow \qquad \liminf_{N\to \infty} \frac{\log J_{N,\,\go}^{\gl, h}}{N} \: \ge \: f(\gl, h)
\end{equation}
for $\bbP$--a.e.~$\go$, so it remains to find a bound for the $\limsup$.

We recall that by hypothesis the density~$f$ of~$X_1$ is supported in the interval~$[-1,+1]$,
and that for some~$n_0 \in \N$ the density~$f_{n_0}$ of~$S_{n_0}$
is bounded, hence we can find two positive
constants~$A,\,M$ such that $f_{n_0}(x) \le A \cdot \ind_{\{|x|\le M\}}$.
Then for~$N\in\N$ by the Markov property we get a first upper bound for~$J$:
\begin{align} \label{eq:first_est}
\begin{split}
    J_{n_0 + N,\,\go} & \;=\; \int_{\R} \dd z \: f_{n_0} (z) \: \bE \big[ \:
    \cG_{n_0,\go} (S_1, ...\, , S_{n_0-1}, z)  \big] \cdot \bE_z \left[ \:
    \cG_{N,\,\theta^{n_0}\go} \, \ind_{\{|S_{N}|\le 1\}} \right] \\
    & \;\le\; A\, K(\go) \cdot \int_{-M}^M \dd z \:  \bE_z \left[ \:
    \cG_{N,\,\theta^{n_0}\go} \cdot \ind_{\{|S_{N}|\le 1\}} \right]
\end{split}
\end{align}
where the constant $K (\go)$ that we have used to bound $\bE [\ldots]$ is simply
\[
K(\go) = K(\gl, h, \{\go_i\}_{1\le i \le n_0}) := \exp \left( \gl \sum_{n=1}^{N_0} (|\go _n| + h) \right)\,.
\]

Next we want to obtain an analogous {\sl lower bound} for~$I$. Observe that by~\eqref{LLTeq}
we can find~$n_1 \in \N$ such that
\begin{equation*}
    f_{n_1} (y) \ge \frac{1}{2} \, \frac{1}{\sqrt{2\pi}} \, \frac{1}{\sqrt{n_1}}
    \qquad \quad \forall y \in [-(M+1), (M+1)]\,.
\end{equation*}
Then for~$N\in\N$ and for all~$x\in[-1,+1]$ we get
\begin{align*}
    & \bE_x \Big[ \cG_{n_1 + N,\go} \, \ind_{\{|S_{n_1 + N}| \le 1\}} \Big] \\
    & \qquad \;=\; \int_{\R} \dd z \: f_{n_1} (z-x) \: \bE \big[ \: \cG_{n_1,\go}
    (S_1, ...\, , S_{n_1-1}, z)  \big] \cdot \bE_z \left[ \: \cG_{N,\,\theta^{n_1}\go} \,
    \ind_{\{|S_{N}|\le 1\}} \right] \\
    & \qquad \;\ge\; \frac{1}{2} \, \frac{1}{\sqrt{2\pi}} \, \frac{1}{\sqrt{n_1}}
    \, K(\go)^{-1} \cdot \int_{-M}^M \dd z \:  \bE_z \left[ \:
    \cG_{N,\,\theta^{n_1}\go} \cdot \ind_{\{|S_{N}|\le 1\}} \right] \,,
\end{align*}
hence
\begin{equation}\label{eq:second_est}
    I_{n_1 + N,\go} \;\ge\; \frac{1}{2} \, \frac{1}{\sqrt{2\pi}} \, \frac{1}{\sqrt{n_1}}
    \, K(\go)^{-1} \cdot \int_{-M}^M \dd z \:  \bE_z \left[ \:
    \cG_{N,\,\theta^{n_1}\go} \cdot \ind_{\{|S_{N}|\le 1\}} \right]\,.
\end{equation}

Combining~\eqref{eq:first_est} with~\eqref{eq:second_est} we get that for all~$N\in\N$
\begin{equation*}
    J_{N,\go} \ \le \ A' \, K(\go)^2 \, I_{N + (n_1 - n_0), \theta^{(n_0 - n_1)} \go}\,,
\end{equation*}
for some positive constant~$A'$ (we recall that we consider two--sided sequence
of charges: $\go=\{\go_n\}_{n\in\Z}$, hence the translations~$\theta^k$ are meaningful
for all~$k\in\Z$). It follows that
\begin{equation*}
    \limsup_{N\to\infty} \frac{\log J_{N,\,\go}^{\gl,h}}{N} \; \le \;
    \limsup_{N\to\infty} \frac{\log I_{N + (n_1-n_0),\, \theta^{(n_0 - n_1)}\go}^{\gl,h}}{N}
    \;=\; f(\gl, h)\,,
\end{equation*}
$\bbP(\dd\go)$--a.s., that is what was to be proven. Notice that the bounds we have obtain yield
easily also the~$\bbL_1(\bbP)$ convergence of $\{\log J_{N,\,\go}^{\gl, h} / N\}_N$
towards~$f(\gl,h)$.

\smallskip

\subsection{Step 3}
\label{sec:c}

Finally we are left with comparing~$J_{N,\go}$ with the original partition
function~$\tilde Z_{N,\go}$, which amounts to removing the restriction $\{|S_N| \le 1\}$.
Observe that by definition $\tilde Z_{N,\go} \ge J_{N,\go}$, hence we can concentrate
on finding a suitable upper bound.

The procedure we follow is very similar to the simple random walk case, cf.~\cite{cf:G}.
The idea is to look at the last point up to epoch~$N$ at which the random walk
changes its sign. More precisely, we define the random variable~$U$ by
\begin{equation*}
    U \;:=\; \min \big\{ k \in \{1, \ldots, N\} :\, \sign(S_k) = \sign(S_{k+1})
    = \ldots = \sign(S_N) \big\}\,,
\end{equation*}
and we disintegrate the partition function according to the range of~$U$. It is convenient
to consider separately the cases~$\{S_N>0\}$ and~$\{S_N <0\}$, that is we split
\begin{equation*}
    \tilde Z_{N,\go} \;=\; \tilde Z_{N,\go}^> \;+\; \tilde Z_{N,\go}^< \;:=\; \bE \big[ \cG_{N,\go}\,,\, S_N > 0  \big]
    \;+\; \bE \big[ \cG_{N,\go}\,,\, S_N < 0  \big]\,.
\end{equation*}
Then we can write
\begin{align*}
    & \tilde Z_{N,\go}^> \;=\; \sum_{k=1}^N \bE \big[ \cG_{N,\go}\,,\, U=k \,,\,S_N>0 \big]\\
    & \quad \;=\; \sum_{k=1}^N \int_{0}^1 \dd z \, f_k(z) \, \bE \big[ \cG_{k,\go} (S_1, \ldots, S_{k-1}, z) \big]
    \, e^{\gl \sum_{i=k+1}^N (\go_i + h)} \, \bP_z \big[ S_1 > 0, \ldots, S_{N-k} > 0 \big]\,.
\end{align*}
Now by the asymptotic behavior in~\eqref{eq:rw_tail_x} and~\eqref{eq:rw_point_x}
it follows that one can find a positive constant~$D$ such that for all~$z \in [0,1]$,
for all~$N\in\N$ and for all~$k \in \{1,\ldots, N\}$ one has
\begin{equation*}
    \bP_z \big[ S_1 > 0, \ldots, S_{N-k} > 0 \big] \;\le\; D\, N\,
    \bP_z \big[ S_1 > 0, \ldots, S_{N-k-1} > 0, S_{N-k} \le 0 \big]\,.
\end{equation*}
Performing this substitution we obtain
\begin{align*}
    \tilde Z_{N,\go}^> \;\le\; e^{2\gl |\go_N|} \, D \, N \, \bE \big[ \cG_{N,\go}\,,\, |S_N| \le 1 \big]
        \;=\; \big( e^{2\gl |\go_N|} \, D \, N \big) \, J_{N,\go} \,.
\end{align*}
As the very same arguments can be performed for~$\tilde Z_{N,\go}^<$, we have definitively shown that
\begin{align*}
    J_{N,\go}^{\gl,h} \;\le\; \tilde Z_{N,\go}^{\gl,h} \;\le\; \big( e^{2\gl |\go_N|}
    \, D' \, N \big) \, J_{N,\go}^{\gl,h} \,,
\end{align*}
for some positive constant~$D'$. From this relation the convergence of~$\{\log \tilde Z_{N,\go}^{\gl,h}/N\}_n$
towards~$f(\gl,h)$ both~$\bbP(\dd\go)$--a.s. and in~$\bbL_1(\bbP)$ follows immediately.


\smallskip
\section{Towards the coarse graining of the free energy}

The coarse graining of the free energy for the copolymer near a selective interface
model is expressed by Theorem~\ref{th:coarse} of Chapter~\ref{ch:first}: it holds
when the underlying random walk is the simple symmetric random walk on~$\Z$, and
the proof of it is the main result of the paper~\cite{cf:BdH} by Bolthausen and
den Hollander. The purpose of this section is to discuss the issue of extending
it to the continuous setting adopted in this chapter.

The idea that lies at the basis of the coarse graining is that when~$\gl \to 0$
the reward to stay close to the interface gets small and consequently
the typical excursions of the polymer away from the interface tend to become
very long. Therefore it
should be possible to approximate both the polymer and the charges by Brownian motions,
and in fact Theorem~\ref{th:coarse} provides a quantitative version of this approximation.

According to this heuristic point of view, the microscopic details of the random walk and
of the charges should not be too relevant, until we work with processes
in the domain of attraction of the Brownian motion.
This is indeed true for the charges, as we already mentioned in Chapter~\ref{ch:first}:
in fact the original proof of Theorem~\ref{th:coarse} in~\cite{cf:BdH}
was given for the Bernoulli case $\bbP(\go_1=+1) = 1-\bbP (\go_1 = -1) = 1/2$,
but it can be easily extended to the general $\go$~case considered here.

On the other hand, the extension to the more general random
walks considered in this chapter appears to be more challenging. In order to outline
the reasons of this fact, we have to look more closely at the original proof
of Theorem~\ref{th:coarse}. Without going into the details, which
are quite long and extremely delicate, we point out that the proof is divided
in four main steps, which we can roughly describe as follows:
\begin{enumerate}
\smallskip
\item\label{point1} first it is shown that when~$\gl$ and~$h$ are small one can safely throw away
the {\sl short excursions} of the walk in the computation of the partition function;
\smallskip
\item\label{point2} then the~$\{\go_n\}$ are replaced by standard Gaussian variables;
\smallskip
\item\label{point3} the law of the (long) excursions under the
rescaled simple random walk measure is then
replaced by the law under the Brownian motion measure, ending up with a Brownian
copolymer model without the short excursions;
\smallskip
\item\label{point4} finally, one reintroduces the short excursion for the Brownian copolymer model.
\smallskip
\end{enumerate}
We observe that in step~\eqref{point2} the random walk plays a minor role
and it is not difficult to adapt the proof to our continuous random walk setting,
while step~\eqref{point4} is a problem involving only the Brownian copolymer model,
hence it requires no change. Therefore the crucial points are step~\eqref{point1} and
step~\eqref{point3}, that will be analyzed separately.

\smallskip
\subsection{Step (1): throwing away the short excursions}

For the first step the original
proof makes use of several peculiar properties enjoyed by the excursions of the
simple random walk, namely:
\begin{itemize}
\smallskip
\item\label{propr:first} there is a complete decoupling between the epochs of the returns to
zero~$\{\tau_k^{\text{SRW}}\}_k$ and the signs~$\{\gs_k^{\text{SRW}}\}_k$ of the excursions:
in fact the sequence~$\{\tau_k^{\text{SRW}}\}_k$ is independent of the
sequence~$\{\gs_k^{\text{SRW}}\}_k$;
\smallskip
\item the signs $\{\gs_k^{\text{SRW}}\}_k$ form an independent sequence of
Bernoulli variables with $\bP (\gs_k^{\text{SRW}} = +1) = \bP (\gs_k^{\text{SRW}} = -1) = 1/2$;
\smallskip
\item\label{propr:last} the zeros $\{\tau_k^{\text{SRW}}\}_k$ form a classical renewal process,
that is the interarrival times $\{\tau_{k+1}^{\text{SRW}} - \tau_{k}^{\text{SRW}}\}_k$
are independent positive random variables.
\smallskip
\end{itemize}

%

The first observation is that the returns to zero are no longer meaningful for a continuous random walk,
since $\bP(S_n = 0\ \text{for some $n$}) = 0$. We point out two possible definitions:
\begin{itemize}
\smallskip
\item[(a)] the epochs at which the random walk crosses the interface:
\begin{equation} \label{eq:cross}
    \tau_0 \;:=\; 0 \qquad \quad \tau_{k+1} \;:=\; \inf \big\{ n > \tau_k :\,
    \sign(S_n) \neq \sign(S_{n-1}) \big\}\,,
\end{equation}
with the signs of the excursions $\{\gs_k\}_{k\ge 1}$ defined by
$\gs_k := \sign(S_{\tau_{k-1}})$;
\smallskip
\item[(b)] the epochs at which the walk gets close to the interface:
\begin{equation} \label{eq:close}
    \tau_0 \;:=\; 0 \qquad \quad \tau_{k+1} \;:=\; \inf \big\{ n > \tau_k :\,
    S_n \in [-1,+1] \big\}\,,
\end{equation}
with the signs of the excursions $\{\gs_k\}_{k\ge 1}$
defined by $\gs_k := \sign(S_{\tau_{k}})$.
\smallskip
\end{itemize}
Notice that with the first definition there is a striking
difference with respect to the simple random
walk case, because the sequence~$\{\gs_k\}_{k\ge 1}$
is almost deterministic: in fact $\bP(\dd S)$--a.s. we have that
$\{\gs_k\}_{k\ge 1} = \sign(S_1) \cdot \{(-1)^k\}_{k \ge 1}$.

\smallskip

In any case none of the above mentioned properties
of~$\{\gs_k^{\text{SRW}}\}_k$ and~$\{\tau_k^{\text{SRW}}\}_k$ holds
anymore for~$\{\gs_k\}_k$ and~$\{\tau_k\}_k$, with any of the two definitions
(a) or~(b).

\smallskip

The most serious problem is that the interarrival times $\{T_k\}_k$ where
$T_k := \tau_{k+1} - \tau_{k}$ are no longer independent. Nevertheless
they enjoy a useful property. Let us introduce the
sequence of random variables~$\{J_k\}_{k\ge 0}$ defined by $J_k := S_{\tau_k}$
(notice that $J_k \in [-1,+1]$ because by hypothesis our random walk has steps bounded
by~$1$ in absolute value). Then it is not difficult to check that the joint process
$\{(J_k,T_k)\}_{k}$ is a {\sl Markov renewal process}~\cite{cf:Asm}, that is a Markov
chain on~$[-1,+1] \times \N$ such that the transition kernel
\begin{equation*}
    \bP \big( J_{k+1} \in \dd y,\, T_{k+1} = n \,\big|\, J_k = x,\, T_k = m  \big)
\end{equation*}
does not depend on~$m$. This implies that, conditionally on~$\{J_k\}_k$, the variables
$\{T_k\}_k$ are independent.
%

Another remarkable fact is that the asymptotic behavior of the probability tail
of the variables~$\{T_k\}_k$ is similar to the simple
random walk case, also conditionally on the~$\{J_k\}_k$:
\begin{equation*}
    \bP \big( T_{k+1} = n \,\big|\, J_k = x \big) \;\sim\; \frac{\text{const.}(x)}{n^{3/2}}
    \qquad \quad (n\to\infty)\,.
\end{equation*}
If one chooses the definition~(a) this relation is just a rephrasing of~\eqref{eq:rw_point_x},
and it is not difficult to check that it holds also with definition~(b).

Thus the situation is not extremely bad. After all,
we have seen that Markov renewal processes have been the fundamental tool
in the study of periodic inhomogeneous polymer models performed in Chapter~\ref{ch:cgz}.
The reason is that a lot of fundamental asymptotic results (renewal theorems)
of classical renewal processes can be extended to the Markov case. We stress
however that the processes of Chapter~\ref{ch:cgz} enjoy the
peculiar property of having a modulating chain~$\{J_k\}_k$ with finite state space
and this is indeed a great simplification, as it is explained in~\cite[Ch.~VII.4]{cf:Asm}.
Dealing with the case when the modulating chain has uncountable state space is
much more delicate and the results are more involved (see for instance~\cite{cf:Als}),
especially in the case of heavy tails.

\smallskip

Up to now we have not succeeded in extending the proof of step~\eqref{point1}
to the continuous random walk setting.
Nevertheless we point out that we are able to prove a weaker form of
step~\eqref{point1} that,
{\sl provided one can extend to the continuous setting step~\eqref{point3}},
is sufficient to yield the first part of~Theorem~\ref{th:coarse},
namely the scaling limit of the free energy expressed by
equation~\eqref{eq:scaling}. This would be an interesting result, but unfortunately
we do not have yet a complete proof of step~\eqref{point3} in the continuous setting.

\smallskip
\subsection{Step (3): from random walk to Brownian motion}

The central point of step~\eqref{point3} is a sharp comparison between the law
of the (long) excursions of the rescaled simple random walk and the law of the excursions
of the Brownian motion. Without getting too much into the details, we mention that
a fundamental estimate of the proof is the following one (cf. equations (4.62) and~(4.66)
in~\cite{cf:BdH}): for~$k,l \in \N$ such that~$k+l \in 2\N$, as~$k,l \to \infty$
jointly we have
\begin{equation} \label{eq:relrel}
    \bP^{(SRW)} \big( S_i \neq 0 \text{ for } k < i < k+l,\, S_{k+l} = 0 \big)
    \;=\; (1+o(1)) \, \frac{2}{\pi} \, \frac{\sqrt k}{(k+l) \sqrt l}\,.
\end{equation}
where~$\bP^{(SRW)}$ is the law of the simple symmetric random walk on~$\Z$.
The proof of this relation is obtained by conditioning on the position
of the walk at epoch~$k$ and then using the reflection principle together
with a strong approximation of the mass function of the simple random walk
by the Gaussian density.

Now let us set~$\ell_N := \max \{ k = 0, \ldots, N: \, \tau_k \le N\}$, where the
$\{\tau_k\}_k$ are defined by \eqref{eq:cross} (but we could also choose definition
\eqref{eq:close}).
Then the continuous analogue of~\eqref{eq:relrel} should be that for~$k,l \in \N$ and
as~$k,l \to \infty$ jointly
\begin{equation} \label{eq:reller}
    \bP \big( \ell_{k+l-1} \le k,\, \ell_{k+l} = k+l \big) \;=\; (1+o(1)) \,
    \frac{1}{\pi} \, \frac{\sqrt k}{(k+l) \sqrt l}
\end{equation}
(the reason for the missing factor~$2$ with respect to~\eqref{eq:relrel}
lies in the periodicity of the returns of the simple random walk).

In order to prove~\eqref{eq:reller}, we start conditioning on the position at
epoch~$k$:
\begin{align*}
    \bP \big( \ell_{k+l-1} \le k,\, \ell_{k+l} = k+l \big) \;=\; & \int_0^{+\infty} \dd x\,
    f_k(x) \, \bP_x \big( S_1 > 0, \ldots, S_{l-1} > 0, S_l \le 0 \big)\\
    & \; + \; \int_{-\infty}^0 \dd x\,
    f_k(x) \, \bP_x \big( S_1 < 0, \ldots, S_{l-1} < 0, S_l \ge 0 \big)\,.
\end{align*}
Let us consider the first integral in the r.h.s. above, the second one being analogous:
when~$k$ is large the asymptotic behavior of $f_k(x)$ is given by the Local Central
Limit Theorem~\eqref{eq:unif_llt} (actually one should use a stronger version valid
in a ratio sense, see~\cite{cf:Ric}). On the other hand the asymptotic behavior
of the term~$\bP_x(S_1 > 0, \ldots, S_{l-1} > 0, S_l \le 0)$ is not immediate: notice in
fact that the relevant values of~$x$ are those of order~$\sqrt k$ and~$k\to\infty$,
hence one cannot use~\eqref{eq:rw_point_x}.

Conditioning on the position at epoch~$l$ we can write
\begin{equation*}
    \bP_x(S_1 > 0, \ldots, S_{l-1} > 0, S_l \le 0) \;=\;
    \int_{0}^{1} \dd y\, \phi_l^{(y)}(x+y) \, f_l(-y-x)
\end{equation*}
where~$\phi_l^{(y)}(z)$ is the value at~$z \in \R$ of the density of the random
variable~$\widehat S_l$ {\sl conditionally on the event}
$\{\widehat S_1 >y, \ldots, \widehat S_l > y\}$, where we have introduced the dual random
walk $\{\widehat S_n\}_n := \{ -S_n\}_n$. Notice however that the value of~$y \in (0,1)$
is actually irrelevant for the asymptotic behavior of~$\phi_l^{(y)}(z)$, because we are
interested in the regime when both~$l$ and~$z$ are large, and it is sufficient to
consider the case~$y=0$.

Therefore an important role is played by
the asymptotic behavior of the density of the variable~$\widehat S_l$
conditionally on the event~$\{\widehat S_1 > 0, \ldots, \widehat S_l > 0\}$, where
$\{\widehat S_n\}_n = \{-S_n\}_n$ is a random walk satisfying the hypothesis stated
in~\S~\ref{sec:deff}. It has been known for a long time~\cite{cf:Bol76} that the only hypothesis
of finite nonzero variance~$\gs$ guarantees the weak convergence
\begin{equation*}
    \frac{\widehat S_l}{\gs \sqrt{l}} \text{ conditionally on } \{\widehat S_1 > 0,
    \ldots, \widehat S_l > 0\} \quad \Rightarrow \quad x e^{-x^2/2} \ind_{(x\ge 0)} \dd x\,.
\end{equation*}
However what we need is rather a local refinement of this weak convergence, exactly as
the Local Limit Theorem~\eqref{eq:unif_llt} is a local refinement of the Central Limit Theorem.

\smallskip

Such a Local Limit Theorem for random walks conditioned to stay positive
does not seem to be known in the literature. We give a proof
in Chapter~\ref{ch:llt} in a very general setting, using the Fluctuation Theory for random
walks. Besides being an interesting result in itself, this theorem is a key step
to prove the asymptotic behavior~\eqref{eq:reller}.
Unfortunately there is still some technical points to be solved
in order to extend the proof of step~\eqref{point3} to the continuous setting,
but we think that a complete solution is not too far.


\chapter[LLT for random walks conditioned to stay positive]
{A local limit theorem for random walks\\ conditioned to stay positive}
\label{ch:llt}


In this chapter we study the asymptotic behavior of random walks conditioned to stay positive.
We consider a real random walk $S_n=X_1+\ldots +X_n$ attracted (without centering) to the normal law: this means that for a suitable norming sequence $a_n$ we have the weak convergence $S_n/a_n \Rightarrow \gp(x)\dd x$, $\gp(x)$ being the standard normal density. A local refinement of this convergence is provided by Gnedenko's and Stone's Local Limit Theorems, in the lattice and nonlattice case respectively.

Now let $\cC_n$ denote the event $(S_1 > 0, \ldots, S_n > 0 )$ and let $S_n^+$ denote the random variable $S_n$ conditioned on $\cC_n$: it is known that $S_n^+/a_n \Rightarrow \gp^+(x) \dd x$, where $\gp^+(x):=x \exp(-x^2/2)\ind_{(x\geq 0)}$. What we are going to establish is an equivalent of Gnedenko's and Stone's Local Limit Theorems for this weak convergence. We also consider the particular case when $X_1$ has an absolutely continuous law: in this case the uniform convergence of the density of $S_n^+/a_n$ towards $\gp^+(x)$ holds under a standard additional hypothesis, in analogy to the classical case. We finally discuss an application of our main results to the asymptotic behavior of the joint renewal measure of the ladder variables process. Unlike the classical proofs of the LLT, we make no use of characteristic functions: our techniques are rather taken from the so--called Fluctuation Theory for random walks.

\smallskip
The article~\cite{cf:C} has been taken from the content of this chapter.

\smallskip
\section{Introduction and results} \label{sec:intro_llt}

\smallskip
\subsection{The nonlattice case} Let $S_n=X_1+\ldots +X_n$ be a real random walk attracted (without centering) to the normal law. This means that $\{X_k\}$ is an IID sequence of real random variables, and for a suitable norming sequence $a_n$ we have the weak convergence
\begin{equation} \label{eq:clt}
    S_n/a_n \Rightarrow \gp(x)\,\dd x\,, \qquad \gp(x) := \frac{1}{\sqrt{2\pi}}\, e^{-x^2/2}\,.
\end{equation}
This is the case for example when $\bE(X_1)=0$ and $\bE(X_1^2) =: \gs^2 \in (0,\infty)$ with $a_n:= \gs \sqrt{n}$, by the Central Limit Theorem.

\smallskip

We recall that, by the standard theory of stability \cite[\S{}IX.8 \& \S{}XVII.5]{cf:Fel2}, for equation \eqref{eq:clt} to hold it is necessary and sufficient that $\bE(X_1)=0$, that the truncated variance $V(t):=\bE(X_1^2 \, \ind_{(|X_1|\leq t)})$ be \textsl{slowly varying} at~$\infty$ (that is $V(ct)/V(t) \to 1$ as $t\to\infty$ for every $c>0$) and that the sequence $a_n$ satisfy the condition $a_n^2 \sim n V(a_n)$ as $n\to\infty$.

\smallskip

For the moment we assume that the law of~$X_1$ is \textsl{nonlattice}, that is not supported in $(b+c\Z)$ for any $b\in\R, c>0$. Then a local refinement of \eqref{eq:clt} is provided by Stone's Local Limit Theorem \cite{cf:Sto65,cf:Sto67}, that in our notations reads as (cf.~\cite[\S{}8.4]{cf:BinGolTeu})
\begin{equation} \label{eq:llt}
    a_n\, \bP \big( S_n \in [x,x+h)\, \big) = h\, \gp(x/a_n) + o(1) \qquad (n\to\infty)\,,
\end{equation}
\textsl{uniformly} for $x\in\R$ and $h$ in compact sets in~$\R^+$.

\smallskip

What we are interested in is the asymptotic behavior of the random walk $\{S_n\}$ conditioned to stay positive. More precisely, let $\cC_n:=(S_1 > 0, \ldots, S_n > 0)$ and let $S_n^+$ denote the random variable $S_n$ under the conditional probability $\bP(\,\cdot\, |\,\cC_n)$: if \eqref{eq:clt} holds then one has an analogous weak convergence result for $S_n^+/a_n$, namely
\begin{equation} \label{eq:pos_clt}
    S_n^+/a_n \Rightarrow \gp^+(x)\,\dd x\,, \qquad \gp^+(x) := x\, e^{-x^2/2} \,\ind_{(x\geq 0)}\,.
\end{equation}
This is an immediate consequence of the fact \cite{cf:Igl74,cf:Bol76,cf:Don85} that, whenever \eqref{eq:clt} holds, the whole process $\{S_{\lfloor nt \rfloor}/a_n\}_{t\in[0,1]}$ under $\bP(\,\cdot\, |\,\cC_n)$ converges weakly as $n\to\infty$ to the standard Brownian meander process $\{B_t^+\}_{t\in[0,1]}$, and $\gp^+(x)\,\dd x$ is the law of $B_1^+$, cf.~\cite{cf:RevYor}.

\medskip

Our main result is an analogue of Stone's LLT for the weak convergence~\eqref{eq:pos_clt}.

\smallskip

\begin{theorem} \label{th:main_llt}
If $X_1$ is nonlattice and \eqref{eq:clt} holds, then
\begin{equation}\label{eq:pos_llt}
a_n\, \bP \big( S_n \in [x,x+h)\, \big| \,\cC_n \big) = h\, \gp^+(x/a_n) + o(1) \qquad (n\to\infty)\,,
\end{equation}
uniformly for $x\in\R$ and $h$ in compact sets in~$\R^+$.
\end{theorem}

\smallskip

The main difficulty with respect to the classical case is given by the fact that under the conditional probability $\bP(\,\cdot\, |\,\cC_n)$ the increments of the walk $\{X_k\}$ are no longer independent. This is a major point in that the standard proof of Stone's LLT relies heavily on characteristic functions methods. As a matter of fact, we make no use of characteristic functions: our methods are rather of combinatorial nature, and we make an essential use of the so--called Fluctuation Theory for random walks. The core of our proof consists in expressing the law of $S_n$ under $\bP(\,\cdot\, |\,\cC_n)$ as a suitable mixture of the laws of $\{S_k\}_{0\leq k \leq n}$ \textsl{under the unconditioned measure}~$\bP$, to which Stone's LLT can be applied. Thus our ``Positive LLT'' is in a sense directly derived from Stone's LLT.

\smallskip

%

We point out that our methods may in principle be applied to the case when the random walk is attracted to a generic stable law (the analogue of~\eqref{eq:pos_clt} in this case is also provided by \cite{cf:Don85}), so that it should be possible to obtain an equivalent of Theorem~\ref{th:main_llt} in this general setting.


\smallskip

\subsection{The lattice case}

Let us consider now the lattice case: we assume that $X_1$ is supported in $(b+c\Z)$, for the least such~$c$. In this case the local version of~\eqref{eq:clt} is given by Gnedenko's Local Limit Theorem \cite[\S{}8.4]{cf:BinGolTeu}, which says that
\begin{equation} \label{eq:lat_llt}
    \frac{a_n}{c}\, \bP\big( S_n=bn+cx \big) = \gp\big((bn+cx)/a_n\big) + o(1) \qquad (n\to\infty)\,,
\end{equation}
uniformly for $x\in\Z$.

\medskip

We can derive the local version of~\eqref{eq:pos_clt} also in this setting.

\begin{theorem} \label{th:lat_main}
If $X_1$ is lattice with span~1 and \eqref{eq:clt} holds, then
\begin{equation*}
\frac{a_n}{c}\, \bP\big( S_n=bn+cx \,\big|\, \cC_n \big) = \gp^+\big((bn+cx)/a_n\big) + o(1) \qquad (n\to\infty)\,,
\end{equation*}
uniformly for $x\in\Z$.
\end{theorem}

The proof is omitted since it can be recovered from the proof of Theorem~1 with only slight modifications (some steps are even simpler).


\smallskip

\subsection{The density case}

When the law of $X_1$ is absolutely continuous with respect to Lebesgue measure and \eqref{eq:clt} holds, one may ask whether the density of $S_n/a_n$ converges to $\gp(x)$ in some pointwise sense. However, it is easy to build examples \cite[\S{}46]{cf:GneKol} satisfying \eqref{eq:clt}, such that for every~$n$ the density of $S_n/a_n$ is unbounded in any neighborhood of~0: therefore without some extra--assumption one cannot hope for convergence to hold at each point. Nevertheless, if one looks for the \textsl{uniform convergence} of the density, then there is a simple condition which turns out to be necessary and sufficient.

\begin{assumption} \label{ass:main}
The law of $X_1$ is absolutely continuous, and for some $k\in\N$ the density $f_k(x)$ of $S_k$ is essentially bounded: $f_k(x) \in L^\infty(\R,\dd x)$.
\end{assumption}

It is easy to see that if this assumption holds, then for large~$n$ the density $f_n(x)$ admits a bounded and continuous version. A proof that Assumption~\ref{ass:main} yields the uniform convergence of the (continuous versions of the) density of $S_n/a_n$ towards~$\gp(x)$, namely
\begin{equation*}
    \sup_{x\in\R} \big| a_n f_n(a_n x) - \gp(x) \big| \to 0 \qquad (n\to\infty)\,,
\end{equation*}
can be found in~\cite[\S{}46]{cf:GneKol}. On the other side, the necessity of Assumption~\ref{ass:main} for the above convergence to hold is evident.

\medskip

We can derive a completely analogous result for $S_n^+$.

\begin{theorem} \label{th:den_main}
Assume that $X_1$ satisfies Assumption~\ref{ass:main}, and that \eqref{eq:clt} holds. Then:
\begin{enumerate}
\item $S_n^+$ has an absolutely continuous law, whose density $f_n^+(x)$ is bounded and continuous (except at $x=0$) for large~$n$;
\item the (continuous version of the) density of $S_n^+/a_n$ converges uniformly to~$\gp^+(x)$:
\begin{equation*}
    \sup_{x\in\R} \big| a_n f_n^+(a_n x) - \gp^+(x) \big| \to 0 \qquad (n\to\infty)\,.
\end{equation*}
\end{enumerate}
\end{theorem}

This Theorem can be proved following very closely the proof of Theorem~1: in fact equation \eqref{eq:main_llt} in Section~\ref{sec:first} provides an explicit expression for $f_n^+(x)$, that can be shown to converge to $\gp^+(x)$ with the very same arguments given in Section~\ref{sec:second}.


\smallskip

\subsection{Asymptotic behavior of the ladder renewal measure}

As a by--product of the Local Limit Theorems described above, we have a result on the asymptotic behavior of the renewal measure of the ladder variables process. For simplicity we take the arithmetic setting, assuming that $X_1$ is supported by~$\Z$ and it is aperiodic, but everything works similarly in the general lattice and nonlattice cases. The renewal mass function $u(n,x)$ of the ladder variables process is defined for $n\in\N,\ x\in\Z$ by
\begin{equation} \label{eq:def_u}
    u(n,x) := \sum_{r=0}^\infty \bP \big( T_r=n, H_r=x \big) = \bP \big( n \text{ is a ladder epoch}, S_n = x \big)\,,
\end{equation}
where $\{(T_k,H_k)\}$ is the (strict, ascending) ladder variables process associated to the random walk (the definitions are given in Section~\ref{sec:fluct}). Generalizing some earlier result of \cite{cf:Kee}, in \cite{cf:AliDon} it has been shown that, for $\{x_n\}$ such that $x_n/a_n\to 0$,
\begin{equation} \label{eq:as_ren}
    u(n,x_n) \sim \frac{1}{\sqrt{2\pi}\, n\, a_n}\, U(x_n-1) \sim \frac 1n\, \bP\big( S_n=x_n \big) \, U(x_n-1) \qquad (n\to\infty)\,,
\end{equation}
where $U(x) := \sum_{r=0}^\infty \bP (H_r \leq x)$ is the distribution function of the renewal measure associated to the ladder heights process (as a matter of fact, the proof of \eqref{eq:as_ren} given in~\cite{cf:AliDon} is carried out under the assumption that $\bE(X_1^2)<\infty$, but it can be easily extended to the general case).

\medskip

With our methods we are able to show that the same relation is valid for $x = O(a_n)$, with no further restriction on~$X_1$ other than the validity of~\eqref{eq:clt}.

\begin{theorem} \label{th:as_ren}
Let $X_1$ be arithmetic with span~1 and such that equation \eqref{eq:clt} holds. Then for $x\in\Z$
\begin{equation} \label{eq:ext_as_ren}
    u(n,x) = \frac 1n\, \bP\big( S_n=x \big) \, U(x-1) \, \big( 1+o(1) \big) \qquad (n\to\infty)\,,
\end{equation}
uniformly for $x/a_n \in [\gep, 1/\gep]$, for every fixed $\gep>0$.
\end{theorem}

The proof of this theorem is a direct consequence of Theorem~\ref{th:lat_main}: the details are worked out in Section~\ref{sec:as_ren}.

\smallskip

Notice that in the r.h.s. of \eqref{eq:ext_as_ren} we could as well write $U(x)$ instead of $U(x-1)$, since $x\to\infty$ as $n\to\infty$. Also observe that putting together equation~\eqref{eq:as_ren} with Theorem~\ref{th:as_ren} one has the stronger result that equation~\eqref{eq:ext_as_ren} holds uniformly for $x/a_n \in [0,K]$, for every fixed $K>0$.

\smallskip

We point out that Theorem~\ref{th:as_ren} has been obtained also in \cite{cf:BryDon}, where the authors study random walks conditioned to stay positive in a different sense.


\smallskip

\subsection{Outline of the exposition}

The exposition is organized as follows: in Section~\ref{sec:fluct} we recall some basic facts on Fluctuation Theory and stable laws, and we set the relative notation; we also give the proof of Theorem~\ref{th:as_ren}. The rest of the chapter is devoted to the proof of Theorem~\ref{th:main_llt}, which has been split in two parts. The first one, in Section~\ref{sec:first}, contains the core of the proof: using Fluctuation Theory we obtain an alternative expression for the law of~$S_n^+$, see equation~\eqref{eq:main_llt}, and we prove a crucial weak convergence result connected to the renewal measure of the ladder variables process. Then in Section~\ref{sec:second} we apply these preliminary results, together with Stone's LLT, to complete the proof. Finally, some minor points have been deferred to the appendix.


\smallskip
\section{Fluctuation Theory and some applications} \label{sec:fluct}

In this section we are going to recall some basic facts about Fluctuation Theory for random walks, especially in connection with the theory of stable laws, and to derive some preliminary results. Standard references on the subject are~\cite{cf:Fel2} and~\cite{cf:BinGolTeu}.


\smallskip

\subsection{Regular variation}
A positive sequence $d_n$ is said to be \textsl{regularly varying} of index $\ga\in\R$ (we denote this by $d_n \in R_{\ga}$) if $d_n \sim L_n\, n^\ga$ as $n\to \infty$, where $L_n$ is \textsl{slowly varying} at~$\infty$ in that $L_{\lfloor tn\rfloor}/L_n\to 1$ as $n\to\infty$, for every $t>0$. If $d_n \in R_\ga$ with $\ga\neq 0$, up to asymptotic equivalence we can (and will) always assume \cite[Th.1.5.3]{cf:BinGolTeu} that $d_n = d(n)$, with $d(\cdot)$ a continuous, strictly monotone function, whose inverse will be denoted by~$d^{-1}(\cdot)$. Observe that if $d_n \in R_\ga$ then $d^{-1}(n) \in R_{1/\ga}$ and $1/d_n \in R_{-\ga}$.

\smallskip

Let us recall two basic facts on regularly varying sequences that will be used a number of times in the sequel. The first one is a uniform convergence property \cite[Th.1.2.1]{cf:BinGolTeu}: if $d_n \in R_\ga$, then
\begin{equation}\label{eq:RV_unif_conv}
    d_{\lfloor tn \rfloor} = t^\ga \,d_n \,\big(1+o(1)\big) \qquad (n\to\infty)\,,
\end{equation}
uniformly for $t\in[\gep, 1/\gep]$, for every fixed $\gep>0$. The second basic fact \cite[Prop.1.5.8]{cf:BinGolTeu} is that if $d_n\in R_\ga$ with $\ga>-1$, then
\begin{equation}\label{eq:RV_asymp}
    \sum_{k=1}^n d_k \sim \frac{n d_n}{\ga+1} \qquad (n\to\infty)\,.
\end{equation}


\smallskip

\subsection{Ladder variables and stability}
The first (strict ascending) \textsl{ladder epoch} $T_1$ of a random walk $S_n = X_1 + \ldots + X_n$ is the first time the random walk enters the positive half line, and the corresponding \textsl{ladder height} $H_1$ is the position of the walk at that time:
\[
    T_1:=\inf \{n> 0: S_n >0\} \qquad H_1:=S_{T_1}\,.
\]
Iterating these definitions one gets the following ladder variables: more precisely, for $k>1$ one defines inductively
\[
    T_{k}:=\inf \{n>T_{k-1}: S_n >H_{k-1}\} \qquad H_k:=S_{T_k}\,,
\]
and for convenience we put $(T_0,H_0):=(0,0)$. The \textsl{weak} ascending ladder variables are defined in a similar way, just replacing $>$ by $\geq$ in the relations $(S_n > 0)$ and $(S_n > H_{k-1})$ above. In the following we will rather consider the weak \textsl{descending} ladder variables $(\oT_k, \oH_k)$, which are by definition the weak ascending ladder variables of the walk $\{-S_n\}$. Observe that, by the strong Markov property, both $\{(T_k, H_k)\}_k$ and $\{(\oT_k, \oH_k)\}_k$ are bidimensional renewal processes, that is random walks on $\R^2$ with step law supported in the first quadrant.

\smallskip

It is known that $X_1$ is in the domain of attraction (without centering) of a stable law if and only if $(T_1,H_1)$ lies in a bivariate domain of attraction, cf. \cite{cf:GreOmeTeu,cf:DonGre,cf:Don95}. This fact will play a fundamental role in our derivation: let us specialize it to our setting. By hypothesis $X_1$ is attracted to the normal law, that is $S_n/a_n \Rightarrow \gp(x) \,\dd x$, so that by the standard theory of stability $a_n \in R_{1/2}$. We define two sequences $b_n,\ c_n$ by
\begin{equation} \label{eq:def_b_c}
    \log \frac{n}{\sqrt 2} = \sum_{m=1}^\infty \frac{\gr_m}{m}\, e^{-\frac{m}{b_n}} \qquad c_n := a(b_n)\,,
\end{equation}
where $\gr_m := \bP(S_m > 0)$: then $b_n \in R_2$, $c_n \in R_1$ and we have the weak convergence
\begin{equation} \label{eq:weak_conv}
    \bigg( \frac{T_n}{b_n}\,, \frac{H_n}{c_n} \bigg) \Rightarrow Z\,, \qquad \bP\big(Z \in (\dd x,\,\dd y)\big) = \frac{e^{-1/2x}}{\sqrt{2\pi}\, x^{3/2}} \ind_{(x\geq 0)} \,\dd x\cdot \, \gd_1(\dd y)\,,
\end{equation}
where $\gd_1(\dd y)$ denotes the Dirac measure at $y=1$.

\smallskip

Thus the first ladder epoch $T_1$ is attracted to the positive stable law of index~$1/2$, as for the simple random walk case:
\begin{equation*}
    \frac{T_n}{b_n} \Rightarrow Y, \qquad \bP \big( Y \in \dd x \big) =  \frac{e^{-1/2x}}{\sqrt{2\pi}\, x^{3/2}} \ind_{(x\geq 0)} \,\dd x\,,
\end{equation*}
while for $\{H_k\}$ one has a generalized law of large numbers, with norming sequence~$c_n$: $H_n/c_n \Rightarrow 1$ (that is $H_1$ is \textsl{relatively stable}, cf.~\cite[\S{}8.8]{cf:BinGolTeu}).

\smallskip

We stress that we choose the sequence $a_n$ to be increasing, and by \eqref{eq:def_b_c} $b_n$ and $c_n$ are increasing too. We also recall that the norming sequence $b_n$ is sharply linked to the probability tail of the random variable~$T_1$, by the relation
\begin{equation}\label{eq:rel_stable}
    \bP\big(T_1 > b_n\big)\sim\sqrt{\frac{2}{\pi}}\; \frac{1}{n}\,.
\end{equation}
In fact, this condition is necessary and sufficient in order that a sequence $b_n$ be such that $T_n/b_n \Rightarrow Y$, cf.~\cite[\S{}XIII.6]{cf:Fel2}.

\begin{rem} \rm
It has already been noticed that when the step~$X_1$ has finite (nonzero) variance and zero mean,
\begin{equation*}
     \bE\big(X_1\big)=0 \qquad \bE\big(X_1^2\big)=:\gs^2 \in\, (0,\infty)\,,
\end{equation*}
by the Central Limit Theorem one can take $a_n = \gs\sqrt{n}$ in order that equation \eqref{eq:clt} holds. In other words, $X_1$ is in the \textsl{normal domain of attraction} of the normal law. In this case the first ladder height~$H_1$ is integrable~\cite{cf:Don80} and the behavior of the tail of~$T_1$ is given by
\begin{equation*}
    \bP\big(T_1 > n\big) \sim \frac{2\, \bE(H_1)}{\gs \sqrt{2\pi}} \frac{1}{\sqrt{n}} \qquad (n\to\infty)\,,
\end{equation*}
cf.~\cite[Th.1 in \S{}XII.7 \& Th.1 in \S{}XVIII.5]{cf:Fel2}.
This means that also $T_1$ and $H_1$ belong to the normal domain of attraction of their respective limit law, and one can take
\begin{equation*}
    b_n = \frac{\bE(H_1)^2}{\gs^2}\,n^2 \qquad c_n = \bE(H_1)\, n
\end{equation*}
in order that \eqref{eq:weak_conv} holds (we have used the law of large numbers for $H_1$ and relation \eqref{eq:rel_stable} for~$T_1$).
\end{rem}


\smallskip

\subsection{An asymptotic result}

As an application of the results exposed so far, we derive the asymptotic behavior of $\bP(\cC_n)$ as $n\to\infty$, which will be needed in the sequel. The connection with Fluctuation Theory is given by the fact that
\begin{equation*}
    \cC_n := \big(S_1>0, \ldots, S_n>0\big) = \big(\oT_1 > n\big)\,.
\end{equation*}
In analogy to what we have seen for~$T_1$, the fact that the random walk is attracted to the normal law implies that $\oT_1$ lies in the domain of attraction of the positive stable law of index~$1/2$. Therefore $\bP(\cC_n) \in R_{-1/2}$, and denoting by $\psi(t):=\bE(\exp (-t\oT_1))$ the Laplace transform of~$\oT_1$, by standard Tauberian theorems \cite[Ex.(c) in \S{}XIII.5]{cf:Fel2} we have that
\begin{equation*}
    \bP(\cC_n) \sim \frac{1}{\sqrt{\pi}}\, \big(1-\psi(1/n)\big) \qquad (n\to\infty)\,.
\end{equation*}
Now, for $\psi(t)$ we have the following explicit expression \cite[Th.1 in \S{}XII.7]{cf:Fel2}:
\begin{equation*}
    -\log (1-\psi(t)) \;=\; \sum_{m=1}^\infty \frac{\overline{\gr}_m}{m} e^{-mt} \;=\; -\log(1-e^{-t}) - \sum_{m=1}^\infty \frac{\gr_m}{m} e^{-mt}\,,
\end{equation*}
where $\overline{\gr}_m := \bP(S_m \leq 0)$. A look to \eqref{eq:def_b_c} then yields the desired asymptotic behavior:
\begin{equation}\label{eq:as_cn}
    \bP(\cC_n) \sim \frac{1}{\sqrt{2\pi}}\, \frac{b^{-1}(n)}{n} \qquad (n\to\infty)\,.
\end{equation}


\smallskip

\subsection{Two combinatorial identities}
The power of Fluctuation Theory for the study of random walks is linked to some fundamental identities. The most famous one is the so-called Duality Lemma \cite[\S{}XII]{cf:Fel2} which can be expressed as
\begin{equation} \label{eq:comb_id1}
\bP\big(n\ \text{is a ladder epoch}, S_n \in \dd x\big) = \bP\big(\cC_n, S_n \in \dd x\big) \,,
\end{equation}
where by $(n\ \text{is a ladder epoch})$ we mean of course the disjoint union $\cup_{k\geq 0}(T_k=n)$, and by $\bP(A,Z\in \dd x)$ we denote the finite measure $B\mapsto \bP(A,Z\in B)$. A second important identity, recently discovered by Alili and Doney \cite{cf:AliDon}, will play a fundamental role for us:
\begin{equation} \label{eq:comb_id2}
\bP \big(T_k=n, H_k \in \dd x\big) = \frac{k}{n}\, \bP \big(H_{k-1} < S_n \leq H_k, S_n \in \dd x\big) \,.
\end{equation}

We point out that both the above identities are of purely combinatorial nature, in the sense that they can be proved by relating the events on the two sides with suitable one to one, measure preserving transformations on the sample paths space.


\smallskip

\subsection{Proof of Theorem~\ref{th:as_ren}} \label{sec:as_ren}

We recall that by hypothesis $\gep$ is a fixed positive number. We start from the definition \eqref{eq:def_u} of $u(n,x)$: applying the Duality Lemma~\eqref{eq:comb_id1} we get
\begin{equation} \label{eq:u_1}
    u(n,x) = \bP \big(\,\cC_n, S_n = x\big) = \bP \big(\,\cC_n\big) \,  \bP \big( S_n = x \,\big|\, \cC_n \big)\,.
\end{equation}
Observe that
\begin{equation*}
    \inf_{z\in[\gep,1/\gep]} \gp^+(z) > 0 \qquad \inf_{z\in[\gep,1/\gep]} \gp(z) > 0\,,
\end{equation*}
which implies that both Theorem~\ref{th:lat_main} and Gnedenko's LLT~\eqref{eq:lat_llt} hold also in a ratio sense, namely
\begin{align*}
    \bP\big( S_n=x \,\big|\, \cC_n \big) &= \frac{1}{a_n}\, \gp^+(x/a_n) \, \big( 1+o(1) \big) \qquad (n\to\infty)\\
    \bP\big( S_n=x  \,\big) &= \frac{1}{a_n}\, \gp\,(x/a_n) \, \big( 1+o(1) \big) \qquad (n\to\infty)\,,
\end{align*}
uniformly for $x/a_n\in[\gep, 1/\gep]$. Since $\gp^+(z) =\sqrt{2\pi}\, z\, \gp(z)$ for $z>0$, it follows that
\begin{equation} \label{eq:u_2}
    \bP\big( S_n=x \,\big|\, \cC_n \big) = \sqrt{2\pi}\, \frac{x}{a_n}\, \bP\big( S_n=x  \,\big) \, \big( 1+o(1) \big)  \qquad (n\to\infty)\,,
\end{equation}
uniformly for $x/a_n\in[\gep, 1/\gep]$.

The asymptotic behavior of $\bP(\cC_n)$ is given by~\eqref{eq:as_cn}, and comparing equation \eqref{eq:ext_as_ren} with \eqref{eq:u_2} and \eqref{eq:u_1} we are left with proving that
\begin{equation*}
    U(x) = x\, \frac{b^{-1}(n)}{a(n)} \, \big( 1+o(1) \big) \qquad (n\to\infty)\,,
\end{equation*}
uniformly for $x/a_n\in[\gep, 1/\gep]$. We recall that $U(x)$ is the distribution function of the renewal measure associated to the ladder height process $\{H_k\}$, which is relatively stable, since $H_n/c_n \Rightarrow 1$ as $n\to\infty$. Then Theorem~8.8.1 in \cite{cf:BinGolTeu} gives that $U(x) \sim c^{-1}(x)$ as $x\to\infty$, hence it rests to show that
\begin{equation*}
    \frac{x}{c^{-1}(x)} \, \frac{b^{-1}(n)}{a(n)} \to 1 \qquad (n\to\infty)\,,
\end{equation*}
uniformly for $x/a_n\in[\gep, 1/\gep]$, or equivalently, setting $x=z\,a_n$, that
\begin{equation*}
    \frac{z\,b^{-1}(n)}{c^{-1}(z\,a(n))} \to 1 \qquad (n\to\infty)\,,
\end{equation*}
uniformly for $z\in[\gep,1/\gep]$. However, as $c^{-1}(\cdot)\in R_1$, by \eqref{eq:RV_unif_conv} we have that
\begin{equation*}
    c^{-1}(z\,a(n)) \sim z\,c^{-1}(a(n)) \qquad (n\to\infty)\,,
\end{equation*}
uniformly for $z\in[\gep,1/\gep]$, and the proof is completed observing that $c^{-1}(a(n)) = b^{-1}(n)$, by the definition \eqref{eq:def_b_c} of~$c_n$.\qed


\smallskip
\section{First part of the proof} \label{sec:first}

\smallskip
\subsection{A fundamental expression}

We are going to use Fluctuation Theory to express the law of $S_n^+$ in a more useful way. For $x>0$ and $n>1$ we have
\begin{align*}
    &n \,\bP \big( \cC_n,\, S_n \in \dd x \big) \overset{\eqref{eq:comb_id1}}{=} n \,\bP \big( n \text{ is a ladder epoch},\, S_n \in \dd x \big)\\
     &\qquad = \sum_{r=1}^\infty n\,\bP \big( T_r=n,\, S_n \in \dd x \big) \overset{\eqref{eq:comb_id2}}{=} \sum_{r=1}^\infty r\,\bP \big( H_{r-1} < x \leq H_r,\, S_n \in \dd x \big)\,,
\end{align*}
where we have used both the combinatorial identities \eqref{eq:comb_id1}, \eqref{eq:comb_id2}. With a simple manipulation we get
\begin{align*}
    & \sum_{r=1}^\infty r\,\bP \big( H_{r-1} < x \leq H_r,\, S_n \in \dd x \big) = \sum_{r=1}^\infty \sum_{k=0}^{r-1} \bP \big( H_{r-1} < x \leq H_r,\, S_n \in \dd x \big) \\
    &\qquad = \sum_{k=0}^\infty \sum_{r=k+1}^{\infty} \bP \big( H_{r-1} < x \leq H_r,\, S_n \in \dd x \big) = \sum_{k=0}^\infty \bP \big( H_{k} < x,\, S_n \in \dd x \big)\,,
\end{align*}
and using the Markov property
\begin{align*}
    \bP \big( H_{k} < x,\, S_n \in \dd x \big) = \sum_{m=0}^{n-1} \int_{[0,x)} \bP \big( T_k=m,\, H_{k} \in \dd z\big) \, \bP \big(S_{n-m} \in \dd x - z \big)\,.
\end{align*}
In conclusion we obtain the following relation (which is essentially the same as equation~(10) in~\cite{cf:AliDon}):
\begin{align}
    & \bP \big( S_n/a_n \in \dd x\,\big|\, \cC_n \big) \nonumber\\
    &\qquad =\; \frac{1}{n \bP(\cC_n)} \, \sum_{m=0}^{n-1} \int_{[0,a_n x)} \Bigg( \sum_{k=0}^\infty \bP \big( T_k=m,\, H_{k} \in \dd z\big) \Bigg) \, \bP \big(S_{n-m} \in a_n\dd x - z \big) \nonumber\\
    &\qquad =\; \frac{b^{-1}(n)}{n \bP(\cC_n)} \: \int_{[0,1)\times[0,x)} \dd\mu_n(\ga,\gb)\; \bP \bigg( \frac{S_{\lfloor n(1-\ga) \rfloor}}{a_n} \in \dd x - \gb \bigg) \label{eq:main_llt} \,,
\end{align}
where $\mu_n$ is the finite measure on $[0,1)\times [0,\infty)$ defined by
\begin{equation}\label{eq:def_mun}
    \mu_n (A) \;:=\; \frac{1}{b^{-1}(n)} \, \sum_{k=0}^\infty \bP \bigg( \bigg( \frac{T_k}{n}, \frac{H_k}{a_n} \bigg) \in A \bigg)\,,
\end{equation}
for $n\in\N$ and for any Borel set $A\subseteq [0,1) \times [0,\infty)$. Notice that $\mu_n$ is nothing but a suitable rescaling of the renewal measure associated to the ladder variables process. Also observe that the sum defining~$\mu_n$ can be stopped at $k=n-1$, since by definition $T_k \geq k$ for every~$k$; hence $\mu_n$ is indeed a finite measure.

\smallskip

Before proceeding, we would like to stress the importance of equation~\eqref{eq:main_llt}, which is in a sense the core of our proof. The reason is that in the r.h.s. the conditioning on~$\cC_n$ has disappeared: we are left with a mixture, governed by the measure $\mu_n$, of the laws of $\{S_{\lfloor n(1- \ga) \rfloor}\}_{\ga\in[0,1)}$ \textsl{without conditioning}, and the asymptotic behavior of these laws can be controlled with Stone's Local Limit Theorem \eqref{eq:llt} (if we exclude the values of~$\ga$ close to~1).

\smallskip

In the following subsection we study the asymptotic behavior of the sequence of measures $\{\mu_n\}$, and in the next section we put together these preliminary results to conclude the proof of Theorem~\ref{th:main_llt}.


\smallskip

\subsection{A weak convergence result} We are going to show that as $n\to\infty$ the sequence of measure $\{\mu_n\}$ converges weakly to the finite measure $\mu$ defined by
\begin{equation}\label{eq:def_mu}
    \mu(A) \;:=\; \int_A \dd\ga\,\dd\gb\, \frac{\gb}{\sqrt{2\pi}\,\ga^{3/2}} \, e^{-\gb^2/2\ga}\,,
\end{equation}
for any Borel set $A\subseteq [0,1)\times[0,\infty)$ (it is easy to check that $\mu$ is really a finite measure, see below). Since we are not dealing with probability measures, we must be most precise: we mean weak convergence with respect to the class $C_b$ of bounded and continuous functions on~$\R^2$: $\mu_n \Rightarrow \mu$ iff $\int h \,\dd\mu_n \to \int h \,\dd\mu$ for every $h\in C_b$. If we introduce the distribution functions $F_n,\, F$ of the measures $\mu_n,\, \mu$:
\begin{equation*}
    F_n(a,b) := \mu_n \big( [0,a] \times [0,b] \big) \qquad F(a,b) := \mu \big( [0,a] \times [0,b] \big)\,,
\end{equation*}
then proving that $\mu_n \Rightarrow \mu$ as $n\to\infty$ is equivalent to showing that $F_n(a,b) \to F(a,b)$ for every $(a,b)\in[0,1]\times[0,\infty]$ (notice that $\infty$ is included, because the total mass of~$\mu_n$ is not fixed).

\begin{proposition} \label{prop:weak_conv}
The sequence of measures $\{\mu_n\}$ converges weakly to the measure~$\mu$.
\end{proposition}

\proof
We start checking the convergence of the total mass:
\begin{equation*}
    F_n(1,\infty) = \frac{1}{b^{-1}(n)} \, \sum_{k=0}^\infty \bP \big( T_k \leq n \big) =: \frac{1}{b^{-1}(n)} \, G(n)\,,
\end{equation*}
where $G(n)$ is the distribution function of the renewal measure associated to the ladder epochs process~$\{T_k\}$. There is a sharp link between the asymptotic behavior as $n\to\infty$ of $G(n)$ and that of $\bP(T_1>n)$, given by \cite[Lem. in \S{}XIV.3]{cf:Fel2}:
\begin{equation} \label{eq:ren_meas}
    G(n) \sim \frac{2}{\pi} \, \frac{1}{\bP (T_1 > n)} \qquad (n\to\infty)\,.
\end{equation}
Since from relation \eqref{eq:rel_stable} we have that
\begin{equation*}
    \bP \big(T_1 > n\big) \sim \sqrt{\frac{2}{\pi}} \, \frac{1}{b^{-1}(n)} \qquad (n\to\infty)\,,
\end{equation*}
it follows that $F_n(1,\infty) \to \sqrt{2/\pi}$ as $n\to\infty$. On the other hand, the check that $F(1,\infty) = \sqrt{2/\pi}$ is immediate:
\begin{equation*}
    F(1,\infty) = \frac{1}{\sqrt{2\pi}} \int_0^1 \dd\ga\, \frac{1}{\ga^{3/2}} \int_0^\infty \dd\gb \,\gb\,e^{-\gb^2/2\ga}  = \frac{1}{\sqrt{2\pi}} \int_0^1 \dd\ga\, \frac{1}{\sqrt{\ga}} = \sqrt{\frac{2}{\pi}}\,.
\end{equation*}

\smallskip

Since the total mass converges, we claim that it suffices to show that
\begin{equation} \label{eq:claim}
    \liminf_{n\to\infty} \mu_n \big( (a_1,a_2] \times (b_1,b_2] \big) \;\geq\; \mu \big( (a_1,a_2] \times (b_1,b_2] \big)
\end{equation}
for all $0<a_1<a_2<1$, $0<b_1<b_2<\infty$, and weak convergence will be proved. The (simple) proof of this claim can be found in \S~\ref{app:claim}.

\smallskip

Directly from the definition of~$\mu_n$ we have
\begin{align*}
    \mu_n\big( (a_1,a_2] \times (b_1,b_2] \big) = \frac{1}{b^{-1}(n)} \sum_{k=0}^\infty \bP \bigg( \frac{T_k}{n} \in (a_1, a_2],\, \frac{H_k}{a_n} \in (b_1,b_2]  \bigg)\,.
\end{align*}
We simply restrict the sum to the set of $k$ such that $k/b^{-1}(n)\in (b_1+\gep, b_2-\gep]$, $\gep$ being a small fixed positive number, getting
\begin{equation} \label{eq:riemann}
    \mu_n\big( (a_1,a_2] \times (b_1,b_2] \big) \;\geq\; \frac{1}{b^{-1}(n)} \sum_{s\in\frac{\Z}{b^{-1}(n)} \cap (b_1+\gep, b_2-\gep]} \xi_n(s)\,,
\end{equation}
where
\begin{equation*}
    \xi_n(s) := \bP \bigg( \frac{T_{\lfloor s b^{-1}(n)\rfloor}}{n} \in (a_1, a_2],\, \frac{H_{\lfloor s b^{-1}(n)\rfloor}}{a_n} \in (b_1,b_2]  \bigg) \,.
\end{equation*}
By the definition \eqref{eq:def_b_c} of $c_n$, we have that $a_n = c(b^{-1}(n))$: then, using the weak convergence \eqref{eq:weak_conv} and the uniform convergence property of regularly varying sequences \eqref{eq:RV_unif_conv}, it is not difficult to check that
\begin{equation*}
    \xi_n(s) \to \bP \bigg( Y \in \bigg( \frac{a_1}{s^2}, \frac{a_2}{s^2}\bigg] \bigg) =: \xi(s) \qquad (n\to\infty)\,,
\end{equation*}
\textsl{uniformly} for $s\in (b_1+\gep, b_2-\gep]$.

Observe that the term in the r.h.s. of \eqref{eq:riemann} is a Riemann sum of the function~$\xi_n(s)$ over the bounded interval $(b_1+\gep, b_2-\gep]$. Since the sequence of functions $\{\xi_n(s)\}$ is clearly equibounded and converges uniformly to~$\xi(s)$, it is immediate to check that the r.h.s. of~\eqref{eq:riemann} does converge to the integral of~$\xi(s)$ over $(b_1+\gep, b_2-\gep]$. Therefore
\begin{align*}
    &\liminf_{n\to\infty} \mu_n\big( (a_1,a_2] \times (b_1,b_2] \big) \;\geq\; \int_{b_1+\gep}^{b_2-\gep} \dd s \; \bP \bigg( Y \in \bigg( \frac{a_1}{s^2}, \frac{a_2}{s^2}\bigg] \bigg)\\
    &\qquad =\; \int_{b_1+\gep}^{b_2-\gep} \dd s \int_{a_1/s^2}^{a_2/s^2} \dd z\; \frac{e^{-1/2z}}{\sqrt{2\pi}\,z^{3/2}} \;=\; \int_{b_1+\gep}^{b_2-\gep} \dd s \int_{a_1}^{a_2} \dd t\; \frac{s\,e^{-s^2/2t}}{\sqrt{2\pi}\,t^{3/2}}\\
    &\qquad =\; \mu \big( (a_1,a_2] \times (b_1+\gep, b_2-\gep] \big)\,,
\end{align*}
and letting $\gep\to 0$ relation \eqref{eq:claim} follows.\qed


\smallskip
\section{Second part of the proof} \label{sec:second}


\smallskip
\subsection{General strategy}

Now we are ready to put together the results obtained in the last section. We start by rephrasing relation~\eqref{eq:pos_llt}, which is our final goal, in terms of $S_n/a_n$, a form that is more convenient for our purposes: we have to prove that
\begin{equation}\label{eq:goal}
    \forall K>0\quad \limsup_{n\to\infty}\, a_n \bigg[\; \sup_{x\in\R^+,\, h \leq K/a_n} \Big|\, \bP \big( S_n/a_n \in x+I_h\, \big| \,\cC_n \big) \;-\; h\, \gp^+(x) \,\Big| \; \bigg] = 0\,,
\end{equation}
where $I_h:=[0,h)$, and $x+I_h := [x,x+h)$.

\medskip

Altough the idea behind the proof is quite simple, our arguments depend on an approximation parameter~$\gep$ and there are a number of somewhat technical points. In order to keep the exposition as transparent as possible, it is convenient to introduce the following notation: given two real functions $f(n,x,h,\gep)$ and $g(n,x,h,\gep)$ of the variables $n\in\N$, $x\in\R^+$, $h\in\R^+$ and $\gep\in(0,1)$, we say that $f \mysim g$ if and only if
\begin{equation*}
     \forall K>0\quad \limsup_{\gep\to 0}\, \limsup_{n\to\infty}\, a_n \bigg[ \sup_{x\in\R^+,\,h\leq K/a_n} \big| f(n,x,h,\gep) - g(n,x,h,\gep) \big| \bigg] = 0\,.
\end{equation*}
With this terminology we can reformulate \eqref{eq:goal} as
\begin{equation} \label{eq:goal1}
    \bP \big( S_n/a_n \in x+I_h\, \big| \,\cC_n \big) \;\mysim\; h\:\gp^+(x)\,.
\end{equation}

\medskip

To obtain a more explicit expression of the l.h.s. of \eqref{eq:goal1}, we resort to equation~\eqref{eq:main_llt}: with an easy integration we get
\begin{equation} \label{eq:start}
\bP \big( S_n/a_n \in x+I_h\, \big| \,\cC_n \big) \;=\; \frac{b^{-1}(n)}{n \bP(\cC_n)} \: \int_{D_1^{x+h}} \dd\mu_n(\ga,\gb)\; \hG^{x,h}_n(\ga,\gb) \,,
\end{equation}
where we have introduced the notation $D_a^b := [0,a)\times[0,b)$, and
\begin{equation} \label{eq:def_Ghat}
    \hG^{x,h}_n(\ga,\gb) \;:=\; \bP \bigg( \frac{S_{\lfloor n(1-\ga)\rfloor}}{a_n} \in \big\{ (x - \gb)+ I_h \big\} \cap [0,\infty) \bigg) \,.
\end{equation}
In order to determine the asymptotic behavior of the r.h.s. of \eqref{eq:start}, we recall that:
\begin{itemize}
\item from~\eqref{eq:as_cn} we have
\begin{equation*}
    \frac{b^{-1}(n)}{n \bP(\cC_n)} \to \sqrt{2\pi}\,;
\end{equation*}
\item from Proposition~\ref{prop:weak_conv} we have that $\mu_n \Rightarrow \mu$;
\item from Stone's LLT \eqref{eq:llt} it follows that, for large~$n$, $\hG^{x,h}_n(\ga,\gb)$ is close to
\begin{equation} \label{eq:def_G}
    \cG^{x,h}(\ga,\gb) \;:=\; h \, \frac{1}{\sqrt{1-\ga}}\; \gp \bigg( \frac{x-\gb}{\sqrt{1-\ga}} \bigg)\,,
\end{equation}
where we have used that $a_{n(1-\ga)} \sim \sqrt{1-\ga}\; a_n$ as $n\to\infty$, by \eqref{eq:RV_unif_conv}.
\end{itemize}
In fact, the rest of this section is devoted to showing that
\begin{equation} \label{eq:diff_goal}
    \bP \big( S_n/a_n \in x+I_h\, \big| \,\cC_n \big) \;\mysim\; \sqrt{2\pi} \int_{D_1^{x}} \dd\mu(\ga,\gb)\; \cG^{x,h}(\ga,\gb) \,.
\end{equation}
It may not be a priori obvious whether this coincides with our goal \eqref{eq:goal1}, that is whether
\begin{equation}\label{eq:conj}
    \gp^+(x) \;=\; \sqrt{2\pi} \int_{D_1^x} \dd\mu(\ga,\gb)\; \frac{1}{\sqrt{1-\ga}}\; \gp \bigg( \frac{x-\gb}{\sqrt{1-\ga}} \bigg)\,.
\end{equation}
Indeed this relation holds true: in fact \eqref{eq:diff_goal} implies the weak convergence of $S_n/a_n$ under $\bP(\,\cdot\, |\, \cC_n)$ towards a limiting law with the r.h.s. of~\eqref{eq:conj} as density, and we already know from \eqref{eq:pos_clt} that $S_n/a_n$ under $\bP(\,\cdot\, |\, \cC_n)$ converges weakly to $\gp^+(x) \,\dd x$. Anyway, a more direct verification of \eqref{eq:conj} is also given in \S~\ref{app:integral}.

\medskip

Thus we are left with proving \eqref{eq:diff_goal}, or equivalently
\begin{equation*}
    \int_{D_1^{x+h}} \dd\mu_n(\ga,\gb)\; \hG^{x,h}_n(\ga,\gb) \;\mysim\; \int_{D_1^x} \dd\mu(\ga,\gb)\; \cG^{x,h}(\ga,\gb)\,.
\end{equation*}
Since $\mysim$ is an equivalence relation, this will be done through a sequence of intermediate equivalences:
\begin{equation*}
    \int_{D_1^{x+h}} \dd\mu_n\; \hG^{x,h}_n  \;\mysim\; \ldots \;\mysim\; \ldots \;\mysim\; \ldots \;\mysim\; \int_{D_1^{x}} \dd\mu\; \cG^{x,h}\,,
\end{equation*}
and for ease of exposition the proof has been accordingly split in four steps. The idea is quite simple: we first restrict the domain from $D_1^{x+h}$ to $D_{1-\gep}^x$ (steps~1--2), then we will be able to apply Stone's LLT and Proposition~\ref{prop:weak_conv} to pass from $(\hG_n^{x,h}, \mu_n)$ to $(\cG^{x,h},\mu)$ (step~3), and finally we come back to the domain~$D_1^x$ (step~4).

\smallskip

Before proceeding, we define a slight variant $G_n^{x,h}$ of $\hG_n^{x,h}$:
\begin{equation}\label{eq:def_Gn}
    G^{x,h}_n(\ga,\gb) \;:=\; \bP \bigg( \frac{S_{\lfloor n(1-\ga)\rfloor}}{a_n} \in (x - \gb)+ I_h \bigg)
\end{equation}
(notice that we have simply removed the set $[0,\infty)$, see~\eqref{eq:def_Ghat}) and we establish a preliminary lemma.

\begin{lemma} \label{lem:prel}
For every $K>0$ there exists a positive constant $\tC = \tC(K)$ such that
\begin{equation*}
    G^{x,h}_n (\ga,\gb) \leq \frac{\tC}{a_{\lfloor(1-\ga)n\rfloor}} \qquad \forall n\in\N,\ \forall x,\gb\in\R,\ \forall \ga\in[0,1),\ \forall h \leq K/a_n\,,
\end{equation*}
and the same relation holds also for~$\hG_n^{x,h}(\ga,\gb)$.
\end{lemma}

\proof
Since by definition $\hG_n^{x,h}(\ga,\gb) \leq G_n^{x,h}(\ga,\gb)$, it suffices to prove the relation for~$G_n^{x,h}$. However, this is a simple consequence of Stone's LLT \eqref{eq:llt}, that we can rewrite in terms of $S_n/a_n$ as
\begin{equation} \label{eq:llt_alt}
    \forall K>0\quad \limsup_{l\to\infty}\, a_l \bigg[\; \sup_{y\in\R,\, h' \leq K/a_l} \Big|\, \bP \big( S_l/a_l \in y+I_{h'}\, \big) \;-\; h'\, \gp(y) \,\Big| \; \bigg] = 0\,.
\end{equation}
In fact from this relation, using the triangle inequality and the fact that $\sup_{x\in\R} |\gp(x)| < \infty$, it follows easily that for every~$K>0$
\begin{equation} \label{eq:lem_step}
    a_l\, \bP \big( S_l/a_l \in y+I_{h'}\, \big) \leq \tC \qquad \forall l\in\N,\ \forall y\in\R,\ \forall h' \leq K/a_l\,,
\end{equation}
for some positive constant $\tC=\tC(K)$. Now it suffices to observe that $G_n^{x,h}$ can be written as
\begin{equation} \label{eq:def_G_alt}
    G^{x,h}_n(\ga,\gb) \;=\; \bP \bigg( \frac{S_{\lfloor n(1-\ga)\rfloor}}{a_{\lfloor n(1-\ga)\rfloor}} \in \frac{a_n}{a_{\lfloor n(1-\ga)\rfloor}} (x - \gb)+ I_{\frac{h\,a_n}{a_{\lfloor n(1-\ga)\rfloor}}} \bigg)\,,
\end{equation}
so that we can apply \eqref{eq:lem_step} with $l=\lfloor n(1-\ga) \rfloor$ and analogous substitutions.\qed

\smallskip


\smallskip

\subsection{First step}

In the first intermediate equivalence we pass from the domain $D_1^{x+h}$ to $D_{1-\gep}^{x+h}$, that is we are going to show that
\begin{equation*}
    \int_{D_1^{x+h}} \dd\mu_n\; \hG^{x,h}_n  \;\mysim\; \int_{D_{1-\gep}^{x+h}} \dd\mu_n\; \hG^{x,h}_n \,.
\end{equation*}
This means by definition that for every $K>0$
\begin{equation}\label{eq:step1_alt}
    \limsup_{\gep\to 0} \limsup_{n\to\infty}\, R_n^\gep = 0 \,,
\end{equation}
where $R_n^\gep := \sup_{\{x\in\R^+, \; h\leq K/a_n\}} r_n^\gep(x,h)$ and
\begin{equation*}
    r_n^\gep(x,h) := a_n \int_{[1-\gep,1)\times[0,x+h)} \dd\mu_n(\ga,\gb)\; \hG^{x,h}_n (\ga,\gb)\,.
\end{equation*}

\smallskip

Applying Lemma~\ref{lem:prel} and enlarging the domain of integration, we get
\begin{equation} \label{eq:intermed}
\begin{split}
    R_n^\gep &\leq \tC\,a_n \int_{[1-\gep,1)\times[0,\infty)} \dd\mu_n(\ga,\gb)\; \frac{1}{a_{\lfloor(1-\ga)n\rfloor}} \\
    &= \tC\,a_n \sum_{m=\lfloor (1-\gep)n\rfloor}^{n-1} \Bigg[ \frac{1}{b^{-1}(n)} \, \sum_{k=0}^\infty \bP \big( T_k = m \big) \Bigg] \frac{1}{a_{n-m}} \\
    &= \tC\, \frac{a_n}{b^{-1}(n)} \sum_{m=\lfloor (1-\gep)n\rfloor}^{n-1} \frac{u(m)}{a_{n-m}}\,,
\end{split}
\end{equation}
where in the second line we have applied the definition \eqref{eq:def_mun} of~$\mu_n$, and in the third line we have introduced $u(m) := \sum_{k=0}^\infty \bP(T_k=m)$, which is the mass function of the renewal measure associated to the ladder epochs process~$\{T_k\}$. In the proof of Proposition~\ref{prop:weak_conv} we have encountered the asymptotic behavior of the distribution function $G(n):=\sum_{m=1}^n u(m)$, see~\eqref{eq:ren_meas}. The corresponding local asymptotic behavior for $u(m)$ follows since the sequence $u(m)$ is decreasing in~$m$ (this is a simple consequence of the Duality Lemma \eqref{eq:comb_id1}, see also~\cite[Th.4]{cf:Don97}): hence
\begin{equation*}
    u(m) \sim \frac{1}{\pi} \, \frac{1}{m\bP (T_1 > m)} \sim \frac{1}{\sqrt{2\pi}} \,\frac{b^{-1}(m)}{m} \qquad (m\to\infty)\,,
\end{equation*}
having used~\eqref{eq:rel_stable}. It follows that $u(m) \leq C_1\, b^{-1}(m) /m$ for every~$m$, for some positive constant~$C_1$. Recalling that $b^{-1}(\cdot)$ is increasing, from \eqref{eq:intermed} we get
\begin{align*}
    R_n^\gep & \;\leq\;  \tC\,C_1 \frac{a_n}{b^{-1}(n)} \sum_{m=\lfloor (1-\gep)n\rfloor}^{n-1} \frac{b^{-1}(m)}{m\,a_{n-m}}\\
    &\;\leq\; \tC C_1\, \frac{a_n}{\lfloor(1-\gep)n\rfloor} \sum_{k=1}^{\lfloor \gep n\rfloor} \frac{1}{a_k} \;\leq\; \tC C_1C_2 \,\frac{\gep}{1-\gep}\,\frac{a_n}{a_{\lfloor \gep n\rfloor}}  \,,
\end{align*}
for some positive constant $C_2$: in the last inequality we have used \eqref{eq:RV_asymp}, since $a_n \in R_{1/2}$. Now from \eqref{eq:RV_unif_conv} we have that $a_n/a_{\lfloor \gep n\rfloor} \to 1/\sqrt{\gep}$ as $n\to\infty$, hence
\begin{equation*}
    \limsup_{n\to\infty}\, R_n^\gep \leq C\, \frac{\sqrt{\gep}}{1-\gep}\,,
\end{equation*}
with~$C:=\tC C_1 C_2$, and \eqref{eq:step1_alt} follows.


\smallskip

\subsection{Second step}

Now we show that we can restrict the domain from $D_{1-\gep}^{x+h}$ to $D_{1-\gep}^{x}$:
\begin{equation*}
    \int_{D_{1-\gep}^{x+h}} \dd\mu_n\; \hG^{x,h}_n  \;\mysim\; \int_{D_{1-\gep}^x} \dd\mu_n\; \hG^{x,h}_n \;=\; \int_{D_{1-\gep}^x} \dd\mu_n\; G^{x,h}_n \,,
\end{equation*}
where the equality simply follows from the fact that by definition (see \eqref{eq:def_Ghat} and \eqref{eq:def_Gn})
\begin{equation*}
    \hG^{x,h}_n(\ga,\gb) = G^{x,h}_n(\ga,\gb) \qquad \text{for } \gb \leq x\,.
\end{equation*}
We have to show that for every $K>0$
\begin{equation}\label{eq:step2_alt}
    \limsup_{\gep\to 0} \,\limsup_{n\to\infty}\, Q_n^\gep = 0\,,
\end{equation}
where $Q_n^\gep := \sup_{\{x\in\R^+, \; h\leq K/a_n\}} q_n^\gep(x,h)$ and
\begin{equation*}
    q_n^\gep(x,h) := a_n \int_{[0,1-\gep)\times[x,x+h)} \dd\mu_n(\ga,\gb)\; \hG^{x,h}_n (\ga,\gb)\,.
\end{equation*}

\smallskip

From Lemma~\ref{lem:prel} and from the fact that $a_n$ is increasing we easily get
\begin{equation*}
    q_n^\gep(x,h) \leq \tC\; \frac{a_n}{a_{\lfloor\gep n\rfloor}} \; \mu_n \big( [0,1-\gep)\times[x,x+h) \big)\,.
\end{equation*}
As $a_n\in R_{1/2}$, we have $a_n/a_{\lfloor\gep n\rfloor}\to 1/\sqrt{\gep}$ as $n\to\infty$ by \eqref{eq:RV_unif_conv}, hence for fixed $\gep>0$ we can find a positive constant $C_1=C_1(\gep)$ such that for all~$n\in\N$
\begin{equation*}
    q_n^\gep(x,h) \leq \tC\, C_1 \; \mu_n \big( [0,1-\gep)\times[x,x+h) \big)\,.
\end{equation*}
However the term in the r.h.s. can be easily estimated: using the definition \eqref{eq:def_mun} of~$\mu_n$, for $h\leq K/a_n$ we get
\begin{align*}
    &\mu_n \big( [0,1-\gep)\times[x,x+h) \big) \;=\; \frac{1}{b^{-1}(n)} \sum_{k=0}^\infty \bP \big( T_k < (1-\gep)n,\, H_k \in [a_n x, a_n x + a_n h) \big) \\
    &\qquad \;\leq\; \frac{1}{b^{-1}(n)} \sum_{k=0}^\infty \bP \big( H_k \in [a_n x, a_n x + K) \big) \;\leq\; \frac{1}{b^{-1}(n)}\; \sup_{z\in\R^+} U\big( [z,z+K) \big)\,,
\end{align*}
where $U(\dd x) := \sum_{k=0}^\infty \bP (H_k \in \dd x)$ is the renewal measure associated to the ladder heights process $\{H_k\}$, that we have already encountered in the proof of Theorem~\ref{th:as_ren}. Notice that
\begin{equation*}
    \forall K>0 \qquad \sup_{z\in\R^+} U \big( [z,z+K)\big) =: C_2 < \infty \,,
\end{equation*}
which holds whenever $\{H_k\}$ is a transient random walk, cf.~\cite[Th.1 in \S{}VI.10]{cf:Fel2}. Thus for every fixed~$\gep > 0$
\begin{equation*}
    Q_n^\gep = \sup_{x\in\R^+,\, h\leq K/a_n} q_n^\gep(x,h) \;\leq\; \tC C_1 C_2\, \frac{1}{b^{-1}(n)} \;\to\; 0 \qquad (n\to\infty)\,,
\end{equation*}
and \eqref{eq:step2_alt} follows.


\smallskip

\subsection{Third step}

This is the central step: we prove that
\begin{equation*}
    \int_{D_{1-\gep}^x} \dd\mu_n\; G^{x,h}_n \;\mysim\; \int_{D_{1-\gep}^x} \dd\mu\; \cG^{x,h}\,,
\end{equation*}
that is for every $K>0$
\begin{equation}\label{eq:step3_alt}
    \limsup_{\gep\to 0} \, \limsup_{n\to\infty} \, \sup_{x\in\R^+,\, h\leq K/a_n} a_n\,\Bigg| \int_{D_{1-\gep}^x} \dd\mu_n\, G^{x,h}_n \;-\; \int_{D_{1-\gep}^x} \dd\mu\, \cG^{x,h} \Bigg| = 0\,.
\end{equation}

By the triangle inequality
\begin{equation}\label{eq:int1}
\begin{split}
    &a_n\,\Bigg| \int_{D_{1-\gep}^x} \dd\mu_n\, G^{x,h}_n \;-\; \int_{D_{1-\gep}^x} \dd\mu\, \cG^{x,h} \Bigg| \\
    &\qquad \leq\; a_n\,\int_{D_{1-\gep}^x} \dd\mu_n\, \big| G^{x,h}_n - \cG^{x,h} \big| \;+\; a_n\,\Bigg| \int_{D_{1-\gep}^x} \dd\mu_n\, \cG^{x,h} \;-\; \int_{D_{1-\gep}^x} \dd\mu\, \cG^{x,h} \Bigg|\,,
\end{split}
\end{equation}
and we study separately the two terms in the r.h.s. above.

\smallskip
\subsubsection{First term}

With a rough estimate we have
\begin{equation} \label{eq:int1.5}
\begin{split}
    &a_n\,\int_{D_{1-\gep}^x} \dd\mu_n\, \big| G^{x,h}_n - \cG^{x,h} \big|\\
    &\qquad \leq\; \bigg[ \sup_{n\in\N} \mu_n\big(D_{1}^\infty\big) \bigg] \, \Bigg( \sup_{(\ga,\gb)\in D_{1-\gep}^\infty} a_n\, \Big| G^{x,h}_n(\ga,\gb) - \cG^{x,h}(\ga,\gb) \Big| \Bigg)\,,
\end{split}
\end{equation}
and notice the prefactor in the r.h.s. is bounded since $\mu_n (D_1^\infty) \to \mu (D_1^\infty)$. For the remaining term, we use the triangle inequality and the definition \eqref{eq:def_G} of~$\cG^{x,h}$, getting
\begin{equation}\label{eq:int2}
\begin{split}
    &a_n\,\Big| G^{x,h}_n(\ga,\gb) - \cG^{x,h}(\ga,\gb) \Big|\\
    &\qquad\quad \;\leq\; \bigg(\frac{a_n}{a_{\lfloor(1-\ga)n\rfloor}}\bigg)\, a_{\lfloor(1-\ga)n\rfloor}\, \bigg| G^{x,h}_n(\ga,\gb) - \frac{h\,a_n}{a_{\lfloor(1-\ga)n\rfloor}} \, \gp \bigg(  \frac{a_n\, (x-\gb)}{a_{\lfloor(1-\ga)n\rfloor}} \bigg)  \bigg|\\
    &\qquad\qquad\ \ +\; (h\,a_n)\, \bigg| \frac{a_n}{a_{\lfloor(1-\ga)n\rfloor}} \, \gp \bigg(  \frac{a_n\, (x-\gb)}{a_{\lfloor(1-\ga)n\rfloor}} \bigg) - \frac{1}{\sqrt{1-\ga}} \, \gp \bigg( \frac{x-\gb}{\sqrt{1-\ga}} \bigg) \bigg|\,.
\end{split}
\end{equation}

Let us look at the first term in the r.h.s. above: by the by the uniform convergence property of regularly varying sequences \eqref{eq:RV_unif_conv} we have
\begin{equation} \label{eq:int3}
    \sup_{\ga\in(0,1-\gep)} \bigg| \frac{a_n}{a_{\lfloor(1-\ga)n\rfloor}} - \frac{1}{\sqrt{1-\ga}}\, \bigg| \to 0 \qquad (n\to\infty)\,,
\end{equation}
hence the prefactor is uniformly bounded. For the remaining part, from the expression \eqref{eq:def_G_alt} for~$G^{x,h}_n$ it is clear that one can apply Stone's LLT, see \eqref{eq:llt_alt}, yielding
\begin{equation*}
    \sup_{(\ga,\gb) \in D_{1-\gep}^\infty,\, x\in\R^+,\, h\leq K/a_n} \; a_{\lfloor(1-\ga)n\rfloor}\, \bigg| G^{x,h}_n(\ga,\gb) - \frac{h\,a_n}{a_{\lfloor(1-\ga)n\rfloor}} \, \gp \bigg(  \frac{a_n\, (x-\gb)}{a_{\lfloor(1-\ga)n\rfloor}} \bigg)  \bigg| \;\to\;0
\end{equation*}
as $n\to\infty$.

For the second term in the r.h.s. of \eqref{eq:int2}, notice that the prefactor $(h\, a_n)$ gives no problem since $h\leq K/a_n$ in our limit. On the other hand, it is easily seen that the absolute value is vanishing as $n\to\infty$, uniformly for $(\ga,\gb) \in D_{1-\gep}^\infty$ and for $x\in\R^+$: this is thanks to relation~\eqref{eq:int3} and to the fact that the function $\gp(x)$ is uniformly continuous. Coming back to equation \eqref{eq:int1.5}, we have shown that
\begin{equation}\label{eq:int4}
    \limsup_{n\to\infty} \, \sup_{x\in\R^+,\,h\leq K/a_n} a_n\,\int_{D_{1-\gep}^x} \dd\mu_n\, \big| G^{x,h}_n - \cG^{x,h} \big| \;=\; 0\,.
\end{equation}

\smallskip
\subsubsection{Second term} Using the definition \eqref{eq:def_G} of~$\cG^{x,h}$, the second term in the r.h.s. of equation \eqref{eq:int1} can be written as
\begin{equation}\label{eq:int5}
\begin{split}
    &a_n\,\Bigg| \int_{D_{1-\gep}^x} \dd\mu_n\, \cG^{x,h} \;-\; \int_{D_{1-\gep}^x} \dd\mu\, \cG^{x,h} \Bigg| \\
    &\qquad =\; (h\,a_n)\, \Bigg| \int_{D_{1-\gep}^\infty} \dd\mu_n\, \Psi(\ga, x-\gb) \;-\; \int_{D_{1-\gep}^\infty} \dd\mu\, \Psi(\ga, x-\gb)\; \Bigg|
\end{split}
\end{equation}
where we have introduced the shorthand
\begin{equation*}
    \Psi(s, t) := \frac{1}{\sqrt{1-s}} \, \gp \bigg( \frac{t}{\sqrt{1-s}} \bigg) \,\ind_{(t\geq 0)}\,.
\end{equation*}

As usual, for us $(h\,a_n)\leq K$ and we can thus concentrate on the absolute value in the r.h.s. of~\eqref{eq:int5}.
Observe that, for fixed~$x\geq 0$, the function $(\ga,\gb) \mapsto \Psi(\ga, x-\gb)$ on the domain $D_{1-\gep}^\infty$
is bounded, and continuous except on the line $\gb=x$: since $\mu_n \Rightarrow \mu$, it follows that for fixed~$x$
the r.h.s. of~\eqref{eq:int5} is vanishing as $n\to\infty$. However, we would like the convergence to be uniform
in~$x\in\R^+$: this stronger result holds true too, as one can verify by approximating $\Psi$ with a sequence of
uniformly continuous functions (the details are carried out in \S~\ref{app:unif}). The net result is
\begin{equation} \label{eq:int6}
    \limsup_{n\to\infty} \, \sup_{x\in\R^+,\,h\leq K/a_n} \, a_n\,\Bigg| \int_{D_{1-\gep}^x} \dd\mu_n\, \cG^{x,h} \;-\; \int_{D_{1-\gep}^x} \dd\mu\, \cG^{x,h} \Bigg| \;=\; 0\,.
\end{equation}

\smallskip

Putting together relations \eqref{eq:int1}, \eqref{eq:int4} and \eqref{eq:int6} it is easily seen that \eqref{eq:step3_alt} holds (even without taking the limit in~$\gep$), and the step is completed.


\smallskip

\subsection{Fourth step}

We finally show that
\begin{equation*}
    \int_{D_{1-\gep}^{x}} \dd\mu\; \cG^{x,h} \;\mysim\; \int_{D_{1}^x} \dd\mu\; \cG^{x,h} \,,
\end{equation*}
that is, for every~$K>0$
\begin{equation} \label{eq:step4_alt}
    \limsup_{\gep\to 0}\, \limsup_{n\to\infty}\, \sup_{x\in\R^+,\, h\leq K/a_n} a_n \int_{[1-\gep,1)\times [0,x)}  \dd\mu (\ga,\gb)\, \cG^{x,h}(\ga,\gb) = 0\,.
\end{equation}

\smallskip

This is very easy: observe that
\begin{equation*}
    \cG^{x,h}(\ga,\gb) \leq \frac{h}{\sqrt{2\pi}\, \sqrt{1-\ga}}\,,
\end{equation*}
as one can check from the explicit expressions for $\cG^{x,h}$~\eqref{eq:def_G} and $\gp(x)$~\eqref{eq:clt}. Hence
\begin{equation*}
    a_n \int_{[1-\gep,1)\times [0,x)}  \dd\mu (\ga,\gb)\, \cG^{x,h}(\ga,\gb) \;\leq\; \frac{(h a_n)}{\sqrt{2\pi}} \int_{[1-\gep,1)\times [0,\infty)}  \dd\mu (\ga,\gb) \frac{1}{\sqrt{1-\ga}} \,,
\end{equation*}
and \eqref{eq:step4_alt} follows, because the function
\begin{equation*}
    \big\{ (\ga,\gb) \mapsto (1-\ga)^{-1/2} \big\} \in \bbL^1\big( D_1^\infty,\,\dd\mu \big)\,,
\end{equation*}
as on can easily verify. This completes the proof of Theorem~\ref{th:main_llt}.


\smallskip
\section{Appendix}

\smallskip
\subsection{An elementary fact} \label{app:claim}

We prove the claim stated in the proof of Proposition~\ref{prop:weak_conv}, in a slightly more general context. Namely, let $\mu_n,\,\mu$ be finite measures on the domain $D:=[0,1)\times[0,\infty)$, with $\mu(\partial D)=0$. Assume that $\mu_n(D) \to \mu(D)$ as $n\to\infty$, and that
\begin{equation} \label{eq:claim1}
    \liminf_{n\to\infty} \mu_n \big( (a_1,a_2] \times (b_1,b_2] \big) \;\geq\; \mu \big( (a_1,a_2] \times (b_1,b_2] \big)\,,
\end{equation}
for all $0<a_1<a_2<1$, $0<b_1<b_2<\infty$. What we are going to show is that
\begin{equation} \label{eq:claim2}
    \exists\ \lim_{n\to\infty} \mu_n \big( (a_1,a_2] \times (b_1,b_2] \big) \;=\; \mu \big( (a_1,a_2] \times (b_1,b_2] \big)\,,
\end{equation}
for all $0<a_1<a_2<1$, $0<b_1<b_2<\infty$, and this implies that $\mu_n \Rightarrow \mu$.

\smallskip

Suppose that \eqref{eq:claim2} does not hold: then for some rectangle $Q:= (x_1,x_2] \times (y_1,y_2]$ contained in the interior of~$D$ and for some $\gep>0$ one has
\begin{equation} \label{eq:absurd}
    \limsup_{n\to\infty} \mu_n(Q) \;\geq\; \mu(Q) + \gep\,.
\end{equation}
We introduce for $\eta \in (0,1/2)$ the rectangle $W:=(\eta,1-\eta] \times (\eta, 1/\eta]$: by choosing~$\eta$ sufficiently small we can assume that $W \supseteq Q$ and that
\begin{equation}\label{eq:app_step1}
    \mu(W) \geq \mu(D)-\gep/2
\end{equation}
(we recall that by hypothesis $\mu(\partial D)=0$). The rectangle $W$ can be easily written as a \textsl{disjoint} union
\begin{equation*}
    W = Q \cup \union_{i=1}^4 Q_i\,,
\end{equation*}
where the rectangles $Q_i$ (whose exact definition however is immaterial) are defined by
\begin{eqnarray*}
    &Q_1 := (\eta,1-\eta]\times(\eta,y_1] \qquad &Q_2 := (\eta, x_1] \times (y_1,y_2]\\
    &Q_3 := (x_2, 1-\eta] \times (y_1,y_2] \qquad &Q_4 := (\eta,1-\eta]\times(y_2,1/\eta]\,.
\end{eqnarray*}
Now, on the one hand we have
\begin{equation*}
    \limsup_{n\to\infty} \mu_n(W) \;\leq\; \limsup_{n\to\infty} \mu_n(D) \;=\; \mu(D)\,,
\end{equation*}
but on the other hand
\begin{align*}
    &\limsup_{n\to\infty} \mu_n(W) \;=\; \limsup_{n\to\infty} \mu_n\bigg(Q \cup \union_{i=1}^4 Q_i \bigg) \;\geq\; \limsup_{n\to\infty} \mu_n(Q) + \liminf_{n\to\infty} \mu_n\bigg(\union_{i=1}^4 Q_i \bigg) \\
    &\qquad \overset{\eqref{eq:absurd}}{\geq} \mu(Q) + \gep + \sum_{i=1}^4 \liminf_{n\to\infty} \mu_n( Q_i) \;\overset{\eqref{eq:claim1}}{\geq} \mu(Q) + \gep + \sum_{i=1}^4 \mu(Q_i) \;=\; \gep + \mu(W) \\
    &\qquad \overset{\eqref{eq:app_step1}}{\geq} \mu(D) + \gep/2\,,
\end{align*}
which evidently is absurd, hence \eqref{eq:claim2} holds true.


\smallskip
\subsection{An integral} \label{app:integral}

We are going to give a more direct proof of relation~\eqref{eq:conj}: substituting the explicit expressions for $\gp(x)$, $\gp^+(x)$, $\mu$ given in equations \eqref{eq:clt}, \eqref{eq:pos_clt}, \eqref{eq:def_mu} and performing an elementary change of variable, we can rewrite it as
\begin{equation} \label{eq:conj1}
    x \, e^{-x^2/2} \;=\; \frac{x^2}{\sqrt{2\pi}} \int_0^1 \dd w \int_0^1 \dd z \; \frac{w}{z^{3/2}(1-z)^{1/2}} \, e^{-\frac{x^2}{2} \big[ \frac{w^2}{z} + \frac{(1-w)^2}{(1-z)} \big]}\,.
\end{equation}

Altough it is possible to perform explicitly the integration in the r.h.s. above, it is easier to proceed in a different way. Let $\{B_t\}$ be a standard Brownian motion and let $T_a := \inf\{t: B_t=a\}$ be its first passage time: then the law of $T_a$ is given by
\begin{equation*}
    \bbP \big( T_a \in \dd t \big) = g(a,t)\, \dd t\,, \qquad g(a,t) := \frac{a}{\sqrt{2\pi}\, t^{3/2}} e^{-a^2/2t}\,.
\end{equation*}
By the strong Markov property, for $x>0$ and $w\in(0,1)$ we have the equality in law $T_x \sim T_{wx} + T_{(1-w)x}$\,, where we mean that $T_{wx}$ and $T_{(1-w)x}$ are independent. Therefore
\begin{equation*}
    g(x,1) = \int_0^1 \dd z\, g\big(wx,z\big)\, g\big((1-w)x,1-z\big)\,,
\end{equation*}
and integrating over $w \in (0,1)$ we get
\begin{equation} \label{eq:app_rel1}
    g(x,1) = \int_0^1 \dd w \int_0^1 \dd z\, g\big(wx,z\big)\, g\big((1-w)x,1-z\big)\,.
\end{equation}

Now observe that relation \eqref{eq:conj1} can be written as
\begin{align*}
    g(x,1) &= \int_0^1 \dd w \int_0^1 \dd z\, \frac{1-z}{1-w}\, g\big(wx,z\big)\, g\big((1-w)x,1-z\big)\\
    &= \int_0^1 \dd w \int_0^1 \dd z\, \frac{z}{w}\, g\big(wx,z\big)\, g\big((1-w)x,1-z\big)\,,
\end{align*}
and comparing with \eqref{eq:app_rel1} we are left with showing that
\begin{equation*}
    \int_0^1 \dd w \int_0^1 \dd z\, \bigg( 1- \frac{z}{w} \bigg)\, g\big(wx,z\big)\, g\big((1-w)x,1-z\big)\;=\; 0 \\
\end{equation*}
However, the l.h.s. above can be decomposed in
\begin{equation*}
    \int_0^1 \dd w \int_w^1 \dd z\; \big( \ldots \big) \;+\; \int_0^1 \dd w \int_0^w \dd z\; \big( \ldots \big) \;=:\; I_1 \;+\; I_2 \,,
\end{equation*}
and with a change of variable one easily verifies that $I_1=-I_2$.



\smallskip
\subsection{A uniformity result} \label{app:unif}

We are going to show that
\begin{equation}\label{eq:app_but}
    \limsup_{n\to\infty} \, \sup_{x\in\R^+} \, \Bigg| \int_{D_{1-\gep}^\infty} \dd\mu_n\, \Psi(\ga, x-\gb) \;-\; \int_{D_{1-\gep}^\infty} \dd\mu\, \Psi(\ga, x-\gb)\; \Bigg| \;=\; 0\,,
\end{equation}
where we recall that $D_a^b := [0,a) \times [0,b)$ and the function $\Psi$ is defined by
\begin{equation*}
    \Psi(s, t) := \frac{1}{\sqrt{1-s}} \, \gp \bigg( \frac{t}{\sqrt{1-s}} \bigg) \,\ind_{(t\geq 0)}\,.
\end{equation*}

\smallskip

Let us consider the fixed domain $T:=[0,1-\gep]\times \R$. Here the function $\Psi$ is bounded, $\|\Psi\|_{\infty,T} = 1/\sqrt{2\pi\gep}$, and continuous except on the line $t=0$. We can easily build a family of approximations $\{\Psi_\gd\}$ of~$\Psi$ that are bounded and uniformly continuous on the whole~$T$, setting for $\gd>0$
\begin{equation*}
    \Psi_\gd (s,t) := \begin{cases}
\Psi(s,t) & t \geq 0\\
\Psi(s,0)\cdot (1+t/\gd) & t \in [-\gd, 0]\\
0 & t \leq -\gd
\end{cases} \,.
\end{equation*}
Notice that $\|\Psi_\gd\|_{\infty,T} = \|\Psi\|_{\infty,T}$, and that for $(s,t)\in T$
\begin{equation} \label{eq:app_int}
    \big|\Psi(s,t) - \Psi_\gd (s,t)\big| \leq \|\Psi\|_{\infty,T} \, \ind_{[-\gd, 0]}(t)\,.
\end{equation}

Let us introduce for short the notation $\Psi^x(\ga,\gb):=\Psi (\ga,x-\gb)$, and analogously for~$\Psi_\gd$. From the triangle inequality we get
\begin{equation} \label{eq:app_step2}
\begin{split}
    \Bigg| \int_{D_{1-\gep}^\infty} \dd\mu_n\, \Psi^x & \;-\; \int_{D_{1-\gep}^\infty} \dd\mu\, \Psi^x\; \Bigg| \;\leq\; \int_{D_{1-\gep}^\infty} \dd\mu_n \, \big| \Psi^x - \Psi^x_\gd \big| \\
    &\;+\; \int_{D_{1-\gep}^\infty} \dd\mu \, \big| \Psi^x - \Psi^x_\gd \big| \;+\; \Bigg| \int_{D_{1-\gep}^\infty} \dd\mu_n\, \Psi^x_\gd \;-\; \int_{D_{1-\gep}^\infty} \dd\mu\, \Psi^x_\gd\; \Bigg|\,.
\end{split}
\end{equation}
Using relation \eqref{eq:app_int}, the first two terms in the r.h.s. above can be estimated by
\begin{equation*}
    \|\Psi\|_{\infty,T} \, \Big( \mu_n \big( [0,1-\gep]\times[x,x+\gd] \big) + \mu \big( [0,1-\gep]\times[x,x+\gd] \big) \Big)\,.
\end{equation*}
Since $\mu$ is an absolutely continuous and finite measure, its distribution function is uniformly continuous: therefore for every $\eta>0$ we can take $\gd_0$ sufficiently small so that
\begin{equation*}
    \sup_{x\in\R^+}\, \mu \big( [0,1-\gep]\times[x,x+\gd_0] \big) \leq \frac{\eta}{4\|\Psi\|_{\infty,T}}\,.
\end{equation*}
On the other hand, we know that for every $x\geq 0$
\begin{equation*}
    \mu_n \big( [0,1-\gep]\times[x,x+\gd_0] \big) \to \mu \big( [0,1-\gep]\times[x,x+\gd_0] \big) \qquad (n\to\infty)\,,
\end{equation*}
and this convergence is uniform for $x\in\R^+$, as it can be easily checked. Hence by the triangle inequality we can choose $n_0$ so large that
\begin{equation*}
    \sup_{n\geq n_0} \, \sup_{x\in\R^+}\, \mu_n \big( [0,1-\gep]\times[x,x+\gd_0] \big) \leq \frac{\eta}{2\|\Psi\|_{\infty,T}}\,.
\end{equation*}
Finally, observe that for fixed $\gd_0$ the family of functions $\{\Psi^x_{\gd_0}\}_{x\in\R^+}$ is equibounded and equicontinuous: since $\mu_n \Rightarrow \mu$, from a classical result \cite[Cor. in \S{}VIII.1]{cf:Fel2} we have that the third term in the r.h.s. of \eqref{eq:app_step2} with $\gd=\gd_0$ is vanishing as $n\to\infty$ uniformly for $x\in\R^+$. Therefore we can assume that $n_0$ has been chosen so large that
\begin{equation*}
    \sup_{n\geq n_0} \,  \sup_{x\in\R^+}\, \Bigg| \int_{D_{1-\gep}^\infty} \dd\mu_n \, \Psi^x_{\gd_0} \;-\; \int_{D_{1-\gep}^\infty} \dd\mu\, \Psi^x_{\gd_0}\; \Bigg| \;\leq\; \frac{\eta}{4}\,.
\end{equation*}

Applying the preceding bounds to equation \eqref{eq:app_step2} with $\gd=\gd_0$, we have shown that for every $\eta>0$ we can find $n_0$ such that for every $n\geq n_0$
\begin{equation*}
    \sup_{x\in\R^+}\, \Bigg| \int_{D_{1-\gep}^\infty} \dd\mu_n\, \Psi^x \;-\; \int_{D_{1-\gep}^\infty} \dd\mu\, \Psi^x\; \Bigg| \;\leq\; \eta\,,
\end{equation*}
and equation \eqref{eq:app_but} is proved.



\end{document}